\documentclass[]{aspm}

\articleinfo{}{}{}

\usepackage{verbatim}
\usepackage{amssymb}
\usepackage{amsbsy}
\usepackage{amscd}
\usepackage{amsmath}
\usepackage{amsthm}
\usepackage[mathscr]{eucal}

\newcommand{\nbiga}{\mathcal{A}}
\newcommand{\nbigb}{\mathcal{B}}
\newcommand{\nbigc}{\mathcal{C}}
\newcommand{\nbigd}{\mathcal{D}}

\newcommand{\nbigh}{\mathcal{H}}
\newcommand{\nbigi}{\mathcal{I}}
\newcommand{\nbigj}{\mathcal{J}}

\newcommand{\nbigl}{\mathcal{L}}
\newcommand{\nbigm}{\mathcal{M}}

\newcommand{\nbigo}{\mathcal{O}}
\newcommand{\nbigp}{\mathcal{P}}

\newcommand{\nbigs}{\mathcal{S}}

\newcommand{\nbigu}{\mathcal{U}}
\newcommand{\nbigv}{\mathcal{V}}

\newcommand{\nbigz}{\mathcal{Z}}

\newcommand{\proj}{\mathbb{P}}
\newcommand{\seisuu}{{\mathbb Z}}

\newcommand{\cnum}{{\mathbb C}}
\newcommand{\real}{{\mathbb R}}
\newcommand{\hyperh}{\mathbb{H}}
\newcommand{\hyperk}{\mathbb{K}}

\newcommand{\gbigb}{\mathfrak B}

\newcommand{\gbigi}{\mathfrak I}

\newcommand{\gbigl}{\mathfrak L}

\newcommand{\gbigv}{\mathfrak V}

\newcommand{\gminia}{\mathfrak a}
\newcommand{\gminib}{\mathfrak b}
\newcommand{\gminic}{\mathfrak c}
\newcommand{\gminid}{\mathfrak d}

\newcommand{\gminim}{\mathfrak m}

\newcommand{\gminis}{\mathfrak s}

\newcommand{\vece}{{\boldsymbol e}}

\newcommand{\vecv}{{\boldsymbol v}}
\newcommand{\vecu}{{\boldsymbol u}}
\newcommand{\vecw}{{\boldsymbol w}}

\newcommand{\veca}{{\boldsymbol a}}
\newcommand{\vecb}{{\boldsymbol b}}

\newcommand{\veck}{{\boldsymbol k}}

\newcommand{\lrarr}{\longrightarrow}

\newcommand{\pf}{{\bf Proof}\hspace{.1in}}

\def\Hom{\mathop{\rm Hom}\nolimits}

\def\End{\mathop{\rm End}\nolimits}

\def\Image{\mathop{\rm Im}\nolimits}

\def\Re{\mathop{\rm Re}\nolimits}

\def\Gr{\mathop{\rm Gr}\nolimits}

\def\SL{\mathop{\rm SL}\nolimits}

\def\rank{\mathop{\rm rank}\nolimits}

\def\Gr{\mathop{\rm Gr}\nolimits}

\def\Res{\mathop{\rm Res}\nolimits}

\def\Tr{\mathop{\rm Tr}\nolimits}
\def\vol{\mathop{\rm dvol}\nolimits}

\def\can{\mathop{\rm can}\nolimits}

\def\id{\mathop{\rm id}\nolimits}

\newcommand{\del}{\partial}
\newcommand{\delbar}{\overline{\del}}

\newcommand{\barz}{\overline{z}}
\newcommand{\zbar}{\barz}
\newcommand{\zetabar}{\overline{\zeta}}

\def\Harm{\mathop{\rm Harm}\nolimits}

\newcommand{\supp}{\gminis}

\newcommand{\openopen}[2]{]#1,#2[}
\newcommand{\closedclosed}[2]{[#1,#2]}

\newcommand{\Etilde}{\widetilde{E}}

\newcommand{\wbar}{\overline{w}}

\newcommand{\htilde}{\widetilde{h}}

\newcommand{\Phitilde}{\widetilde{\Phi}}
\newcommand{\vtilde}{\widetilde{v}}

\newcommand{\ftilde}{\widetilde{f}}

\newcommand{\Itilde}{\widetilde{I}}

\newcommand{\Utilde}{\widetilde{U}}
\newcommand{\Dtilde}{\widetilde{D}}
\newcommand{\Xtilde}{\widetilde{X}}

\newcommand{\Ltilde}{\widetilde{L}}
\newcommand{\nbigmtilde}{\widetilde{\nbigm}}

\def\mero{\mathop{\rm mero}\nolimits}

\def\ess{\mathop{\rm ess}\nolimits}
\def\rc{\mathop{\rm c}\nolimits}
\def\an{\mathop{\rm an}\nolimits}
\def\Ai{\mathop{\rm Ai}\nolimits}
\def\Toda{\mathop{\rm Toda}\nolimits}

\newcommand{\Ubar}{\overline{U}}
\newcommand{\Dbar}{\overline{D}}

\newcommand{\Ibar}{\overline{I}}
\newcommand{\Wtilde}{\widetilde{W}}

\newcommand{\gtilde}{\widetilde{g}}

\newcommand{\Ctilde}{\widetilde{C}}

\newcommand{\Stilde}{\widetilde{S}}

\newcommand{\gbigltilde}{\widetilde{\gbigl}}

\newcommand{\Ytilde}{\widetilde{Y}}
\newcommand{\Ybar}{\overline{Y}}

\newcommand{\ztilde}{\widetilde{z}}

\newcommand{\nbigdtilde}{\widetilde{\nbigd}}

\newcommand{\gbigvtilde}{\widetilde{\gbigv}}

\newcommand{\taubar}{\overline{\tau}}

\newcommand{\Ktilde}{\widetilde{K}}

\newcommand{\cnumtilde}{\widetilde{\cnum}}

\newcommand{\Wbar}{\overline{W}}
\newcommand{\Nbar}{\overline{N}}

\newcommand{\Xbar}{\overline{X}}

\newcommand{\projtilde}{\widetilde{\proj}}

\newcommand{\hyperhbar}{\overline{\hyperh}}
\newcommand{\zetatilde}{\widetilde{\zeta}}
\newcommand{\nbigubar}{\overline{\nbigu}}

\newcommand{\ttG}{{\texttt G}}

\newtheorem{thm}{Theorem}[section]
\newtheorem{cor}[thm]{Corollary}

\newtheorem{rem}[thm]{Remark}
\newtheorem{lem}[thm]{Lemma}
\newtheorem{prop}[thm]{Proposition}
\newtheorem{df}[thm]{Definition}

\newtheorem{condition}[thm]{Condition}

\newtheorem{notation}[thm]{Notation}

\title[Isolated singularities]{Isolated singularities of Toda equations and cyclic Higgs bundles}

\author[Q. Li \and T. Mochizuki]{Qiongling Li \and Takuro Mochizuki}

\address{{\textrm Qiongling Li:}\\
Chern Institute of Mathematics and LPMC, Nankai University,\\
Tianjin 300071, China}

\email{qiongling.li@nankai.edu.cn}

\address{{\textrm Takuro Mochizuki:}\\
Research Institute for Mathematical Sciences, Kyoto University, \\
Kyoto 606-8512, Japan}

\email{takuro@kurims.kyoto-u.ac.jp}

\rcvdate{}
\rvsdate{}



\begin{document}

\begin{abstract}
This paper is the second part of our study on the Toda equations
and the cyclic Higgs bundles
associated with $r$-differentials
over non-compact Riemann surfaces.
We classify all the solutions up to boundedness
around the isolated singularity of an $r$-differential
under the assumption that
the $r$-differential is meromorphic or
has some type of essential singularity.
As a result, for example,
we classify all the solutions on
${\mathbb C}$ if the $r$-differential is
a finite sum of the exponential of polynomials.

\end{abstract}

\maketitle

\tableofcontents

\section{Introduction}

\subsection{Harmonic bundles and Toda equations}

\subsubsection{Higgs bundles associated with $r$-differentials
   and harmonic metrics}
 
Let $X$ be any Riemann surface.
We fix a line bundle $K_X^{1/2}$
with an isomorphism
$K_X^{1/2}\otimes K_X^{1/2}\simeq K_X$.
Let $r$ be a positive integer.
We set $\hyperk_{X,r}:=\bigoplus_{i=1}^r K_X^{(r+1-2i)/2}$.
We define the actions of $G_r=\{a\in\cnum\,|\,a^r=1\}$
on $K_{X}^{(r+1-2i)/2}$
by $a\bullet v=a^{i}v$.
They induce a $G_r$-action on $\hyperk_{X,r}$.
For any $r$-differential $q\in H^0(X, K_X^r)$,
let $\theta(q)$ denote the Higgs field
of $\hyperk_{X,r}$ induced by
the identity map $\theta(q)_i: K_X^{(r+1-2i)/2}\rightarrow
K_X^{(r+1-2(i+1))/2}\otimes K_X$
$(i=1,\cdots, r-1)$ and
$\theta(q)_r=q: K_X^{(1-r)/2}\rightarrow K_X^{(r-1)/2}\otimes K_X$. 
Note that
$\theta(q)$ is homogeneous with respect to the $G_r$-action.
It implies that
if $h$ is a harmonic metric of
$(\hyperk_{X,r},\theta(q))$
then $a^{\ast}(h)$ is again a harmonic metric of
$(\hyperk_{X,r},\theta(q))$.

Let $\Harm(q)$ denote the set of $G_r$-invariant
harmonic metrics $h$ of
$(\hyperk_{X,r},\theta(q))$
such that $\det(h)=1$. 
By the $G_r$-invariance,
the decomposition
$\hyperk_{X,r}=\bigoplus_{i=1}^rK_X^{(r+1-2i)/2}$
is orthogonal with respect to
any $h\in \Harm(q)$,
and hence
we obtain the decomposition
$h=\bigoplus h_{|K_X^{(r+1-2i)/2}}$.
Note that
$K_X^{(r+1-2i)/2}$
and $K_X^{(r+1-2(r+1-i))/2}=K_X^{(-r-1+2i)/2}$
are mutually dual.
We say that $h\in\Harm(q)$ is real
if
$h_{|K_X^{(r+1-2i)/2}}$
and
$h_{|K_X^{(-r-1+2i)/2}}$
are mutually dual.
Let $\Harm^{\real}(q)$ denote the subset of
$h\in\Harm(q)$ which are real.

Recall that if $X$ is compact
then the classification of harmonic metrics of
$(\hyperk_{X,r},\theta(q))$ is well known.
Indeed,
if $X$ is hyperbolic,
the Higgs bundle
$(\hyperk_{X,r}(q),\theta(q))$ is stable for any $q$.
Hence, according to the Kobayashi-Hitchin correspondence
for Higgs bundles
due to Hitchin and Simpson
(\cite{Hitchin87}, \cite{s1}),
$(\hyperk_{X,r},\theta(q))$ has a unique harmonic metric $h$
such that $\det(h)=1$.
Moreover,
as observed by Baraglia \cite{Baraglia},
the uniqueness implies that 
$h$ is $G_r$-invariant and real.
In other words,
for a compact hyperbolic Riemann surface $X$,
$\Harm(q)=\Harm^{\real}(q)$
consists of a unique element.
If $X$ is an elliptic curve,
it is easy to see that
$\Harm(q)$ consists of a unique element if $q\neq 0$,
and that $\Harm(0)$ is empty.
If $X$ is $\proj^1$,
it is easy to see that
there is no non-zero holomorphic $r$-differential,
and that $\Harm(0)$ is empty.
Note that if $X$ is compact hyperbolic
and if $q=0$,
as observed by Hitchin and Simpson,
$(\hyperk_{X,r},\theta(0))$ with the harmonic metric
is naturally a polarized variation of Hodge structure,
which particularly implies the $G_r$-invariance.

If $X$ is non-compact,
the uniqueness of harmonic metrics of
the Higgs bundle $(\hyperk_{X,r},\theta(q))$
does not necessarily hold,
and a harmonic metric of unit determinant is
not necessarily $G_r$-invariant.
But, if $X$ is the complement of a finite subset
in a compact Riemann surface $\Xbar$,
and if $q$ is meromorphic on $\Xbar$,
we may still apply 
the Kobayashi-Hitchin correspondence
for singular Higgs bundles
due to Simpson \cite{s2} in the tame case
and Biquard-Boalch \cite{BiquardBoalch}
in the wild case.
(See also \cite{Mochizuki-wild}
for the extension of a wild harmonic bundle
to a filtered Higgs bundle.)
Indeed,
in \cite{GuestItsLin, GuestLin, Toda-lattice, Toda-latticeII},
$\Harm(q)$ was classified
in the case $X=\cnum^{\ast}$ with  $q=z^m(dz)^r$
motivated by the relation with the $tt^{\ast}$-geometry \cite{CV},
and the method in \cite{Toda-lattice, Toda-latticeII}
owes on the Kobayashi-Hitchin correspondence.
(See \cite{GuestItsLin, GuestLin}
for a different approach to the same issue.)

\vspace{.1in}
In this paper,
as a continuation of \cite{Note0},
we investigate a classification of
$G_r$-invariant harmonic metrics on $(\hyperk_{X,r},\theta(q))$
over a more general non-compact Riemann surface $X$
with a more general $r$-differential $q$,
which is not covered by the theory of wild harmonic bundles.
In particular, we shall closely study
the case where $q$ has some type of essential singularity.
 
\subsubsection{Toda equations associated with $r$-differentials}
 
Let $g$ be any K\"{a}hler metric of $X$.
It induces a $G_r$-invariant Hermitian metric $h^{(1)}(g)$
of $\hyperk_{X,r}$.
For any other $G_r$-invariant Hermitian metric $h$
such that $\det(h)=1$,
we obtain a tuple of $\real$-valued functions
$\vecw=(w_1,\ldots,w_r)$ such that $\sum w_i=0$
by the relation
\[
 h_{|K_X^{(r+1-2i)/2}}
=e^{w_i}h^{(1)}(g)_{|K_X^{(r+1-2i)/2}}
=e^{w_i}g^{-(r+1-2i)/2}.
\]
Then, $h$ is contained in $\Harm(q)$ 
if and only if
the following type of Toda equation is satisfied:
{\small
\begin{equation}
\label{eq;20.7.2.1}
\left\{
\begin{array}{l}
\sqrt{-1}\Lambda_g\del\delbar
 w_1=e^{-w_r+w_1}|q|_g^2-e^{-w_1+w_2}
 -\frac{r-1}{2}k_g
  \\
\sqrt{-1}\Lambda_g\del\delbar
  w_i=
  e^{-w_{i-1}+w_i}-e^{-w_i+w_{i+1}}
 -\frac{r+1-2i}{2}k_g\quad (i=2,\ldots,r-1)\\
\sqrt{-1}\Lambda_g\del\delbar
  w_r=
  e^{-w_{r-1}+w_r}-e^{-w_r+w_1}|q|_g^2
 -\frac{1-r}{2}k_g
\end{array}
\right.
\end{equation}
}
Here,
$\Lambda_g$ denotes the adjoint
of the multiplication of the associated K\"ahler form,
and 
$k_g=\sqrt{-1}\Lambda_g R(g)$ is the Gaussian curvature of $g$.
Let $\Toda_1(q,g)$ denote the set of solutions
$\vecw$ of
(\ref{eq;20.7.2.1})
satisfying $\sum w_i=0$.
A solution $\vecw\in\Toda_1(q,g)$ is called real
if $w_i+w_{r+1-i}=0$.
Let $\Toda_1^{\real}(q,g)$ denote the set of
real solutions of (\ref{eq;20.7.2.1}).
As explained, there is a natural bijection
$\Harm(q)\simeq\Toda_1(q,g)$,
which induces
$\Harm^{\real}(q)\simeq\Toda_1^{\real}(q,g)$.

\begin{rem}
The equation {\rm(\ref{eq;20.7.2.1})} is equivalent to
the Toda equation in {\rm\cite{Note0}}.
See {\rm\S\ref{subsection;20.10.10.2}} below.
\end{rem}

\begin{rem}
If we change the sign in the system {\rm(\ref{eq;20.7.2.1})},
we obtain the classical Toda equation studied extensively
in integrable system, e.g. see {\rm\cite{BPW, BP}}.
Geometrically, the classical Toda equation gives rise
to harmonic maps from surface to compact flag manifolds.
\end{rem}

\begin{rem}
A solution of the Toda equation gives rise to
an equivariant harmonic map
$f$ from a universal covering of $X$
to the symmetric space $\SL(r,\cnum)/{\mathop{\rm SU}}(r)$
such that $\Tr(\partial f^{\otimes i})=0$ $(i=1,\cdots,r-1)$
 except for $\Tr(\partial f^{\otimes r})$
 is a nonzero constant multiple of $q$.
In lower rank,
the Toda equation is encoded with much richer geometry. 
If $r=2$, the Toda equation coincides with
the Bochner equation for harmonic maps between surfaces
for a given Hopf differential,
e.g., see {\rm\cite{SchoenYau, WAN, AuWAN, Wolf}}.
For the study of solutions for
given meromorphic quadratic differentials,
one can check {\rm\cite{Gupta, HTTW, Wolf91b}}.
If $r=3$, the Toda equation for a real solution
coincides with Wang's equation for hyperbolic affine spheres
in $\mathbb R^3$ for a given Pick differential,
 e.g., see {\rm\cite{BH1, Labourie, Loftin0, Wang}}.
For the study of solutions for
given meromorphic cubic differentials,
one can check
 {\rm\cite{BH,DumasWolf,Loftin1, Loftin2, Nie}}.
If $r=4$, the Toda equation for a real solution coincides
with the Gauss-Ricci equation for maximal surfaces
in $\mathbb H^{2,2}$, e.g., see {\rm\cite{CTT, TW}}. 
\end{rem}

\subsubsection{Existence and uniqueness of complete solutions}
\label{subsection;20.10.6.30}

A solution $\vecw\in\Toda_1(q,g)$ is called complete
if the metrics $e^{w_{i+1}-w_i}g$ $(i=1,\ldots,r-1)$
are complete.
In terms of harmonic metrics,
it is equivalent to the condition that
the K\"ahler metrics
$g(h)_i$ $(i=1,\ldots,r-1)$
induced by
$h_{|K_X^{(r+1-2(i+1))/2}}\otimes
 h_{|K_X^{(r+1-2i)/2}}^{-1}$
are complete on $X$,
where we naturally identify the tangent bundle $K_X^{-1}$ of $X$
with
$K_X^{(r+1-2(i+1))/2}\otimes
 (K_X^{(r+1-2i)/2})^{-1}$.
The following fundamental theorem is proved in \cite{Note0}.
(See \cite{Note0} for more detailed properties of
complete solutions.)
\begin{thm}[\cite{Note0}]
 \label{thm;20.9.23.100}
Let $X$ be a non-compact Riemann surface
with a holomorphic $r$-differential $q$.
Assume either
(i) $X$ is hyperbolic,
or (ii) $X$ is parabolic and $q\neq 0$.
Here, we say that $X$ is hyperbolic (resp. parabolic)
 if a universal covering of $X$ is the upper half plane
(resp. $\cnum$).
 Then, there uniquely exists a complete solution
 $\vecw^{\rc}\in\Toda_1(q,g)$.
It is real, i.e., $\vecw^{\rc}\in\Toda_1^{\real}(q,g)$.
 Moreover, there exists a complete K\"ahler metric $\gtilde$ of $X$
 such that $e^{w_{i+1}-w_i}g$ $(i=1,\ldots,r-1)$ are mutually bounded
 with $\gtilde$,
 and that $|q|_{\gtilde}$ is bounded.
 \hfill\qed
\end{thm}
When we emphasize the dependence on $(q,g)$,
we use the notation $\vecw^{\rc}(q,g)$.
Let $h^{\rc}$ denote the corresponding
harmonic metric,
which is independent of the choice of $g$.
We shall also use the notation $h^{\rc}(q)$
when we emphasize the dependence on $q$.

\subsubsection{Uniqueness and non-uniqueness of general solutions}

Suppose that $q$ has finitely many zeros.
In this case,
the results in \cite{Note0}
clarify whether general solutions of (\ref{eq;20.7.2.1})
are uniquely determined or not.
Let $N$ be any relatively compact open neighbourhood of
the zero set of $q$.
On $X\setminus N$,
we obtain the K\"ahler metric $|q|^{2/r}:=(q\cdot\overline{q})^{1/r}$.
We proved the following proposition in \cite{Note0}
on the uniqueness.
Note that such a uniqueness was first proved
in \cite{Li}
for $X=\cnum$ with a polynomial $r$-differential $(r=2,3)$.
\begin{prop}[\cite{Note0}]
If $|q|^{2/r}$ induces a complete distance on $X\setminus N$,
then the equation {\rm(\ref{eq;20.7.2.1})} has a unique solution,
i.e.,
$\Toda_1(q,g)=\{\vecw^{\rc}(q,g)\}$
and $\Harm(q)=\{h^{\rc}(q)\}$.
\hfill\qed
\end{prop}

\begin{rem}
In fact, the condition $|q|^{2/r}$ induces a complete distance on $X\setminus N$ is
equivalent to that $X=\bar X\setminus D$ where $\bar X$ is compact and $D$ is finite, 
$q$ is meromorphic on $\bar X$ with poles at each point of $D$ of pole order at least $r$.  
Here is a brief argument. Since $q$ has finitely many zeros, 
one can easily extend $|q|^{2/r}$ to a smooth metric on $X$, which is obviously of finite total curvature.
By using Huber's theorem on complete surfaces with finite total curvature, 
we obtain $X=\bar X\setminus D$ where $\bar X$ is compact and $D$ is finite. 
The rest is proven using a similar argument as in \cite[Lemma B.3]{Nie}.
\end{rem}

As for the non-uniqueness,
the following proposition is proved in \cite{Note0}.

\begin{prop}[\cite{Note0}]
\label{prop;20.9.23.101}
If $q$ has finitely many zeros,
there exists a solution
$\vecw\in \Toda_1^{\real}(q,g)$
such that outside a compact subset $K$, 
\[
 w_i\sim-\frac{r+1-2i}{r}\log|q|_g,\quad 1\leq i\leq r.
\]
In particular,
solutions of {\rm (\ref{eq;20.7.2.1})} are not unique
if $|q|^{2/r}$ does not induce a complete distance of $X\setminus N$.
\hfill\qed
\end{prop}

When the uniqueness does not hold,
it is natural to study
a classification of general solutions of {\rm(\ref{eq;20.7.2.1})}
under a mild assumption for $q$.
For example,
let $X=\cnum^{\ast}$
and $q=z^m(dz/z)^r$ $(m>0)$.
According to
{\rm\cite{GuestItsLin, GuestLin, Toda-lattice, Toda-latticeII}},
for any $h\in \Harm(q)$,
we obtain a tuple of real numbers
$\vecb(h)=(b(h)_1,\ldots,b(h)_r)$
determined by
\[
b(h)_i:=\inf\Bigl\{
c\in\real\,\Big|\,
 |z|^{c+i}|(dz)^{(r+1-2i)/2}|_h\,\,
 \mbox{ is bounded around $0$}
 \Bigr\}.
\]
The correspondence $h\longmapsto \vecb(h)$ induces a bijection
\begin{multline}
 \Harm(q)\simeq \\
 \left\{(b_1,\ldots,b_r)\,\left|\,
 b_1\geq b_2\geq\cdots\geq b_r\geq b_1-m,\,\,
 \sum b_i=-r(r+1)/2
 \right.
 \right\}.
\end{multline}

It is our purpose in this paper to 
pursue a similar classification of $\Harm(q)$
in a more general situation,
which we shall explain in the following subsections.

\subsubsection{Remark on the difference of the conventions
for the induced Hermitian metrics}
\label{subsection;20.10.10.2}

We use a different convention from {\rm\cite{Note0}}
for the induced Hermitian metric on $K_X^{j/2}$.
For a K\"{a}hler metric $g=g_0\,dz\otimes d\zbar$,
 we set $|(dz)^{j/2}|^2_g=(\frac{g_0}{2})^{-j/2}$ in this paper.
(For example, see {\rm\cite{Griffiths-Harris}}.)
In {\rm\cite{Note0}},
we used the norm
 $(|(dz)^{j/2}|'_g)^2=g_0^{-j/2}$.
 In particular,
 $|q|_g^2=2^r(|q|'_g)^2$.
Let $h^{(0)}(g)$ denote the induced Hermitian metric of
$\hyperk_{X,r}$ obtained by the latter convention.
Then, for an $\real^r$-valued function $\vecu$,
a Hermitian metric
$\bigoplus_{i=1}^r e^{u_i}h^{(0)}(g)_{|K_X^{(r+1-2i)/2}}$
is a harmonic metric of
$(\hyperk_{X,r},\theta(q))$
if and only if the following equation is satisfied,
which is studied in \cite{Note0}:
{\small
\begin{equation}
\label{eq;20.10.10.1}
\left\{
\begin{array}{l}
\frac{1}{2}
\sqrt{-1}\Lambda_g\del\delbar
 u_1=e^{-u_r+u_1}|q|_g^{\prime\,2}-e^{-u_1+u_2}
 -\frac{r-1}{4}k_g
 \\
\frac{1}{2}
\sqrt{-1}\Lambda_g\del\delbar
  u_i=
  e^{-u_{i-1}+u_i}-e^{-u_i+u_{i+1}}
 -\frac{r+1-2i}{4}k_g \quad (i=2,\ldots,r-1)\\
\frac{1}{2}
\sqrt{-1}\Lambda_g\del\delbar
  u_r=
  e^{-u_{r-1}+u_r}-e^{-u_r+u_1}|q|_g^{\prime\,2}
 -\frac{1-r}{4}k_g
\end{array}
\right.
 \end{equation}
}
 An $\real^r$-valued function $\vecw$ satisfies {\rm(\ref{eq;20.7.2.1})}
 if and only if
the tuple
$u_i=w_i+\frac{r+1-2i}{2}\log 2$ $(i=1,\ldots,r)$
satisfies {\rm(\ref{eq;20.10.10.1})}.

\subsection{Isolated singularities}

Let $D\subset X$ be a finite subset.
Let us consider the case where
$q$ is a holomorphic $r$-differential on $X\setminus D$,
which is not constantly $0$.
We shall study the classification of
the behaviour of $h\in\Harm(q)$ up to boundedness
around each point $P$ of $D$.
Let $(X_P,z_P)$ denote 
a holomorphic coordinate neighbourhood
around $P$ with $z_P(P)=0$.
We set $X_P^{\ast}:=X_P\setminus\{P\}$.

\subsubsection{The case of poles}
\label{subsection;20.7.2.20}

Let us recall that
if $P$ is a pole of $q$,
a general theory of harmonic bundles
\cite{s1,Mochizuki-wild, Decouple}
allows us to obtain the classification
in terms of parabolic structure up to boundedness.
We have the expression
\[
 q_{|X_P^{\ast}}
 =z_P^{m_P}\alpha_P\cdot(dz_P/z_P)^r,
\]
where $m_P$ denotes an integer,
and $\alpha_P$ induces a holomorphic function
on $X_P$ such that $\alpha_P(P)\neq 0$.
If $X_P$ is sufficiently small,
$\alpha_P$ is nowhere vanishing.

If $m_P\leq 0$,
we obtain the following estimate:
\[
\log
 \bigl|
 (dz_P)^{(r+1-2i)/2}
 \bigr|_h
+\frac{r+1-2i}{2r}\log|z_P|^{m_P-r}
=O(1).
\]
In particular, any $h_1$ and $h_2\in\Harm(q)$
are mutually bounded around $P$.
(See \S\ref{subsection;20.9.24.10} for more details.)

If $m_P>0$,
there exists
$\vecb_P(h)=(b_{P,1}(h),\ldots,b_{P,r}(h))\in\real^r$
satisfying
\begin{equation}
\label{eq;20.10.20.2}
 \left\{
\begin{array}{l}
 b_{P,1}(h)\geq b_{P,2}(h)\geq\cdots\geq b_{P,r}(h)
 \geq b_{P,1}(h)-m_P,
 \\
 \sum_{i=1}^r b_{P,i}(h)=-\frac{r(r+1)}{2},
\end{array}
 \right.
\end{equation}
for which the following estimates hold:
\begin{multline}
\label{eq;20.7.2.2}
 \log
 \bigl|
 (dz_P)^{(r+1-2i)/2}
 \bigr|_h
+(b_{P,i}(h)+i)\log|z_P|
-\frac{k_i(\vecb_P(h))}{2}\log\bigl(-\log|z_P|\bigr)
\\
 =O(1).
\end{multline}
Here,
$\veck(\vecb_P(h))=
(k_1(\vecb_P(h)),\ldots,k_r(\vecb_P(h)))$
denotes the tuple of integers determined by
$\vecb_P(h)$
as in \S\ref{subsection;20.7.2.1}.
It particularly implies that
$h_i\in\Harm(q)$ $(i=1,2)$ are mutually bounded around $P$
if and only if
$\vecb_P(h_1)=\vecb_P(h_2)$.
(See \S\ref{subsection;20.7.2.1} for more details.)

\subsubsection{Holomorphic functions with multiple growth orders}
\label{subsection;20.8.6.1}

To study more general cases,
we introduce a reasonable subclass of essential singularities
of holomorphic functions.
Let $f$ denote a holomorphic function on $X_P^{\ast}$.

Let $\varpi:\Xtilde_D\lrarr X$ denote the oriented real blowing up.
Let $Q\in \varpi^{-1}(P)$.
We fix a branch of $\log z_P$ around $Q$.
\begin{itemize}
 \item
We say that $f$ is regularly bounded around $Q$
if there exist a neighbourhood $\nbigu$ of $Q$ in $\Xtilde_D$,
non-zero complex numbers $\alpha_j$ $(j=1,\ldots,m)$,
and mutually distinct real numbers $c_j$ $(j=1,\ldots,m)$
for which
$\bigl|
f-\sum_{j=1}^m\alpha_jz_P^{\sqrt{-1}c_j}\bigr|
\to 0$ holds as $|z_P|\to 0$ in
$\nbigu\setminus\varpi^{-1}(P)$.
\item
We say that $f$ has a single growth order at $Q$
if there exists
$\gminia(f,Q)\in\bigoplus_{b>0}\cnum z_P^{-b}$,
$a(f,Q)\in\real$
and $j(f,Q)\in\seisuu_{\geq 0}$
such that
$e^{-\gminia(f,Q)}z_P^{-a(f,Q)}(\log z_P)^{-j(f,Q)}f$
is regularly bounded around $Q$.
\item      
We say that $f$ has multiple growth orders at $Q$
if $f$ is expressed as
a finite sum of holomorphic functions $f_i$
$(i=1,\ldots,m)$
with a single growth order at $Q$
such that
$\gminia(f_i,Q)$ $(i=1,\ldots,m)$
are mutually distinct.
\item
We say that $f$ has multiple growth orders at $P$
if $f$ has multiple growth orders at any $Q\in\varpi^{-1}(P)$.
\end{itemize}

\begin{rem}
If $f$ satisfies a linear differential equation
\[
 \del_{z_P}^nf+\sum_{j=0}^{n-1}\gamma_j(z_P)\del_{z_P}^jf=0
\]
for meromorphic functions $\gamma_j$ on $(X_P,P)$,
then $f$ has multiple growth orders at $P$.
(For example, see {\rm\cite[\S II.1]{Majima}}.) 
\end{rem}

Note that
when $f$ has single growth order at $Q\in\varpi^{-1}(P)$,
\[
|\gminia(f,Q)|^{-1}\Re(\gminia(f,Q))
\]
induces a continuous function on
a neighbourhood of $Q$.
The function is also denoted by
$|\gminia(f,Q)|^{-1}\Re(\gminia(f,Q))$.
As the value of the function at $Q$,
we obtain a real number
$|\gminia(f,Q)|^{-1}\Re(\gminia(f,Q))(Q)$.

If $f$ has multiple growth orders at $P$,
there exists a finite subset
$\nbigz(f)\subset \varpi^{-1}(P)$
such that the following holds
at $Q\in\varpi^{-1}(P)\setminus \nbigz(f)$.

\begin{itemize}
 \item $f$ has a single growth order at $Q$.
 \item $|\gminia(f,Q)|^{-1}\Re(\gminia(f,Q))(Q)\neq 0$,
\end{itemize}

Note that $\varpi^{-1}(P)$ $(P\in D)$
are identified with $S^1$.
An interval $I$ of $\varpi^{-1}(P)$
is called special with respect to $f$
if there exists $\alpha_I\in\cnum^{\ast}$
and $\rho_I>0$ such that the following holds.
\begin{itemize}
 \item The length of $I$ is $\pi/\rho_I$.
 \item For any $Q\in I\setminus\nbigz(f)$,
       we obtain that
\[
       \gminia(f,Q)-\alpha_I z_P^{-\rho_I}\in
       \bigoplus_{0<b<\rho_I}\cnum z_P^{-b},
\]
       and that
       $|\alpha_I z_P^{-\rho_I}|^{-1}\Re(\alpha_I z_P^{-\rho_I})(Q)<0$.
       Here, we assume that the branches of $\log z_P$ at $Q$
       are analytically continued along $I$.
\end{itemize}

\subsubsection{The case where $q$ is not meromorphic but has multiple growth orders}

For $P\in D$,
we have the expression $q=f_P(dz_P)^r$,
where $f_P$ is a holomorphic function on $X_P^{\ast}$.
We say that $q$ has multiple growth orders at $P$,
if $f_P$ has multiple growth orders at $P$.

Assume that $q$ is not meromorphic
but has multiple growth orders at $P$.
Let $\nbigs(q,P)$ denote the (possibly empty) set of intervals
in $\varpi^{-1}(P)$
which are special with respect to $f_P$.
Let $\nbigp$ denote the set of
$\veca\in\real^r$ such that
\[
 a_1\geq a_2\geq \cdots\geq a_r\geq a_1-1,
 \quad
 \sum a_i=0.
\]

\begin{thm}[Theorem \ref{thm;20.6.9.30}]
\mbox{{}}  \label{thm;20.8.6.2}
\begin{itemize}
 \item 
For any $h\in\Harm(q)$,
there exist $\veca_I(h)\in\nbigp$ $(I\in\nbigs(q,P))$
and $\epsilon>0$      
such that the following estimates hold
as $|z_P|\to 0$ on
$\bigl\{
 |\arg(\alpha_Iz_P^{-\rho(I)})-\pi|<(1-\delta)\pi/2
\bigr\}$ for any $\delta>0$:
\begin{equation}
\label{eq;20.7.5.10}
 \log\bigl|
 (dz_P)^{(r+1-2i)/2}
 \bigr|_{h}
 +a_i(h)
 \Re\bigl(
 \alpha_Iz_P^{-\rho_I}
 \bigr)
=O\bigl(|z_P|^{-\rho_I+\epsilon}\bigr).
\end{equation}
\item
     If $\veca_I(h_1)=\veca_I(h_2)$
     for any $I\in \nbigs(q,P)$,
     there exists a relatively compact neighbourhood $X_P$
     of $P$
     such that
     $h_1$ and $h_2$ are mutually bounded
     on $X_P^{\ast}$.
 \item
      For any $\veca_I\in\nbigp$ $(I\in \nbigs(q,P))$,
      there exists $h\in\Harm(q)$ such that
      $\veca_I(h)=\veca_I$ $(I\in\nbigs(q,P))$.
\end{itemize}
\end{thm}

Let us explain rough ideas for the proof.
Let $h\in\Harm(q)$.
If $f_P$ has a single growth order at $Q$
with
\[
 |\gminia(f_P,Q)|^{-1}\Re(\gminia(f_P,Q))(Q)>0,
\]
then 
$h$ should be close to a canonical harmonic metric $h_{\can}$
around $Q$
outside of some neighbourhoods of the zeroes of $q$,
which follows from variants of Simpson's main estimate.
(See Proposition \ref{prop;20.6.12.20}.
See \S\ref{subsection;20.9.25.10} for $h_{\can}$.)
Therefore, around such $Q$,
we may apply the subharmonicity of
the difference of two harmonic metrics
(see \S\ref{subsection;20.8.16.10})
to obtain that two harmonic metrics are
mutually bounded.
Let $I$ be an interval of $\varpi^{-1}(P)$
such that
\[
 |\gminia(f_P,Q)|^{-1}\Re(\gminia(f_P,Q))(Q)<0
\]
for any $Q\in I\setminus \nbigz(f_P)$.
Assume that $I$ is maximal among such intervals.
There exist $\rho_I\in\real_{>0}$ and $\alpha_I\in\cnum^{\ast}$
such that
$\gminia(f_P,Q)-\alpha_Iz_P^{-\rho_I}
\in\bigoplus_{0<b<\rho_I}\cnum z_P^{-b}$
for any $Q\in I\setminus \nbigz(f)$.
Moreover, the length of $I$ is not strictly greater than
$\pi/\rho_I$.
If the length of $I$ is strictly smaller than
$\pi/\rho_I$,
we can apply Phragm\'{e}n-Lindel\"{o}f theorem
(Proposition \ref{prop;20.4.27.10})
to obtain that any $h_i\in\Harm(q)$ $(i=1,2)$
are mutually bounded around
the closure $\Ibar$ of $I$
(Theorem \ref{thm;20.6.9.21}).
If the length of $I$ is $\pi/\rho$,
i.e.,
if $I$ is special,
then the Nevanlinna formula (Proposition \ref{prop;20.4.22.11})
says that there exists the tuple $\veca_I(h)$
for any $h\in\Harm(q)$
such that the estimate (\ref{eq;20.7.5.10}) holds.
(See Theorem \ref{thm;20.6.14.10}.)
We can also apply Phragm\'{e}n-Lindel\"{o}f theorem
to obtain that any $h_i\in\Harm(q)$ $(i=1,2)$
are mutually bounded around $\Ibar$
if $\veca_I(h_1)=\veca_I(h_2)$.

For the construction of $h\in\Harm(q)$
such that $\veca_I(h)=\veca_I$ $(I\in\nbigs(q,P))$,
we shall find a neighbourhood $\nbigu_I$ of
$I\in\nbigs(q,P)$,
a $G_r$-invariant Hermitian metric $h_{\veca}$
of $\hyperk_{X\setminus D,r}$ with $\det(h_{\veca})=1$,
and an increasing sequence $Y_i$ of
relatively compact open subsets of $X\setminus D$
such that the following holds.
\begin{itemize} 
 \item The estimate (\ref{eq;20.7.5.10}) holds
       for $h_{\veca}$ and each $I\in \nbigs(q,P)$.
 \item The boundary $\del Y_i$ are smooth,
       and $\bigcup Y_i=X\setminus D$.
 \item For $h_i\in\Harm(q_{|Y_i})$ such that
       $h_{i|\del Y_i}=h_{\veca|\del Y_i}$,
       we obtain
       $\log\Tr(h_{i}\cdot h_{\veca|Y_i}^{-1})
       \leq C|z_P|^{-\rho_I+\epsilon}$
       on $Y_i\cap(\nbigu_I\setminus\varpi^{-1}(P))$,
       where $C>0$ and $0<\epsilon<\rho_I$ are independent of $i$.       
\end{itemize}
Note that there always exists
$h_i\in\Harm(q_{|Y_i})$
such that 
$h_{i|\del Y_i}=h_{\veca|\del Y_i}$
by a theorem of Donaldson
\cite{Donaldson-boundary-value}
(see Proposition \ref{prop;20.5.29.20}).
For any compact subset $K\subset X\setminus D$,
there exists $C_K>0$ such that
$\Tr(h_{i}\cdot h_{\veca|Y_i}^{-1})<C_K$
on $K$ for any $i$ such that $K\subset Y_i$.
By taking a subsequence,
we may assume that
the sequence $h_i$ is convergent to
$h_{\infty}\in\Harm(q)$.
(See Proposition \ref{prop;20.6.15.30}.)
By the uniform estimate 
$\log\Tr(h_{i}h_{\veca|Y_i}^{-1})
\leq C|z_P|^{-\rho_I+\epsilon}$
on $Y_i\cap(\nbigu_I\setminus\varpi^{-1}(P))$,
we obtain the estimate
$\log\Tr(h_{\infty}h_{\veca}^{-1})
\leq C|z_P|^{-\rho_I+\epsilon}$
on $\nbigu_I\setminus\varpi^{-1}(P)$,
and hence
$\veca_I(h_{\infty})=\veca_I$.

\subsection{Global classification
   in the case where $q$ has multiple growth orders}
\label{subsection;20.9.24.2}
   
Suppose that
$q$ has multiple growth orders at each point of $D$.
Let $D_{\mero}$ denote the set of $P\in D$
such that $q$ is meromorphic at $P$.
We divide $D_{\mero}=D_{>0}\sqcup D_{\leq 0}$,
where $D_{>0}$ denotes the set of $P\in D$
such that $m_P>0$,
and $D_{\leq 0}:=D\setminus D_{>0}$.
For $P\in D_{>0}$,
let $\nbigp(q,P)$ denote the set of
$\vecb\in\real^r$
satisfying the condition (\ref{eq;20.10.20.2}).
We set
$D_{\ess}:=D\setminus D_{\mero}$.
We set
$\nbigs(q):=\coprod_{P\in D_{\ess}}\nbigs(q,P)$.
Then, we obtain the map
\begin{equation}
\label{eq;20.7.2.3}
 \Harm(q)\lrarr \prod_{P\in D_{>0}}\nbigp(q,P)
 \times\prod_{I\in \nbigs(q)}\nbigp.
\end{equation}
If $D_{>0}$ is empty,
$\prod_{P\in D_{>0}}\nbigp(q,P)$
denotes a set which consists of one element.
Similarly, if $\nbigs(q)$ is empty,
$\prod_{I\in \nbigs(q)}\nbigp$
denotes a set which consists of one element.

Let $\nbigp^{\real}(q,P)$
denote the set of $\vecb\in\nbigp(q,P)$
such that $b_i+b_{r+1-i}=-r-1$.
Let $\nbigp^{\real}\subset\nbigp$
denote the set of $\veca\in\nbigp$
such that $a_i+a_{r+1-i}=0$.
The map (\ref{eq;20.7.2.3}) induces the following:
\begin{equation}
\label{eq;20.7.2.4}
 \Harm^{\real}(q)\lrarr
 \prod_{P\in D_{>0}}\nbigp^{\real}(q,P)
 \times
 \prod_{I\in\nbigs(q)}\nbigp^{\real}.
\end{equation}
The following theorem
is a special case of
Theorem \ref{thm;20.7.6.100}
below.
\begin{thm}
\label{thm;20.7.6.5}
 If $X$ is compact,
 then the maps
 {\rm(\ref{eq;20.7.2.3})} and
 {\rm(\ref{eq;20.7.2.4})} are
 bijective.
\end{thm}

As a special case of Theorem \ref{thm;20.7.6.5},
we obtain the following corollary
in the meromorphic case.

\begin{cor}
Suppose that $X$ is compact,
and that $q$ is meromorphic on $(X,D)$.
Then, there exists a natural bijection
 $\Harm(q)\simeq \prod_{P\in D_{>0}}\nbigp(q,P)$.
If moreover $m_P\leq 0$ at each $P\in D$,
we obtain $\Harm(q)=\{h^{\rc}(q)\}$.
\hfill\qed
\end{cor}

If $X$ is not compact,
it is reasonable to impose an additional boundary condition
on the behaviour of harmonic metrics
around the infinity of $X$.
Recall that $D$ is assumed to be finite.
Let $\Harm(q;D,\rc)$ denote the set of
$h\in\Harm(q)$
such that
for any relatively compact open neighbourhood $N$ of $D$,
the K\"{a}hler metrics
$g(h)_{i|X\setminus N}$ $(i=1,\ldots,r-1)$
are complete.
(See \S\ref{subsection;20.10.6.30}
for $g(h)_i$.)
It turns out that
$h\in\Harm(q)$ is contained in
$\Harm(q;D,\rc)$
if and only if
for any relatively compact open neighbourhood $N$ of $D$,
$h_{|X\setminus N}$
and $h^{\rc}(q)_{|X\setminus N}$
are mutually bounded.
We set
$\Harm^{\real}(q;D,\rc):=
\Harm^{\real}(q)
\cap
\Harm(q;D,\rc)$.

We obtain the following map
as the restriction of (\ref{eq;20.7.2.3}):
\begin{equation}
\label{eq;20.7.5.1}
 \Harm(q;D,\rc)
\lrarr \prod_{P\in D_{>0}}\nbigp(q,P)
 \times\prod_{I\in \nbigs(q)}\nbigp.
\end{equation}
Similarly,
we obtain the following map
as the restriction of (\ref{eq;20.7.2.4}):
\begin{equation}
\label{eq;20.7.5.2}
 \Harm^{\real}(q;D,\rc)\lrarr
 \prod_{P\in D_{>0}}\nbigp^{\real}(q,P)
 \times
 \prod_{I\in\nbigs(q)}\nbigp^{\real}.
\end{equation}

\begin{thm}[Theorem \ref{thm;20.6.26.10}]
\label{thm;20.7.6.100}
 The maps {\rm(\ref{eq;20.7.5.1})}
 and {\rm(\ref{eq;20.7.5.2})}
 are bijective.
\end{thm}

Let us explain rough ideas
for the proof of Theorem \ref{thm;20.7.6.100}.
For the uniqueness,
the argument in \cite{Note0} is available
even in this situation.
According to \cite[Lemma 3.1]{s1}
(see \S\ref{subsection;20.8.16.10}),
for any two $h_j\in\Harm(q)$ $(j=1,2)$,
we obtain the inequality (\ref{eq;20.7.5.22})
on $X\setminus D$.
If
$\vecb_P(h_1)=\vecb_P(h_2)$ $(P\in D_{>0})$
and
$\veca_I(h_1)=\veca_I(h_2)$ $(I\in \nbigs(q))$,
we obtain that
$h_1$ and $h_2$ are mutually bounded,
and hence
\cite[Lemma 2.2]{s2} implies that
the inequality (\ref{eq;20.7.5.22})
weakly holds on $X$.
Then, as in \cite{Note0},
we may apply Omori-Yau maximum principle
to obtain $h_1=h_2$.
Note that the argument can be simplified
if $X$ is compact.
Indeed,
because $X\setminus D$ is potential theoretically parabolic
in this case,
we obtain $h_1=h_2$
immediately from the mutual boundedness of $h_1$ and $h_2$.
(See Corollary \ref{cor;20.7.6.1} and Remark \ref{rem;20.7.10.1}.)

As for the existence,
for any $(\vecb_P)_{P\in D_{>0}}$
and $(\veca_I)_{I\in\nbigs(q)}$,
we can construct a $G_r$-invariant
Hermitian metric $h_0$ of $\hyperk_{X\setminus D,r}$
such that the following holds.
\begin{itemize}
 \item There exists a relatively compact neighbourhood $N_1$
       of $D$ such that
       $h_{0|X\setminus N_1}=h^{\rc}(q)_{|X\setminus N_1}$.
 \item There exists a relatively compact neighbourhood $N_2\subset N_1$
       such that
       $h_{0|N_2\setminus D}\in\Harm(q_{|N_2\setminus D})$
       and that
       \[
       \vecb_P(h_{0|N_2\setminus D})=\vecb_P\quad (P\in D_{>0}),
       \quad\quad
       \veca_I(h_{0|N_2\setminus D})=\veca_I\quad (I\in\nbigs(q)).
       \]
\end{itemize}
Note that there exists a compact subset $K\subset X\setminus D$
such that $h_{0|X\setminus K}\in\Harm(q_{|X\setminus K})$.
We develop a variant of Kobayashi-Hitchin correspondence
for harmonic bundles
(Proposition \ref{prop;20.6.11.10},
 Proposition \ref{prop;20.6.29.11}),
which allows us to obtain $h\in\Harm(q)$ such that
$h$ and $h_0$ are mutually bounded.
We remark that if $X$ is compact,
the argument can be simplified, again.
Indeed, we immediately obtain the desired metric
by applying a theorem of Simpson \cite{s1}
to $h_0$.

\subsection{Comparison of the methods and the results}

In \cite{Note0} and this paper,
we apply rather different two methods.
The both methods have their own advantage.
In fact, one of the goals of this work is to present readers
as many as tools which can be useful in dealing with
cyclic Higgs bundles over non-compact Riemann surfaces.

On one hand,
the study in \cite{Note0} is more p.d.e theoretic,
and it is available in a quite general setting.
The p.d.e tools like Omori-Yau and  Cheng-Yau maximum principles work well
for complete Riemann manifolds with bounded curvature.
The method is particularly powerful
in analyzing real solutions of (\ref{eq;20.7.2.1}),
which allows us to obtain precise estimates of complete solutions.
Therefore, we manage to use the maximum principles
together with the method of super-subsolution
to construct a complete solution for any Riemann surface
and any holomorphic $r$-differential.

On the other hand,
the study in this paper
heavily owes to the techniques and tools
developed in the theory of
Kobayashi-Hitchin correspondence for Higgs bundles,
pioneered by Donaldson, Hitchin and Simpson.
It is more efficient if the global information on $q$ is provided.
In the case of $(\hyperk_{X,r},\theta(q))$,
the Higgs field is generically regular semisimple,
and the spectral curve is easily described as
the set of $r$-th roots of $q$.
It allows us to use efficiently
variants of Simpson's main estimate
(see \S\ref{subsection;20.4.16.10})
in controlling the behaviour of
harmonic metrics on $(\hyperk_{X,r},\theta(q))$
around the isolated singularities
by using the information of $q$.
Therefore, for instance,
under the assumptions that $X$ is compact,
that $D$ is finite,
and that $q$ has at most multiple growth orders at each point of $D$,
we can classify the set of solutions of the Toda equation
associated with $q$.

\vspace{.1in}

In the case that both methods can be applied,
we obtain different features of solutions
from these two methods.
Suppose that $X$ is compact, and
that $q$ has multiple growth orders
at each point of $D$.
Using Theorem \ref{thm;20.7.6.5},
we obtain the full set $\Harm(q)$ of solutions
$h^{\mathcal B,\mathcal A}$,
where $(\nbigb,\nbiga)=\big((\vecb_P)_{P\in D_{>0}},
 (\veca_I)_{I\in \nbigs(q)}\big)\in \prod_{P\in D_{>0}}\nbigp(q,P)
 \times\prod_{I\in \nbigs(q)}\nbigp$ satisfies
 \begin{eqnarray*}
  &&\sum_i b_{P,i}=-\frac{r(r+1)}{2},
   \quad b_{P,i}\geq b_{P,i+1}\quad(i=1,\cdots, r-1),\\
  && b_{P,r}\geq b_{P,1}-m_P,
  \quad\text{for $P\in D_{>0}$},\\
&&\sum_i a_{I,i}=0,\quad a_{I, 1}\geq a_{I,2}\geq \cdots\geq a_{I,r}\geq a_{I,1}-1, \quad\text{for $I\in \mathcal S(q)$}.
 \end{eqnarray*}
Suppose in addition that $q$ has finitely many zeros.
>From Theorem \ref{thm;20.9.23.100}
and Proposition \ref{prop;20.9.23.101},
two special solutions are singled out.
One is the complete solution in Theorem \ref{thm;20.9.23.100},
which corresponds to
\begin{eqnarray*}
&&b_{P,i}=-\frac{r+1}{2},\quad  
\text{for $P\in D_{> 0}$},\\
 &&a_{I,i}=0, \quad\text{for $I\in \mathcal S(q)$}.
\end{eqnarray*}
The other one is the real solution
purely controlled by $q$
which corresponds to 
\begin{eqnarray*}
 &&b_{P,i}=\frac{r+1-2i}{2r}m_P-\frac{r+1}{2},\quad \text{for $P\in D_{>0}$},\\
 &&a_{I, i}=\frac{r+1-2i}{2r}, \quad\text{for $I\in \mathcal S(q)$}.
\end{eqnarray*}

\subsection{A generalization}
\label{subsection;20.8.5.100}

Assume that $X\setminus D$ is non-compact,
i.e.,
$X$ is non-compact,
or $X$ is compact and $D$ is non-empty.
(See Remark \ref{rem;20.8.6.10}
for the case that $X\setminus D$ is compact,
i.e., $X$ is compact and $D$ is empty.)
Let $L_i$ $(i=1,\ldots,r)$ be 
holomorphic line bundles on $X\setminus D$.
Let
$\psi_i:L_i\lrarr
L_{i+1}\otimes K_{X\setminus D}$ $(i=1,\cdots,r-1)$
and 
$\psi_r:L_r\lrarr
 L_1\otimes K_{X\setminus D}$
be non-zero morphisms.
We set $E:=\bigoplus L_i$.
Let $\theta$ be the cyclic Higgs field 
of $E$ induced by $\psi_i$
$(i=1,\ldots,r)$.

We obtain a holomorphic section  
$q:=\psi_{r}\circ\cdots\circ\psi_1$
of $K_{X\setminus D}^{\otimes r}$,
and a holomorphic section
$q_{\leq r-1}:=
 \psi_{r-1}\circ\cdots\circ\psi_1$
of $\Hom(L_1,L_r)\otimes K_{X\setminus D}^{\otimes (r-1)}$.
We assume the following.
\begin{itemize}
\item
     The zero set of $q_{\leq r-1}$ is finite.
 \item $q$ is not constantly $0$.
       Moreover, $q$ has multiple growth orders at each $P\in D$.
\end{itemize}
Let $D_{\mero}$ denote the set of $P\in D$
such that $q$ is meromorphic at $P$.
We put $D_{\ess}:=D\setminus D_{\mero}$.
For $P\in D_{\mero}$,
we describe
$q_{|X_P^{\ast}}=z_P^{m_P}\beta_P(dz_P/z_P)^r$,
where $\beta_P$ is a nowhere vanishing holomorphic function on $X_P$.
Let $D_{>0}$ denote the set of $P\in D_{\mero}$
such that $m_P>0$.
We set $D_{\leq 0}:=D_{\mero}\setminus D_{>0}$.

Let $h_{\det(E)}$ be a flat metric of $\det(E)$.
For each $P\in D$,
there exist a tuple of frames
$\vecv_P=(v_{P,i}\,|\,i=1,\ldots,r)$ of $L_{i|X_P^{\ast}}$
and a real number $c(\vecv_P)$
such that
$\theta(v_{P,i})=v_{P,i+1}dz_P/z_P$
$(i=1,\ldots,r-1)$
and
\[
 \bigl|
v_{P,1}\wedge\cdots\wedge v_{P,r}
\bigr|_{h_{\det(E)}}=|z_P|^{-c(\vecv_P)}.
\]
Let $\nbigp(q,P,\vecv_P)$
be the set of $\vecb=(b_i)\in\real^r$
such that 
\[
 \sum_{i=1}^r b_i=c(\vecv_P),\quad
 b_i\geq b_{i+1}\,\,(i=1,\ldots,r-1),
\quad
 b_r\geq b_{1}-m_{P}.
\]

\begin{rem}
If $(E,\theta)=(\hyperk_{X\setminus D,r},\theta(q))$,
we canonically choose
\[
 v_{P,i}=z_{P}^{i}(dz_P)^{(r+1-2i)/2}.
\]
Because $\det(\hyperk_{X\setminus D,r})=\nbigo_{X\setminus D}$,
we choose $h_{\det(E)}=1$.
Then, $c(\vecv_P)=-r(r+1)/2$.
\end{rem}

We consider the $G_r$-action on $L_i$
by $a\bullet u_i=a^iu_i$,
which induces a $G_r$-action on $E$.
Let $\Harm^{G_r}(E,\theta,h_{\det(E)})$
denote the set of $G_r$-invariant harmonic metrics $h$
of $(E,\theta)$
such that $\det(h)=h_{\det(E)}$.
We recall that
$D$ is assumed to be a finite subset of $X$.

\begin{prop}[Proposition
 \ref{prop;20.7.13.41},
 Proposition \ref{prop;20.7.13.110}]
\label{prop;20.8.5.110}
The following holds for any $h\in \Harm^{G_r}(E,\theta,h_{\det(E)})$.
\begin{itemize}
 \item For any $P\in D_{>0}$,
       there exists $\vecb_P(h)\in\nbigp(q,P,\vecv_P)$
       determined by the following condition.
       \[
       b_{P,i}(h)=
       \inf\Bigl\{
       b\in\real\,\Big|\,
      \mbox{\rm $|z_P|^b|v_{P,i}|_h$ is bounded}
       \Bigr\}.
       \]
 \item For any $P\in D_{\ess}$ and any $I\in \nbigs(q,P)$,
       there exist $\veca_I(h)\in\nbigp$ and $\epsilon>0$
       such that the following estimates hold
       as $|z_P|\to 0$ on
       $\bigl\{
        |\arg(\alpha_Iz_P^{-\rho(I)})-\pi|<(1-\delta)\pi/2
       \bigr\}$ for any $\delta>0$:
       \[
       \log|v_i|_h+a_{I,i}(h)\Re(\alpha_Iz_P^{-\rho(I)})
       =O\bigl(|z_P|^{-\rho(I)+\epsilon}\bigr)
       \]
       (See {\rm\S\ref{subsection;20.8.6.1}}
       for $\alpha_I$ and $\rho(I)$ for a special interval $I$.)
\end{itemize}
Moreover, if $h_i\in\Harm^{G_r}(E,\theta,h_{\det(E)})$
 satisfy
 $\vecb_P(h_1)=\vecb_P(h_2)$ $(P\in D_{>0})$
 and
 $\veca_I(h_1)=\veca_I(h_2)$
 $(I\in \coprod_{P\in D_{\ess}}\nbigs(q,P))$,
 then $h_1$ and $h_2$
 are mutually bounded
 on $N\setminus D$
 for any relatively compact neighbourhood $N$ of $D$.
\end{prop}

Let $Z(q_{\leq r-1})$ denote the zero set of $q_{\leq r-1}$,
which is assumed to be finite.
We set $\Dtilde=D\cup Z(q_{\leq r-1})$.
Note that $\psi_i$ $(i=1,\ldots,r-1)$ induce
isomorphisms
$K_{X\setminus\Dtilde}^{-1}
\simeq
 (L_{i+1}/L_i)_{|X\setminus\Dtilde}$.
Hence, for any $h\in\Harm^{G_r}(E,\theta,h_{\det(E)})$,
we obtain K\"ahler metrics $g(h)_i$ $(i=1,\ldots,r-1)$
induced by the restrictions of $h$
to $L_{i|X\setminus\Dtilde}$
and $L_{i+1|X\setminus\Dtilde}$.
Let
\[
\Harm^{G_r}(E,\theta,h_{\det(E)};D,\rc)
\]
denote the set of
$h\in\Harm^{G_r}(E,\theta,h_{\det(E)})$
such that
for any relatively compact open neighbourhood $N$
of $\Dtilde$,
the metrics $g(h)_{i|X\setminus N}$ $(i=1,\ldots,r-1)$
are complete.

By Proposition \ref{prop;20.8.5.110},
we obtain the map
\begin{equation}
\label{eq;20.8.5.120}
 \Harm^{G_r}(E,\theta,h_{\det(E)}; D,\rc)
 \lrarr
  \prod_{P\in D_{>0}}\nbigp(q,P,\vecv_P)
  \times
  \prod_{P\in D_{\ess}}
  \prod_{I\in\nbigs(q,P)} \nbigp.
\end{equation}

 \begin{thm}[Theorem
  \ref{thm;20.8.5.121}]
 The map {\rm(\ref{eq;20.8.5.120})}
 is a bijection.
 \end{thm}

\subsection{Examples on $\cnum$}
      
Let $\gamma(z)$ be any non-zero polynomial.
We set $q=\gamma(z)e^{\sqrt{-1}z}(dz)^r$.
It is easy to see that
$\{\infty\}\times\{0<\arg(z)<\pi\}$
is the unique special interval with respect to $q$
in this case.
Hence, we obtain the following proposition
as a corollary of Theorem \ref{thm;20.8.6.2}.
 \begin{prop}[A special case of
  Proposition \ref{prop;20.7.6.2}]
 We have the bijection
 $\Harm(q)\simeq \nbigp$.
\hfill\qed
 \end{prop}

More generally,
let $\gamma_i(z)$ $(i=1,\ldots,m)$ be non-zero polynomials.
Let $\alpha_i$ $(i=1,\ldots,m)$ be
mutually distinct complex numbers.
We set
$q=\sum_{i=1}^m\gamma_i(z)e^{\sqrt{-1}\alpha_i z}(dz)^r$.

\begin{prop}[A special case of 
Proposition \ref{prop;20.7.6.3}]
\label{prop;20.9.24.20}
 If there exists a non-zero complex number $\alpha$
 such that $\alpha_i/\alpha\in\real_{>0}$ $(i=1,\ldots,m)$,
 then there exists a bijection
 $\Harm(q)\simeq \nbigp$.
 Otherwise,
 $\Harm(q)=\{h^{\rc}\}$. 
\hfill\qed
\end{prop}

\begin{cor}
For any non-zero polynomial $\gamma(z)$,
we obtain that
$\Harm\bigl(\gamma(z)\cos(z)(dz)^r\bigr)=\{h^{\rc}\}$.
\hfill\qed
\end{cor}

We may apply Theorem \ref{thm;20.7.6.5}
to the case $q=f(dz)^r$
in which $f$ is the product of a non-zero polynomial
and a non-zero solution of
a linear differential equation with polynomial coefficients.
For instance,
let $\Ai(z)$ be the Airy function given as
in (\ref{eq;20.7.2.40}) below,
which is a solution of
the differential equation
$\del_z^2u-zu=0$.
Let $\gamma(z)$ be any non-zero polynomial.
For the $r$-differential $q=\gamma(z)\Ai(z)(dz)^r$,
$\{-\pi/3<\arg(z)<\pi/3\}$
is the unique special interval in this case.
(See (\ref{eq;20.6.26.10}).
See \cite[\S23]{Wasow} for a more detailed
asymptotic expansion of $\Ai(z)$.)
Hence, we obtain the following.

\begin{prop}[Proposition \ref{prop;20.7.6.4}]
If $q=\gamma(z)\Ai(z)(dz)^r$,
we have a bijection
$\Harm(q)\simeq \nbigp$. 
\hfill\qed
\end{prop}

\subsection{Acknowledgement} 
This study starts from a discussion during the workshop
``Higgs bundles and related topics''
in University of Nice, in 2017.
The first author is partially supported by the National Key 
R$\&$D Program of China No. 2022YFA1006600, 
the Fundamental Research
Funds for the Central Universities and Nankai Zhide foundation.
The second author is grateful to Martin Guest and Claus Hertling
for discussions on harmonic bundles and Toda equations.
The second author is partially supported by
the Grant-in-Aid for Scientific Research (S) (No. 16H06335),
the Grant-in-Aid for Scientific Research (S) (No. 17H06127),
the Grant-in-Aid for Scientific Research (A) (No. 21H04429),
the Grant-in-Aid for Scientific Research (A) (No. 22H00094),
the Grant-in-Aid for Scientific Research (A) (No. 23H00083),
the Grant-in-Aid for Scientific Research (C) (No. 15K04843),
and the Grant-in-Aid for Scientific Research (C) (No. 20K03609),
Japan Society for the Promotion of Science.
The second author is also partially supported by
the Research Institute for Mathematical Sciences,
an International Joint Usage/Research Center located in Kyoto University.

\section{Some existence results of harmonic metrics}

\subsection{Dirichlet problem}

Let $X$ be any Riemann surface.
Let $(E,\delbar_E,\theta)$ be a Higgs bundle on $X$.
Let $Y\subset X$ be a relatively compact connected open subset
with smooth boundary $\del Y$.
Assume that $\del Y$ is non-empty.
Let $h_{\del Y}$ be any Hermitian metric of
$E_{|\del Y}$.
The following proposition is essentially due to Donaldson
\cite{Donaldson-boundary-value}.

\begin{prop}[Donaldson]
\label{prop;20.5.29.20}
 There exists a unique harmonic metric $h$
 of $(E,\delbar_E,\theta)$
 such that $h_{|\del Y}=h_{\del Y}$.
\end{prop}
\pf
Donaldson proved this theorem in the case where $Y$ is a disc
in \cite{Donaldson-boundary-value}.
The general case is similar.
We shall give a proof for the convenience of the readers.
We may assume that $X$ is an open Riemann surface.
According to \cite{Gunning-Narasimhan},
there exists a nowhere vanishing
holomorphic $1$-form $\tau$ om $X$.
Let $f$ be the automorphism of $E$
determined by $\theta=f\,\tau$.
We consider the K\"ahler metric $g_X=\tau\,\taubar$ of $X$.

 Let $\Gamma$ be a lattice of $\cnum$
 and let $T$ be a real $2$-dimensional torus
 obtained as $\cnum/\Gamma$.
 We set $g_T=dz\,d\zbar$.
 We set $\Xtilde:=X\times T$
 with the projection $p:\Xtilde\lrarr X$.
 It is equipped with the flat K\"ahler metric
 $g_{\Xtilde}$ induced by $g_T$ and $g_X$.
 We set $\Ytilde:=p^{-1}(Y)$.

 Let $\Etilde$ be the pull back of $E$
 with the holomorphic structure
 $p^{\ast}(\delbar_{E})+p^{\ast}(f)\,d\zbar$.
According to the dimensional reduction of Hitchin,
a Hermitian metric $h$ of $E_{|Y}$ is
a harmonic metric of $(E,\delbar_E,\theta)_{|Y}$
if and only if
$\Lambda_{\Ytilde}R(p^{\ast}h)=0$.
According to a theorem of
Donaldson \cite{Donaldson-boundary-value},
there exists a unique Hermitian metric $\htilde$
of $\Etilde$
such that
$\Lambda_{\Ytilde}R(\htilde)=0$
and that
$\htilde_{|\del Y}=p^{\ast}(h_{\del Y})$.
By the uniqueness,
$\htilde$ is $T$-invariant.
Hence, there uniquely exists a harmonic metric $h$
of $(E,\delbar_{E},\theta)_{|Y}$
which induces $\htilde$.
It satisfies
$h_{|\del Y}=h_{\del Y}$.
\hfill\qed

\vspace{.1in}
Let $h_0$ be a Hermitian metric of $E$.
Assume that $\det(h_0)$ is flat.
 \begin{cor}
There exists a unique harmonic metric $h$ of
  $E_{|Y}$ such that
  $h_{|\del Y}=h_{0|\del Y}$
  and that
  $\det(h)=\det(h_0)_{|Y}$.
 \end{cor}
 \pf
 There exists a unique harmonic metric $h$
 such that $h_{|\del Y}=h_{0|\del Y}$.
 We obtain
 $\det(h)_{|\del Y}=\det(h_0)_{|\del Y}$.
 Note that 
 both $\det(h)$ and $\det(h_0)_{|Y}$ are flat.
 By the uniqueness in Proposition \ref{prop;20.5.29.20},
 we obtain $\det(h)=\det(h_0)_{|Y}$. 
 \hfill\qed
 
\subsubsection{Homogeneous case}

Let $G$ be a compact Lie group.
Let $\kappa:G\lrarr S^1$ be a character.
Suppose that $X$ is equipped with a $G$-action,
which can be trivial.
A Higgs bundle $(E,\delbar_E,\theta)$
is called $(G,\kappa)$-homogeneous
if $(E,\delbar_E)$ is $G$-equivariant,
and $g^{\ast}(\theta)=\kappa(g)\theta$ for any $g\in G$.

\begin{lem}
 If $Y$ and $h_{\del Y}$ are $G$-invariant
 in Proposition {\rm\ref{prop;20.5.29.20}},
 the harmonic metric $h$ is also $G$-invariant.
\end{lem}
\pf
For any $g\in G$,
$g^{\ast}h$ is a harmonic metric of
$(E,\delbar_E,\kappa(g)\theta)$
such that $(g^{\ast}h)_{|\del Y}=h_{\del Y}$.
Because $|\kappa(g)|=1$,
$g^{\ast}(h)$ is also a harmonic metric of
$(E,\delbar_E,\theta)$.
By the uniqueness,
we obtain $g^{\ast}(h)=h$.
\hfill\qed

\subsubsection{Donaldson functional}

Let $\nbigh(E_{|Y},h_{\del Y})$
be the space of $C^{\infty}$-Hermitian
metrics $h$ of $E_{|Y}$
such that $h_{|\del Y}=h_{\del Y}$.
For two metrics $h_1,h_2\in \nbigh(E_{|Y},h_{\del Y})$,
let $u$ be the automorphism of $E$
which is self-adjoint with respect to both $h_i$,
determined by $h_2=h_1e^u$.
Then, we define
\[
 M(h_1,h_2):=
 \sqrt{-1}\int_Y\Tr(u\Lambda F(h_1))
 +\int_Y
 h_1\Bigl(
 \Psi(u)\bigl((\delbar_E+\theta)u\bigr),
 (\delbar_E+\theta)u
  \Bigr).
\]
(See \S\ref{subsection;20.8.16.10} below
for $F(h_1)$.
See \cite[\S4]{s1} for $\Psi(u)$.)

 \begin{lem}
  Let $h_0$ be any element of $\nbigh(E_{|Y},h_{\del Y})$.
  Let $h$ be the harmonic metric as in
  Proposition {\rm\ref{prop;20.5.29.20}}.
  Then, we obtain $M(h_0,h)\leq 0$.
 \end{lem}
 \pf
 By the dimensional reduction as in the proof of
 Proposition \ref{prop;20.5.29.20},
 the claim is reduced to
 \cite[Proposition 2.18]{KH-infinite}.
 \hfill\qed

\subsection{Convergence}

Let $X$ be an open Riemann surface.
Let $G$ be a compact Lie group acting on $X$.
Let $\kappa:G\lrarr S^1$ be a character.
Let $(E,\delbar_E,\theta)$ be a $(G,\kappa)$-homogeneous
Higgs bundle on $X$.
Let $h_0$ be any $G$-invariant Hermitian metric of $E$.

\begin{df}
 An exhaustive family  $\{X_i\}$ of a Riemann surface $X$
 means an increasing sequence of relatively compact
 $G$-invariant open subsets
 $X_1\subset X_2\subset\cdots$ of $X$
 such that $X=\bigcup X_i$.
 The family is called smooth if
 $\del X_i$ are smooth.
\end{df}

Let $\{X_i\}$ be a smooth exhaustive family of $X$.
The restriction $h_{0|X_i}$ is denoted by $h_{0,i}$.
Let $h_i$ $(i=1,2,\ldots)$
be $G$-invariant harmonic metrics of
$(E,\delbar_{E},\theta)_{|X_i}$.
Let $s_i$ be the automorphism of
$E_{|X_i}$ determined by
$h_i=h_{0,i}\cdot s_i$.
Let $f$ be an $\real_{>0}$-valued function on $X$
such that each $f_{|X_i}$ is bounded.
\begin{prop}
\label{prop;20.5.29.1}
Assume that 
 $|s_i|_{h_{0,i}}+|s^{-1}_i|_{h_{0,i}}\leq f_{|X_i}$
 for any $i$.
 Then, there exists a subsequence
 $s_{i(j)}$ which is convergent
 to a $G$-invariant automorphism $s_{\infty}$ of $E$
 on any relatively compact subset of $X$
 in the $C^{\infty}$-sense.
 As a result, we obtain a $G$-invariant harmonic metric
 $h_{\infty}=h_0s_{\infty}$ of $(E,\delbar_E,\theta)$
 as the limit of the subsequence $h_{i(j)}$.
 Moreover, we obtain
 $|s_{\infty}|_{h_0}+|s_{\infty}^{-1}|_{h_0}\leq f$.
 In particular, if $f$ is bounded, $h_0$ and $h_{\infty}$
 are mutually bounded.
\end{prop}
\pf
We explain an outline of the proof.
Let $g_X$ be a K\"ahler metric of $X$.
According to a general formula (\ref{eq;20.8.16.3}) below,
the following holds on any $X_i$:
\begin{equation}
\label{eq;20.6.30.1}
\sqrt{-1}\Lambda\delbar\del\Tr(s_i)
=
-\Tr\bigl(
s_i\Lambda F(h_{0,i})
\bigr)
-\bigl|(\delbar+\theta)(s_i)\cdot s_i^{-1/2}\bigr|^2_{h_{0,i},g}.
\end{equation}

Let $K$ be any compact subset of $X$.
Let $N$ be a relatively compact neighbourhood of $K$ in $X$.
Let $\chi:X\lrarr\real_{\geq 0}$ be a $C^{\infty}$-function
such that (i) $\chi_{|K}=1$,
(ii) $\chi_{|X\setminus N}=0$,
(iii) $\chi^{-1/2}\del\chi$ and $\chi^{-1/2}\delbar\chi$
on $\{P\in X\,|\,\chi(P)>0\}$ induces a $C^{\infty}$-function on $X$.

There exist $i_0$ such that
$N$ is a relatively compact open subset of $X_i$
for any $i\geq i_0$.
We obtain the following:
\begin{multline}
\label{eq;20.6.30.2}
 \sqrt{-1}\Lambda\delbar\del(\chi\Tr(s_i))
 =
 \chi\sqrt{-1}\Lambda\delbar\del\Tr(s_i)
 +(\sqrt{-1}\Lambda\delbar\del\chi)\cdot\Tr(s_i)
\\
 +\sqrt{-1}\Lambda(\delbar\chi\del\Tr(s_i))
 -\sqrt{-1}\Lambda(\del\chi\delbar\Tr(s_i)).
\end{multline}
Note that
$|\delbar_Es_i|_{h_i,g_X}
=|\del_{E,h_i}s_i|_{h_i,g_X}$,
and that
\begin{equation}
\label{eq;20.6.30.3}
\left|
 \int_X\sqrt{-1}\Lambda(\delbar\chi\del\Tr(s_i))
 \right|
 \leq
 \left(
 \int_X|\chi^{-1/2}\delbar\chi|^2
 \right)^{1/2}
 \left(
 \int_X \chi|\del_{E,h_i} s_i|^2_{h_{0},g_X}
 \right)^{1/2}\!\!\!\!.
\end{equation}
Note that there exists $C_0>0$
such that $|s_i|_{h_0}+|s_i^{-1}|_{h_0}\leq C_0$ on $N$
for any $i$.
By (\ref{eq;20.6.30.1}), (\ref{eq;20.6.30.2})
and (\ref{eq;20.6.30.3}),
there exist $C_j>0$ $(j=1,2)$ such that
the following holds for any sufficiently large $i$:
\begin{multline}
 \int\chi
 \bigl|
 \delbar_Es_i
 \bigr|^2_{h_{0},g_X}
 +\int\chi\bigl|[\theta, s_i]\bigr|^2_{h_{0},g_X}
 \leq
 \\
 C_1+
 C_2\left(
 \int
 \chi\bigl|
 \delbar_Es_i
 \bigr|^2_{h_{0},g_X}
 +\int\chi\bigl|[\theta, s_i]\bigr|^2_{h_{0},g_X}
 \right)^{1/2}.
\end{multline}
Therefore, there exists $C_3>0$ such that
the following holds for any sufficiently large $i$:
\[
\int\chi
 \bigl|
 \delbar_Es_i
 \bigr|^2_{h_{0},g_X}
 +\int\chi\bigl|[\theta, s_i]\bigr|^2_{h_{0},g_X}
 \leq C_3.
\]
We obtain the boundedness of the $L^2$-norms of
$\delbar_Es_i$
and $\del_{E,h_i}s_i$ $(i\geq i_0)$
on $K$
with respect to $h_0$ and $g_X$.
By a variant of Simpson's main estimate
(see \cite[Proposition 2.1]{Decouple}),
we obtain the boundedness of the sup norms of
$\theta$ on $N$ with respect to
$h_i$ and $g_X$.
By the Hitchin equation,
we obtain the boundedness of
the sup norms of 
$\delbar_E(s_i^{-1}\del_{E,h_i} s_i)$ on $N$
with respect to $h_i$ and $g_X$.
By using the elliptic regularity,
we obtain that
the $L_1^p$-norms of
$s_i^{-1}\del_{E,h_i}(s_i)$ on
a relatively compact neighbourhood of $K$
are bounded
for any $p>1$.
It follows that
$L_2^p$-norms of $s_i$ on
a relatively compact neighbourhood of $K$
are bounded for any $p$.
Hence,
a subsequence of $s_i$
is weakly convergent in $L_2^p$
on a relatively compact neighbourhood of $K$.
By the bootstrapping argument
using a general formula (\ref{eq;20.8.16.5}) below,
we obtain that the sequence is convergent
on a relatively compact neighbourhood of $K$
in the $C^{\infty}$-sense.
By using the diagonal argument,
we obtain that 
a subsequence of $s_i$ is weakly convergent
in $C^{\infty}$-sense on any compact subset.
\hfill\qed

\vspace{.1in}

Let us give a complement to Proposition \ref{prop;20.5.29.1}.
 \begin{prop}
  Suppose that $h_{i|\del X_i}=h_{0|\del X_i}$
  and that $\det(h_0)$ is flat.
  Then,
  in Proposition {\rm\ref{prop;20.5.29.1}},
  we obtain  $\det(s_{\infty})=1$.
  Moreover, if $\int_X|F(h_0)|_{h_0,g_X}<\infty$,
  and if $f$ is bounded,
  then $|(\delbar_E+\theta)(s_{\infty})|_{h_0,g_X}$
  is $L^2$,
  where $g_X$ denotes any K\"ahler metric of $X$.
  \end{prop}
  \pf
  Because $\det(s_i)=1$,
  the first claim is clear.
  Let $g_X$ be any K\"ahler metric of $X$. 
Note that $\Tr(s_i)\geq \rank (E)$.
We obtain the following by using Green formula:
\begin{equation}
\label{eq;20.6.30.4}
 \int_{X_i}
  \sqrt{-1}\Lambda\delbar\del\Tr(s_i)\geq 0.
\end{equation}
We obtain the following from
(\ref{eq;20.6.30.1})
and (\ref{eq;20.6.30.4}):
\[
 \int_{X_i}
   \bigl|(\delbar_E+\theta)(s_i)\cdot s_i^{-1/2}\bigr|^2_{h_0,g_X}
   \leq
 -\int_{X_i} \sqrt{-1}\Tr\bigl(s_i\Lambda F(h_0)\bigr).
\]
We obtain the claim of the corollary by Fatou's lemma. 
\hfill\qed

\begin{rem}
In {\rm\cite{Ni1}} and {\rm\cite{KH-infinite}},
a Hermitian-Einstein metric of a holomorphic vector bundle
$(E,\delbar_E)$
on a K\"ahler manifold $Y$
is constructed
as a limit of Hermitian-Einstein metrics
$(E,\delbar_E)_{|Y_i}$
for an exhaustive family $Y_i$ of $Y$.   
For the proof of Proposition {\rm\ref{prop;20.5.29.1}},
we may also apply an argument in
{\rm\cite{Ni1}} and {\rm\cite[\S2.8]{KH-infinite}}
by using the dimensional reduction.
\end{rem}

\subsection{A Kobayashi-Hitchin correspondence}
\label{subsection;20.7.1.2}

Let $X$ be an open Riemann surface.
Let $g_X$ be any K\"ahler metric of $X$.
Let $(E,\delbar_E,\theta)$ be a
$(G,\kappa)$-homogeneous Higgs bundle on $X$.
Let $h_0$ be a $G$-invariant Hermitian metric of $E$
such that $\det(h_0)$ is flat.

 \begin{condition}
\label{condition;24.1.6.30}
  The support of $F(h_0)$ is compact.
 \end{condition}

The condition clearly implies that
$|F(h_0)|$ is integrable on $X$.
For any Higgs subbundle $E'$ of $E$,
let $h_0'$ be the induced metric of $E'$,
and we set
\[
 \deg^{\an}(E',h_0):=
 \int\sqrt{-1}\Lambda\Tr F(h_0').
\]
Recall the Chern-Weil formula \cite[Lemma 3.2]{s1}:
\[
 \deg^{\an}(E',h_0)=
  \int\sqrt{-1}\Tr(p_{E'} \Lambda F(h_0))
 -\int \bigl|(\delbar_E+\theta)p_{E'}\bigr|^2_{h_0,g_X},
\]
where $p_{E'}$ denote the orthogonal projection of $E$
onto $E'$.
By Condition \ref{condition;24.1.6.30},
$\deg^{\an}(E',h_0)$ is well defined in
$\real\cup\{-\infty\}$.
Because $\det(h_0)$ is assumed to be flat,
we obtain $\Tr(F(h_0))=0$,
and hence $\deg^{\an}(E,h_0)=0$.

\begin{df}
 $(E,\delbar_E,\theta,h_0)$ is analytically stable
 with respect to the $G$-action
 if $\deg^{\an}(E',h_0)<0$
 for any proper $G$-equivariant Higgs subbundle $E'$.
\end{df}

\begin{prop}
\label{prop;20.6.11.10}
 If $(E,\delbar_E,\theta,h_0)$ is analytically stable
 with respect to the $G$-action,
 there exists a $G$-invariant harmonic metric $h$ of
 $(E,\delbar_E,\theta)$
 such that
 (i) $h_0$ and $h$ are mutually bounded,
 (ii) $\det(h)=\det(h_0)$,
 (iii) $(\delbar_E+\theta)(h\cdot h_0^{-1})$
 is $L^2$.
\end{prop}
\pf
Let $N_0$ be a relatively compact neighbourhood of
the support of $F(h_0)$.
Let $N_1$ be a relatively compact neighbourhood of $\Nbar_0$.
Let $X_1\subset X_2\subset\cdots$
be a smooth exhaustive sequence of $X$
such that $\Nbar_1\subset X_1$.
We set $h_{0,i}:=h_{0|X_i}$.

There exists a harmonic metric $h_i$
of $(E,\delbar_{E},\theta)_{|X_i}$
such that
$h_{i|\del X_i}=h_{0|\del X_i}$.
Let $s_i$ be the automorphism of $E_{|X_i}$
determined by
$h_{i}=h_{0,i}\cdot s_i$.
Note that $\det(s_i)=1$.
According to the inequality (\ref{eq;20.8.16.6}) below,
we obtain
\begin{equation}
\label{eq;20.5.29.11}
 \sqrt{-1}\Lambda\delbar\del
 \log\Tr(s_i)
 \leq
  \bigl|
   F(h_{0,i})
   \bigr|_{h_{0},g_X}.    
\end{equation}
Note that $\Tr s_{i|\del X_i}=\rank (E)$
and that $F(h_{0,i})=0$ on $X_i\setminus N_0$.
By (\ref{eq;20.5.29.11}),
we obtain the following
\[
 \max_{X_i}\log\Tr(s_i)
 \leq
  \max\Bigl\{
   \max_{N_0}\log\Tr(s_i),
   \rank(E)
   \Bigr\}.
\]
Let $\|\log\Tr(s_i)\|_{h_0,N_1}$ denote the $L^1$-norm of
$\log\Tr(s_i)_{|N_1}$ with respect to $h_0$ and $g_X$.
By using the argument in the proof of \cite[Proposition 2.1]{s1}
together with (\ref{eq;20.5.29.11}),
we obtain $C_j>0$ $(j=0,1)$
depending only on $|F(h_0)|_{h_0,g_X}$
such that
\[
 \max_{N_0}|\log\Tr(s_i)|_{h_0}
 \leq
  C_0\|\log\Tr(s_i)\|_{h_0,N_1}
  +C_1.
\]
Let $u_i$ be the endomorphism of $E_{|X_i}$
which is self-adjoint with respect to $h_{0,i}$ and $h_i$
determined by $s_i=e^{u_i}$.
There exists $C_j$ $(j=2,3)$ depending only on $\rank(E)$
such  that
\[
 |u_i|_{h_0}\leq C_2(\log|s_i|_{h_0}+1),
 \quad\quad
 \log|s_i|_{h_0}\leq C_3(|u_i|_{h_0}+1).
\]
Let $\|u_i\|_{h_0,N_1}$ denote the $L^1$-norm of
$u_{i|N_1}$ with respect to $h_0$ and $g_X$.
There exists $C_j$ $(j=4,5)$,
depending only on $\rank(E)$ and $|F(h_0)|$
such that
\[
 \max_{X_i}|u_i|_{h_0}
\leq
 C_4\|u_i\|_{h_0,N_1}
 +C_5.
\]

The rest is almost the same as 
the proof of \cite[Theorem 2.5]{KH-infinite},
particularly
\cite[Proposition 2.32]{KH-infinite}
and \cite[\S2.7.4]{KH-infinite}.
We explain only an outline.

Suppose that
$\sup_{X_i}|s_i|_{h_{0,i}}\to\infty$ as $i\to\infty$.
It implies
$\|u_i\|_{h_{0},N_1}\to\infty$.
We set
$v_i:=u_i/\|u_i\|_{h_{0},N_1}$.
They are endomorphisms of $E_{|X_i}$
which are self adjoint with respect to
$h_i$ and $h_{0,i}$.
We can prove the following lemma
by using the argument in \cite[Lemma 5.4]{s1}.

\begin{lem}
\label{lem;20.6.30.10}
 There exists a $G$-invariant
 $L_1^2$-section $v_{\infty}$ of $\End(E)$ on $X$
 such that the following holds.
 \begin{itemize}
  \item $v_{\infty}\neq 0$.
  \item A subsequence of $\{v_i\}$ is weakly convergent to
	  $v_{\infty}$ in $L_1^2$ on any compact subset of $X$.
  \item Let $\Phi:\real\times\real\lrarr\real_{>0}$
	be a $C^{\infty}$-function such that
	$\Phi(y_1,y_2)<(y_1-y_2)^{-1}$ if $y_1>y_2$.
	Then, the following holds:
\[
	\sqrt{-1}\int_X\Tr\bigl(
	 v_{\infty}\Lambda F(h_0)
	\bigr)
	+\int_X h_0\Bigl(
	\Phi(v_{\infty})\bigl(
	(\delbar_E+\theta)v_{\infty}\bigr),
	(\delbar_E+\theta)v_{\infty}
	 \Bigr)\leq 0.
\]
	(See {\rm\cite[\S4]{s1}} for $\Phi(v_{\infty})$.)
\hfill\qed
 \end{itemize}
\end{lem}

By applying Lemma \ref{lem;20.6.30.10}
with the argument
in the proof of \cite[Lemma 5.5]{s1},
we obtain that
the eigenvalues of $v_{\infty}$ are constant.
Let $\lambda_1<\cdots<\lambda_{\ell}$ be
the eigenvalues of $v_{\infty}$.
As the orthogonal projections
onto the sum of the eigen spaces associated with
$\lambda_1,\ldots,\lambda_i$,
we obtain an $L_1^2$-section $\pi_i$ of $\End(E)$
for which $\pi_i^2=\pi_i$
and $(\id_E-\pi_i)\circ(\delbar_E+\theta)\pi_i=0$.
(See \cite[\S4]{s1} for more precise construction of $\pi_i$.)
According to the regularity of Higgs $L_1^2$-subbundle
\cite[Proposition 5.8]{s1},
there exists a Higgs subbundle $E_i\subset E$
such that the orthogonal projection onto $E_i$ is equal to $\pi_i$.
By the argument in the proof of \cite[Proposition 5.3]{s1},
we obtain $\deg(E_i,h_0)>0$ for one of $i$,
which contradicts with the stability condition.
Hence,
there exists $C>0$, which is independent of $i$,
such that $|s_i|_{h_{0,i}}<C$.
Then, the claim of Proposition \ref{prop;20.6.11.10}
follows from Proposition \ref{prop;20.5.29.1}.
\hfill\qed

\subsection{Appendix}
\label{subsection;20.8.16.10}

We recall some fundamental formulas due to
Simpson \cite[Lemma 3.1]{s1}
for the convenience of the readers.
Let $(E,\delbar_E,\theta)$ be a Higgs bundle
on a Riemann surface $X$.
For a Hermitian metric $h$ of $E$,
we obtain the Chern connection
$\nabla_h=\delbar_E+\del_{E,h}$
of $(E,\delbar_E,h)$.
The curvature of $\nabla_h$
is denoted by $R(h)$.
We also obtain the adjoint
$\theta^{\dagger}_h$ of $\theta$
with respect to $h$.
The curvature of
$\nabla_h+\theta+\theta_h^{\dagger}$
is denoted by $F(h)$,
i.e.,
$F(h)=R(h)+[\theta,\theta^{\dagger}_h]$.

Let $h_i$ $(i=1,2)$ be Hermitian metrics of $E$.
We obtain the automorphism $s$ of $E$ determined by
$h_2=h_1\cdot s$.
Let $g$ be a K\"ahler metric of $X$,
let $\Lambda$ denote the adjoint of the multiplication of
the associated K\"ahler form.
Then, according to \cite[Lemma 3.1 (a)]{s1},
we obtain the following on $X$:
\begin{multline}
\label{eq;20.8.16.5}
 \sqrt{-1}\Lambda
 \bigl(\delbar_E+\theta\bigr)
 \circ
 \bigl(\del_{E,h_1}
 +\theta^{\dagger}_{h_1}\bigr)s
=
\\
 s\sqrt{-1}\Lambda\bigl(
 F(h_2)-F(h_1)
 \bigr)
+\sqrt{-1}\Lambda
 \Bigl(
 \bigl(\delbar_{E}+\theta\bigr)(s)
  s^{-1}
  \bigl(\del_{E,h_1}
  +\theta^{\dagger}_{h_1}\bigr)(s)
 \Bigr).
\end{multline}
By taking the trace,
and by using \cite[Lemma 3.1 (b)]{s1},
we obtain
\begin{equation}
\label{eq;20.8.16.3}
 \sqrt{-1}\Lambda\delbar\del
 \Tr(s)
=
 \sqrt{-1}
 \Tr\Bigl(
  s\Lambda\bigl(F(h_2)-F(h_1)\bigr)
 \Bigr)
-\Bigl|
 \bigl(\delbar_{E}+\theta\bigr)(s)
  s^{-1/2}
 \Bigr|^2_{h_1,g}.
\end{equation}
Note that
$(\delbar_{E}+\theta)(s)
=\delbar_{E}(s)
+[\theta,s]$.
Moreover,
$\delbar_E(s)$
is a $(0,1)$-form,
and
$[\theta,s]$
is a $(1,0)$-form.
Hence, (\ref{eq;20.8.16.3}) is also rewritten as follows:
\begin{multline}
\label{eq;20.8.16.2}
 \sqrt{-1}\Lambda\delbar\del\Tr(s)
 =
 \sqrt{-1}
 \Tr\Bigl(
  s\Lambda\bigl(F(h_2)-F(h_1)\bigr)
 \Bigr)
\\
 -\bigl|
[\theta,s]s^{-1/2}
\bigr|^2_{h_1,g}
-\bigl|
\delbar_E(s)s^{-1/2}
\bigr|_{h_1,g}.
\end{multline}
We also recall
the following inequality \cite[Lemma 3.1 (d)]{s1}:
\begin{equation}
\label{eq;20.8.16.6}
 \sqrt{-1}\Lambda\delbar\del\log\Tr(s)
  \leq
  \bigl|\Lambda F(h_1)\bigr|_{h_1}
 +\bigl|\Lambda F(h_2)\bigr|_{h_2}.
\end{equation}
In particular,
if both $h_i$ are harmonic,
the functions $\Tr(s)$
and $\log\Tr(s)$ are subharmonic:
\begin{equation}
\label{eq;20.8.17.1}
 \sqrt{-1}\Lambda\delbar\del\Tr(s)
  =-\bigl|
  (\delbar_E+\theta)(s)s^{-1/2}
  \bigr|^2_{h,g}\leq 0,
  \quad
  \sqrt{-1} \Lambda\delbar\del\log\Tr(s)
  \leq 0.
\end{equation}

\section{Preliminaries for
 harmonic metrics of cyclic Higgs bundles}

\subsection{Cyclic Higgs bundles associated with $r$-differentials}

Let $X$ be any Riemann surface
equipped with a line bundle $K_X^{1/2}$
and an isomorphism
$(K_{X}^{1/2})^2=K_{X}$.
We set
$\hyperk_{X,r}:=\bigoplus_{i=1}^rK_{X}^{(r+1-2i)/2}$.
We set $G_r:=\{a\in\cnum\,|\,a^r=1\}$.
We define the $G_r$-actions on $K_{X}^{(r+1-2i)/2}$
by $a\bullet v=a^{i}v$.
They induce a $G_r$-action on $\hyperk_{X,r}$.

For any holomorphic $r$-differential $q$ on $X$,
Let $\theta(q)$ be obtained as
$\theta(q)=\sum_{i=1}^r
 \theta(q)_i$,
 where 
$\theta(q)_i:K_X^{(r+1-2i)/2}\lrarr K_X^{(r+1-2(i+1))/2}\otimes K_X$
$(i=1,\ldots,r-1)$
are induced by the identity,
and
$\theta(q)_r:K_X^{(-r+1)/2}\lrarr K_X^{(r-1)/2}\otimes K_X$
is induced by the multiplication of $q$.
Let $\Harm(q)$ denote the set of $G_r$-invariant
harmonic metrics $h$ of
$(\hyperk_{X,r},\theta(q))$
such that $\det(h)=1$.
By the $G_r$-invariance,
the decomposition
$\hyperk_{X,r}=\bigoplus_{i=1}^rK_X^{(r+1-2i)/2}$
is orthogonal with respect to
any $h\in \Harm(q)$,
and hence
we obtain the decomposition
$h=\bigoplus h_{|K_X^{(r+1-2i)/2}}$.

Note that
$K_X^{(r+1-2i)/2}$
and $K_X^{(r+1-2(r+1-i))/2}=K_X^{(-r-1+2i)/2}$
are mutually dual.
We say that $h\in\Harm(q)$ is real
if
$h_{|K_X^{(r+1-2i)/2}}$
and
$h_{|K_X^{(-r-1+2i)/2}}$
are mutually dual.
Let $\Harm^{\real}(q)$ denote the subset of
 $h\in\Harm(q)$ which are real.

Note that
$\bigl(
K_X^{(r+1-2i)/2}\bigr)^{-1}
\otimes
K_X^{(r+1-2(i+1))/2}$
is naturally isomorphic to the tangent bundle of $X$.
Hence,
\begin{equation}
\label{eq;20.6.30.20}
 h_{|K_X^{(r+1-2i)/2}}^{-1}
 \otimes
 h_{|K_X^{(r+1-2(i+1))/2}}
 \quad (i=1,\ldots,r-1)
\end{equation}
induce K\"ahler metrics of $X$.
If the metrics (\ref{eq;20.6.30.20})
 induce complete distances on $X$,
 $h$ is called complete.

\begin{rem}
According to {\rm\cite{Note0}},
there uniquely exists $h^{\rc}\in\Harm(q)$
which is real and complete.
If $q$ is nowhere vanishing on $X$,
there exists a harmonic metric $h_{\can}$
 as in {\rm\cite{Note0}}
 (see also {\rm\S\ref{subsection;20.4.16.10}} below).
\end{rem}

\subsection{Hermitian metrics on a vector space with a cyclic automorphism}
\label{subsection;20.6.12.1}

Let $r$ be a positive integer.
Let $V$ be an $r$-dimensional $\cnum$-vector space
with a base
$\vece=(e_0,\ldots,e_{r-1})$.
Let $\alpha$ be a non-zero complex number.
Let $f$ be the automorphism of $V$
determined by
$f(e_i)=e_{i+1}$ $(i=0,1,\ldots,r-2)$
and $f(e_{r-1})=\alpha^r e_0$.
We put $G_r:=\{a\in\cnum\,|\,a^r=1\}$.
We consider the action of $G_r$ on $V$
determined by
$a\bullet e_i=a^ie_i$.

Let $\nbigh(f,\vece)$ be the set of
$G_r$-invariant Hermitian metrics $h$ of $V$
such that
$\prod_{i=0}^{r-1}h(e_i,e_i)=1$.
The $G_r$-invariance is equivalent to
the orthogonality $h(e_i,e_j)=0$ $(i\neq j)$.

We put $\tau:=\exp(2\pi\sqrt{-1}/r)\in G_r$.
We set
$v_i:=\sum_{j=0}^{r-1}\tau^{-ij}\alpha^{-j}e_j$
$(i=0,1,\ldots,r-1)$.
Then, $\vecv=(v_0,\ldots,v_{r-1})$ is a base of $V$
such that $f(v_i)=\tau^{i}\alpha v_i$.
Note that
$\tau\bullet v_i=v_{i-1}$ $(i=1,\ldots,r-1)$
and $\tau \bullet v_0=v_{r-1}$.
Let us remark the following well known lemma.
\begin{lem}
\label{lem;20.8.17.10}
 Let $V'\subset V$ be a $G_r$-invariant $\cnum$-subspace
 such that $f(V')\subset V'$.
Then, either $V'=0$ or $V'=V$ hold.
\end{lem}
\pf
Because $f(V')\subset V'$,
there is a subset $I\subset\{0,\ldots,r-1\}$
such that $V'=\bigoplus_{i\in I}\cnum v_i$.
Because $V'$ is $G_r$-invariant,
and because $G_r$ acts on $\{v_0,\ldots,v_{r-1}\}$
in a cyclic way,
we obtain $I=\emptyset$ or $I=\{0,\ldots,r-1\}$.
\hfill\qed

\vspace{.1in}

Let $h\in\nbigh(f,\vece)$.
By the $G_r$-invariance,
$b(h):=h(v_i,v_i)$
is independent of $i$. 
We obtain
\begin{equation}
\label{eq;20.4.22.1}
 h(e_j,e_j)
 =\frac{b(h)}{r}|\alpha|^{2j}
 +\frac{1}{r^2}
 \sum_{i\neq \ell}
  |\alpha|^{2j} \tau^{j(i-\ell)}h(v_i,v_{\ell}).
 \end{equation}

\subsubsection{Canonical metric}
\label{subsection;20.9.25.10}

Let $h_{\can}\in\nbigh(f,\vece)$ be determined by
  \begin{equation}
   \label{eq;20.4.15.1}
    h_{\can}(e_j,e_j)=|\alpha|^{-(r-1)+2j}
    \quad (j=0,\ldots,r-1).
  \end{equation}

\begin{lem}
\label{lem;20.4.15.2}
$h_{\can}$ is the unique Hermitian metric
contained in
$\nbigh(f,\vece)$
such that the base $\vecv$ is orthogonal
with respect to $h_{\can}$.
\end{lem}
\pf
We can check the orthogonality of $\vecv$
by a direct computation.
Suppose that $\vecv$
is orthogonal with respect to
$h_1\in\nbigh(f,\vece)$.
By (\ref{eq;20.4.22.1}),
we obtain
\[
h_1(e_j,e_j)=\frac{b(h_1)}{r}|\alpha|^{2j}.
\]
By the condition
$1=\prod_{j=0}^{r-1}h_1(e_j,e_j)$,
we obtain $b(h_1)=r|\alpha|^{-(r-1)}$,
and hence $h_1=h_{\can}$.
\hfill\qed

\subsubsection{$\epsilon$-orthogonality}
 
Suppose $r>1$,
and let $0\leq \epsilon\leq 1$.
Suppose that the following holds
for a metric $h\in\nbigh(f,\vece)$:
\[
  |h(v_i,v_j)|\leq
  \epsilon b(h)
  =\epsilon\cdot|v_i|_h\cdot |v_j|_h
  \quad(i\neq j).
\]
The second term in the right hand side of (\ref{eq;20.4.22.1})
is dominated by
\[
\frac{b(h)|\alpha|^{2j}}{r}(r-1)\epsilon.
\]
Hence, we obtain the following description
\[
 h(e_j,e_j)=\frac{b(h)|\alpha|^{2j}}{r}(1+c_j),
\]
where
$|c_j|\leq (r-1)\epsilon$ and $1+c_j>0$.
By the condition
$\prod h(e_j,e_j)=1$,
we obtain
\[
 1=\left(\frac{b(h)}{r}\right)^r|\alpha|^{r(r-1)}
  \prod_{j=0}^{r-1}(1+c_j).
\]
In particular,
we obtain
$1=\left(\frac{b(h_{\can})}{r}\right)^r|\alpha|^{r(r-1)}$.
Hence, we obtain
\[
 b(h_{\can})=
 b(h)\cdot
   \prod_{j=0}^{r-1}(1+c_j)^{1/r}.
\]
For $0<\delta<1$,
we set
$C_{\delta}:=(1-\delta)^{-1}(r-1)$.
If  $(r-1)\epsilon<\delta$,
we obtain
$|\log (1+c_j)|\leq C_{\delta}\epsilon$.
For such $\epsilon$,
we obtain
\[
\left|
 \log\bigl(
  b(h)/b(h_{\can})
  \bigr)
  \right|
=\left|
  \frac{1}{r}\sum_{j=0}^{r-1}
 \log(1+c_j)
\right|
 \leq C_{\delta}\epsilon.
\]
We also obtain
\begin{equation}
 \label{eq;20.4.16.4}
\left|
 \log\Bigl(
  h(e_j,e_j)/h_{\can}(e_j,e_j)
  \Bigr)
\right|
=\Bigl|
\log\bigl(b(h)/b(h_{\can})\bigr)
+\log(1+c_j)
 \Bigr|
 \leq 2C_{\delta}\epsilon.
\end{equation}

\subsubsection{Norm of the automorphism}

For any $h\in\nbigh(f,\vece)$,
we obtain $|f|_h\geq \sqrt{r}|\alpha|$.
Let
$C$ be any positive number such that
$C\geq \sqrt{r}|\alpha|$.

\begin{lem}
\label{lem;20.4.21.1}
 For any $h\in \nbigh(f,\vece)$
 such that $|f|_h\leq C$,
 we obtain
 \[
 |e_i|_h\leq |\alpha|^{-r}C^{(r+1)/2+i}
 \quad (i=0,\ldots,r-2),
 \quad
 |e_{r-1}|_h\leq C^{(r-1)/2}.
 \]
\end{lem}
\pf
Because $e_{i+1}=f(e_i)$,
we obtain $|e_{i+1}|_h\leq C|e_i|_h$.
We obtain
$|e_{r-1}|_h\leq C^{r-1-i}|e_i|_h$.
Because $\prod |e_i|_h=1$,
we obtain
\[
 |e_{r-1}|_h^rC^{-r(r-1)/2}\leq 1,
\quad
 \mbox{\rm i.e., }
 |e_{r-1}|_h\leq C^{(r-1)/2}.
\]
Because $f(e_{r-1})=\alpha^r e_0$,
we obtain
$|e_0|_h\leq C|\alpha|^{-r}|e_{r-1}|_h$,
and hence
$|e_0|_h\leq |\alpha|^{-r}C^{(r+1)/2}$.
Then, we obtain
$|e_i|_h\leq
|\alpha|^{-r}C^{(r+1)/2+i}$.
\hfill\qed

\vspace{.1in}
\begin{cor}
\label{cor;20.4.21.2}
 We obtain
 \[
 |\alpha|^r(C+1)^{-2r}\leq |e_i|_h
  \leq |\alpha|^{-r}(C+1)^{2r}.
 \]
\end{cor}
\pf
We obtain the second inequality
from Lemma \ref{lem;20.4.21.1}.
By using the duality,
we obtain the first inequality.

\hfill\qed

\begin{cor}
 \label{cor;20.4.17.10}
 For any $h_i\in \nbigh(f,\vece)$ $(i=1,2)$
such that $|f|_{h_i}\leq C$,
we obtain
\[
 |\alpha|^{2r}(C+1)^{-4r}
 \leq 
 \frac{|e_i|_{h_1}}{|e_i|_{h_2}}
 \leq
 |\alpha|^{-2r}(C+1)^{4r}.
\]
\hfill\qed
\end{cor}

\subsection{Harmonic metrics of a cyclic Higgs bundle
  on a disc}
\label{subsection;20.4.16.10}

We shall recall variants of Simpson's main estimate
for harmonic bundles,
which was pioneered in \cite{s2},
and further developed in \cite{Mochizuki-wild}
and \cite{Decouple}.

We set $U(R):=\bigl\{z\in\cnum\,\big|\,|z|<R\bigr\}$.
We set
$E:=\bigoplus_{i=0}^{r-1}\nbigo_{U(R)}\,e_i$.
Let $\beta$ be a holomorphic function on $U(R)$.
Let $f$ be the endomorphism of $E$
determined by
$f(e_i)=e_{i+1}$ $(i=0,\ldots,r-2)$
and $f(e_{r-1})=\beta e_0$.
We set $\theta:=f\,dz$.
We set $G_r:=\{a\in\cnum\,|\,a^r=1\}$
as in \S\ref{subsection;20.6.12.1}.
We define the $G_r$-action on $E$
by $a\bullet e_i=a^ie_i$.
Note that $a^{\ast}(\theta)=a^{-1}\theta$.
Let $\Harm(\beta)$ denote the set of
$G_r$-invariant harmonic metrics of
the Higgs bundle $(E,\theta)$ such that
$\prod_{i=0}^{r-1}h(e_i,e_i)=1$.

\begin{prop}[\mbox{\cite[Proposition 2.1]{Decouple}}]
\label{prop;20.4.15.11}
 Fix $0<R_1<R$.
 Then, there exist $C_1,C_2>0$
 depending only on $r$, $R_1$, and $R$
 such that
\[
 \sup_{U(R_1)}|f|_h\leq C_1+C_2\max_{U(R)}|\beta|^{1/r}
\]
 for any $h\in\Harm(\beta)$.
\hfill\qed
\end{prop}

 \begin{cor}
  \label{cor;20.6.11.1}
  Let $h_0$ be a $G_r$-invariant Hermitian metric
  of $E$ such that $\prod_{i=0}^{r-1} h_0(e_i,e_i)=1$.
  For $h\in\Harm(\beta)$,
  we obtain the automorphism $s$ of $E$
  which is self-adjoint with respect to both $h$ and $h_0$,
  determined by $h=h_0s$.
  Let $C_1$ and $C_2$ be constants
  as in 
  Proposition {\rm\ref{prop;20.4.15.11}}.
  Let $C_{3}>0$ be a positive constant such that
  $|f|_{h_0}\leq C_{3}$.  
  We set
  \[
  A_1:=
  \max\Bigl\{
  C_{3},
  C_1+C_2\max_{z\in U(R)}|\beta(z)|^{1/r}
  \Bigr\}.
  \]
Then, the following holds
at $z\in U(R_1)$
such that $\beta(z)\neq 0$:
\[
  \bigl|
  \log |s|_{h_0}(z)
  \bigr|
  \leq
  \log\Bigl(
  \sqrt{r}|\beta(z)|^{-2}(A_1+1)^{4r}
  \Bigr).
 \]
\end{cor}
\pf
It follows from
Corollary \ref{cor;20.4.17.10} and
Proposition \ref{prop;20.4.15.11}.
\hfill\qed

\vspace{.1in}

Suppose moreover that $\beta$ is nowhere vanishing on $U(R)$.
We fix an $r$-th root $\alpha$ of $\beta$,
i.e., $\alpha^r=\beta$.
We put $\tau:=\exp(2\pi\sqrt{-1}/r)$.
We set
$v_i:=\sum_{j=0}^{r-1}\tau^{-ij}\alpha^{-j}e_j$
$(i=0,1,\ldots,r-1)$.
We obtain the decomposition
$E=\bigoplus_{i=0}^{r-1}\nbigo_{U(R)}v_i$
and $f(v_i)=\tau^i\alpha v_i$.

\begin{prop}[\mbox{\cite[Corollary 2.6]{Decouple}}]
\label{prop;20.4.15.10}
 Assume that $|\alpha(0)-\alpha(z)|\leq \frac{1}{100}|\alpha(0)|$
 on $U(R)$.
 Fix $0<R_2<R_1$.
 Then, there exist $C_{10}>0$ and $\epsilon_{10}>0$
 depending only on $r$, $R$, $R_1$  and $R_2$
 such that the following holds
 on $U(R_2)$
 for any $h\in \Harm(\beta)$
 and for $i\neq j$:
 \[
 |h(v_i,v_j)|
 \leq C_{10}\cdot
 \exp\Bigl(-\epsilon_{10}|\alpha(0)|\Bigr)\cdot
  |v_i|_h\cdot |v_j|_h.
 \]
 Note that
 $|v_i|_h=|v_j|_h$
 by the $G_r$-invariance of $h$.
 \hfill\qed
\end{prop}

Let $h_{\can}$ be the $G_r$-invariant Hermitian metric determined by
\[
h_{\can}(e_j,e_j)=|\alpha|^{-(r-1)+2j}.
\]
Then, the frame $(v_0,\ldots,v_{r-1})$ is orthogonal,
and we have
\[
h_{\can}(v_i,v_i)=r|\alpha|^{-(r-1)}.
\]
Hence,
the curvature of the Chern connection of $h_{\can}$ is $0$.
Because $[\theta,\theta^{\dagger}_{h_{\can}}]=0$,
$h_{\can}$ is a harmonic metric of the Higgs bundle $(E,\theta)$.

 \begin{cor}
   \label{cor;20.4.16.11}
  Suppose that
  the assumption of Proposition {\rm\ref{prop;20.4.15.10}}
  is satisfied.
  Moreover, we assume that
  \[
(r-1)C_{10}\cdot
 \exp\Bigl(-\epsilon_{10}|\alpha(0)|\Bigr)\leq 1/2.
  \]
There exist $\epsilon_i>0$ and $C_i>0$ $(i=11,12,13)$
  depending only on $r$, $R$, $R_1$ and $R_2$ such that
  the following holds for any $h\in\Harm(\beta)$
  on $U(R_2)$:
\begin{equation}
\label{eq;20.4.16.1}
  \Bigl|
 \log\Bigl(
 h(e_j,e_j)/h_{\can}(e_j,e_j)
  \Bigr)
  \Bigr|
  \leq
  C_{11}
  \exp\bigl(-\epsilon_{11}|\alpha(0)|\bigr). 
\end{equation}
\begin{equation}
\label{eq;20.4.16.2}
  \Bigl|
 \del_z\log\Bigl(
 h(e_j,e_j)/h_{\can}(e_j,e_j)
  \Bigr)
  \Bigr|
  \leq
  C_{12}
  \exp\bigl(-\epsilon_{12}|\alpha(0)|\bigr). 
\end{equation}
\begin{equation}
\label{eq;20.4.16.3}
  \Bigl|
 \del_{\zbar}\del_z\log\Bigl(
 h(e_j,e_j)/h_{\can}(e_j,e_j)
  \Bigr)
  \Bigr|
  \leq
  C_{13}
  \exp\bigl(-\epsilon_{13}|\alpha(0)|\bigr). 
\end{equation}
 \end{cor}
\pf
We obtain (\ref{eq;20.4.16.1})
from (\ref{eq;20.4.16.4}) and Proposition \ref{prop;20.4.15.10}.
Note that the Chern connection of the curvature of $h_{\can}$ is $0$,
and that the decomposition $E=\bigoplus \nbigo e_i$
is orthogonal with respect to $h$ and $h_{\can}$.
Then, we obtain (\ref{eq;20.4.16.3})
from Proposition \ref{prop;20.4.15.11} and Proposition \ref{prop;20.4.15.10}.
We obtain (\ref{eq;20.4.16.2})
from (\ref{eq;20.4.16.1}) and (\ref{eq;20.4.16.3})
by using the elliptic regularity.
\hfill\qed

\subsubsection{Estimates near boundary}

Let $Y$ be an open subset in $U(R)$.
Let $\Ybar$ denote the closure of $Y$ in $\Ubar(R)$.
Let $h\in\Harm(\beta_{|Y})$
such that $h$ extends to a continuous metric on $\Ybar$.

\begin{prop}
\label{prop;20.6.13.32}
Fix $0<R_{10}<R$.
There exists $C_{i}$ $(i=1,2,3)$,
depending only on $R_{10}$, $R$ and $\rank(E)$,
such that the following holds on $U(R_{10})\cap Y$:
\[
 |f|_h\leq
 C_1+C_2\max_{z\in U(R)}|\beta(z)|^{1/r}
 +C_3\max_{U(R)\cap \del Y}|f|_h.
\]
 \end{prop}
\pf
It follows from Proposition \ref{prop;20.6.13.40} below.
\hfill\qed

 \begin{prop}
  Let $h_0$ be a $G_r$-invariant Hermitian metric of $E$
  such that $\prod_{i=0}^{r-1}h_0(e_i,e_i)=1$.
  For $h\in\Harm(\beta)$,
  let $s(h)$ be the automorphism of $E_{|Y}$
  which is self-adjoint with respect to $h$ and $h_0$,
  determined by $h=h_0\cdot s(h)$.
  Let $C_i$ $(i=1,2,3)$ be
  as in Proposition {\rm\ref{prop;20.6.13.32}},
  and let $C_4$ be a positive constant
  such that $|f|_{h_0}<C_4$.
  We set
 \[
  A_1:=
  \max\Bigl\{
  C_4,\,
  C_1+C_2\max_{z\in U(R)}|\beta(z)|^{1/r}
  +C_3\max_{U(R)\cap\del Y}|f|_h
  \Bigr\}.
 \]
  Then, the following holds
  at any $z\in U(R_{10})\cap Y$
  such that $\beta(z)\neq 0$:
\[
  \bigl|
  \log |s(h)|_{h_0}(z)
  \bigr|
  \leq
  \log\Bigl(
  \sqrt{r}
  |\beta(z)|^{-2}
  (A_1+1)^{4r}
  \Bigr).
\]
\end{prop}
\pf
It follows from 
Proposition \ref{prop;20.6.13.32}
and Corollary \ref{cor;20.4.17.10}.
\hfill\qed
 
\subsubsection{Appendix: A variant of Simpson's main estimate near boundary}

Let $(E,\delbar_E,\theta)$ be a Higgs bundle
defined on a neighbourhood of $\Ubar(R)$ in $\cnum$.
Let $Y$ be an open subset in $U(R)$.
Let $\Ybar$ denote the closure of $Y$ in $\Ubar(R)$.
Let $h$ be a harmonic metric of
$(E,\delbar_E,\theta)_{|Y}$
which extends to a continuous metric of
$E_{|\Ybar}$.

Let $f$ be the endomorphism of $E$
determined by $\theta=f\,dz$.
Let $M$ be a positive constant
such that
the following holds.
\begin{itemize}
 \item Let $\alpha$ be any eigenvalue of $f_{|P}$
       for any $P\in \Ubar(R)$.
       Then, $|\alpha|\leq M$.
\end{itemize}

\begin{prop}
\label{prop;20.6.13.40}
Let $R_1<R$.
There exist positive constant
$C_i$ $(i=1,2,3)$,
depending only on $R$, $R_1$ and $\rank(E)$,
such that
the following holds on
$Y\cap U(R_1)$
\[
 |f|_h\leq
 C_1+C_2M+C_3\max_{U(R)\cap \del Y}|f|_h.
\]
\end{prop}
\pf
Recall the following inequality on $Y$:
\[
 -\del_z\del_{\zbar}\log|f|^2_h
 \leq
 -\frac{\bigl|[f,f^{\dagger}_h]\bigr|^2}{|f|_h^2}.
\]
For any $P\in \Ubar(R)$,
let $\alpha_1,\ldots,\alpha_{\rank E}$
be the eigenvalues of $f_{|P}$,
and we set
$g(P):=\sum |\alpha_i|^2$.
There exists $C_{10}>0$ such that
the following holds for any $P\in Y$:
\[
 |[f_{|P},f_{|P}^{\dagger}]|_h
 \geq
 C_{10}(|f_{|P}|_h^2-g(P)).
\]
We obtain the following on $Y$:
\[
  -\del_z\del_{\zbar}\log|f|^2_h(P)
 \leq
 -C_{10}^2
 \frac{(|f_{|P}|_h^2-g(P))^2}{|f_{|P}|_h^2}.
\]
If $|f_{|P}|^2\geq 2g(P)$,
then
we obtain
\begin{equation}
\label{eq;20.6.13.30}
 -\del_z\del_{\zbar}\log|f|^2_h(P)
 \leq
 -\frac{C_{10}^2}{4}
 |f_{|P}|^2_h.
\end{equation}

We consider the following function on $U(R)$:
\[
 A(z):=\frac{2R^2}{(R^2-|z|^2)^2}
 \geq 2R^{-2}.
\]
Note that
$-\del_z\del_{\zbar}\log A
 =-A$.
Let $C_{11}$ be a constant such that
$C_{11}>(C_{10}^2/4)^{-1}$.
We set
\[
 C_{12}:=C_{11}+R^2\sup_{U(R)\cap \del Y}|f|_h^2+10\rank(E)R^2M^2.
\]
Note that the following holds on $U(R)$:
\[
 -\del_z\del_{\zbar}\log (C_{12}A)
=-A
=-C_{12}^{-1}(C_{12}A)
\geq -\frac{C_{10}^2}{4}(C_{12}A).
\]

Let $Z$ be the set of $P\in Y$
such that $|f_{|P}|_h^2>C_{12}A(P)$.
We assume that $Z$ is non-empty,
and we shall derive a contradiction.
Note that
$|f|_h^2<C_{12}A$ on $\del Y\cap U(R)$.
We also note that
$|f|_h$ is bounded
around any point of $\del Y\cap \del U(R)$,
and that $A(z)\to\infty$ as $|z|\to R$.
Hence,
$Z$ is a relatively compact open subset of $Y$,
which implies $|f|_h^2=C_{12}A$ on $\del Z$.
For $P\in Z$,
we obtain
$|f_{|P}|_h^2\geq 2g(P)$,
and hence (\ref{eq;20.6.13.30}) holds at $P$.
We obtain the following at $P\in Z$:
\begin{equation}
\label{eq;20.6.13.31}
 -\del_z\del_{\zbar}
 \Bigl(
\log|f|_h^2
-\log(C_{12}A)
 \Bigr)(P)
\leq
 -\frac{C_{10}^2}{4}(|f_{|P}|^2-C_{12}A(P))
 <0. 
\end{equation}
Together with $|f|_h^2=C_{12}A$ on $\del Z$,
we obtain $|f|_h^2\leq C_{12}A$ on $Z$,
which contradicts the construction of $Z$.
Hence, we obtain that $Z=\emptyset$.
\hfill\qed

\subsection{Some existence results
  of harmonic metrics on cyclic Higgs bundles}

Let $X$ be a Riemann surface.
Let $E=\bigoplus_{i=0}^{r-1} E_i$
be a graded holomorphic vector bundle on $X$.
Let $\theta$ be a Higgs field of $E$
such that
$\theta(E_i)\subset E_{i+1}\otimes K_X$
$(i=0,\ldots,r-2)$
and $\theta(E_{r-1})\subset E_{0}\otimes K_X$.
We assume that $\det(\theta)$ is not constantly $0$.
Assume that there exists a flat metric
$h_{\det(E)}$ of $\det(E)$.
We consider the $G_r$-action on $E_i$
given by $a\bullet v_i=a^{i}v_i$,
which induces a $G_r$-action on $E$.
For any open subset $Y$,
let $\Harm^{G_r}((E,\theta)_{|Y},h_{\det(E)})$
denote the set of
$G_r$-invariant harmonic metrics $h$
of $(E,\theta)_{|Y}$
such that $\det(h)=h_{\det(E)|Y}$.

\begin{rem}
 For any $r$-differential $q$ on $X$,
 we obtain
 $\Harm(q)=
 \Harm^{G_r}\bigl(
 (\hyperk_{X,r},\theta(q)),1
 \bigr)$
 by definition.
\end{rem}

\subsubsection{Convergence}
 
 Let $X_1\subset X_2\subset\cdots$
 be a smooth exhaustive family of $X$.
  Suppose that $h_i\in \Harm^{G_r}((E,\theta)_{|X_i},h_{\det(E)})$.
Let $h_0$ be any $G_r$-invariant Hermitian metric of $E$
such that $\det(h_0)=h_{\det(E)}$.
 Let $s_i$ be the automorphism of $E_{|X_i}$
 determined by $h_i=h_{0|X_i}s_i$.

\begin{prop}
\label{prop;20.6.15.30}
There exists a locally bounded function $f$ on $X$
 such that
 $|s_i|_{h_{0,i}}+|s_i^{-1}|_{h_{0,i}}\leq f_{|X_i}$
 for any $i$.
 As a result,
 there exists a convergent subsequence of
 the sequence $\{s_i\}$ with the limit $s_{\infty}$,
 and we obtain $h_{\infty}=h_0s_{\infty}\in
  \Harm^{G_r}((E,\theta),h_{\det(E)})$.
\end{prop}
\pf
Let $Z$ denote the zero set of $\det(\theta)$.
It is discrete in $X$.
By Corollary \ref{cor;20.6.11.1},
for any $P\in X\setminus Z$,
there exists a relatively compact neighbourhood $X_P$
of $P$ in $X\setminus Z$
with a constant $C_P>1$ such that
$C_P^{-1}
\leq
|s_i|_{h_0}\leq C_P$ for any large $i$.
By using the subharmonicity of $\Tr(s_i)$
(see \S\ref{subsection;20.8.16.10}),
we also obtain a similar estimate
locally around any point of $Z$.
Then, the claim follows
from Proposition \ref{prop;20.5.29.1}.
\hfill\qed

\begin{cor}
\label{cor;20.7.1.1}
 Suppose that there exists a subset $Y\subset X$
 and an $\real_{>0}$-valued
 locally bounded function $f_Y$ on $Y$
 such that
 $|s_{i|Y\cap X_i}|_{h_0}\leq f_{Y|Y\cap X_i}$.
 Then, we obtain $|s_{\infty|Y}|_{h_0}\leq f_Y$
 for $s_{\infty}$ in
 Proposition {\rm\ref{prop;20.6.15.30}}.
 \hfill\qed
\end{cor}

\subsubsection{Control of growth order}

Let $g_X$ be any K\"ahler metric of $X$.
Let $h_0$ be a $G_r$-invariant Hermitian metric of $E$
such that $\det(h_0)=h_{\det(E)}$.
Let $\{X_i\}$ be a smooth exhaustive family of $X$.
Let $h_i\in\Harm^{G_r}((E,\theta)_{|X_i},h_{\det(E)})$
such that $h_{i|\del X_i}=h_{0|\del X_i}$.

Let $Y$ be an open subset of $X$.
Suppose the following.
\begin{itemize}
 \item There exists an $\real_{\geq 0}$-valued function $A_0$
       on a neighbourhood $N$ of $\Ybar$
       such that the following holds on $N$:
\[
 \sqrt{-1}\Lambda\delbar\del A_0\geq |F(h_0)|_{h_0,g_X}.
\]
 \item There exists an $\real_{\geq 0}$-valued
       harmonic function $A_1$ on $N$
       such that
       $\log\Tr(h_i\cdot h_0^{-1})\leq A_1+\log r$
       on $X_i\cap\del Y$.
\end{itemize}
\begin{prop}
 There exists $h\in\Harm^{G_r}((E,\theta),h_{\det(E)})$
 such that
 $\log\Tr(h\cdot h_0^{-1})\leq A_0+A_1+\log r$ on $Y$.
\end{prop}
\pf
By Proposition \ref{prop;20.6.15.30},
we may assume that $\{h_i\}$ is convergent to
$h_{\infty}$.
Let $s_i$ be the automorphism of
$E_{|X_i}$ obtained as
$s_i=h_i\cdot h_{0|X_i}^{-1}$.
By our choice,
we obtain
\[
 \sqrt{-1}\Lambda\delbar\del
 \bigl(
  \log\Tr(s_i)
  -A_0-A_1-\log r
  \bigr)
  \leq 0.
\]
On $Y\cap\del X_i$,
we obtain
\[
 \log\Tr(s_i)-A_0-A_1-\log r
 \leq
  \log\Tr(s_i)-\log r= 0.
\]
On $X_i\cap \del Y$, we obtain
\[
 \log\Tr(s_i)-A_0-A_1-\log r
 \leq
 \log\Tr(s_i)-A_1-\log r\leq 0.
\]
Hence, we obtain
$\log\Tr(s_i)-A_0-A_1-\log r\leq 0$
on $Y\cap X_i$.
We obtain
$\log\Tr(h_{\infty}h_0^{-1})-A_0-A_1-\log r\leq 0$
on $Y$
as in Corollary \ref{cor;20.7.1.1}.
\hfill\qed

\subsubsection{Compact subsets}

\begin{prop}
\label{prop;20.6.29.11}
Let $K\subset X$ be a compact subset.
For any
$h_{1}\in\Harm^{G_r}((E,\theta)_{|X\setminus K},h_{\det(E)})$,
there exists $h\in\Harm^{G_r}((E,\theta),h_{\det(E)})$
such that $h_{|X\setminus N}$ and $h_{1|X\setminus N}$
are mutually bounded
for any relatively compact neighbourhood $N$ of $K$ in $X$.
\end{prop}
\pf
Let $h_2$ be a $G_r$-invariant Hermitian metric of $E$
such that
$\det(h_2)=h_{\det(E)}$
and that
$h_{2|X\setminus N}=h_{1|X\setminus N}$
for a relatively compact neighbourhood $N$ of $K$.
Note that $(E,\theta,h_2)$
is analytically stable with respect to the $G_r$-action
in the sense of \S\ref{subsection;20.7.1.2}
because there is no proper $G_r$-invariant Higgs subbundle
by Lemma \ref{lem;20.8.17.10}.
(Recall that $\det(\theta)$ is not constantly $0$.)
Then, the claim follows from
Proposition \ref{prop;20.6.11.10}.
\hfill\qed

\subsection{Asymptotic behaviour around poles}

Let $U$ be a neighbourhood of $0$ in $\cnum$.
We set $U^{\circ}:=U\setminus\{0\}$.
Let $q$ be a nowhere vanishing holomorphic $r$-differential on $U^{\circ}$
which is meromorphic at $0$.
We have the expression
$q=z^m\,\alpha\,(dz/z)^r$,
where $\alpha$ induces a nowhere vanishing holomorphic
function on $U$.

\subsubsection{The case $m\leq 0$}
\label{subsection;20.9.24.10}

For a positive number $\epsilon>0$,
we set
\[
 \rho_{\epsilon}(z):=
 \left\{
\begin{array}{ll}
 \exp(-\epsilon|z|^{m/r})& (m<0)\\
 |z|^{\epsilon} & (m=0).
\end{array}
 \right.
\]
There exists a relatively compact
neighbourhood $U_1$ of $0$ in $U$
such that 
$q$ is nowhere vanishing on $U_1^{\circ}:=U_1\cap U^{\circ}$.
Recall that there exists
the canonical harmonic metric $h_{\can}\in\Harm(q_{|U_1^{\circ}})$.

\begin{prop}
\label{prop;20.7.1.25}
 Any $h\in\Harm(q)$ is mutually bounded
 with $h_{\can}$ on $U_1^{\circ}$.
 More strongly,
 there exist $C>0$ and $\epsilon>0$
 such that the following holds on $U_1^{\circ}$:
\[
 \left|
 \log |(dz)^{(r+1-2i)/2}|_{h}
-\log |(dz)^{(r+1-2i)/2}|_{h_{\can}}
 \right|
 \leq
  C\rho_{\epsilon}.
\]
\end{prop}
\pf
We may assume that
$\{|z|\leq 1\}\subset U_1$.
Let $\Psi:\cnum\lrarr \cnum^{\ast}$
be the map given by
$\Psi(\zeta)=e^{\sqrt{-1}\zeta}$.
We have the expression
\[
\Psi^{\ast}(q)=
e^{\sqrt{-1}m\zeta}
\Psi^{\ast}(\alpha)(\sqrt{-1}d\zeta)^r.
\]
There exists
$C_1>1$ such that
$C_1^{-1}\leq |\Psi^{\ast}(\alpha)|\leq C_1$
on $\{\Image(\zeta)\geq 0\}$.

Let us consider the case $m=0$.
For $\zeta_0$ with $y_0=\Image(\zeta_0)>0$,
let $F_{\zeta_0}:\{|w|<1\}\lrarr \{\Image(\zeta)>0\}$
given by
$F_{\zeta_0}(w)=y_0(w+\sqrt{-1})$.
We set
$\htilde:=(\Psi\circ F_{\zeta_0})^{\ast}(h)$
and
$\htilde_{\can}:=(\Psi\circ F_{\zeta_0})^{\ast}(h_{\can})$.
By Proposition \ref{prop;20.4.15.11}
and Corollary \ref{cor;20.4.16.11},
there exists $C_2>0$ and $\epsilon_2>0$
such that
\[
 \Bigl|
  \log|(dw)^{(r+1-2i)/2}|_{\htilde}
- \log|(dw)^{(r+1-2i)/2}|_{\htilde_{\can}}
  \Bigr|
\leq C_2\exp(-\epsilon_2 y_0).
\]
It implies the desired estimate
in the case $m=0$.

Let us study the case $m<0$.
Take a large $T>0$.
For $\zeta_0$ with $y_0=\Image(\zeta_0)>0$,
we consider
$F_{\zeta_0}:\{|w|<1\}\lrarr \{\Image(\zeta)>0\}$
determined by
$F_{\zeta_0}(w)=T(w+\zeta_0)$.
By applying Proposition \ref{prop;20.4.15.11}
and Corollary \ref{cor;20.4.16.11}
to the pull back by $\Psi\circ F_{\zeta_0}$,
we obtain the desired estimate in the case $m<0$.
\hfill\qed

\begin{rem}
 Proposition {\rm\ref{prop;20.7.1.25}}
 can be also obtained 
 as a consequence of
 the classification of $G_r$-equivariant
 good filtered Higgs bundles
 over a cyclic Higgs bundle
 on a punctured disc.
 See {\rm\cite[\S3.1.3]{Toda-lattice}}.
 \end{rem}

\subsubsection{The case $m>0$}
\label{subsection;20.7.2.1}

Let $\nbigp_{r,m}$ denote the set of
$\vecb=(b_i)\in\real^r$ satisfying
\[
 b_1\geq b_2\geq \cdots\geq b_r\geq b_1-m,
 \quad
 \sum b_i=-\frac{r(r+1)}{2}.
\]
Let $\nbigp^{\real}_{r,m}\subset\nbigp_{r,m}$
denote the subset of
$\vecb\in\nbigp_{r,m}$ such that
$b_i+b_{r+1-i}=-(r+1)$.
We obtain the following proposition from \cite{s2}
(see also \cite{Toda-lattice}).

\begin{prop}\mbox{{}}
\label{prop;20.7.1.20}
 \begin{itemize}
   \item 
 For any $h\in\Harm(q)$,
 there exist
 $\vecb(h)=(b_{i}(h))
 \in\nbigp_{r,m}$
 determined by the following condition
 around $z=0$:
 \begin{multline}
\label{eq;20.7.13.10}
 \log|(dz)^{(r+1-2i)/2}|_h
  +(b_{i}(h)+i)\log|z|
  \\
  = 
  O\Bigl(
 \log\bigl(
 -\log |z|
 \bigr)
 \Bigr)
 \quad
  (i=1,\ldots,r).
 \end{multline}
If $h\in\Harm^{\real}(q)$,
 then $\vecb(h)\in\nbigp^{\real}_{r,m}$.
 \item
 Let $h_1,h_2\in\Harm(q)$ such that
 $\vecb(h_1)=\vecb(h_2)$.
 Then, for any relatively compact neighbourhood
 $U'$ of $0$ in $U$,
 the restrictions
 $h_{1|U'\setminus\{0\}}$
and 
 $h_{2|U'\setminus\{0\}}$
 are mutually bounded.
 \end{itemize}
\end{prop}
\pf
Let us explain an outline of the proof
(see also \cite[\S3.1.2, \S3.2.8]{Toda-lattice}).
Let $h\in\Harm(q)$,
i.e.,
$h$ is a $G_r$-invariant harmonic metric
of $(\hyperk_{U^{\circ},r},\theta(q))$
such that $\det(h)=1$.
Let $f$ be the endomorphism of
$\hyperk_{U^{\circ},r}$
determined by
$\theta(q)=f\,dz/z$.
Because $m>0$, the harmonic bundle
$(\hyperk_{U^{\circ},r},\theta(q),h)$ is tame.
For any $a\in\real$,
and for any open subset $\nbigu\ni 0$,
let $\nbigp^h_a\hyperk_{U^{\circ},r}(\nbigu)$
denote the space of sections $s$ of
$\nbigu\setminus\{0\}$
such that
$|s|_h=O\bigl(|z|^{-a-\epsilon}\bigr)$
for any $\epsilon>0$.
For any open subset $\nbigu\not\ni 0$,
let
$\nbigp^h_a\hyperk_{U^{\circ},r}(\nbigu)$
denote the space of holomorphic sections of
$\hyperk_{U^{\circ},r}(\nbigu)$.
Thus, we obtain $\nbigo_{U}$-modules
$\nbigp^h_a\hyperk_{U^{\circ},r}$ $(a\in\real)$.
According to \cite[Theorem 2]{s2},
$\nbigp^h_a\hyperk_{U^{\circ},r}$
are locally free $\nbigo_U$-modules,
and
$f(\nbigp^h_{a}\hyperk_{U^{\circ},r})
\subset\nbigp^h_a\hyperk_{U^{\circ},r}$ for any $a$.
Because $h$ is $G_r$-invariant,
$\nbigp^h_a\hyperk_{U^{\circ},r}$ is naturally
$G_r$-equivariant.
Hence, we obtain the decomposition
\begin{equation}
\label{eq;20.7.13.11}
 \nbigp^h_a\hyperk_{U^{\circ},r}
 =\bigoplus_{i=1}^r
 \nbigp^h_aK_{U^{\circ}}^{(r+1-2i)/2}.
\end{equation}
Here, the locally free $\nbigo_U$-modules
$\nbigp^h_aK_{U^{\circ}}^{(r+1-2i)/2}$
are obtained from
$K_{U^{\circ}}^{(r+1-2i)/2}$ with $h$
as above.

Because $m>0$,
the eigenvalues of $f_{|z}$ goes to $0$ as $|z|\to 0$.
Hence, according to Simpson's main estimate
\cite[Theorem 1]{s2},
we obtain $|f|_h=O\bigl((-\log|z|)^{-1}\bigr)$.

We set $v_i:=z^{i}(dz)^{(r+1-2i)/2}$.
Then, we obtain $f(v_i)=v_{i+1}$ $(i=1,\ldots,r-1)$,
and hence
$|v_{i+1}|_h\leq |v_{i}|_h(-\log|z|)^{-1}$.
Because $\det(h)=1$,
we obtain
\[
|v_r|^{r}\cdot (-\log|z|)^{r(r-1)/2}
\leq
\prod_{i=1}^h |v_i|_h=|z|^{r(r+1)/2}.
\]
Hence, there exists $a\in\real$
such that $v_r$ is a section of
$\nbigp^h_{(r+1)/2}\hyperk_{U^{\circ},r}$.
Because $f(v_r)=z^m\alpha v_1$,
we obtain
$|z|^m\cdot |\alpha|\cdot |v_1|_h\leq |v_r|_h(-\log |z|)^{-1}$.
Therefore,
there exists $b\in\real$
such that
$v_i\in\nbigp^h_b\hyperk_{U^{\circ},r}$
$(i=1,\ldots,r)$.
We set
\[
 b_i(h):=\min
  \bigl\{b\in\real\,\big|\,
   v_i\in\nbigp^h_{b}\hyperk_{U^{\circ},r}\bigr\}.
\]
Then,
by the norm estimate \cite[\S7]{s2}
(Proposition \ref{prop;20.7.1.26} below),
we obtain (\ref{eq;20.7.13.10}).
Because $f(v_i)=v_{i+1}$
and $f(v_{r})=z^m\alpha v_1$,
we obtain $b_i\geq b_{i+1}$
and $b_r\geq b_{1}-m$.
We also remark
the filtered bundle
$\nbigp^h_{\ast}\hyperk_{U^{\circ},r}$
is uniquely determined by
the numbers
$\vecb(h)=(b_i(h))$
because of the decomposition (\ref{eq;20.7.13.11}).

Let $h_i\in\Harm(q)$ $(i=1,2)$
such that $\vecb(h_1)=\vecb(h_2)$.
We obtain
$\nbigp^{h_1}_{\ast}\hyperk_{U^{\circ},r}
=\nbigp^{h_2}_{\ast}\hyperk_{U^{\circ},r}$.
Hence,
according to \cite[Corollary 4.3]{s2},
$h_1$ and $h_2$ are mutually bounded.
\hfill\qed

\vspace{.1in}

Let us refine the estimate (\ref{eq;20.7.13.10}).
For $\vecb\in\nbigp_{r,m}$,
we introduce a non-negative integer $\ell$
and integers
$\nu_0$, $\nu_1,\ldots,\nu_{\ell},\nu_{\ell+1}$
as follows.
If $b_1=b_r$, we set
$\ell=0$, $\nu_0=0$ and $\nu_1=r$.
If $b_1>b_r$,
we obtain the numbers
$\ell\in\seisuu_{\geq 1}$
and $1\leq\nu_1<\cdots<\nu_{\ell}<r$
determined by the following condition.
\begin{itemize}
 \item $b_1=b_{\nu_1}$,
       $b_{\nu_j}>b_{\nu_{j}+1}=b_{\nu_{j+1}}$ $(j=1,\ldots,\ell-1)$,
       and $b_{\nu_{\ell}}>b_{\nu_{\ell}+1}=b_r$.
\end{itemize}
Moreover, we set $\nu_{0}=0$ and $\nu_{\ell+1}=r$.
For $1\leq i\leq r$,
we obtain the number $j(i)$ determined by
$\nu_{j(i)}< i\leq \nu_{j(i)+1}$.
If $j(i)\neq 0,\ell$,
we set
\begin{equation}
\label{eq;20.7.1.21}
 k_i(\vecb):=
  \nu_{j(i)+1}-\nu_{j(i)}+1
  -2(i-\nu_{j(i)})
  =\nu_{j(i)+1}+\nu_{j(i)}+1-2i.
\end{equation}
If $b_1-m\neq b_{r}$,
we define $k_i(\vecb)$
for $i$ such that
$j(i)=0,\ell$
by the same formula (\ref{eq;20.7.1.21}).
If $b_1-m=b_{r}$,
we define $k_i(\vecb)$ for $i$
such that $j(i)=0$ or $j(i)=\ell$
as follows:
\begin{equation}
\label{eq;20.7.1.22}
 k_i(\vecb):=
\left\{
\begin{array}{ll}
 \nu_1+\nu_{\ell}-r+1-2i & (j(i)=0)\\
 \nu_1+\nu_{\ell}+r+1-2i & (j(i)=\ell).
\end{array}
\right.
\end{equation}
Thus, we obtain
$\veck(\vecb)\in\seisuu^r$.
The following is a consequence of
the norm estimate of Simpson in \cite{s2}.
\begin{prop}
\label{prop;20.7.1.26}
 For $h\in\Harm(q)$,
 we set $\veck(h):=\veck(\vecb(h))$,
 where 
 $\vecb(h)\in\nbigp_{r,m}$
 is given as in Proposition {\rm\ref{prop;20.7.1.20}}.
Then,
we obtain the following refined estimates for $i=1,\ldots,r$:
\begin{equation}
\label{eq;20.7.13.20}
 \log|(dz)^{(r+1-2i)/2}|_h
 +(b_{i}(h)+i)\log|z|
 -\frac{k_{i}(h)}{2}\log\bigl(-\log|z|\bigr)
 =O(1).  
\end{equation}
\end{prop}
\pf
We explain an outline of the proof
(see also \cite{Toda-lattice}).
We use the notation in the proof of
Proposition \ref{prop;20.7.1.20}.
We set
$\nbigp^h_{<a}(\hyperk_{U^{\circ},r})$.
For $-1<a\leq 0$,
we obtain the following finite dimensional
complex vector space
\[
\Gr^F_{a}(\nbigp^h_0\hyperk_{U^{\circ},r|0}):=
\nbigp^h_a(\hyperk_{U^{\circ},r})\big/
\nbigp^h_{<a}(\hyperk_{U^{\circ},r}).
\]
Let $W$ denote the filtration on
$\Gr^F_a(\nbigp^h_0\hyperk_{U^{\circ},r|0})$
obtained as the monodromy weight filtration 
of the nilpotent endomorphism $\Res(\theta)$.

Let $n_i(h)$ denote the integer determined by
$-1<b_i(h)+n_i(h)\leq 0$.
We set $c_i(h):=b_i(h)+n_i(h)$.
Then,
$\vtilde_i:=z^{-n_i(h)}v_i$ is a section of
$\nbigp^h_{c_i(h)}(\hyperk_{U^{\circ},r})$.
Let $[\vtilde_i]$ denote the induced element of
$\Gr^F_{c_i(h)}(\nbigp^h_0\hyperk_{U^{\circ},r|0})$.
Then, the tuple 
$([\vtilde_i]\,|\,c_i(h)=a)$
is a base of 
$\Gr^F_{a}(\nbigp_0\hyperk_{U^{\circ},r|0})$.

We obtain the following for any $i=1,\ldots,r-1$:
\begin{equation}
\label{eq;20.7.13.30}
 \Res(\theta)[\vtilde_i]
 =\left\{
 \begin{array}{ll}
  [\vtilde_{i+1}]& (b_i(h)=b_{i+1}(h))\\
  0   & (b_i(h)>b_{i+1}(h))
 \end{array}
  \right.
\end{equation}
We also obtain the following:
\[
 \Res(\theta)[\vtilde_r]=
 \left\{
 \begin{array}{ll}
 \alpha(0)[\vtilde_1]  & (b_r(h)=b_1(h)-m)\\
  0 & \mbox{\rm (otherwise)}
 \end{array}
 \right.
\]
Hence, we obtain
\begin{equation}
\label{eq;20.7.13.21}
 W_k\Gr^F_{a}(\nbigp^h_0\hyperk_{U^{\circ},r|0})
 =\bigoplus_{\substack{
  c_i(h)=a\\
  k_i(h)\leq k}}
  \cnum [\vtilde_i].
\end{equation}
According to the norm estimate for tame harmonic bundles
\cite[\S7]{s2},
there exists $C>1$ such that
\[
 C^{-1}|z|^{-c_i(h)}(-\log|z|)^{k_i(h)/2}
 \leq
  |\vtilde_i|_h
 \leq
 C|z|^{-c_i(h)}(-\log|z|)^{k_i(h)/2}.
\]
Thus, we obtain (\ref{eq;20.7.13.20}).
\hfill\qed

\vspace{.1in}

We also remark the following existence.
\begin{prop}
\label{prop;20.7.6.20}
 For any $\vecb\in\nbigp_{r,m}$,
 there exist a neighbourhood $U_1$ of $0$ in $U$
 and $h\in \Harm(q_{|U_1\cap U^{\circ}})$
 such that
 $\vecb(h)=\vecb$.
\end{prop}
\pf
We may assume that
$q$ is nowhere vanishing on $U^{\circ}$.
\begin{lem}
There exist a neighbourhood $U_1$ of $0$
and a holomorphic function $\zeta$
on $U_1$ such that
 $\del_z\zeta(0)=1$
 and that 
$q=\alpha(0)\zeta^m(d\zeta/\zeta)^r$. 
The germ of $\zeta$ at $0$ is unique.
\end{lem}
\pf
We set $\beta:=\alpha(0)^{-1}\alpha$
which is nowhere vanishing on $U$.
We have the Taylor expansion
$\beta=z+\sum_{j=1}^{\infty} \beta_jz^j$.
First,
let us prove
that there exists a unique
formal power series
$\zeta^{(f)}=1+\sum_{j=2}^{\infty} \zeta^{(f)}_jz^j$
satisfying 
$\beta=(\zeta^{(f)}/z)^m(\del_z\zeta^{(f)})^r$.
The coefficient of
$z^{j}$ $(j>1)$
in $(\zeta^{(f)}/z)^m(\del_z\zeta^{(f)})^r$
is described as
the sum of
$(m+r(j+1))\zeta^{(f)}_{j+1}$
and a polynomial of
$\zeta^{(f)}_i$ $(i=2,\ldots,j)$.
Hence, $\zeta^{(f)}_{j+1}$ are uniquely determined
by an easy induction.
It particularly
implies that the germ of $\zeta$
is uniquely determined.
Moreover,
if there exists $n\in\seisuu_{>0}$
such that
$\beta_j=0$ unless $j\in n\seisuu_{\geq 0}$,
then we obtain
$\zeta^{(f)}_{j+1}=0$ unless $j\in n\seisuu_{\geq 0}$.

Let us prove the existence of a convergent solution
in the case where
$m_1:=m/r\in\seisuu_{>0}$.
Let $\beta^{1/r}$ be the holomorphic function on $U$
determined by the conditions
$(\beta^{1/r})^r=\beta$
and $\beta^{1/r}(0)=1$.
For any $z\in U$,
we take a path $\gamma(z)$ connecting $0$ and $z$ in $U$,
and we set
\[
F(z)=\int_{\gamma(z)}z^{m_1-1}\beta^{1/r}dz.
\]
Then, $F(z)$ is a holomorphic function on $U$
satisfying
\[
dF(z)=z^{m_1-1}\beta^{1/r}dz,
\quad
F(0)=0.
\]
Note that
$F(z)/z^{m_1}$ is holomorphic at $z=0$,
and
$(F(z)/z^{m_1})(0)=m_1^{-1}$.
There exist a neighbourhood $U_1$ of $0$ in $U$
and a holomorphic function $\zeta$ on $U_1$
such that
$F(z)=\frac{1}{m_1}\zeta^{m_1}$
and $\del_z\zeta(0)=1$.
Then, we obtain
$\zeta^{m_1-1}d\zeta=z^{m_1-1}\beta^{1/r}\,dz$,
which implies
$\zeta^m\alpha(0)(d\zeta/\zeta)^r
=z^m\alpha\cdot(dz/z)^r$.

Let us study the existence in the general case.
Let $\varphi_r:\cnum\lrarr \cnum$ be determined by
$\varphi_r(w)=w^r$.
We set $\Utilde:=\varphi_r^{-1}(U)$.
There exists a neighbourhood $U_1$
of $0$ in $U$
and a holomorphic function $\zetatilde$ on $\varphi_r^{-1}(U_1)$
such that
\[
\varphi_r^{\ast}(q)
=r^r\alpha(0)\zetatilde^{mr}(d\zetatilde/\zetatilde)^r
=\alpha(0)(\zetatilde^r)^m
\cdot (d\zetatilde^r/\zetatilde^r)^r.
\]
By the consideration in the first paragraph
of this proof,
we obtain that
$\zetatilde^r$ is a convergent power series of
$w^r$.
Hence, there exists
a holomorphic function $\zeta$ on $U_1$
such that $\varphi^{\ast}(\zeta)=\zetatilde^r$,
for which
$q=\zeta^m\alpha(0)\cdot (d\zeta/\zeta)^r$ holds.

\hfill\qed

\vspace{.1in}

We obtain
an embedding $\zeta:U_1\lrarr \cnum\subset\proj^1$.
Let $w$ be the standard coordinate on $\cnum$.
We set $q_0=\alpha(0)w^m(dw/w)^r$.
We obtain
$\zeta^{\ast}(q_0)=q$.
According to \cite{Toda-lattice},
there exists $h_1\in \Harm(q_0)$
such that
$\vecb(h_1)=\vecb$.
Then, we obtain
$h=\zeta^{\ast}(h_1)\in\Harm(q)$
which satisfies
$\vecb(h)=\vecb$.
\hfill\qed

\subsection{Examples on $\cnum$}

We introduce some examples
as a preparation for the proof of
Theorem \ref{thm;20.6.14.10}.
We shall see more examples later
in \S\ref{subsection;20.7.1.10}.

\subsubsection{Preliminary}
\label{subsection;20.6.24.200}
Let $\nbigp$ denote the set of
$\veca=(a_1,\ldots,a_r)\in\real^r$
satisfying
\[
 a_1\geq a_2\geq\cdots\geq a_r\geq a_1-1,
 \quad
 \sum a_i=0.
\]
For any $\veca\in\nbigp$,
we introduce a non-negative integer $\ell$
and integers $\nu_0,\nu_1,\ldots,\nu_{\ell},\nu_{\ell+1}$
as follows.
If $a_1=a_r$, we set
$\ell=0$, $\nu_0=0$ and $\nu_1=r$.
If $a_1>a_r$,
we obtain positive integers $\ell$
and $1\leq\nu_1<\cdots<\nu_{\ell}<r$
by the following condition.
\begin{itemize}
 \item $a_{\nu_1}=a_1$,
       $a_{\nu_j}>a_{\nu_j+1}$ $(j=1,\ldots,\ell-1)$
       and
       $a_{\nu_{\ell}}>a_{\nu_{\ell}+1}=a_r$.
\end{itemize}
Moreover, we set $\nu_{0}=0$
and $\nu_{\ell+1}=r$.
For $1\leq i\leq r$,
we obtain the number $j(i)$ determined by
$\nu_{j(i)}< i\leq \nu_{j(i)+1}$.
If $j(i)\neq 0,\ell$,
we set
\begin{equation}
\label{eq;20.4.27.52}
 k_i(\veca):=
  \nu_{j(i)+1}-\nu_{j(i)}+1
  -2(i-\nu_{j(i)})
  =\nu_{j(i)+1}+\nu_{j(i)}+1-2i.
\end{equation}
If $a_1-1\neq a_{r}$,
we define $k_i(\veca)$
for $i$ such that
$j(i)=0,\ell$
by the same formula (\ref{eq;20.4.27.52}).
If $a_1-1=a_{r}$,
we define $k_i(\veca)$ for $i$
such that $j(i)=0$ or $j(i)=\ell$
by
\begin{equation}
\label{eq;20.7.1.23}
 k_i(\veca):=
\left\{
\begin{array}{ll}
 \nu_1+\nu_{\ell}-r+1-2i & (j(i)=0)\\
 \nu_1+\nu_{\ell}+r+1-2i & (j(i)=\ell).
\end{array}
\right.
\end{equation}

\subsubsection{Examples}

Let $z$ be the standard coordinate of $\cnum$.
Let $(x,y)$ be the real coordinate system
obtained as $z=x+\sqrt{-1}y$.
We set $q:=\alpha e^{\sqrt{-1}z}\,(dz)^r$,
where $\alpha$ is a non-zero complex number.

\begin{prop}
 \label{prop;20.6.24.200}
For any $\veca\in\nbigp$,
there exists $h_{\veca}\in\Harm(q)$
 such that
 the following estimate holds
 on $\{y\geq 0\}$:
\[
 \log\bigl|
  (dz)^{(r+1-2i)/2}
 \bigr|_{h_{\veca}}
 -a_iy
 -\frac{k_i(\veca)}{2}\log (y+2)
 =O(1).
\]
\end{prop}
\pf
We apply the Kobayashi-Hitchin correspondence
on $\cnum^{\ast}$
to construct such $h\in\Harm(q)$
in \cite{Toda-lattice}.
Let $\varphi:\cnum\lrarr \cnum^{\ast}$ be determined by
$\varphi(z)=-\sqrt{-1}re^{\sqrt{-1}z/r}$.
We obtain
$\varphi^{\ast}((dw)^r)=e^{\sqrt{-1}z}(dz)^r$.

We consider the $r$-differential $q_1:=\alpha (dw)^r$ on $\cnum^{\ast}$.
For any $\veca\in\nbigp$,
we define $\vecb\in\nbigp_{r,r}$
by the following relation:
\[
 a_i=\frac{1}{r}\left(b_i+\frac{r+1}{2}\right).
\]
Note that $\veck(\veca)=\veck(\vecb)$.
According to \cite{Toda-lattice},
there exists $h_{1}\in \Harm(q_1)$
such that 
\[
 \log |(dw)^{(r+1-2i)/2}|_{h_{1}}
 +(b_i+i)\log |w|
 -\frac{k_i(\vecb)}{2}\log\bigl(-\log|w|\bigr)
 =O(1)
\]
around $w=0$.
We set $h_{\veca}:=\varphi^{-1}(h_{1})\in\Harm(q)$.
Note that
\[
 \varphi^{\ast}((dw)^{(r+1-2i)/2})
 =e^{\frac{r+1-2i}{2r}\sqrt{-1}z}(dz)^{(r+1-2i)/2}.
\]
Hence, we obtain the following estimates for $i=0,\ldots,r-1$
as $y\to\infty$:
\[
 \log|(dz)^{(r+1-2i)/2}|_{h_{\veca}}
 -\left(b_i+\frac{r+1}{2}\right)\frac{y}{r}
 -\frac{k_i(\veca)}{2}\log (y+2)
 =O(1).
\]
Because $h_{\veca}$ is invariant
under the translation by $(x,y)\longmapsto (x+2\pi,y)$,
we obtain the desired estimate for $h_{\veca}$.
\hfill\qed

\begin{rem}
Later in {\rm\S\ref{subsection;20.7.1.10}},
we shall see that 
the harmonic metric $h_{\veca}$ is uniquely determined
by the condition
\[
    \log\bigl|
  (dz)^{(r+1-2i)/2}
 \bigr|_{h_{\veca}}
-a_iy
 =O\bigl(\log(2+y)\bigr)
\]
 on $\{y\geq 0\}$.
Moreover,
  for any $h\in\Harm(q)$,
  there exists $\veca\in\Harm(q)$
  such that $h=h_{\veca}$.
\end{rem}

\begin{rem}
 The examples in Proposition {\rm\ref{prop;20.6.24.200}}
 will be used
 in the proof of Lemma {\rm\ref{lem;20.6.15.100}}. 
\end{rem}

\section{Preliminary from the classical analysis}

\subsection{Subharmonic functions on sectors}

\subsubsection{Phragm\'{e}n-Lindel\"{o}f theorem for subharmonic functions}

For any $1/2<a$ and $0\leq R$,
we set
$D_a(R):=\bigl\{
z\in\cnum\,\big|\,|z|>R,\,|\arg(z)|<\pi/2a
\bigr\}$.
Let $\Dbar_a(R)$ denote the closure in $\cnum$.

\begin{prop}[\mbox{\cite[\S7.3, Theorem 3]{Entire-functions}}]
\label{prop;20.4.27.10}
 Let $f$ be a continuous function on $\Dbar_a(0)$
 which is subharmonic in $D_a(0)$.
 Assume that there exist $C>0$, $0<\rho<a$, $M\in\real$
 and an exhaustive family $\{K_i\}$ of $D_a(0)$
 such that the following holds.
 \begin{itemize}
  \item
  $f\leq C(|z|^{\rho}+1)$ on $\del K_i$.
  \item 
  $f\leq M$ on $\del D_a(0)$.
 \end{itemize}
 Then, $f\leq M$ on $D_a(0)$.
\end{prop}
\pf
There exists $C_1>0$ such that
$f\leq C_1\Re\bigl(z^{\rho}+1\bigr)$
on $\del K_i$.
Because $f$ is assumed to be subharmonic,
we obtain that
$f\leq
 C_1\Re\bigl(z^{\rho}+1\bigr)$ on $K_i$,
and hence on $D_a(0)$.
There exists $C_2>0$ such that
$f\leq C_2(|z|^{\rho}+1)$ on $D_a(0)$.
Then, the claim follows from
Phragm\'{e}n-Lindel\"{o}f theorem for subharmonic functions
(see \cite[\S7.3, Theorem 3]{Entire-functions}).
\hfill\qed

\begin{cor}
\label{cor;20.5.2.1}
 Let $f$ be a subharmonic $C^{\infty}$-function on
 a neighbourhood of $\Dbar_a(R)$.
 Assume the following.
\begin{itemize}
 \item $f$ is bounded from above
	on $\del \Dbar_a(R)$.
 \item	There exist $C>0$, $0<\rho<a$ and an exhaustive
	family $\{K_i\}$ of $D_a(R)$ such that
	$f\leq C(|z|^{\rho}+1)$ on $\del K_i$.
\end{itemize}
Then, $f$ is bounded from above on $\Dbar_a(R)$.
\end{cor}
\pf
There exists an $\real$-valued $C^{\infty}$-function $g$
on a neighbourhood of
$\Dbar_a(0)$
such that $g=f$
on $\Dbar_a(R+1)$.
Let $\chi$ be a $C^{\infty}$-function on $\Dbar_a(0)$
with compact support such that
(i) $0\leq \chi\leq 1$,
(ii) $\chi(z)=1$ on
$\bigl\{z\in \Dbar_a(0)\,\big|\,|z|\leq R+1\bigr\}$.
Note that $(1-\chi)\triangle g\geq 0$ on $\Dbar_a(0)$,
where $\triangle=4\del_z\del_{\zbar}$.
There exists a continuous bounded function $G$
on $\Dbar_a(0)$
such that
$\triangle G=\chi\triangle g$
and that $G_{|\del D_a(0)}=0$.
There exists $C_1>0$ such that
$g-G\leq C_1$ on $\del \Dbar_a(0)$.
There exists $C_2>0$ such that
$g-G\leq C_2\Re\bigl(z^{\rho}+1\bigr)$
on $\del K_i$.
Because $g-G$ is subharmonic,
we obtain that
$g-G_1\leq C_2\Re\bigl(z^{\rho}+1\bigr)$
on $K_i$,
and hence on $D_a(R)$.
Because $D_a(0)\setminus D_a(R)$ is relatively compact,
there exists $C_3>0$ such that
$g-G_1\leq C_3(|z|^{\rho}+1)$ on $D_a(0)$.
By the Phragm\'{e}n-Lindel\"{o}f theorem,
we obtain the desired boundedness.
\hfill\qed

\vspace{.1in}
Let us give a variant.
Let $b_+$, $b_-$ and $\kappa$ be non-negative numbers
such that $b:=\max\{b_+,b_-\}\leq \kappa$.

\begin{cor}
\label{cor;20.6.12.21}
Let $f$ be a subharmonic function on $D_a(R)$
which extends to a continuous function on $\Dbar_a(R)$
such that the following holds.
\begin{itemize}
 \item For any $\delta>0$, there exists $M_{-,\delta}>0$
       such that 
\[
       f\leq M_{-,\delta}(|z|^{b_-}+1)
\]
       on
       $\bigl\{z\in D_a(R)\,\big|-\pi/2a< \arg(z)< -\delta\bigr\}$.
 \item For any $\delta>0$, there exists $M_{+,\delta}>0$
       such that 
\[
        f\leq M_{+,\delta}(|z|^{b_+}+1)
\]
       on $\bigl\{z\in D_a(R)\,\big|\,\delta< \arg(z)< \pi/2a\bigr\}$.
 \item There exist an exhaustive family $\{K_i\}$ of $D_a(R)$
       and $C>0$ such that
\[
 f\leq C(|z|^{\kappa}+1)
\]
       on $\del K_i$.
\end{itemize}
Then, there exists $C_1>0$ such that 
  $f\leq C_1(|z|^{b}+1)$ on $D_a(R)$.
 \end{cor}
\pf
We explain the case $0\leq b_{-}\leq b_+$.
The other case can be discussed similarly.
For any $a_1>a$,
we take a decreasing sequence $c_i$ in $\{a_1<c<2a_1\}$
such that $\lim c_i=a_1$.
We set
$K_{1,i}=K_i\cap D_{c_i}(R)$.
Then, $\{K_{1,i}\}$ is an exhaustive family of $D_{a_1}(R)$,
and there exists $C_2>0$ such that
$f\leq C_2(|z|^{\kappa}+1)$ on $\del K_{1,i}$.
Hence, by making $a$ larger,
we may assume that $\kappa<a$ from the beginning.

There exists $C_3>0$ such that
$f-C_3\Re\bigl(z^{b_+}+1\bigr)\leq 0$
 on
$\bigl\{z\in D_a(R)\,\big|\,
 \pi/4a<|\arg(z)|<\pi/2a
\bigr\}$. 
Then, we obtain the claim of the corollary
from Proposition \ref{prop;20.4.27.10}.
\hfill\qed

\subsubsection{Nevanlinna formula}

We set $\hyperh:=\{z\in\cnum\,|\,\Image(z)> 0\}$.
For any $R\geq 0$,
we set
$\hyperh(R):=\bigl\{z\in\hyperh\,\big|\,|z|>R\bigr\}$.
Let $\hyperhbar$ and $\hyperhbar(R)$
denote the closure of $\hyperh$ and $\hyperhbar(R)$,
respectively.
Let $(x,y)$ be the real coordinate system on $\cnum$
determined by $z=x+\sqrt{-1}y$.

Let $f$ be an $\real$-valued
$C^{\infty}$-function on $\hyperhbar(R)$
for some $R\geq 0$.
Suppose that there exist 
$C>0$, $0<\rho<1$ and $\delta\geq 0$
such that the following conditions are satisfied.
\begin{itemize}
 \item $f$ is $C^{\infty}$ in $\hyperh(R)$,
       and $|\triangle (f)|\leq C(1+y^2)^{-1-\delta}$.
       Here, $\triangle=\del_x^2+\del_y^2$.
 \item There exists an exhaustive family $\{K_i\}$ of $\hyperh(R)$
       such that $|f|\leq y+C(|z|^{\rho}+1)$
       on $\del K_i$ $(i=1,2,\ldots)$.
\end{itemize}

\begin{prop}
\label{prop;20.4.22.11}
 There exist $-1\leq a(f)\leq 1$
 and a constant $C_1>0$
 such that
 the following holds on
 $\hyperh(R)$:
 \[
 |f-a(f)y|\leq C_1(|z|^{\rho}+1).
 \]
If $|f|$ is bounded on $\del \hyperh(R)$,
 we obtain the following stronger estimate
 on $\hyperh(R)$
 for a positive constant $C_2>0$:
 \[
 |f-a(f)y|\leq C_2\log (2+y).
 \] 
 If $|f|$ is bounded on $\del \hyperh(R)$,
 and if $\delta>0$,
 then $|f-a(f)y|$ is bounded on $\hyperh(R)$.
\end{prop}
\pf
There exists a $C^{\infty}$-function
$\ftilde$ on $\hyperhbar$
such that
$\ftilde=f$ on $\hyperhbar(R+1)$.
There exists $\Ctilde>0$ such that the following holds.
\begin{itemize}
 \item $|\triangle\ftilde|\leq \Ctilde(1+y^2)^{-1-\delta}$.      
 \item There exists an exhaustive family
$\{\Ktilde_i\}$ such that
$|\ftilde|\leq y+\Ctilde(1+|z|^{\rho})$.
\end{itemize}
If $|f|$ is bounded on $\del \hyperh(R)$,
then $|\ftilde|$ is bounded on $\del\hyperh$.
It is enough to obtain the estimate for $\ftilde$.
Therefore,
we may assume $R=0$ from the beginning.

Let $\zeta$ be a complex variable,
and let $(\xi,\eta)$ be the real variables
determined by $\zeta=\xi+\sqrt{-1}\eta$.
On $\{y>0\}$,
we consider the following function
\[
 F_{1}(z):=
 \frac{1}{4\pi}
 \int_{-\infty}^{\infty}d\xi
 \int_{0}^{\infty}\,d\eta
 \left(
 \triangle f(\zeta)
 \cdot
  \log\left(
  \frac{(x-\xi)^2+(y-\eta)^2}{(x-\xi)^2+(y+\eta)^2}
  \right)
  \right).
\]
By Lemma \ref{lem;20.4.21.3} below,
the integral is convergent,
and we obtain the estimate
\[
F_{1}(z)=
\left\{
\begin{array}{ll}
 O\bigl(\log(1+y)\bigr) & (\delta=0) \\
 \\
O\bigl(y(1+y)^{-1}\bigr) & (\delta>0)
\end{array}
\right.
\]
Moreover, we obtain
$\triangle F_{1}=\triangle f$,
and hence
$F_{2}:=f-F_{1}$
is a harmonic function on $\hyperh$,
and continuous on $\hyperhbar$.

We set
$\phi:=\Re\bigl(e^{-\pi\rho\sqrt{-1}/2}(z+\sqrt{-1})^{\rho}\bigr)$,
where we consider the branch of
$(z+\sqrt{-1})^{\rho}$
such that
$(b\sqrt{-1}+\sqrt{-1})^{\rho}=(b+1)^{\rho}e^{\pi\rho\sqrt{-1}/2}$
for $b>0$.
It is a harmonic function,
and there exists $\epsilon>0$
such that $\epsilon(|z|^{\rho}+1)\leq \phi$
on $\hyperhbar$.
By the assumption,
there exist $C_{10}>0$
such that
$|F_2|\leq C_{10}\bigl(y+\phi\bigr)$
on $\del K_i$ $(i=1,2,\ldots)$.
Because $F_2$ is harmonic,
we obtain
$|F_2|\leq
C_{10}\bigl(y+\phi\bigr)$
on $K_i$,
and hence on $\hyperhbar$.
In particular,
we obtain
$|F_2|\leq C_{10}\phi$
 on $\del\hyperh$,
 and 
 $F_2+C_{10}\bigl(y+\phi\bigr)
 \geq 0$
 on $\hyperh$.
Then, according to the Nevanlinna formula
\cite[\S14, Theorem 1]{Entire-functions}
(see also \cite[Remark in page 101]{Entire-functions}),
there exists a real number $c$ such that
the following holds for $y>0$:
\[
 F_2+C_{10}\bigl(y+\phi\bigr)
 =cy+
 \frac{y}{\pi}
 \int_{-\infty}^{\infty}
  \frac{F_{2}(t,0)+C_{10}\phi(t,0)\,dt}{(t-x)^2+y^2}.
\]
Note that 
\[
\phi(x,y)=\frac{y}{\pi}
 \int_{-\infty}^{\infty}
 \frac{\phi(t,0)}{(t-x)^2+y^2}\,dt.
\]
We obtain
\[
f=(c-C_{10})y+F_1
+\frac{y}{\pi}\int_{-\infty}^{\infty}
 \frac{f(t,0)\,dt}{(t-x)^2+y^2}.
\]
Note that
$|f|=|F_2|\leq C_{10}\phi$ on $\del\hyperh$.
We obtain that
$\Bigl|
 \frac{y}{\pi}\int_{-\infty}^{\infty}
 \frac{f(t,0)\,dt}{(t-x)^2+y^2}
 \Bigr|
 \leq C_{10}\phi$ on $\hyperh$.
Then, we can easily obtain the claims of the proposition.
\hfill\qed

\vspace{.1in}
In the proof, we have used the following lemma.

\begin{lem}
\label{lem;20.4.21.3}
 There exists a constant $C_0>0$
 such that the following holds
 for any $z\in\hyperh$:
\begin{equation}
\label{eq;20.4.16.20}
0\leq \int_{-\infty}^{\infty}\,d\xi
  \int_{0}^{\infty}
   \frac{d\eta}{1+\eta^2}
  \left(
 \log\left|
  \frac{\zeta-\zbar}{\zeta-z}
	  \right|^2
       \right)
  \leq C_0\log (1+y).
\end{equation}
 For any $\delta>0$,
there exists a constant $C_{\delta}>0$
 such that the following holds
 for any $z$ with $y>0$:
\begin{equation}
\label{eq;20.4.24.1}
0\leq \int_{-\infty}^{\infty}\,d\xi
  \int_{0}^{\infty}
   \frac{d\eta}{(1+\eta^2)^{1+\delta}}
  \left(
 \log\left|
  \frac{\zeta-\zbar}{\zeta-z}
	  \right|^2
       \right)
  \leq C_{\delta}\frac{y}{1+y}.
\end{equation}
\end{lem}
\pf
Let us study the case $\delta=0$.
We set
$Z_1:=\{(\xi,\eta)\,|\,\eta\geq 2y\}$,
$Z_2:=\{(\xi,\eta)\,|\,y\leq 2\eta\leq 4y\}$,
and 
$Z_3:=\{(\xi,\eta)\,|\,0<2\eta\leq y\}$.
We obtain
\begin{multline}
\label{eq;20.4.27.1}
 \int_{Z_1}
\frac{d\xi d\eta}{1+\eta^2}
 \log\left|
  \frac{\zeta-\zbar}{\zeta-z}
	  \right|^2
= \int_{Z_1}
\frac{d\xi d\eta}{1+\eta^2}
 \log\left(
  \frac{\xi^2+(\eta+y)^2}{\xi^2+(\eta-y)^2}
  \right)
\\
 =\int_{Z_1}
 \frac{d\xi\,d\eta}{1+\eta^2}
 \log\left(
1+
 \frac{4y\eta}{\xi^2+(\eta-y)^2}
\right)
 \\
 \leq
 \int_{-\infty}^{\infty}\,d\xi
 \int_{2}^{\infty}
 \frac{y^2}{1+y^2\eta^2}
 \frac{4\eta d\eta}{\bigl(\xi^2+(\eta-1)^2\bigr)}.
\end{multline}
For $y\geq 1$,
(\ref{eq;20.4.27.1}) is dominated by
\[
  \int_{-\infty}^{\infty}\,d\xi
 \int_{2}^{\infty}
 \frac{4d\eta}{\eta\bigl(\xi^2+(\eta-1)^2\bigr)}
 <\infty.
\]
For $0<y\leq 1$,
(\ref{eq;20.4.27.1}) is dominated by
\[
 C_1\int_{2}^{\infty}
 \frac{y^2\,d\eta}{1+y^2\eta^2}
 \leq
 C_1\int_{-\infty}^{\infty}
 \frac{y^2\,d\eta}{1+y^2\eta^2}
 \leq
  C_2y.
\]

Similarly,
we obtain
\begin{multline}
\label{eq;20.4.27.2}
 \int_{Z_2}
\frac{d\xi\, d\eta}{1+\eta^2}
 \log\left|
  \frac{\zeta-\zbar}{\zeta-z}
	  \right|^2
=\int_{Z_2}
 \frac{d\xi\,d\eta}{1+\eta^2}
 \log\left(
1+
 \frac{4y\eta}{\xi^2+(y-\eta)^2}
\right)	  
 \\
 \leq
 \int_{|\xi|\geq 1}
 d\xi
 \int_{\frac{1}{2}}^2
 \frac{y^2}{1+y^2\eta^2}
 \frac{4\eta d\eta}{\bigl(\xi^2+(\eta-1)^2\bigr)}
\\
 +\int_{|\xi\leq 1}d\xi
 \int_{\frac{1}{2}}^{2}
 \frac{y^2\,d\eta}{1+y^2\eta^2}
 \log\left(
1+
 \frac{\eta}{\xi^2+(\eta-1)^2}
 \right).
\end{multline}
If $y\geq 1$,
(\ref{eq;20.4.27.2}) is dominated by
\begin{multline}
\int_{|\xi|\geq 1}
 d\xi
 \int_{\frac{1}{2}}^2
 \frac{4d\eta}{\eta\bigl(\xi^2+(\eta-1)^2\bigr)}
 +
\\ \int_{|\xi\leq 1}d\xi
 \int_{\frac{1}{2}}^{2}
  \frac{d\eta}{\eta^2}
 \log\left(
1+
 \frac{\eta}{\xi^2+(\eta-1)^2}
 \right)
 <\infty
\end{multline}
If $y\leq 1$,
(\ref{eq;20.4.27.2}) is dominated by
\[
 8y^2
 \int_{|\xi|\geq 1}\frac{d\xi}{\xi^2}
 \int_{\frac{1}{2}}^2d\eta
 +y^2\int_{|\xi|\leq 1}d\xi
 \int_{\frac{1}{2}}^2
 d\eta
  \log\left(
1+
 \frac{\eta}{\xi^2+(\eta-1)^2}
 \right)
 \leq C_3y^2.
\]

As for the integral over $Z_3$,
we obtain
\begin{multline}
 \int_{Z_3}
\frac{d\xi\, d\eta}{1+\eta^2}
 \log\left|
  \frac{\zeta-\zbar}{\zeta-z}
	  \right|^2
\leq
  \int_{Z_3}
\frac{d\xi\, d\eta}{1+\eta^2}
\frac{4y\eta}{\xi^2+(\eta-y)^2}
\\
 \leq
  C_{10}\int_{0}^{y/2}
 \frac{d\eta}{1+\eta^2}
 \frac{y\eta}{y-\eta}
 \leq
  2C_{10}\int_0^{y/2}
 \frac{\eta d\eta}{1+\eta^2}
 \leq C_{11}\log(1+y).
\end{multline}
Thus, we obtain (\ref{eq;20.4.16.20})
in the case $\delta=0$.

As for the case $\delta>0$,
the integrals over $Z_i$ $(i=1,2)$
are dominated as in the case of (\ref{eq;20.4.16.20}).
The integral over $Z_3$ is dominated by
$\int_0^{y/2}(1+\eta^2)^{-1-\delta}\,\eta\,d\eta
<C_{20}(\delta)y(1+y)^{-1}$.
Thus, we obtain (\ref{eq;20.4.24.1}).
\hfill\qed

\subsubsection{Holomorphic line bundles with a Hermitian metric
on upper half plane}

We state a consequence of
the Nevanlinna formula (Proposition \ref{prop;20.4.22.11})
on a holomorphic line bundle $L$ with a Hermitian metric $h$
on $\hyperh$,
which is fundamental in the classification of
solutions of Toda equations in terms of parabolic structures.
Let $R(h)$ denote the curvature of the Chern connection
associated with $(L,h)$.
By using the standard Euclidean metric $g=dz\,d\zbar$,
we assume the following condition.
\begin{itemize}
 \item $|R(h)|_{h,g}=O\bigl((1+y^2)^{-1}\bigr)$.
\end{itemize}
Note that this condition is analogue to the acceptability condition
 in {\rm\cite{s1, s2}}.

 \begin{lem}
  \label{lem;20.6.14.101}
Let $v$ be a global frame of $L$.
  If there exist $C>0$, $0<\rho<1$,
 and an exhaustive family $\{K_i\}$ of $\hyperh$
 such that $\bigl|\log|v|_h\bigr|\leq y+C(|z|^{\rho}+1)$ 
 on $\del K_i$,
then there exists $-1\leq a(h,v)\leq 1$
such that
$\log|v|_h-a(h,v)y=O(|z|^{\rho}+1)$
on $\hyperh$. 
If $\log|v|_h$ is bounded on $\del \hyperh$,
then we obtain
$\log|v|_h-a(h,v)y=O\bigl(\log(y+2)\bigr)$
on $\hyperh$.
\end{lem}
\pf
Because $R(h)=\frac{1}{4}\triangle(\log|v|^2_h)\,d\zbar\,dz$,
we obtain $\triangle(\log|v|_h^2)=O\bigl((1+y^2)^{-1}\bigr)$.
Then, the claim follows from Proposition \ref{prop;20.4.22.11}.
\hfill\qed

\begin{df}
 The number $a(h,v)$ is called the parabolic order of
 $v$ with respect to $h$.
\end{df}

\subsection{Finite exponential sums and their perturbation}

Let us recall some results in \cite{Tamarkin, Wilder},
which are summarized in \cite{Langer}.

\subsubsection{Finite exponential sums}
Let $c_0< c_1<c_2<\cdots<c_n$ be real numbers.
Let $a_i$ $(i=0,\ldots,n)$ be non-zero complex numbers.
We consider the entire function
\[
 F(\zeta)=\sum_{i=0}^n a_ie^{\sqrt{-1}c_i\zeta}.
\]
We set $Z(F):=\{\zeta\in\cnum\,|\,F(\zeta)=0\}$.
It is easy to see that there exists $L>0$ such that
any $\zeta\in Z(F)$ satisfies
$|\Image(\zeta)|<L$.
For any $x_1<x_2$,
we set
$Y(x_1,x_2):=
\bigl\{\zeta\in\cnum\,\big|\,
|\Image(\zeta)|\leq L,\,\,x_1\leq \Re(\zeta)\leq x_2
\bigr\}$.
\begin{prop}[\cite{Langer, Tamarkin, Wilder}]
\label{prop;20.6.20.2}
 If $Z(F)\cap\bigcup\{\Re(\zeta)=x_i\}=\emptyset$,
we obtain
\[
 -n+\frac{c_n-c_0}{2\pi}(x_2-x_1)
 \leq
 \#\Bigl(
  Y(x_1,x_2)\cap Z(F)
  \Bigr)
  \leq
  n+\frac{c_n-c_0}{2\pi}(x_2-x_1).
\]
It particularly implies that
the order of the zero of $F$ at any point is not larger than $n$. 
\hfill\qed 
\end{prop}

Note that the proposition implies that
$Z(F)\cap \{\Re(\zeta)=x_i\}\leq n$.
Hence, 
 we obtain
 the following for any $x_1<x_2$:
 \[
 -3n+\frac{c_n-c_0}{2\pi}(x_2-x_1)
 \leq
 \#\Bigl(
  Y(x_1,x_2)\cap Z(F)
  \Bigr)
  \leq
  3n+\frac{c_n-c_0}{2\pi}(x_2-x_1).
\] 

\subsubsection{Perturbation}

Let $G$ be a holomorphic function defined on
$W(C_1,C_2):=
\{\zeta\in\cnum\,|
\,\Re(\zeta)> C_1,\,\,|\Image(\zeta)|<C_2\}$
such that
\[
\lim_{\Re(\zeta)\to\infty}|F(\zeta)-G(\zeta)|=0.
\]
We set $Z(G):=\{\zeta\in W(C_1,C_2)\,|\,G(\zeta)=0\}$.
For any closed subset $K\subset\cnum$ and $\zeta\in\cnum$,
let $d(\zeta,K)$ denote the Euclidean distance between $\zeta$ and $A$,
i.e.,
$d(\zeta,A)=\min\bigl\{|\zeta'-\zeta|\,\big|\,\zeta'\in A\bigr\}$.

\begin{prop}
 \label{prop;20.6.21.20}
Take $0<C_2'<C_2$.
 \begin{itemize}
  \item 
	There exist $C_3>0$ such that
	the following holds for any $C_1+1<x_1<x_2$:
\[
 \#\Bigl(
 Y(x_1,x_2)\cap
 W(C_1,C_2')
 \cap Z(G)
 \Bigr)
 \leq
 C_3(1+|x_2-x_1|).
\]
  \item
There exist $C_1'>C_1$ and $C_4>0$
such that 
$|G(\zeta)|\geq C_4\delta^{n}$ holds
for any $0<\delta\leq 1$
and
any $\zeta\in W(C_1',C_2')$ satisfying
       $d(\zeta,Z(G))\geq \delta$.
 \end{itemize}
\end{prop}
\pf
We apply the argument in \cite{Wilder}.
Let us prove the first claim in the case $|x_2-x_1|=1$.
Assume that there exists a sequence
$x_j$ such that
\[
  \#\Bigl(
 Y(x_j,x_j+1)\cap
 W(C_1,C_2')
 \cap Z(G)
 \Bigr)
 \to \infty.
\]
Let
\begin{multline}
\Psi_j:
\{\zeta\in\cnum\,|\,
-\epsilon<\Re(\zeta)<1+\epsilon,\,|\Image(\zeta)|<C_2\}
\simeq \\
\{\zeta\in\cnum\,|\,
x_j-\epsilon<\Re(\zeta)<x_j+1+\epsilon,|\Image(\zeta)|<C_2\}
\end{multline}
be the isomorphism defined by $\Psi_j(\zeta)=\zeta+x_j$.
Going into a subsequence,
we may assume that
the sequence
$\Psi_j^{\ast}(F)$ is convergent to
$\sum b_ie^{\sqrt{-1}c_i\zeta}$,
where $b_i$ are complex numbers
such that $|b_i|=|a_i|$.
It implies that
the sequence $\Psi_j^{\ast}(G)$
is convergent to
$\sum b_ie^{\sqrt{-1}c_i\zeta}$.
Then,
we obtain the contradiction
by Proposition \ref{prop;20.6.20.2}
and Lemma \ref{lem;20.6.20.1} below.
Hence, there exists a constant $C_{10}>0$
such that
\[
  \#\Bigl(
 Y(x,x+1)\cap
 W(C_1,C_2')
 \cap Z(G)
 \Bigr)
 <C_{10}
\]
for any $x>C_1+1$.

Take any $C_1+1<x_1<x_2$.
We set $m_0:=\min\{m\in\seisuu\,|\,x_2-x_1\leq m\}\geq 1$.
Then, we obtain
\begin{multline}
  \#\Bigl(
 Y(x_1,x_2)\cap
 W(C_1,C_2')
 \cap Z(G)
 \Bigr)
 \leq \\
   \#\Bigl(
 Y(x_1,x_1+m_0)\cap
 W(C_1,C_2')
 \cap Z(G)
 \Bigr) 
 \leq m_0C_{10} \\
 \leq (1+|x_2-x_0|)C_{10}.
\end{multline}
Thus, we obtain the first claim.

To prove the second claim,
we prepare the following lemma.
\begin{lem}
\label{lem;24.1.5.2}
There exists $C_1''>C_1$
such that
the order of the zero of $G$ at any point of
$\overline{W(C_1'',C'_2)}\cap Z(G)$
is not larger than $n$.
\end{lem}
\pf
Suppose the contrary.
Let $\zeta_j\in \overline{W(C_1,C_2')}\cap Z(G)$
such that the order of zero of $G$ at $\zeta_j$
is strictly larger than $n$.
Let
\begin{multline}
\label{eq;24.1.5.1}
 \Psi_j:
 \bigl\{\zeta\in\cnum\,\big|\,
  -1<\Re(\zeta)<1,\,|\Image(\zeta)|<C_2\bigr\}
\simeq \\
 \bigl\{\zeta\in\cnum\,\big|\,
 -1<\Re(\zeta)-\Re(\zeta_j)<1,\,|\Image(\zeta)|<C_2
  \bigr\}
\end{multline}
be the isomorphism determined by
$\Psi_j(\zeta)=\zeta+\Re(\zeta_j)$.
We may assume that
the sequence $\Psi_j^{\ast}(F)$ is convergent to
a holomorphic function of the form
$\sum b_ie^{\sqrt{-1}c_i\zeta}$,
where $b_i$ are complex numbers such that
$|b_i|=|a_i|$.
It particularly implies that
$\Psi_j^{\ast}(G)$ is convergent to
$\sum b_ie^{\sqrt{-1}c_i\zeta}$.
We may also assume that
$\sqrt{-1}\Image(\zeta_j)$ is convergent
to $\sqrt{-1}\beta$ for some $-C_2'\leq \beta\leq C_2'$.
Because
the order of zero of $\Psi_j^{\ast}(G)$
at $\sqrt{-1}\Image(\zeta_j)$
are strictly larger than $n$,
we obtain that
the order of $\sum b_ie^{\sqrt{-1}c_i\zeta}$
is strictly larger than $n$.
But, it contradicts 
Proposition \ref{prop;20.6.20.2}.
\hfill\qed

\vspace{.1in}

Let us prove the second claim.
Suppose the contrary.
There exist sequences
$\epsilon_j>0$,
$0<\delta_j\leq 1$
and $\zeta_j\in W(C_1'',C_2')$
satisfying
$\lim\epsilon_j=0$,
$d(\zeta_j,Z(G))\geq \delta_j$
and
$|G(\zeta)|\leq \epsilon_j\delta_j^n$.
We shall deduce a contradiction.

Let us consider the case that the sequence $\zeta_j$ contains
a bounded sequence.
By going to a subsequence,
the sequence $\zeta_j$ is convergent to
$\zeta_{\infty}\in \overline{W(C_1'',C_2')}$.
It is easy to observe that
$\zeta_{\infty}\in Z(G)$,
and the order of zero at $\zeta_{\infty}$ is
strictly larger than $n$.
Hence, we obtain a contradiction
by Lemma \ref{lem;24.1.5.2}.

Let us consider the case that
the sequence $\zeta_j$ does not contain
a bounded sequence.
In particular, $\Re(\zeta_j)\to\infty$ as $j\to\infty$.
We consider the maps $\Psi_j$ as in (\ref{eq;24.1.5.1}).
We may assume that
the sequence $\Psi_j^{\ast}(F)$ is convergent to
a holomorphic function of the form
$\sum b_ie^{\sqrt{-1}c_i\zeta}$,
where $b_i$ are complex numbers such that
$|b_i|=|a_i|$.
It particularly implies that
$\Psi_j^{\ast}(G)$ is convergent to
$\sum b_ie^{\sqrt{-1}c_i\zeta}$.
Then, we obtain a contradiction
by Proposition \ref{prop;20.6.20.2}
and Lemma \ref{lem;20.6.18.20} below.
Thus, we obtain the second claim of the proposition.
\hfill\qed

\begin{cor}
\label{cor;24.1.5.10}
For any $0<\delta\leq 1$,
there exists $C_5>0$ such that
we obtain $|G(\zeta)|\geq C_4\delta^n$
for any $\zeta\in W(C_1,C_2')$
satisfying $d(\zeta,Z(G))\geq \delta$.
\hfill\qed
\end{cor}

\subsubsection{Appendix}

Let $Y\subset\cnum$ be a connected open subset.
Let $Y_1$ be a relatively compact open subset of $Y$.
Let $F_i:Y\lrarr \cnum$ $(i=1,2,\ldots)$
be holomorphic functions
which uniformly converge to
a holomorphic function $F_{\infty}:Y\lrarr \cnum$.
Note that any derivatives
$\del_z^jF_i$ converges to
$\del_z^jF_{\infty}$.
Assume that $F_{\infty}$ is not constantly $0$.

For $i=1,2,\ldots,\infty$,
we set
$Z(F_i):=\{\zeta\in Y\,|\,F_i(\zeta)=0\}$.
Let $\gminim_i(\zeta)\in\seisuu_{>0}$
denote the order of the zero of $F_i$ at $\zeta$.

\begin{lem}
\label{lem;20.6.20.1}
 There exists $i_0$ such that
 the following holds for any $i\geq i_0$:
\[
 \sum_{\zeta\in Z(F_i)\cap\Ybar_1}\gminim_i(\zeta)
 \leq
 \sum_{\zeta\in Z(F_{\infty})\cap\Ybar_1}\gminim_{\infty}(\zeta).
\]
\end{lem}
\pf
We may assume that $Y_1$ is simply connected.
There exists a relatively compact open subset
$Y_2\subset Y$
such that
(i) $\del Y_2$ is smooth,
(ii) $Y_1\subset Y_2$,
(iii) $Z(F_{\infty})\cap \Ybar_1=Z(F_{\infty})\cap \Ybar_2$.
The condition (iii) implies that
$F_{\infty}(\zeta)\neq 0$ for any $\zeta\in \del Y_2$.
There exists $\epsilon>0$
such that $\min_{\zeta\in\del Y_2}|F_{\infty}(\zeta)|>2\epsilon$.
There exists $i_0$ such that
$|F_i(\zeta)-F_{\infty}(\zeta)|\leq\epsilon/2$
for any $\zeta\in \del Y_2$ and $i\geq i_0$.
Then, the claim holds according to Rouch\'{e}'s theorem.
\hfill\qed

\vspace{.1in}

Let $n$ be a positive integer such that
$\gminim_{\infty}(\zeta)\leq n$
for any $\zeta\in \Ybar_1$.
For any closed subset $K\subset \cnum$,
let $d(\zeta,K)$ denote the Euclidean distance between
$\zeta$ and $K$.
There exists $\delta_0>0$
such that
$\{\zeta\in\cnum\,|\,d(\zeta,\Ybar_1)\leq 2\delta_0\}
\subset Y$.
For any $0<\delta<\delta_0$,
we set
$W_i(\delta):=\{\zeta\in Y_1\,|\,
d(\zeta,Z(F_i))\geq \delta\}$.

\begin{lem}
\label{lem;20.6.18.20}
 There exist $C>0$ and $i_0$
 such that
 we obtain
 $|F_i|\geq C\delta^n$
 on $W_i(\delta)$
 for any $i\geq i_0$ and any $0<\delta<\delta_0$.
\end{lem}
\pf
Assume the contrary.
Then, 
there exist
a sequence of positive numbers
$C_1>C_2>\cdots$,
a sequence $i(j)$,
positive numbers $\delta_j>0$
and $\zeta_j\in W_{i(j)}(\delta_j)$
such that
$\lim C_j=0$
and that
$|F_{i(j)}(\zeta_j)|<C_j\delta_j^n$.
Going to a subsequence,
we may assume that $i(j)=j$.
We may assume $\delta_j= d(\zeta_j,Z(F_j))$.
Going to a subsequence,
we may assume that
the sequence $\zeta_j$ is convergent to
$\zeta_{\infty}\in \Ybar_1$.
Because $F_{\infty}(\zeta_{\infty})=\lim F_j(\zeta_j)$,
we obtain $F_{\infty}(\zeta_{\infty})=0$,
and hence
$\zeta_{\infty}\in Z(F_{\infty})\cap \Ybar_1$.
For any $\epsilon>0$,
there exists $j_0$
such that $d(\zeta_{\infty},Z(F_j))<\epsilon$
for any $j\geq j_0$.
Hence, we obtain $\lim \delta_j=0$.
Going to a subsequence,
we may assume that the sequence $\delta_j$ is decreasing.

Fix a relatively compact open
neighbourhood $U$ of $\zeta_{\infty}$
in $Y$
such that
$\Ubar\cap Z(F_{\infty})=\{\zeta_{\infty}\}$.
We set $\ell:=\gminim_{\infty}(\zeta_{\infty})\leq n$.
Going to a subsequence,
we may assume that
$F_{j}$ is nowhere vanishing on $\del U$
for any $j$.
We set $S_j:=U\cap Z(F_j)$.
We obtain
$\sum_{a\in S_j}\gminim_j(a)=\ell$.
We set
\[
 G_{U,j}(\zeta):=
 F_{j|U}(\zeta)\cdot \prod_{a\in S_j}(\zeta-a)^{-\gminim_j(a)},
\]
which is nowhere vanishing on $U$ for any $j$.
Because the sequence $F_{j}$ is convergent to $F_{\infty}$,
the sequence
$\prod_{a\in S_j}(\zeta-a)^{\gminim_j(a)}$
is convergent to
$(\zeta-\zeta_{\infty})^{\ell}$.
The sequence
$G_{U,j}$ is convergent to
$G_{U,\infty}(\zeta):=
F_{\infty|U}(\zeta)(\zeta-\zeta_{\infty})^{-\ell}$,
which is also nowhere vanishing
on $U$.

By the assumption,
$|\zeta_j-a|\geq \delta_j$
for any $a\in S_j$.
Hence,
$\prod_{a\in S_j}(\zeta_j-a)^{\gminim_j(a)}
\geq \delta_j^{\ell}$.
Then, there exists $A>0$ such that
the following holds for any $j$:
\[
 |F_{j}(\zeta_j)|
 \geq A\cdot \delta_j^{\ell}.
\]
Hence, we obtain
$A\delta_j^{\ell}
 \leq
 |F_{j}(\zeta_j)|
 \leq C_j\delta_j^{n}$.
It implies
 $0<A\leq C_j\delta_j^{n-\ell}\to 0$,
 which is a contradiction.
\hfill\qed

\subsection{Holomorphic functions with multiple growth orders}
\label{subsection;20.7.2.11}

\subsubsection{Growth orders}

Let $\varpi:\cnumtilde\lrarr \cnum$
be the oriented real blowing up
at $0$.
We set
$\cnum^{\ast}:=\cnum\setminus\{0\}$.
Let $\iota:\cnum^{\ast}\lrarr\cnumtilde$
denote the natural inclusion.
Let $z$ denote the standard coordinate of $\cnum$.
We identify $\cnumtilde$ with $\real_{\geq 0}\times S^1$
by the polar decomposition of $z$.

 \begin{notation}
  Let $Q=(0,e^{\sqrt{-1}\theta_0})\in\varpi^{-1}(0)$.
Let $\iota_{\ast}(\nbigo_{\cnum^{\ast}})_Q$
denote the stalk of the sheaf
$\iota_{\ast}(\nbigo_{\cnum^{\ast}})$ at $Q$.
  We fix a branch of $\log z$ around $Q$.
\begin{itemize}
 \item
 Let $\gbigi(Q)\subset \iota_{\ast}(\nbigo_{\cnum^{\ast}})_Q$
 be the $\cnum$-linear subspace generated by
 $z^{-a}$ $(a>0)$,
 i.e.,
$\gbigi(Q)=\bigoplus_{a>0}\cnum\,z^{-a}$.
 \item
 Note that for any $\gminia\in\gbigi(Q)\setminus\{0\}$,
 $|\gminia|^{-1}\Re(\gminia)$
 induces a continuous function
 on a neighbourhood of $Q$.
 The continuous function is also denoted by
 $|\gminia|^{-1}\Re(\gminia)$.
\item
For any $\gminia\in\gbigi(Q)\neq 0$,
we have the expression
$\gminia=\alpha z^{-\rho}+\sum_{0<c<\rho}\alpha_cz^{-c}$
where $\alpha\neq 0$.
We set
$\deg(\gminia):=\rho$.
\item
For $\gminia,\gminib\in\gbigi(Q)$,
we say
$\gminia\prec_Q\gminib$
     if $|\gminib-\gminia|^{-1}\Re(\gminib-\gminia)(Q)>0$
     or $\gminia=\gminib$.
It defines a partial order $\prec_Q$ on $\gbigi(Q)$.
\end{itemize}
 \end{notation}

\begin{rem}
 Note that the orders $\prec_{Q}$ are opposite
 to the order used in {\rm\cite{Mochizuki-wild}}.
\end{rem}

\subsubsection{Regularly bounded holomorphic functions}

\begin{df}
 We say that $f\in \iota_{\ast}(\nbigo_{\cnum^{\ast}})_Q$
 is regularly bounded
 if there exist a neighbourhood $\nbigu$ of $Q$ in $\cnumtilde$,
 non-zero complex numbers
 $\alpha_1,\ldots,\alpha_m$,
 and 
 mutually distinct real numbers $c_1,c_2,\ldots,c_m$
 such that the following holds.
 \begin{itemize}
  \item $f$ induces a holomorphic function on
	$\nbigu\setminus\varpi^{-1}(0)$.
  \item $|f(z)-\sum_{i=1}^m\alpha_i z^{\sqrt{-1}c_i}| \to 0$
	as $|z|\to 0$ in $\nbigu\setminus\varpi^{-1}(0)$.
 \end{itemize}
\end{df}

\subsubsection{Holomorphic functions with single growth order}

 \begin{df}
\label{df;20.9.11.10}
We say that $f\in\iota_{\ast}(\nbigo_{\cnum^{\ast}})_Q$
has a single growth order
if there exist
 $\gminia(f,Q)\in \gbigi(Q)$,
 $a(f,Q)\in\real$,
 $j(f,Q)\in\seisuu_{\geq 0}$
 such that
 $e^{-\gminia(f,Q)}z^{-a(f,Q)}(\log z)^{-j(f,Q)} f$
 is regularly bounded.
 \end{df}
 
\begin{df}
 Suppose that $f\in\iota_{\ast}(\nbigo_{\cnum^{\ast}})_Q$
 has single growth order at $Q$.
\begin{itemize}
 \item
       We say that $f$ is simply positive (resp. negative) at $Q$
      if
\[
      |\gminia(f,Q)|^{-1}\Re(\gminia(f,Q))(Q)>0\quad
      \mbox{(resp. $|\gminia(f,Q)|^{-1}\Re(\gminia(f,Q))(Q)<0$)}.
\]
  \item If $\gminia(f,Q)=0$,
	then we say that $f$ is neutral at $Q$.
 \item If $\gminia(f,Q)\neq 0$ but
       $|\gminia(f,Q)|^{-1}\Re(\gminia(f,Q))(Q)=0$,
       we say that $f$ is turning at $Q$.
\end{itemize}
 \end{df}
 
The following lemma is easy to see.
\begin{lem}
\label{lem;20.9.11.2}  
 If $f\in\iota_{\ast}(\nbigo_{\cnum^{\ast}})_Q$
 has a single growth order,
 there exists a neighbourhood $\nbigu$ of $Q$
 such that
 (i) $\nbigu\cap\varpi^{-1}(0)$ is connected,
 (ii) $f$ induces a holomorphic function on
 $\nbigu\setminus\varpi^{-1}(0)$,
 (iii) $f$ has a single growth order at
       any point $Q'\in\nbigu'\cap\varpi^{-1}(0)$,
 (iv) $\gminia(f,Q')=\gminia(f,Q)$.
 \hfill\qed
 \end{lem} 
Note that (iv) implies that
the conditions
 ``simply positive'',
 ``simply negative''
 and  ``neutral''
 are preserved
 if $\nbigu$ is sufficiently small.
If $f$ is turning at $Q$,
$|\gminia(f,Q)|^{-1}\Re(\gminia(f,Q))$
is positive on a connected component of
$(\nbigu\cap\varpi^{-1}(0))\setminus\{Q\}$
and negative on the other connected component
of $(\nbigu\cap\varpi^{-1}(0))\setminus\{Q\}$.

\begin{lem}
\label{lem;20.9.11.3}
Let $f$ be a section of
$\iota_{\ast}\nbigo_{\cnum^{\ast}}$
on an open subset $V\subset\widetilde{\cnum}$.
Let $I\subset V\cap\varpi^{-1}(0)$ be an interval.
If $f$ has single growth order at each $Q\in I$,
then we obtain $\gminia(f,Q_1)=\gminia(f,Q_2)$
for any $Q_1,Q_2\in I$.
\end{lem}
\pf
It follows from Lemma \ref{lem;20.9.11.2}.
\hfill\qed

\subsubsection{Holomorphic functions with multiple growth orders}
 
\begin{df}
\label{df;20.8.17.10}
 We say that $f\in\iota_{\ast}(\nbigo_{\cnum^{\ast}})_Q$
 has multiple growth orders
 if $f$ is expressed as a finite sum
 $\sum_{i=1}^N f_i$
 such that 
 (i) each $f_i$ has a single growth order,
 (ii) $\gminia(f_i,Q)\neq \gminia(f_j,Q)$ $(i\neq j)$.
\end{df}

\begin{lem}
\label{lem;20.8.18.1}
Suppose that $f\in\iota_{\ast}(\nbigo_{\cnum^{\ast}})_Q$
has multiple growth orders.
\begin{itemize}
 \item There exists a finite subset $\nbigi\subset\gbigi(Q)$
       and an expression $f=\sum_{\gminia\in\nbigi}f_{\gminia}$
       such that
  (i) for two distinct elements
  $\gminia,\gminib\in\nbigi$,
  neither $\gminia\prec_Q\gminib$ nor $\gminib\prec_Q\gminia$
  holds,
  (ii) each $f_{\gminia}$  has single growth order
  with $\gminia(f_{\gminia},Q)=\gminia$.
  Such an expression is called reduced.
 \item
  There exists a neighbourhood $\nbigu$ of $Q$
  such that
  (i) $\nbigu\cap\varpi^{-1}(0)$ is connected,
  (ii) $f$ induces a holomorphic function on
  $\nbigu\setminus\varpi^{-1}(0)$,
  (iii) $f$ has single growth order at
      any point $Q'\in\bigl(\nbigu\cap\varpi^{-1}(0)\bigr)\setminus\{Q\}$,
      and
      $\gminia(f,Q')$ is the unique maximal element of
      the partially ordered set
      $(\nbigi,\prec_{Q'})$.
\end{itemize}
 \end{lem}
\pf
There exists an expression
$f=\sum_{i=1}^mf_i$
as in Definition \ref{df;20.8.17.10}.
Let $\nbigi$ denote the set of the maximal elements
in the partially ordered set
$\{\gminia(f_1,Q),\ldots,\gminia(f_m,Q)\}$
with $\prec_Q$.
There exists a decomposition
\[
 \{\gminia(f_1,Q),\ldots,\gminia(f_m,Q)\}
=\coprod_{\gminia\in\nbigi}\nbigj_{\gminia}
\]
such that any $\gminib\in \nbigj_{\gminia}$
satisfy $\gminib\prec_Q\gminia$.
We set
$f_{\gminia}:=\sum_{\gminia(f_i,Q)\in\nbigj_{\gminia}}f_i$.
We obtain the expression
$f=\sum_{\gminia\in\nbigi}f_{\gminia}$.
It is easy to see that
$f_{\gminia}$ has a single growth order
with $\gminia(f_{\gminia},Q)=\gminia$,
and the first claim is proved.
The second claim is clear.
\hfill\qed

\vspace{.1in}

Note that a reduced expression $f=\sum_{\gminia\in\nbigi}f_{\gminia}$
is not uniquely determined.

 \begin{notation}
 Let $V$ be an open subset of $\cnumtilde$.
 Let $f$ be a holomorphic function on $V\setminus\varpi^{-1}(0)$
 which has multiple growth orders at any point of
  $V\cap\varpi^{-1}(0)$.
  Then, let $\nbigz(f)$ denote the set of
  the points $Q\in V\cap\varpi^{-1}(0)$
  such that one of the following condition is satisfied.
  \begin{itemize}
   \item $f$ does not have single growth order at $Q$.
   \item $f$ has single growth order and turning
	 at $Q$.
  \end{itemize}
Note that $\nbigz(f)$ is discrete in
$V\cap\varpi^{-1}(0)$.
\end{notation}

 \begin{rem}
  Let $U$ be a neighbourhood of $0$ in $\cnum$.
  Suppose that a non-zero holomorphic function $f$
  on $U\setminus\{0\}$
 satisfies a linear differential equation
 $\del_z^nf+\sum_{j=0}^{n-1} a_j(z) \del_z^jf=0$,
 where $a_j$ are meromorphic functions on $(U,0)$.
  Then, $f$ has multiple growth orders
  at any point of $\varpi^{-1}(0)$,
  which is a consequence of
  the classical asymptotic analysis.
 (For example, see {\rm\cite[\S II.1]{Majima}}.)
 \end{rem}

\subsubsection{Coordinate change}
 
Let $U$ be an open neighbourhood of $0$ in $\cnum$.
Let $\gamma:U\lrarr \cnum$ be a holomorphic function
such that $\gamma(0)=0$ and $\del_z\gamma(0)=1$.
It induces an automorphism $\gamma^{\ast}$
of $\iota_{\ast}(\nbigo_{\cnum^{\ast}})_Q$.

For any $z^{-a}$ $(a>0)$,
we have the expansion
\[
\gamma^{\ast}(z^{-a})
=z^{-a}(\gamma^{\ast}(z)/z)^{-a}
=z^{-a}+\sum_{b>-a}\alpha_{a,b}z^{b}.
\]
We define
\[
\gamma^{\ast}_-(z^{-a}):=
z^{-a}+\sum_{-a<b<0}\alpha_{a,b}z^{b}.
\]
It induces an injection
$\gamma_-^{\ast}:\gbigi(Q)\lrarr\gbigi(Q)$.
The following lemma is easy to see.
\begin{lem}
\label{lem;20.6.22.10}
$f\in\iota_{\ast}(\nbigo_{\cnum^{\ast}})_Q$
has a single growth order
if and only if
$\gamma^{\ast}(f)$ has
a single growth order.
In that case,
we have
$\gminia(\gamma^{\ast}(f),Q)
=\gamma_-^{\ast}(\gminia(f,Q))$.
As a result  
  $f\in\iota_{\ast}(\nbigo_{\cnum^{\ast}})_Q$
 has multiple growth orders
  if and only if
 $\gamma^{\ast}(f)$ has
 multiple growth orders.
\hfill\qed
\end{lem}

Let $\gbigb_{\cnumtilde,Q}$ denote the set of
$f\in\iota_{\ast}(\nbigo_{\cnum^{\ast}})_Q$
with multiple growth orders.
It is independent of the choice of
a coordinate $z$
by Lemma \ref{lem;20.6.22.10}.

\subsubsection{Sheaf of holomorphic functions
   with multiple growth orders}

Let $X$ be any Riemann surface
with a discrete subset $D$.
Let $\varpi_{X,D}:\Xtilde_D\lrarr X$
denote the oriented real blowing up.
Let $\iota_{X\setminus D}:X\setminus D\lrarr \Xtilde_D$
denote the inclusion.
 
 \begin{notation}
  Let $\gbigb_{\Xtilde_D}\subset
  \iota_{X\setminus D,\ast}(\nbigo_{X\setminus D})$
  denote the subsheaf determined as follows.
  \begin{itemize}
   \item $\gbigb_{\Xtilde_D|X\setminus D}=\nbigo_{X\setminus D}$.
   \item For any $P\in D$,
	 we take a holomorphic embedding
	 $z_P:X_P\lrarr \cnum$ of a neighbourhood $X_P$ of $P$
	 in $X$ such  that $z_P(P)=0$.
	 Note that it induces the isomorphisms
	 $\ztilde_P:\varpi_{X,D}^{-1}(P)\simeq\varpi^{-1}(0)$
	 and
	 $\ztilde_P^{\ast}:
	 \iota_{\ast}(\nbigo_{\cnum^{\ast}})_{\ztilde_P(Q)}
	 \simeq
	 \iota_{X\setminus D,\ast}(\nbigo_{X\setminus D})_Q$
	 for any $Q\in\varpi_{X,D}^{-1}(P)$.
	 Then,
	 $\gbigb_{\Xtilde_D,Q}=
	 \ztilde_P^{-1}\bigl(\gbigb_{\cnumtilde,\ztilde_P(Q)}\bigr)$
	 holds.
  \end{itemize}
 \end{notation}

 If $f$ is a section of $\gbigb_{\Xtilde_D}$
 on an open subset $V\subset\Xtilde_D$,
 $\alpha f$ $(\alpha\in\cnum^{\ast})$ are also sections of
 $\gbigb_{\Xtilde_D}$ on $V$.
 But, even if $f_i$ $(i=1,2)$ are sections of
 $\gbigb_{\Xtilde_D}$ on $V$,
 $f_1+f_2$ is not necessarily
 a section of $\gbigb_{\Xtilde_D}$ on $V$.
 For example,
 we set
 $f_1=e^{-z^{-1}}$
 and $f_2=-e^{-z^{-1}}+e^{-e^{z^{-1}}}$
 around $\arg(z)=0$.
 Then, $f_i$ are sections of $\gbigv_{\Xtilde_D}$,
 but $f_1+f_2=e^{-e^{z^{-1}}}$ is not.

\subsubsection{Intervals with some properties}

Let $f$ be a section of
$\gbigb_{\cnumtilde}$ on an open subset $V\subset\cnumtilde$.

\begin{df}
\label{df;20.8.18.10}
 An open interval $I\subset\varpi^{-1}(0)\cap V$
 is called positive (resp. negative) with respect to $f$
if $f$ is simply positive (resp. negative) at each point of
$I\setminus\nbigz(f)$. 
It is called maximal if moreover
there does not exist an interval
$I_1\subset\varpi^{-1}(0)\cap V$
such that
(i) $I_1$ is positive (resp. negative) with respect to $f$,
(ii) $I\subsetneq I_1$.
\end{df}

\begin{df}
\label{df;20.8.18.11}
 An open interval $I\subset\varpi^{-1}(0)\cap V$
 is called neutral with respect to $f$
if $f$ is neutral at each point of
$I\setminus\nbigz(f)$.
 It is called maximal if moreover
 there does not exist an interval
 $I_1\subset\varpi^{-1}(0)\cap V$
 such that
 (i) $I_1$ is neutral with respect to $f$,
 (ii) $I\subsetneq I_1$.
\end{df}

\begin{df}
\label{df;20.8.18.12}
 An open interval $I\subset\varpi^{-1}(0)\cap V$
is called non-positive with respect to $f$
if $f$ is neutral or simply negative
at each point of $I\setminus\nbigz(f)$.
It is called maximal if moreover
there does not exist an interval
$I_1\subset\varpi^{-1}(0)\cap V$
such that
(i) $I_1$ is non-positive with respect to $f$,
(ii) $I\subsetneq I_1$.
\end{df}

In Definitions {\rm\ref{df;20.8.18.10}}--{\rm\ref{df;20.8.18.12}},
if $\Ibar$ is contained in $V$
and if $I$ is maximal,
then $\del I\subset\nbigz(f)$. 

\begin{lem}
\label{lem;20.8.18.20}
 If an interval $I$ is non-positive with respect to $f$,
 then $I$ is either negative
 or neutral with respect to $f$.
 Moreover, if $I$ is neutral with respect to $f$,
 then $I\cap\nbigz(f)=\emptyset$. 
\end{lem}
\pf
Let $I$ be a non-positive interval.
Let $I\setminus\nbigz(f)=\coprod I_i$
be the decomposition into the connected components.
The numbering $I_i$ is given in the counter-clockwise way.
If $f$ is negative (neutral) at a point $Q\in I_i$,
then $f$ is negative at any point of $I_i$.

Let $Q\in I\cap \nbigz(f)$
such that $Q=\Ibar_i\cap \Ibar_{i+1}$.
There exists a reduced expression
$f=\sum_{\gminia\in\nbigi}f_{\gminia}$ at $Q$
as in Lemma \ref{lem;20.8.18.1}.
Take $Q_i\in I_{i}$ and $Q_{i+1}\in I_{i+1}$.
Suppose that $f$ is negative at $Q_i$
and neutral at $Q_{i+1}$,
and we shall deduce a contradiction.
Note that
$\gminia(f,Q_i)=\max(\nbigi,\prec_{Q_i})=:\gminib$
and
$\gminia(f,Q_{i+1})=\max(\nbigi,\prec_{Q_{i+1}})=0$.
We obtain
$|\gminib|^{-1}\Re(\gminib)(Q)=0$.
If $Q_{i+1}$ is sufficiently close to $Q$,
we obtain
$|\gminib|^{-1}\Re(\gminib)(Q_{i+1})>0$,
i.e., $0\prec_{Q_{i+1}}\gminib$.
It contradicts $\gminib\prec_{Q_{i+1}}0$.
Therefore,
if $f$ is negative at $Q_{i}$,
then $f$ is negative at $Q_{i+1}$.
If $f$ is neutral at $Q_i$,
then $f$ is neutral at $Q_{i+1}$,
and moreover we obtain $\nbigi=\{0\}$.
Then, we easily obtain the claim of the lemma.
\hfill\qed

\begin{lem}
\label{lem;20.6.22.3}
Suppose that $I$ is negative with respect to $f$
such that $I\cap\nbigz(f)$ consists of one point $Q$.
Choose $Q_1$ and $Q_2$ from the two connected components
of $I\setminus\{Q\}$.
 Let $f=\sum_{\gminia\in\nbigi} f_{\gminia}$
 be a reduced expression at $Q$.
Then, the following holds.
\begin{itemize}
 \item $\deg\gminia(f,Q_1)=\deg\gminia(f,Q_2)$.
 \item For any $\gminia\in\nbigi$,
       we obtain
       $\deg(\gminia)=\deg(\gminia(f,Q_i))$.
       Moreover,
       $|\gminia|^{-1}\Re(\gminia)(Q)<0$.
\end{itemize}
\end{lem}
\pf
Note that $\gminia(f,Q_i)\in\nbigi$.
Let us prove
\begin{equation}
\label{eq;24.1.6.40}
 \deg\gminia(f,Q_1)\geq \deg\gminia(f,Q_2).
\end{equation}
Suppose
$\deg\gminia(f,Q_1)<\deg\gminia(f,Q_2)$,
and we shall deduce a contradiction.
If $|\gminia(f,Q_2)|^{-1}\Re(\gminia(f,Q_2))(Q)<0$,
we obtain
$\gminia(f,Q_2)\prec_Q\gminia(f,Q_1)$,
and hence
$\gminia(f,Q_2)\prec_{Q_2}\gminia(f,Q_1)$
holds,
which contradicts
$\gminia(f,Q_2)=\max(\nbigi,\prec_{Q_2})$.
If $|\gminia(f,Q_2)|^{-1}\Re(\gminia(f,Q_2))(Q)\geq 0$,
because
$\gminia(f,Q_1)\prec_{Q_2}\gminia(f,Q_2)$,
we obtain
$|\gminia(f,Q_2)|^{-1}\Re(\gminia(f,Q_2))(Q_2)>0$,
which contradicts that $f$ is negative at $Q_2$.
Thus, we obtain (\ref{eq;24.1.6.40}).
Similarly,
we obtain $\deg\gminia(f,Q_1)\leq \deg\gminia(f,Q_2)$,
and hence
$\deg\gminia(f,Q_1)=\deg\gminia(f,Q_2)$.

Take $\gminia\in\nbigi$.
Let us prove $\deg(\gminia)=\deg(\gminia(f,Q_i))$.
Suppose $\deg\gminia>\deg\gminia(f,Q_i)$.
If $|\gminia|^{-1}\Re(\gminia)(Q)<0$,
we obtain $\gminia\prec_Q\gminia(f,Q_i)$,
which contradicts that
the expression
$f=\sum_{\gminia\in\nbigi} f_{\gminia}$
is reduced.
If $|\gminia|^{-1}\Re(\gminia)(Q)\geq 0$,
then either
$|\gminia|^{-1}\Re(\gminia)(Q_1)>0$
or
$|\gminia|^{-1}\Re(\gminia)(Q_2)>0$ holds,
which contradicts that $f$ is negative at $Q_i$.
Suppose $\deg\gminia<\deg\gminia(f,Q_i)$.
Because
$|\gminia(f,Q_i)|^{-1}\Re(\gminia(f,Q_i))(Q_i)<0$,
we obtain
$\gminia(f,Q_i)\prec_{Q_i}\gminia$,
which contradicts
$\gminia(f,Q_i)=\max(\nbigi,\prec_{Q_i})$.
In all,
we obtain $\deg(\gminia)=\deg(\gminia(f,Q_i))$.
\hfill\qed

\vspace{.1in}
Recall that we identify 
$\cnumtilde\simeq\real_{\geq 0}\times S^1$
by the polar coordinate system.
For any $a<b$,
we set $\openopen{a}{b}=\{\theta\in\real\,|\,a<\theta<b\}$.
If $b-a<2\pi$,
we may naturally regard it as an interval in $S^1=\real/2\pi\seisuu$.

\begin{df}
 An open interval $I\subset \varpi^{-1}(0)\cap V$ is called special
 with respect to $f$
if the following holds.
\begin{itemize}
 \item The length of $I$ is $\pi/\rho$ for $\rho>1/2$.
       Namely, $I=\{0\}\times\openopen{\theta_1}{\theta_1+\pi/\rho}$
       for some $\theta_1$.
 \item There exists $\alpha\in\cnum\setminus\{0\}$ such that
       (i) $\deg(\gminia(f,Q)-\alpha z^{-\rho})<
       \deg(\gminia(f,Q))=\rho$
       for any $Q\in I\setminus \nbigz(f)$,
       (ii) $\Re(\alpha e^{-\rho\sqrt{-1}\theta})<0$ on $I$,
       which particularly implies that 
\[
       \Re(\alpha e^{-\rho\sqrt{-1}\theta_1})
       =\Re(\alpha e^{-\rho\sqrt{-1}(\theta_1+\pi/\rho)})=0.
\]
\end{itemize}
\end{df}

\begin{lem}
 Let $I$ be an interval which is
 negative with respect to $f$.
 Note that $\rho:=\deg(\gminia(f,Q))$ $(Q\in I\setminus \nbigz(f))$
 is well defined.
 Then, either one of the following holds.
 \begin{itemize}
  \item $I$ is special with respect to $f$.
  \item The length of $I$ is strictly smaller than $\pi/\rho$.
 \end{itemize}
\end{lem}
\pf
There exists an interval
$\Itilde=\{\psi_0\leq \theta\leq\psi_1\}\subset\real$
such that
the map $\theta\longmapsto e^{\sqrt{-1}\theta}$
induces a diffeomorphism $\Itilde\simeq I$.
Let $\psi_0<\theta_1<\cdots<\theta_{\ell}<\psi_1$
be determined by
$\{(0,e^{\sqrt{-1}\theta_i})\}=I\cap\nbigz(f)$.
We set
$\theta_0:=\psi_0$
and $\theta_{\ell+1}:=\psi_1$.
We set $\Itilde_i=\{0\}\times\openopen{\theta_{i-1}}{\theta_i}
\subset \Itilde$
$(i=1,\ldots,\ell+1)$.
We obtain the corresponding subsets $I_i\subset I$.
Choose $Q_i\in I_i$.
There exist $\alpha_i\in\cnum^{\ast}$
such that
$\deg(\gminia(f,Q_i)-\alpha_iz^{-\rho})
<\deg(\gminia(f,Q_i))$.
There exist intervals
$J_i=\openopen{\varphi_i}{\varphi_i+\pi/\rho}$
such that
(i) $\varphi_i\leq \theta_{i-1}<\theta_i\leq\varphi_i+\pi/\rho$,
(ii) $\Re(\alpha_ie^{-\rho\sqrt{-1}\theta})<0$ on $J_i$.

For $1\leq i\leq \ell$,
we obtain
$\Re(\alpha_ie^{-\rho\sqrt{-1}\theta_i})
=\Re(\alpha_{i+1}e^{-\rho\sqrt{-1}\theta_{i}})<0$.
We also obtain
\[
\Re(\alpha_ie^{-\rho\sqrt{-1}(\theta_i+\epsilon)})
<\Re(\alpha_{i+1}e^{-\rho\sqrt{-1}(\theta_{i}+\epsilon)}),
\]
\[
\Re(\alpha_ie^{-\rho\sqrt{-1}(\theta_i-\epsilon)})
>\Re(\alpha_{i+1}e^{-\rho\sqrt{-1}(\theta_{i}-\epsilon)})
\]
for any sufficiently small positive number $\epsilon$.
Therefore, we obtain
$\varphi_{i+1}\leq \varphi_i$.
Moreover, if $\varphi_{i+1}=\varphi_i$,
there exists $a>0$ such that
$\alpha_{i+1}=a\alpha_i$.
Because 
$\Re(\alpha_ie^{-\rho\sqrt{-1}\theta_i})
=\Re(\alpha_{i+1}e^{-\rho\sqrt{-1}\theta_{i}})<0$,
we obtain $\alpha_{i+1}=\alpha_{i}$.
Therefore, we obtain
\[
 \theta_{\ell+1}\leq \varphi_{\ell+1}+\pi/\rho
 \leq\cdots
 \leq \varphi_1+\pi/\rho\leq\theta_0+\pi/\rho.
\]
Namely,
we obtain
$\theta_{\ell+1}-\theta_0\leq \pi/\rho$.
If $\theta_{\ell+1}=\theta_0+\pi/\rho$,
we obtain
$\alpha_1=\alpha_2=\cdots=\alpha_{\ell+1}$.
\hfill\qed

\vspace{.1in}

We also obtain the following lemma.

\begin{lem}
 Suppose that $V$ is a neighbourhood of
 $\varpi^{-1}(0)$,
 and that
 $f$ is non-positive at each point of
 $\varpi^{-1}(0)\setminus\nbigz(f)$.
 Then, $f$ is meromorphic at $0$.
\end{lem}
\pf
If $f$ is neutral at a point of
$\varpi^{-1}(0)\setminus\nbigz(f)$,
we obtain that
$f$ is neutral at any point of 
$\varpi^{-1}(0)\setminus\nbigz(f)$
according to Lemma \ref{lem;20.8.18.20},
and hence
we obtain that $f$ is meromorphic.
Suppose that $f$ is simply negative at any point of
$\varpi^{-1}(0)\setminus\nbigz(f)$.
By Lemma \ref{lem;20.6.22.3},
there exist $\epsilon>0$ and $\rho>0$
such that
$|f|=O\bigl(e^{-\epsilon|z|^{-\rho}}\bigr)$
on $V\setminus\varpi^{-1}(0)$,
which implies $f=0$.
\hfill\qed

\subsubsection{Estimates in a single case}

For any $R>0$ and $0<L<\pi$,
we set
$W(R,L):=
\bigl\{w\in\cnum\,\big|\,|w|>R,\,|\arg(w)|<L\bigr\}$.
Let $\projtilde^1_{\infty}\lrarr\proj^1$
be the oriented real blowing up
at $\infty$.
We regard $W(R,L)$ as an open subset of $\projtilde^1_{\infty}$.
Let $\Wbar(R,L)$ denote the closure of
$W(R,L)$ in $\projtilde^1_{\infty}$.
Let $Q_0\in\Wbar(R,L)$
denote the point corresponding to
$+\infty$,
i.e.,
$Q_0$ is the limit of
any sequence of positive numbers $t_i$ in $\projtilde^1_{\infty}$
with $t_i\to\infty$.
Let $\Phi_{\beta,b}:W(R,L)\lrarr \cnum^{\ast}$
be given by
$\Phi_{\beta,b}(w)=\beta w^{-b}$
for a positive number $b$
and a non-zero complex number $\beta$.
It induces $\Phitilde_{\beta,b}:\Wbar(R,L)\lrarr \cnumtilde$.

Let $f$ be a section of $\gbigb_{\cnumtilde}$
on an open subset $V\subset\cnumtilde$.
Suppose that $f$ has a single growth order
at $Q\in V\cap\varpi^{-1}(0)$.
Let $0<L<\Ltilde<\pi$.
By choosing $\beta\in\cnum^{\ast}$,
$b>0$ and $R>0$ appropriately,
we obtain
$\Phitilde_{\beta,b}:W(R,\Ltilde)\to \cnumtilde$
such that
$\Phitilde_{\beta,b}(Q_0)=Q$
and 
$\Phitilde_{\beta,b}(\Wbar(R,\Ltilde))
\subset V$.
We obtain the holomorphic function
$\Phi_{b,\beta}^{\ast}(f)$ on $W(R,\Ltilde)$
which induces a section of
$\gbigb_{\projtilde^1_{\infty}}$
on $\Wbar(R,\Ltilde)$.
If $b$ is sufficiently small
and $R$ is sufficiently large,
\[
 e^{-\gminia(f,Q)}z^{-a(f,Q)}(\log z)^{-j(f,Q)}f
\]
is regularly bounded
on $\Phitilde_{\beta,b}(\Wbar(R,\Ltilde))\setminus\varpi^{-1}(0)$.

\begin{lem}
\label{lem;20.6.22.40}
For any $0<L_2<L_1<L$,
there exist
  $R_2>R_1>R$ and
 subsets 
 $\nbigc_1\subset\nbigc_2\subset W(R,L)$
such that the following holds.
\begin{itemize}
  \item There exists $C>0$ such that
	\[
	\bigl|
	 \Phi_{\beta,b}^{\ast}(f)
	\bigr|
	\geq
	Ce^{\Re\Phi_{\beta,b}^{\ast}\gminia(f,Q)}
	|w|^{ba(f,Q)}(\log |w|)^{j(f,Q)}
	\]
	on $W(R,L)\setminus \nbigc_1$.
  \item Let $\nbigd$ be any connected component of
	$\nbigc_2$ such that
	\[
	 \nbigd\cap W(R_2,L_2)\neq\emptyset.
	\]
	Then, $\nbigd$ is relatively compact in
	$W(R_1,L_1)$.

  \item For any $w_0\in
	W(R_1,L_1)\setminus\nbigc_2$,
	we obtain
	$\{|w-w_0|<1\}\subset
	W(R,L)\setminus\nbigc_1$.
 \end{itemize}
 \end{lem}
\pf
It is enough to prove the case
where 
$\gminia(f,Q)=0$,
$a(f,Q)=0$
and $j(f,Q)=0$,
i.e.,
$f$ is regularly bounded at $Q$.

For any $B_1\in\real$ and $B_2>0$,
we set
$\Wtilde(B_1,B_2):=
\{\zeta\in\cnum\,|\,\Re(\zeta)>B_1,\,|\Image(\zeta)|<B_2\}$.
Let $\Phi_1:\cnum\lrarr \cnum$
be given by
$\Phi_1(\zeta)=e^{\zeta}$.
It induces
$\Wtilde(B_1,B_2)\simeq W(e^{B_1},B_2)$
for any $B_1,B_2>0$.
We set
$\ftilde:=\Phi_1^{\ast}\Phi_{\beta,b}^{\ast}(f)$
on $W(\log R,\Ltilde)$.
Let $Z(\ftilde)$
denote the zero set of $\ftilde$.
By Proposition \ref{prop;20.6.21.20},
there exists $N>0$
such that the following holds for any
$x_1>\log R$:
\[
\#\Bigl(
 Z(\ftilde)
 \cap
 \bigl\{\zeta\in\cnum\,\big|\,
 x_1< \Re(\zeta)<x_1+L,\,\,
 |\Image(\zeta)|<L
 \bigr\}
 \Bigr)
 <N.
\]
Take
\[
0< \delta<
\frac{1}{10^3N}
 \min\bigl\{
  |L-L_1|,|L_1-L_2|,|L_2|
 \bigr\}.
\]
We set
\[
 \nbigb_1:=\bigcup_{\zeta_1\in Z(\ftilde)}
 \{|\zeta-\zeta_1|< \delta\}
 \cap \Wtilde(\log R,L),
\]
\[
  \nbigb_2:=\bigcup_{\zeta_1\in Z(\ftilde)}
  \{|\zeta-\zeta_1|<3\delta\}
  \cap \Wtilde(\log R,L).
\]
There exists $C>0$ such that
$|\ftilde|>C$
on $\Wtilde(\log R,L)\setminus \nbigb_1$
by Corollary \ref{cor;24.1.5.10}.

We set $\nbigc_i=\Phi_1(\nbigb_i)\subset W(R,L)$.
Then, the first condition is satisfied.
Take
$R_1>\max\{Re^L,10\delta^{-1}\}$.
For any $\zeta_1,\zeta_2\in \Wtilde(\log R)$,
we have
\[
|\zeta_1-\zeta_2|\leq
 |e^{\zeta_1}-e^{\zeta_2}|
 \max\{e^{-\Re(\zeta_1)},e^{-\Re(\zeta_2)}\}.
\]
Then, the third condition is satisfied.
Take $R_2>R_1e^L$.
Let $\nbigdtilde$ be a connected component
of $\nbigb_2$
such that
$\nbigdtilde\cap\Wtilde(\log R_2,L_2)\neq\emptyset$.
Then,
there exists
a subset $Z_{\nbigdtilde}(\ftilde)\subset Z(\ftilde)$
such that
$\nbigdtilde$ is the union of
$\{|\zeta-\zeta_j|<3\delta\}$ $(\zeta_j\in Z_{\nbigdtilde}(\ftilde))$.
Let us observe that 
$\# Z_{\nbigdtilde}(\ftilde)\leq N$.
Take $\zeta_0\in Z_{\nbigdtilde}(\ftilde)$.
If $\# Z_{\nbigdtilde}(\ftilde)> N$,
there exists
$\eta\in Z_{\nbigdtilde}(\ftilde)$
such that
$|\Re(\eta-\zeta_0)|>L/2$.
There exists a sequence
$\zeta_1,\ldots,\zeta_m,\zeta_{m+1}=\eta
\in Z_{\nbigdtilde}(\ftilde)$
such that
$d(\zeta_i,\zeta_{i+1})<3\delta$
for $i=0,\ldots,m$.
There exists $i_0$ such that
$|\Re(\zeta_0-\zeta_i)|\leq L/2$ for $i<i_0$
and $|\Re(\zeta_0-\zeta_{i_0})|>L/2$.
Note that $i_0\leq N$.
Because of the choice of $\delta$,
we obtain $d(\zeta_0,\zeta_{i_0})<L/2$.
Thus, we arrive at a contradiction,
and hence we obtain
$\# Z_{\nbigdtilde}(\ftilde)\leq N$.
Then, we obtain that
$\nbigdtilde$ is relatively compact in $\Wtilde(\log R_1,L_1)$,
and hence the second condition is satisfied.
\hfill\qed

\vspace{.1in}

We set $g=w+\alpha w^c$ on $W(R,\Ltilde)$
for $\alpha\in\cnum\setminus\{0\}$ and $0<c<1$.
For any $B^{(1)}\in\real$ and $0<L^{(1)}<L$,
there exists $R^{(1)}\geq R$
such that
\begin{multline}
S(g,R^{(1)},B^{(1)},L^{(1)})= \\
\bigl\{w\in W(R,\Ltilde)\,\big|\,
|g(w)|>R^{(1)},\,\,\Image(g(w))\geq B^{(1)},\,\,
\arg(g(w))<L^{(1)}
\bigr\} \\
\subset W(R,L).
\end{multline}
Let $Z(\Phi_{b,\beta}^{\ast}(f))$ denote the zero set of
$\Phi_{b,\beta}^{\ast}(f)$.

\begin{lem}
\label{lem;20.6.22.30}
 There exists $N>0$ such that
 the following holds for any $M>R^{(1)}e^L$:
\[
\#\Bigl(
 Z(\Phi_{b,\beta}^{\ast}(f))
 \cap
 \bigl\{
 w\in
 S(g,R^{(1)},B^{(1)},L^{(1)})\,\big|\,
 M<|g(w)|<Me^L
 \bigr\}
 \Bigr)
<N.
\]
\end{lem}
\pf
There exist $M_0>R$ and $\epsilon>0$ such that
the following holds for any $M>M_0$:
\begin{multline}
  \bigl\{
 w\in
 S(g,R^{(1)},B^{(1)},L^{(1)})\,\big|\,
 M<|g(w)|<Me^L
 \bigr\}
 \subset
 \\
 \bigl\{
 w\in W(R,L)\,\big|\,
  Me^{-\epsilon}<|w|<Me^{L+\epsilon}
 \bigr\}.
\end{multline}
Then, the claim follows from Proposition \ref{prop;20.6.21.20}.
\hfill\qed

\begin{lem}
\label{lem;20.6.22.31}
For any $\rho>0$ and $\kappa>0$,
 there exist
$R_0^{(1)}>0$,
$C_1>0$ and $\kappa_1>0$
such that
\[
 e^{-\Re\Phi_{\beta,b}^{\ast}\gminia(f,Q)}
 |w|^{-ba(f,Q)}
 (\log |w|)^{-j(f,Q)}
 \bigl|
 \Phi_{\beta,b}^{\ast}(f)(w)
 \bigr|\geq
 C_{1}e^{-\kappa_1|w|^{\rho}}
\]
 for any $w\in S(g,R_0^{(1)},B^{(1)},L^{(1)})$ such that
 $d(w,Z(\Phi_{\beta,b}^{\ast}(f)))\geq
  \frac{1}{2}e^{-\kappa|w|^{\rho}}$.
\end{lem}
\pf
We have only to study the case where
$\gminia(f,Q)=0$,
$a(f,Q)=0$
and $j(f,Q)=0$.
We use the notation in the proof of
Lemma \ref{lem;20.6.22.40}.

For $\zeta_1,\zeta_2\in\cnum$
with $|\zeta_1-\zeta_2|\leq 1$,
we have
$|e^{\zeta_1}-e^{\zeta_2}|
 \leq
 e^{\Re(\zeta_1)+1}|\zeta_1-\zeta_2|$.
Hence, there exists $C_{2}>0$ such that
the following holds for any 
$w\in S(g,R^{(1)},B^{(1)},L^{(1)})$
satisfying 
$d(w,Z(\Phi_{\beta,b}^{\ast}(f)))\geq\frac{1}{2}e^{-\kappa|w|^{\rho}}$.
\begin{itemize}
 \item Take $\zeta\in \Wtilde(\log R,L)$ such that
       $e^{\zeta}=w$.
       Then, we obtain
       $d(\zeta,Z(\ftilde))>C_{2}|w|^{-1}e^{-\kappa|w|^{\rho}}$.
\end{itemize}
By Corollary \ref{cor;24.1.5.10},
there exist constants 
$C_i>0$ $(i=3,4,5,6)$
such that
the following holds for any
$w\in S(g,R^{(1)},B^{(1)},L^{(1)})$
satisfying 
$d(w,Z(\Phi_{\beta,b}^{\ast}(f)))\geq\frac{1}{2}e^{-\kappa|w|^{\rho}}$:
\[
\bigl|
 \Phi_{\beta,b}^{\ast}(f)(w)
 \bigr|
 \geq C_3(|w|^{-1}e^{-\kappa|w|^{\rho}})^{C_4}
 \geq C_5e^{-C_6\kappa|w|^{\rho}}.
\]
Thus, we obtain the claim of the lemma.
\hfill\qed

\begin{lem}
\label{lem;20.6.24.1}
 For any $B^{(1)}_1>B^{(1)}$,
 $0<L^{(1)}_1<L^{(1)}$
 and  $\rho>0$,
 there exist
 $R^{(1)}_1>R^{(1)}$,
 and a subset
 $\nbigc\subset S(g,R^{(1)},B^{(1)},L^{(1)})$
 such that the following holds:
 \begin{itemize}
  \item There exist $C_i>0$ $(i=1,2)$ such that
\begin{multline}
   \label{eq;20.8.4.10}
        e^{-\Re(\Phi_{\beta,b}^{\ast}\gminia(f,Q))}
	|w|^{-ba(f,Q)}
	(\log|w|)^{-j(f,Q)}
	\bigl|\Phi_{\beta,b}^{\ast}(f)(w)\bigr|
	\geq \\
	C_1\exp(-C_2|w|^{\rho}) 
\end{multline}
	on $S(g,R^{(1)},B^{(1)},L^{(1)})\setminus\nbigc$.
  \item
       Let $\nbigd$ be a connected component of
       $\nbigc$ such that
\[
       \nbigd\cap S(g,R_1^{(1)},B^{(1)}_1,L^{(1)}_1)\neq\emptyset.
\]
       Then, $\nbigd$ is relatively compact in
       $S(g,R^{(1)},B^{(1)},L^{(1)})$.
  \item There exist an increasing sequence
	 of positive numbers $T_i$
	such that (i) $\lim T_i=\infty$,
	(ii) $\nbigc\cap \{|g|=T_i\}\cap
         S(g,R^{(1)},B^{(1)},L^{(1)})=\emptyset$.       
\end{itemize}
\end{lem}
\pf
We have only to study the case where
$\gminia(f,Q)=0$,
$a(f,Q)=0$
and $j(f,Q)=0$.
Let $N$ be as in Lemma \ref{lem;20.6.22.30}.
Take $\kappa>0$ satisfying
{\small
\begin{multline}
 10^3N
 e^{-\kappa (R^{(1)})^{\rho}}
 <\\
 \min\Bigl\{
 |B^{(1)}_1-B^{(1)}|,
 e^{R^{(1)}}\!\!\!\sin\bigl((L^{(1)}-L_1^{(1)})/2\bigr),
  e^{R^{(1)}}\!\!\!\sin(L_1^{(1)}/2),
  R^{(1)}(e^L-1)
 \Bigr\}.
\end{multline}
}
We set
\[
 \nbigc:=\bigcup_{w_1\in Z(\Phi_{\beta,b}^{\ast}(f))}
 \bigl\{
 |w-w_1|<e^{-\kappa|w_1|^{\rho}}
 \bigr\}
 \cap
 S(g,R^{(1)},B^{(1)},L^{(1)}).
\]

For $0<\nu<1$
and for $w,a\in\cnum$ with $|w|>1$ and $|a|<\nu$,
we have
\[
\bigl|
 |w+a|^{\rho}-|w|^{\rho}
\bigr|
\leq
 \rho\nu\left(
 |w|+\nu
 \right)^{\rho-1}
+\rho\nu\left(
 |w|-\nu
 \right)^{\rho-1}.
\]

Let $R_0^{(1)}\geq R^{(1)}$ be as in Lemma \ref{lem;20.6.22.31}.
There exists $R_2^{(1)}>R_0^{(1)}e^L+1$
such that
the following holds
for any $w\in\cnum$ with $|w|>R_1^{(1)}$:
\[
 \frac{1}{2}
 \exp\Bigl(
 \kappa \rho
 e^{-\kappa |w|^{\rho}}
 \bigl(
(|w|+1)^{\rho-1}
+(|w|-1)^{\rho-1}
+2|w|^{\rho-1}
 \bigr)
 \Bigr)
 \leq 3/4.
\]
For $w_2,w_3\in\cnum$ satisfying $|w_i|>R_2^{(1)}$
satisfying
$|w_2-w_3|
 \leq
  \frac{1}{2}
 e^{-\kappa|w_2|^{\rho}}$,
we obtain
$|w_2-w_3|
\leq
 (3/4)e^{-\kappa|w_3|^{\rho}}$.
Hence,
by Lemma \ref{lem;20.6.22.31},
there exist $C'_i>0$ $(i=1,2)$ such that
\begin{multline}
        e^{-\Re(\Phi_{\beta,b}^{\ast}\gminia(f,Q))}
	|w|^{-ba(f,Q)}
	(\log|w|)^{-j(f,Q)}
	\bigl|\Phi_{\beta,b}^{\ast}(f)(w)\bigr|
	\geq \\
	C'_1\exp(-C_2'|w|^{\rho}) 
\end{multline}
on $S(g,R_2^{(1)},B^{(1)},L^{(1)})\setminus\nbigc$.
Because
\[
S(g,R^{(1)},B^{(1)},L^{(1)})
\setminus
S(g,R^{(1)}_1,B^{(1)},L^{(1)})
\]
is relatively compact,
we obtain the first claim.
Take $R^{(1)}_1>R^{(1)}e^L$.
We can check the second claim
by using the argument in the proof of
Lemma \ref{lem;20.6.22.40}.
By Lemma \ref{lem;20.6.22.30},
there exists an infinite sequence $T_i$
as desired.
\hfill\qed

\subsubsection{Estimates in a multiple case}

Let $f$ be a section of $\gbigb_{\cnumtilde}$
on an open subset $V\subset\cnumtilde$.
Let $Q\in V\cap\varpi^{-1}(0)$.
\begin{lem}
\label{lem;20.6.23.30}
 There exist
 $\gminia\in\gbigi(Q)=\bigoplus_{b>0}\cnum z^{-b}$,
 positive constants $C_i$ $(i=1,2)$
 and
 a neighbourhood $\nbigu$ of $Q$ in $V$
 such that 
\[
 |f|\leq
 C_1|z|^{-C_2}e^{\Re\gminia}
\]
on $\nbigu\setminus\varpi^{-1}(0)$.
If $f$ has an expression $f=\sum f_j$
such that each $f_j$ is simply negative,
then we may assume $\gminia\prec_Q0$. 
\end{lem}
\pf
Let $f=\sum f_{j}$ be an expression
such that each $f_j$ has a single growth order at $Q$.
There exists $\gminia\in\gbigi(Q)$
such that
$\gminia(f_j,Q)\prec_Q\gminia$ for any $j$.
There exists $C_2>0$ such that
$C_2>\max\{\gminia(f_j,Q)\}$ for any $j$.
(See Definition \ref{df;20.9.11.10}
for $\gminia(f_j,Q)$.)
Then, there exist $C_3>0$
and a neighbourhood $\nbigu$ of $Q$ in $V$
such that
$|f_j|\leq C_3|z|^{-C_2}e^{\Re(\gminia)}$
on $\nbigu\setminus\varpi^{-1}(0)$
for any $j$.
Then, we obtain the claim of the lemma.
\hfill\qed

\vspace{.1in}

Let us state a lower bound of $f$ around $Q$
outside an exceptional subset.
We identify $\cnumtilde=\real_{\geq 0}\times S^1$
by the polar decomposition
$z=re^{\sqrt{-1}\theta}$.
We choose $\theta_0$ such that
$Q$ is expressed as $(0,e^{\sqrt{-1}\theta_0})$.
Suppose that
$f$ is expressed as a sum $\sum_{i=1}^m f_i$ on 
\[
 V_0:=
 \Bigl\{(r,e^{\sqrt{-1}\theta})\,\Big|\,
 |\theta-\theta_0|< \frac{\pi}{2}\kappa_0,\,
 \,\,0<r< r_0
  \Bigr\}\subset V,
\]
where
$f_i$ are holomorphic functions
with a single growth order,
and $\kappa_0$ and $r_0$ are positive constants.

\begin{lem}
\label{lem;20.6.23.31}
 Take $0<\kappa<\kappa_0$ such that
 $\max\{\kappa \deg\gminia(f_i,Q)\}<1$.
 There exist a neighbourhood $\nbigu$ of $Q$ in $\cnumtilde$
 and a subset $Z\subset \nbigu\setminus\varpi^{-1}(0)$
 such that the following holds.
\begin{itemize}
 \item There exists a subset $Z_1\subset \real_{>0}$
       with $\int_{Z_1}dt/t<\infty$
       such that
       $Z\subset Z_1\times S^1$.
 \item For any $\delta>0$, there exists $C>0$ such that
       the following holds on $\nbigu\setminus (\varpi^{-1}(0)\cup Z)$:
\[
       |f|\geq
       C\exp\Bigl(-\delta|z|^{-1/\kappa}\Bigr).
\]       
\end{itemize}
\end{lem}
\pf
We may assume that $\theta_0=0$.
We set
\[
 W_2:=\{w\in\cnum\,|\,\Re(w)\geq 0\}.
\]
We define the map
$\Psi:W_2\lrarr \cnum$
by $\Psi(w)=(w+C_1)^{-\kappa}$
for some sufficiently large $C_1>0$
such that the image of $\Psi$ is contained in $V_0$.
We obtain the expression
\[
 \Psi^{\ast}(f)
 =\sum_{i=1}^m
  \Psi^{\ast}(f_i)
\]
on $W_2$.
There exist $C_2>0$ such that
\[
 \bigl|
  \Psi^{\ast}(f_i)
 \bigr|
  \leq
  C_2
  e^{\Re(\Psi^{\ast}\gminia(f_i,Q))}
  |w+C_1|^{\kappa\cdot a(f_i,Q)}
  (\log|w+C_1|)^{j(f_i,Q)}
\]
on $W_2$.
Take $b$ such that
$\max\{\kappa\deg(\gminia(f_i,Q))\}
<b<1$.
There exists $A>0$ such that
the following holds on $W_2$:
\begin{multline}
\Re(\Psi^{\ast}\gminia(f_i,Q))
+\log\bigl(
  |w+C_1|^{\kappa\cdot a(f_i,Q)}
  (\log|w+C_1|)^{j(f_i,Q)}
  \bigr)
\\
 \leq
 A\Re\bigl((w+C_1)^b\bigr). 
\end{multline}
Hence, there exists $B>0$ such that
\begin{equation}
\label{eq;24.1.5.11}
 \log |\Psi^{\ast}(f)|\leq
   A\Re\Bigl(
   (w+C_1)^{b}
   \Bigr)
   +B.
\end{equation}
Let $a_1,a_2,\ldots$  be the zeroes of
$\Psi^{\ast}(f)$ in $W_2$.
According to \cite[\S14.2, Theorem 2]{Entire-functions},
we obtain
\[
 \sum_{j=1}^{\infty}
 \frac{\Re(a_j)}{1+|a_j|^2}<\infty.
\]
According to \cite[\S14.2, Theorem 3]{Entire-functions},
we obtain the following description of
$\log |\Psi^{\ast}(f)|$:
\begin{multline}
 \log |\Psi^{\ast}(f)(w)|
 =\sum_{j=1}^{\infty}\log\left|
  \frac{w-a_j}{w+\overline{a}_j}
 \right|
 \\
  +
  \frac{\Re(w)}{\pi}
  \int_{-\infty}^{\infty}
   \frac{d\nu_1(t)}{|t-\sqrt{-1}w|^2}
 -\frac{\Re(w)}{\pi}
   \int_{-\infty}^{\infty}
   \frac{d\nu_2(t)}{|t-\sqrt{-1}w|^2}
 +\sigma\Re(w).
\end{multline}
Here, $\sigma$ is a real number,
and $\nu_i$ are non-negative measures on $\real$
such that
\[
 \int_{-\infty}^{\infty}
  \frac{d\nu_i(t)}{1+t^2}<\infty.
\]
According to the Hayman theorem
\cite[\S15, Theorem 1]{Entire-functions},
there exist
a sequence $w_1,w_2,\ldots$
in $\{w\in\cnum\,|\,\Re(w)>0\}$
and a sequence of positive numbers
$\rho_1,\rho_2,\ldots$
such that $\sum \rho_j/|w_j|<\infty$
and that
\begin{multline}
\left|
\sum_{j=1}^{\infty}\log\left|
  \frac{w-a_j}{w+\overline{a}_j}
  \right|
 \right|
  +
\left|
  \frac{\Re(w)}{\pi}
  \int_{-\infty}^{\infty}
   \frac{d\nu_1(t)}{|t-\sqrt{-1}w|^2}
 \right|
\\ 
+\left|
\frac{\Re(w)}{\pi}
   \int_{-\infty}^{\infty}
   \frac{d\nu_2(t)}{|t-\sqrt{-1}w|^2}
   \right|
=o(|w|)
\end{multline}
outside of $\bigcup\bigl\{w\,\big|\,|w-w_j|<\rho_j\bigr\}$.
By (\ref{eq;24.1.5.11}),
we obtain that $\sigma=0$,
and hence
$\log\bigl|\Psi^{\ast}(f)\bigr|=o\bigl(|w|\bigr)$
outside of
$\bigcup_{j=1}^{\infty}\bigl\{w\,\big|\,|w-w_j|<\rho_j\bigr\}$.
We note that
$\{|w-w_j|<\rho_j\}
\subset
\{|w_j|-\rho_j<|w|<|w_j|+\rho_j\}$.
Because $\sum \rho_j/|w_j|<\infty$,
we have
$\rho_j/|w_j|\to 0$ as $j\to\infty$,
and there exists $A>0$ such that
\[
 \int_{|w_j|-\rho_j}^{|w_j|+\rho_j}dt/t
 =\log\left(
 \frac{|w_j|+\rho_j}{|w_j|-\rho_j}
 \right)
 \leq A\frac{\rho_j}{|w_j|}.
\]
Then, the claim of the lemma follows. 
\hfill\qed

\section{Cyclic Higgs bundles with multiple growth orders
 on sectors}

\subsection{Statements}
\label{subsection;20.6.14.1} 

Let $U$ be a neighbourhood of $0$ in $\cnum$.
Let $\varpi:\Utilde\lrarr U$ be the oriented real blowing up.
Let $V$ be an open subset of $\Utilde$.
Let $f$ be a section of $\gbigb_{\Utilde}$ on $V$.
We set $q=f\cdot (dz)^r$.

\subsubsection{Positive intervals}

Let $I\subset\varpi^{-1}(0)\cap V$
be a positive interval with respect to $f$.
\begin{prop}
 \label{prop;20.6.12.20}
 Let $h_1,h_2\in\Harm(q)$.
 For any relatively compact subset $K\subset I$,
 there exists a neighbourhood $\nbigu_K$ of $K$
 in $V$ such that
 $h_1$  and $h_2$ are mutually bounded
 on $\nbigu_K\setminus\varpi^{-1}(0)$.
\end{prop}

\subsubsection{Maximal non-positive but non-special intervals}

Suppose that there exists an interval $I$
in $\varpi^{-1}(0)\cap V$
satisfying the following conditions.
\begin{itemize}
 \item $I$ is maximally non-positive but not special
       with respect to $f$.
 \item $\Ibar\subset V$.
\end{itemize}
 
\begin{thm}
\label{thm;20.6.9.21}
Let $h_1,h_2\in\Harm(q)$.
There exists a neighbourhood $\nbigu$ of $\Ibar$ in $V$
such that 
 $h_1$ and $h_2$ are mutually bounded
 on $\nbigu\setminus\varpi^{-1}(0)$.
\end{thm}

\subsubsection{Special intervals}

Suppose that there exists an interval $I$
in $\varpi^{-1}(0)\cap V$
satisfying the following conditions.
\begin{itemize}
 \item $\Ibar\subset V$.
 \item $I$ is special with respect to $f$.
\end{itemize}

There exist a non-zero complex number $\alpha$
and $\rho>0$
such that
$\deg(\gminia(f,Q)-\alpha z^{-\rho})<\rho$
for any $Q\in I\setminus\nbigz(f)$.
Let $\nbigp$ be the set of the tuples
$\veca=(a_1,\ldots,a_r)\in\real^r$
such that
\[
 a_1\geq a_2\geq \cdots\geq a_r\geq a_1-1,
 \quad\quad
 \sum a_i=0.
\]
\begin{thm}
\label{thm;20.6.14.10}  \mbox{{}}
\begin{itemize}
 \item 
For any $h\in\Harm(q)$,
there exist
$\veca_I(h)=(a_{I,i}(h))\in\nbigp$
and $\epsilon>0$
such that the following holds on 
$\bigl\{|\arg(\alpha z^{-\rho})-\pi|<(1-\delta)\pi/2\bigr\}$
for any $\delta>0$:
\[
 \log\bigl|
 (dz)^{(r+1-2i)/2}
 \bigr|_h
 +a_{I,i}(h)\Re\bigl(
 \alpha z^{-\rho}
 \bigr)
 =O\bigl(|z|^{-\rho+\epsilon}\bigr).
\]
The tuple $\veca_I(h)$ is uniquely determined
by the condition.
\item
 For $h_1,h_2\in\Harm(q)$  such that
 $\veca_I(h_1)=\veca_I(h_2)$,
 there exists a neighbourhood $\nbigu$ of $\Ibar$
 in $V$
 such that
 $h_1$ and $h_2$ are mutually bounded
 on $\nbigu\setminus\varpi^{-1}(0)$.
 \item For any $\veca\in\nbigp$,
       there exists $h\in\Harm(q)$ such that
       $\veca_I(h)=\veca$.      
\end{itemize}
\end{thm}

We shall also prove the following auxiliary statement.

\begin{prop}
\label{prop;20.6.14.11}
 There exist relatively compact
 open neighbourhoods $\nbigu_{I,i}$ $(i=1,2)$ of $I$ in $V$,
 an $\real_{\geq 0}$-valued harmonic function
 $\phi_I$ on $\nbigu_{I,1}\setminus\varpi^{-1}(0)$,
 and a smooth exhaustive family
 $\{K_i\}$ of $V\setminus\varpi^{-1}(0)$
 such that the following holds.
\begin{itemize}
 \item $\nbigu_{I,2}\subset\nbigu_{I,1}$.
 \item $\nbigu_{I,i}\cap\varpi^{-1}(0)=I$ holds,
       and $\del(\nbigu_{I,i}\setminus\varpi^{-1}(0))$ is smooth
       in $V\setminus\varpi^{-1}(0)$.
 \item There exists $\epsilon>0$
       such that
       $\phi_I=O\bigl(|z|^{-\rho+\epsilon}\bigr)$
       on
       $\nbigu_{I,1}\setminus\varpi^{-1}(0)$.
 \item Take $h\in\Harm(q)$ and
       $h_i\in \Harm(q_{|K_i})$ such that
       $h_{i|\nbigu_{I,1}\cap\del K_i}=h_{|\nbigu_{I,1}\cap\del K_i}$.
       Let $s_i$ be the automorphism of
       $\hyperk_{\nbigu_{I,1}\cap K_i,r}$
       determined by
       $h_{i|\nbigu_{I,1}\cap K_i}
       =h_{|\nbigu_{I,1}\cap K_i}\cdot s_i$.
       Then, we obtain
\[
       \log\Tr(s_i)\leq \phi_I
\]
       on $\nbigu_{I,2}\cap K_i$.
 \end{itemize}
\end{prop}

\subsection{Some estimates on sectors around infinity}

\subsubsection{Positive and negative cases}

\label{subsection;20.6.23.10}

For any $R>0$ and $0<L<\pi$,
we set
$W(R,L):=
\bigl\{
w\in\cnum\,\big|\,
|w|>R,\,\,|\arg(w)|<L
\bigr\}$.
Let $f$ be a holomorphic function on
a sector $W(R,L)$.
We set $q=f\,(dw)^r$.
Let $Z(f)$ denote the set of the zeroes of $f$.
Let $\rho$ be a positive number such that
$\rho<\min\{1,\pi/2L\}$.
Let $\alpha$ be a non-zero complex number.
Let $0<\kappa_1<\cdots<\kappa_m<\rho$
and $\alpha_i\in\cnum^{\ast}$ $(i=1,\ldots,m)$.
We set 
\[
 \gminia:=\alpha w^{\rho}+\sum_{i=1}^m \alpha_{i}w^{\kappa_i}.
\]
Let $a\in\real$ and $n\in\seisuu_{\geq 0}$.
We impose the following condition on $f$
in this subsection.
\begin{condition}
 \label{condition;20.4.28.100}
 There exists $C_0>0$ such that
 the following holds on $W(R,L)$:
\begin{equation}
\label{eq;20.8.4.2}
	e^{-\Re\gminia(w)}
	|w|^{-a}|\log w|^{-n}
	|f|
	\leq C_0.
\end{equation}
 Moreover, for any $L_2<L_1<L$,
 there exist
 $R_2>R_1>R$
 and subsets
 $\nbigc_1\subset\nbigc_2\subset
 W(R,L)$
 such that the following holds.
 \begin{itemize}
  \item
	Let $\nbigd$ be any connected component of
	$\nbigc_2$ such that
\[
        \nbigd\cap W(R_2,L_2)\neq\emptyset.
\]
	Then, $\nbigd$
	is relatively compact in
	$W(R_1,L_1)$.
  \item
	For any
	$w_0\in W(R_1,L_1)\setminus\nbigc_2$,
	$\{|w-w_0|<1\}$ is contained in
	$W(R,L)\setminus \nbigc_1$.
  \item
        There exists
	$C_1>0$ such that
	the following holds on
	$W(R,L)\setminus\nbigc_1$:
	\begin{equation}
	\label{eq;20.8.4.1}
	e^{-\Re\gminia(w)}
	|w|^{-a}|\log w|^{-n}
	|f|\geq C_1.
	\end{equation}
\end{itemize}
\end{condition}

\begin{lem}
   \label{lem;20.4.28.101}
Suppose that 
$\Re(\alpha e^{\sqrt{-1}\rho\theta})>0$ for $|\theta|\leq L$.
Then, any $h_1,h_2\in\Harm(q)$ are mutually bounded
 on $W(R',L')$
 for any $R'>R$ and $0<L'<L$.
  \end{lem}
\pf
Let $s$ be the automorphism of
$\hyperk_{W(R,L),r}$
determined by
$h_2=h_1\cdot s$.
We set $L_2=L'$.
We take $L_1$ such that
$L_2<L_1<L$.
We take $R_2>R_1>R$
and $\nbigc_1\subset\nbigc_2$
as in Condition {\rm\ref{condition;20.4.28.100}}.
Let $w_0\in W(R_1,L_1)\setminus\nbigc_2$.
On $\{|w-w_0|<1\}$,
the estimates
(\ref{eq;20.8.4.2}) and (\ref{eq;20.8.4.1}) hold,
and we may apply Corollary \ref{cor;20.4.16.11}.
Hence, we obtain the boundedness of
$\Tr(s)$ on $W(R_1,L_1)\setminus\nbigc_2$.
Let $\nbigd$ be a connected component of
$\nbigc_2$ with $\nbigd\cap W(R_2,L_2)\neq\emptyset$.
Then, it is relatively compact in
$W(R_1,L_1)$.
Recall that
$\Tr(s)$ is subharmonic
(see \S\ref{subsection;20.8.16.10}).
By the maximum principle,
we obtain that
$\sup_{\nbigd}\Tr(s)
\leq
\sup_{W(R_1,L_1)\setminus\nbigc_2}\Tr(s)$.
Hence, we obtain the boundedness of
$\Tr(s)$ on $W(R_2,L_2)$.
For any $0<R'<R_2$,
$W(R_2,L_2)\setminus W(R',L_2)$
is relatively compact
in $W(R,L)$.
Hence, we obtain the boundedness of
$\Tr(s)$ on $W(R',L_2)$.
\hfill\qed

 \begin{lem}
\label{lem;20.6.15.3}
Suppose that 
$\Re(\alpha e^{\sqrt{-1}\rho\theta})<0$ for $|\theta|\leq L$.
Then, for any $h\in\Harm(q)$,
and for any $R'>R$ and $0<L'<L$,
there exists $C>0$
such that the following holds
on $W(R',L')$:
 \[
 \Bigl|
 \log
 \bigl|
 (dw)^{(r+1-2i)/2}
 \bigr|_h
 \Bigr|
  \leq
  C\bigl|\Re(\alpha w^{\rho})\bigr|
  =-C\Re(\alpha w^{\rho}).
 \]
 \end{lem}
\pf
Let $\theta(q)$ denote the Higgs field
of $\hyperk_{W(R,L),r}$
associated with $q$.
Take $L'<L''<L$.
There exist $\epsilon>0$, $R_{10}>0$
and the map
$\Psi:W(R_{10},(1+\epsilon)\pi/2)
\lrarr W(R,L)$
defined by
$\Psi(\zeta)=\zeta^{2L''/\pi}$.

For $\zeta_1\in W(2R_{10},\pi/2)$,
we may assume that
there exists
$\{|\zeta-\zeta_1|<1\}\subset
\Psi^{-1}(W(R,L))$.
By Proposition \ref{prop;20.4.15.11},
there exists $C_{10}>0$,
which is independent of $\zeta_1$,
such that 
$|\Psi^{-1}(\theta(q))|_{\Psi^{-1}(h)}\leq C_{10}$
on $\{|\zeta-\zeta_1|<1/2\}$.
Hence, we obtain
$|\Psi^{-1}(\theta(q))|_{\Psi^{-1}(h)}\leq C_{10}$
on $W(2R_{10},\pi/2)$.

For $\zeta_1\in W(R_{10},\pi/2)$
such that $\Re(\zeta_1)>2R_{10}$,
$\bigl\{|\zeta-\zeta_1|<\Re(\zeta_1)-R_{10}
\bigr\}$ is contained in 
$W(R_{10},\pi/2)$.
There exists the isomorphism
$\Phi_{\zeta_1}:\{|\eta|<1\}\simeq
\bigl\{|\zeta-\zeta_1|<\Re(\zeta_1)-R_{10}
\bigr\}$
given by
$\Phi_{\zeta_1}(\eta)=\zeta_1+(\Re\zeta_1-R_{10})\eta$.
By applying Proposition \ref{prop;20.4.15.11} to
$\Phi_{\zeta_1}^{\ast}
\Psi^{\ast}(\hyperk_{X(R,L),r},\theta(q),h)$,
we obtain that there exists
$C_{11}>0$, which is independent of $\zeta_1$,
such that 
$\bigl|
\Phi_{\zeta_1}^{\ast}\Psi^{\ast}(\theta(q))
\bigr|_{\Phi_{\zeta_1}^{\ast}\Psi^{\ast}(h)}
\leq C_{11}$
on $\{|\eta|<1/2\}$.
Because
$\Phi_{\zeta_1}^{\ast}(d\zeta)
=(\Re\zeta_1-R_{10})d\eta$,
we obtain
\[
\Psi^{\ast}(\theta(q))
 =O\bigl((|\Re(\zeta)|+1)^{-1}\bigr)\,d\zeta 
\]
with respect to $\Psi^{-1}(h)$
on $W(R_{10},\pi/2)$.
Recall that $R(\Psi^{\ast}(h))$
denotes the curvature of the Chern connection of
$\Psi^{\ast}(\hyperk_{W(R,L),r},h)$.
The Hitchin equation implies
$R(\Psi^{\ast}(h))
=O\bigl((|\Re(\zeta)|+1)^{-2}\bigr)\,d\zeta\,d\zetabar$
on $W(R_{10},\pi/2)$.
It implies the following estimates with respect to
$\Psi^{-1}(h)$
on $W(R_{10},(1-\kappa)\pi/2)$
for any $\kappa>0$:
\[
\Psi^{\ast}(\theta(q))=O(|\zeta|^{-1})\,d\zeta,
\quad
R(\Psi^{\ast}(h))=O(|\zeta|^{-2})\,d\zeta\,d\zetabar.
\]

Take $L'<L'''<L''$.
Recall that $R(h)$ denotes the curvature of
the Chern connection of
$(\hyperk_{W(R,L),r},h)$.
By the previous consideration,
and by $w=\zeta^{2L''/\pi}$,
we obtain the estimates
$\theta(q)=O(|w|^{-1})\,dw$
and 
$R(h)=O(|w|^{-2})\,dw\,d\wbar$
with respect to $h$
on $W(R,L''')$.
 
We set $L_2:=L'$
and take $L_2<L_1<L'''$.
Take $R_2>R_1>R$
and $\nbigc_1\subset\nbigc_2$
as in Condition {\rm\ref{condition;20.4.28.100}}.
Note that
$\theta(q)=O(|w|^{-1})\,dw$ with respect to $\Psi^{-1}(h)$
on $W(R_1,L_1)$.
At any point of $W(R_1,L_1)\setminus\nbigc_1$,
the estimate (\ref{eq;20.8.4.1}) holds.
By Corollary \ref{cor;20.4.21.2},
there exists $C_{12}>0$
such that the following holds for $i=1,\ldots,r$
on $W(R_1,L_1)\setminus\nbigc_1$:
 \[
 \Bigl|
 \log
 \bigl|
 (dw)^{(r+1-2i)/2}
 \bigr|_h
 \Bigr|
 \leq
 C_{12}\bigl|
  \Re(\alpha w^{\rho})
  \bigr|.
 \]
Because $R(h)=O(|w|^{-2})dw\,d\wbar$
on $W(R_1,L''')$,
we obtain the following estimate on $W(R,L''')$:
\[
 \del_w\del_{\wbar}
\log \bigl|
 (dw)^{(r+1-2i)/2}
 \bigr|_h
 =O(|w|^{-2}).
\]
Hence,
by Lemma \ref{lem;20.4.21.3},
we can prove that
there exist functions
$\beta_i$
on $W(R_1,L_1)$
such that
(i) the functions $\log \bigl|
 (dw)^{(r+1-2i)/2}
 \bigr|_h-\beta_i$ are harmonic,
(ii) $|\beta_i|=O\bigl(\log(|w|+1)\bigr)$.
 There exists $C_{13}>0$ such that
 the following holds on $W(R_1,L_1)\setminus\nbigc_1$:
 \begin{equation}
\label{eq;20.6.23.20}
 -C_{13}\bigl|
 \Re(\alpha w^{\rho})
 \bigr|
\leq
 \log
 \bigl|
 (dw)^{(r+1-2i)/2}
 \bigr|_h
 -\beta_i
 \leq
 C_{13}
 \bigl|\Re(\alpha w^{\rho})
 \bigr|.
 \end{equation}
 Let $\nbigd$ be any connected component of $\nbigc_2$
 such that $\nbigd\cap W(R_2,L_2)\neq\emptyset$.
Because it is relatively compact in $W(R_1,L_1)$,
the inequalities (\ref{eq;20.6.23.20}) hold
on $\nbigd$.
Hence, (\ref{eq;20.6.23.20})
holds on $W(R_2,L_2)$.
Then, we can deduce the claim of the lemma easily.
\hfill\qed

\subsubsection{Neutral case}

We set
\[
 \Wtilde(R,L):=\bigl\{\zeta\in\cnum\,
\big|\,\Re(\zeta)>\log R,\,\,|\Image(\zeta)|< L
\bigr\}.
\]
Let $\Phi_1:\cnum\lrarr \cnum^{\ast}$
be given by $\Phi_1(\zeta)=e^{\zeta}$.
It induces
$\Wtilde(R,L)\simeq W(R,L)$.
Let $a\in\real$ and $n\in\seisuu_{\geq 0}$.
We impose the following condition to $f$.
\begin{condition}
\label{condition;20.6.23.51}
 For any $L_2<L_1<L$,
 there exist
 $R_2>R_1>R$
 and subsets
 $\nbigb_1\subset\nbigb_2\subset
 \Wtilde(R,L)$
 such that the following holds.
 \begin{itemize}
  \item Let $\nbigd$ be any connected component of
	$\nbigb_2$ such that
\[
 	\nbigd\cap \Wtilde(R_2,L_2)\neq\emptyset.
\]
	Then, $\nbigd$
	is relatively compact in
	$\Wtilde(R_1,L_1)$.
  \item For any
	$\zeta_0\in \Wtilde(R_1,L_1)\setminus\nbigb_2$,
	we obtain
	$\{|\zeta-\zeta_0|<1\}\subset
	\Wtilde(R,L)\setminus \nbigb_1$.
  \item There exists
	$C>1$ such that
	the following holds on
	$\Wtilde(R,L)\setminus\nbigb_1$:
	\begin{equation}
	 \label{eq;20.6.23.50}
	C^{-1}\leq
        e^{-a\Re(\zeta)}|\zeta|^{-n}
	|\Phi_1^{\ast}(f)(\zeta)|
	\leq C.
	\end{equation}
 \end{itemize}
\end{condition}

\begin{lem}
 \label{lem;20.6.23.100}
If Condition {\rm\ref{condition;20.6.23.51}} is satisfied,
the following holds.
 \begin{itemize}
  \item If $a\geq -r$,
	any $h_1,h_2\in\Harm(q)$
	are mutually bounded on
	$W(R',L')$ for any $R'>R$ and $0<L'<L$.
  \item Suppose $a<-r$.
	Then, 
	for any $h\in\Harm(q)$,
  and for any $R'>R$ and $0<L'<L$,
there exists $C>0$
	such that the following holds
	for $i=1,\ldots,r$ on $W(R',L')$:
  \begin{equation}
   \label{eq;20.6.23.60}
 \Bigl|
 \log
 \bigl|
 (dw)^{(r+1-2i)/2}
 \bigr|_h
 \Bigr|
	\leq
	C\log|w|.	  
  \end{equation}
 \end{itemize}
\end{lem}
\pf
Let us study the first claim.
We note that
$\Phi_1^{\ast}\bigl(
 f\cdot(dw)^r
 \bigr)
 =\Phi_1^{\ast}(f)\cdot
  e^{r\zeta}(d\zeta)^r$.
The inequalities (\ref{eq;20.6.23.50})
is rewritten as follows:
\[
C^{-1}\leq
 e^{-(a+r)\Re(\zeta)}|\zeta|^{-n}
\bigl|\Phi_1^{\ast}(f)(\zeta)\cdot e^{r\zeta}\bigr|
 \leq
 C.
\]
For $h_1,h_2\in\Harm(q)$,
let $s$ be the automorphism of
$\hyperk_{\Wtilde(R,L),r}$ determined by
$h_2=h_1\cdot s$.
We set $L_2:=L'$,
and take $L_2<L_1<L$.
Let $R_2>R_1>R$
and $\nbigb_1\subset\nbigb_2$
as in Condition {\rm\ref{condition;20.6.23.51}}.
\begin{itemize}
 \item
If $(a+r,n)\in
(\real_{\geq 0}\times\seisuu_{\geq 0})
\setminus\{(0,0)\}$,
we obtain that $\Tr(s)$ is bounded
on $\Wtilde(R_1,L_1)\setminus\nbigb_2$
by Corollary \ref{cor;20.4.16.11}.
 \item
      If $(a+r,n)=(0,0)$,
      we obtain the boundedness of
$\Tr(s)$ on $\Wtilde(R_1,L_1)\setminus\nbigb_2$
from Corollary \ref{cor;20.6.11.1}.
\end{itemize}
By applying the argument in the proof of
Lemma \ref{lem;20.4.28.101},
we obtain that $\Tr(s)$ is bounded
on $\Wtilde(R_2,L_2)$.
Then, we obtain the first claim immediately.

Let us study the second claim.
Let $h\in\Harm(q)$.
By using the argument in the proof of 
Lemma \ref{lem;20.6.15.3},
we can prove that there exists $C_1>0$
such that
the following holds on $\Wtilde(R',L')$:
\begin{equation}
\label{eq;20.6.23.61}
 \Bigl|
 \log
 \bigl|
 (d\zeta)^{(r+1-2i)/2}
 \bigr|_{\Phi_1^{\ast}(h)}
 \Bigr|
 \leq
 C_1\bigl(
  |\Re(\zeta)|+1
  \bigr)   
\end{equation}
Because $\Phi_1^{\ast}(dw/w)=d\zeta$,
we obtain (\ref{eq;20.6.23.60})
from (\ref{eq;20.6.23.61}).
\hfill\qed

\subsubsection{Estimates on sectors around infinity
   in the multiple case}

Let us continue to use the notation in
\S\ref{subsection;20.6.23.10}.

\begin{lem}
 \label{lem;20.6.14.2}
 Suppose that $f$ satisfies the following conditions.
 \begin{itemize}
  \item There exist $C_0>0$ and $\kappa_0>0$ 
       such that $|f|\leq C_0\exp(\kappa_0 |w|^{\rho})$
	on $W(R,L)$.
  \item There exist $C_1>0$, $\kappa_1>0$,
	and a subset $Z_1\subset\real_{>0}$ with
	$\int_{Z_1}dt/t<\infty$,
	such that $|f|\geq C_1\exp(-\kappa_1|w|^{\rho})$
	on $\bigl\{w\in W(R,L)\,\big|\,|w|\not\in Z_1\bigr\}$.
 \end{itemize}
Let $h_1,h_2\in\Harm(q)$.
Let $s$ be the automorphism of
$\hyperk_{W(R,L),r}$
determined by $h_2=h_1s$.
 Then, for any $R'>R$ and $0<L'<L$,
 there exist
 $C>0$ such that
 $\bigl|
 \log\Tr(s)\bigr|\leq C|w|^{\rho}$
 on $\bigl\{w\in W(R',L')\,\big|\,|w|\not\in Z_1\bigr\}$.
 \end{lem}
\pf
By Proposition \ref{prop;20.4.15.11},
there exists a constant $C_{10}>0$
and $\kappa_{10}>0$
such that
the following holds on $W(R',L')$:
\[
 |\theta(q)|_{h_j}\leq
 C_{10}\exp(\kappa_{10}|w|^{\rho}).
\]
By Corollary \ref{cor;20.4.17.10},
there exist
$C_{11}>0$ and  $\kappa_{11}>0$
such that the following holds on
$\bigl\{w\in W(R',L')\,\big|\,|w|\not\in Z_1\bigr\}$:
\[
 |s|_{h_1}\leq
  C_{11}\exp(\kappa_{11}|w|^{\rho}).
\]
Thus, we obtain the claim of the lemma.
\hfill\qed

\subsection{Proof of Proposition \ref{prop;20.6.12.20}
and Theorem \ref{thm;20.6.9.21}}

\subsubsection{Proof of Proposition \ref{prop;20.6.12.20}}

It is enough to prove the following claim
for any $Q\in I$.
\begin{itemize}
 \item
 There exists a relatively compact neighbourhood
 $\nbigu$ of $Q$ in $V$
 such that
 $h_1$ and $h_2$ are mutually bounded
 on $\nbigu\setminus\varpi^{-1}(0)$.
\end{itemize}
If $Q\in I\setminus\nbigz(f)$,
the claim follows from 
Lemma \ref{lem;20.4.28.101}.
Let us study the case where $Q\in I\cap\nbigz(f)$.
We may assume that $I\cap\nbigz(f)=\{Q\}$.
Let $s$ be determined by $h_2=h_1\cdot s$.
We have already proved that
for any $Q'\in I\setminus\{Q\}$
there exists a neighbourhood $\nbigu_{Q'}$
in $V$ such that
$s$ is bounded on $\nbigu_{Q'}\setminus\varpi^{-1}(0)$.
By Lemma \ref{lem;20.6.23.30},
Lemma \ref{lem;20.6.23.31}
and Lemma \ref{lem;20.6.14.2},
there exist $N>0$,
a neighbourhood of $\nbigu$ of $Q$ in $V$,
and a subset
$Z_1\subset\real_{>0}$
with $\int_{Z_1}dt/t<\infty$
such that
$\log\Tr(s)=O(|z|^{-N})$
on $\nbigu\setminus
\bigl(
\varpi^{-1}(0)
\cup(Z_1\times S^1)
\bigr)$.
Then, we obtain the boundedness of
$\log\Tr(s)$
on $\nbigu\setminus\varpi^{-1}(0)$
by Corollary \ref{cor;20.6.12.21}.
\hfill\qed

\subsubsection{Proof of Theorem \ref{thm;20.6.9.21}}
\label{subsection;20.6.24.100}

Let $Q_1,Q_2$ denote the points of $\del I$.
Let $I_i$ be a neighbourhood of $Q_i$
in $\varpi^{-1}(0)\cap V$
such that $I_i\cap\nbigz(f)=\{Q_i\}$.
Then, $f$ is simply positive at each
point of $I_i\setminus \Ibar$.

Let us study the case where
$I$ is negative with respect to $f$.
Then,
there exists $\rho>0$
such that
$\deg(\gminia(f,Q))=\rho$
for any $Q\in I\setminus\nbigz(f)$.
Because $I$ is not special with respect to $f$,
the length of $I$ is strictly less than $\pi/\rho$.

Let $s$ be the automorphism determined by
$h_2=h_1\cdot s$.
According to Proposition \ref{prop;20.6.12.20},
$\Tr(s)$ is bounded
around any point of 
$(I_1\cup I_2)\setminus \Ibar$.
According to Lemma \ref{lem;20.6.15.3},
we obtain 
$\log\Tr(s)=O\Bigl(|z|^{-\rho}\Bigr)$
around any point of
$I\setminus\nbigz(f)$.
Moreover, according to Lemma \ref{lem;20.6.14.2},
there exist $N>0$,
a neighbourhood $\nbigu_0$
of $\nbigz(f)\cap \Ibar$ in $V$,
and a subset $Z_1\subset\real_{>0}$
with $\int_{Z_1}dt/t<\infty$,
such that $\log\Tr(s)=O\Bigl(|z|^{-N}\Bigr)$
on $\nbigu\setminus
\bigl(
\varpi^{-1}(0)\cup
(Z_1\times S^1)
\bigr)$.
By Corollary \ref{cor;20.6.12.21},
we obtain that
$\log\Tr(s)=O(|z|^{-\rho})$
on $\nbigu\setminus\varpi^{-1}(0)$.
Because the length of $I$ is strictly smaller than
$\pi/\rho$,
the Phragm\'{e}n-Lindel\"{o}f theorem
(see Corollary \ref{cor;20.5.2.1})
implies that $\log\Tr(s)$ is bounded.
Thus, we obtain the claim of the proposition
in the case where $I$ is negative with respect to $f$.

We can prove the claim by a similar and easier argument
in the case where $I$ is neutral with respect to $f$
by using Lemma \ref{lem;20.6.23.100}
instead of Lemma \ref{lem;20.6.15.3}.
\hfill\qed

\subsection{Outline of the proof of Theorem \ref{thm;20.6.14.10}
and Proposition \ref{prop;20.6.14.11}}

\subsubsection{Setting in a simple case}
\label{subsection;24.1.5.20}

Let $\varpi_{\infty}:\projtilde^1_{\infty}\lrarr \proj^1$
denote the oriented real blow up at $\infty$.
We identify $\varpi_{\infty}^{-1}(\infty)\simeq S^1$
by the polar decomposition $w=|w|e^{\sqrt{-1}\theta}$.

For $R>0$ and $\delta>0$,
we set
\[
X(R,\delta):=
\bigl\{w\in\cnum\,\big|\,|w|>R,\,|\arg(w)-\pi/2|<
 \pi/2+\delta
\bigr\}.
\]
We regard $X(R,\delta)$ as an open subset of
$\projtilde^1_{\infty}$.
Let $\Xbar(R,\delta)$ denote the closure of
$X(R,\delta)$ in $\projtilde^1_{\infty}$.

For the proof of Theorem \ref{thm;20.6.14.10}
and Proposition \ref{prop;20.6.14.11},
we introduce the coordinate $w$
determined by $\alpha z^{-\rho}=\sqrt{-1}w$,
and we study
harmonic metrics of
the Higgs bundle associated with $q=f(dw)^r$
on $X(R,\delta)$
such that
$\deg(\gminia(f,\theta)-\sqrt{-1}w)<\deg(\gminia(f,\theta))$
for any $0<\theta<\pi$.

In the following of this subsection,
$C_i$ will denote positive constants.

As an outline of the proof  
of Theorem \ref{thm;20.6.14.10}
and Proposition \ref{prop;20.6.14.11},
we explain rather detail of our arguments
by assuming that $f$ satisfies
the following stronger condition.
 \begin{itemize}
  \item There exist $C_i>0$ $(i=1,2,3)$ such that
the following inequalities 
are satisfied on $\{\Image(w)>R\}$:
	\begin{equation}
\label{eq;20.8.6.110}
 |f|\leq C_1e^{-\Image(w)}.
	\end{equation}
\begin{equation}
\label{eq;20.8.6.111}
 |f|\geq C_2e^{-\Image(w)-C_3|w|^{1-\epsilon}}.
\end{equation}
  \item There exist $C_{i}>0$ $(i=4,5,6)$
	such that
\begin{equation}
\label{eq;20.8.6.130}
 C_4e^{-C_4|w|^{C_4}}\leq
 |f|
 \leq
 C_5e^{C_5|w|^{C_5}}
\end{equation}
on
$\bigl\{
|\arg(w)|<C_{6}
\bigr\}
\cup
\bigl\{
|\arg(w)-\pi|<C_{6}
\bigr\}$.
 \end{itemize}
Note that
the interval
$\{\infty\}\times\openopen{0}{\pi}$
is special with respect to $f$,
and hence $f$ is simply positive at
any points of
$\{\infty\}\times
\bigl(
\openopen{-\delta}{0}
\cup
\openopen{\pi}{\pi+\delta}
\bigr)$.

\subsubsection{Parabolic structure}

Let us explain an outline of the proof of
the first claim of Theorem \ref{thm;20.6.14.10}
under the simplified setting in \S\ref{subsection;24.1.5.20}.
(We shall study a general case in \S\ref{subsection;20.8.6.500}.)

Let $h\in\Harm(q)$.
For any $w_1$ with $\Image(w_1)>R$,
there exists an isomorphism
\[
 \Phi_{w_1}:\{|\eta|<1\}\simeq \{|w-w_1|<\Image(w_1)-R\}
\]
defined by
$\Phi_{w_1}(\eta)=w_1+(\Image(w_1)-R)\eta$.
Note that
$\Phi_{w_1}^{\ast}(f(dw)^r)
=\Phi_{w_1}^{\ast}(f)\cdot (\Image(w_1)-R)^r\,(d\eta)^r$.
By (\ref{eq;20.8.6.110}),
there exists $C^{(1)}_1>0$,
which is independent of $w_1$
such that the following holds
on $\{|\eta|<1\}$:
\[
|\Phi_{w_1}^{\ast}(f)|(\Image(w_1)-R)^r\leq C^{(1)}_1.
\]
Hence, by Proposition \ref{prop;20.4.15.11},
there exists $C^{(1)}_2>0$,
which is independent of $w_1$,
such that the following holds
on $\{|\eta|<1/2\}$:
\[
 |\Phi_{w_1}^{\ast}(\theta(q))|_{\Phi_{w_1}^{\ast}(h)}
 \leq C^{(1)}_2.
\]
Because $\Phi_{w_1}^{\ast}(dw)=(\Image(w_1)-R)\,d\eta$,
there exists $C^{(1)}_3>0$
such that the following holds
on $\{\Image(w)>R\}$:
\begin{equation}
\label{eq;20.8.6.100}
 |\theta(q)|_h
 \leq C^{(1)}_3 (\Image w-R)^{-1}.
\end{equation}
By the Hitchin equation,
there exists $C^{(1)}_4>0$
such that the following holds
on $\{\Image(w)>R\}$:
\begin{equation}
\label{eq;20.8.7.1}
 |R(h)|_h
 \leq C^{(1)}_4(\Image w-R)^{-2}.
\end{equation}
By (\ref{eq;20.8.7.1}),
there exists $C^{(1)}_5>0$
such that the following holds
on $\{\Image(w)>R\}$:
\begin{equation}
\label{eq;20.8.6.113}
\Bigl|
 \del_z\del_{\zbar}\log\bigl|(dw)^{(r+1-2i)/2}\bigr|_h
 \Bigr|
 \leq
 C^{(1)}_5(\Image w-R)^{-2}.
\end{equation}

By (\ref{eq;20.8.6.111}),
(\ref{eq;20.8.6.100})
and Corollary \ref{cor;20.4.21.2},
there exists
$C^{(1)}_6$ such that
\begin{equation}
 \label{eq;20.8.6.120}
\Bigl|
 \log
 \bigl|
 (dw)^{(r+1-2i)/2}
 \bigr|_h
 \Bigr|
 \leq C^{(1)}_6(\Image(w)+|w|^{1-\epsilon}).
\end{equation}
By (\ref{eq;20.8.6.113}),
(\ref{eq;20.8.6.120})
and
the Nevanlinna formula (Proposition \ref{prop;20.4.22.11}),
there exists $a_i(h)$ such that
the following holds on
$\{\Image(w)> R\}$:
\begin{equation}
\label{eq;20.8.6.121}
 \log|(dw)^{(r+1-2i)/2}|_h
 -a_i(h)\Image(w)
=O(|w|^{1-\epsilon}).
\end{equation}
Because of (\ref{eq;20.8.6.100})
and $\theta(q)(dw)^{(r+1-2i)/2}
=(dw)^{(r+1-2(i+1))/2}$,
we obtain
$a_i(h)\geq a_{i+1}(h)$.
By (\ref{eq;20.8.6.111}),
(\ref{eq;20.8.6.100}),
and the relation
$\theta(q)(dw)^{(-r+1)/2}
=f(dw)^{(r-1)/2}$,
we also obtain
$a_r(h)\geq a_1(h)-1$.
In this way,
we obtain the parabolic structure of $h$.

\subsubsection{Mutually boundedness}

Let us explain an outline of the proof of
the second claim of Theorem \ref{thm;20.6.14.10}
under the simplified setting in \S\ref{subsection;24.1.5.20}.
(We shall study a general case in \S\ref{subsection;20.8.6.501}.)

Suppose that
$h_i\in\Harm(q)$ $(i=1,2)$ satisfy
$\veca(h_1)=\veca(h_2)$.
Let $s$ be the automorphism of
$\hyperk_{X(R,\delta),r}$
determined by
$h_2=h_1\cdot s$.
By (\ref{eq;20.8.16.6}),
$\log\Tr(s)$ is subharmonic.

Note that
$f$ is simply positive at
any points of
$\{\infty\}\times
\bigl(
\openopen{-\delta}{0}
\cup
\openopen{\pi}{\pi+\delta}
\bigr)$.
Hence, by Proposition \ref{prop;20.6.12.20},
$h_1$ and $h_2$
are mutually bounded
on $X(R_1,\delta_1)\setminus X(R_1,\delta_2)$
for any $R_1>R$
and $0<\delta_2<\delta_1<\delta$.

By the assumption $\veca(h_1)=\veca(h_2)$,
and by (\ref{eq;20.8.6.121}),
we obtain
\[
 \log\Tr(s)=O(|w|^{1-\epsilon})
\]
on $\{\Image(w)>R\}$.
By (\ref{eq;20.8.6.130})
and Proposition \ref{prop;20.4.15.11},
there exists $C^{(2)}_{1}$ 
such that
the following holds
on $\{|\arg(w)|<C_6/2 \}
\cup \{|\arg(w)-\pi|<C_6/2\}$:
\begin{equation}
\label{eq;20.8.6.140}
 |\theta(q)|_{h_j}\leq
 C^{(2)}_1
 e^{C^{(2)}_1|w|^{C_1^{(2)}}}
 \quad(j=1,2).
\end{equation}
By (\ref{eq;20.8.6.130}),
(\ref{eq;20.8.6.140})
and Corollary \ref{cor;20.4.17.10},
there exists $C^{(2)}_2$
such that the following holds
on $\{|\arg(w)|<C_6/2 \}
\cup \{|\arg(w)-\pi|<C_6/2\}$:
\begin{equation}
 \log\Tr(s)\leq
  C^{(2)}_2\bigl(
  1+|w|^{C_2^{(2)}}
  \bigr).
\end{equation}
By using
Phragm\'{e}n-Lindel\"{o}f theorem
(Corollary \ref{cor;20.6.12.21})
on small sectors around 
$\arg(w)=0$ and $\arg(w)=\pi$,
we obtain that there exists
$C_3^{(2)}$ such that
the following holds
on $X(R_1,\delta_1)$:
\[
 \log\Tr(s)\leq
 C_3^{(2)}(1+|w|^{1-\epsilon}).
\]
By using 
Phragm\'{e}n-Lindel\"{o}f theorem
(Corollary \ref{cor;20.5.2.1}) again,
we obtain that
$\log\Tr(s)$ is bounded
on $X(R_1,\delta_3)$
if $\delta_3>0$ is sufficiently small.
Then, we obtain that
$\log\Tr(s)$ is bounded
on $X(R_1,\delta_4)$
for any $0<\delta_4<\delta$,
i.e.
$h_1$ and $h_2$ are mutually bounded
on $X(R_1,\delta_4)$.

\subsubsection{Auxiliary metrics}

We introduce auxiliary metrics
under the simplified setting in \S\ref{subsection;24.1.5.20}.
(We shall study a general
case in \S\ref{subsection;20.8.6.502}.)

We set $q_0=e^{\sqrt{-1}w}(dw)^r$ on $\cnum$.
For any $\veca\in\nbigp$,
there exists $h_0\in\Harm(q_0)$
such that $\veca(h_0)=\veca$.
(See Proposition \ref{prop;20.6.24.200}.)
There exists
$C^{(3)}_1$
such that the following holds
on $\{\Image(w)>R\}$:
\begin{equation}
\label{eq;20.8.6.210}
 |\theta(q_0)|_{h_0}
 \leq C^{(3)}_1 \Image(w)^{-1}.
\end{equation}
Note that
\[
\theta(q)(dw)^{(r+1-2i)/2}
=(dw)^{(r+1-2(i+1))/2}\,dw
=\theta(q_0)(dw)^{(r+1-2i)/2}
\]
for $i=1,\ldots,r-1$.
We also note that
\[
\theta(q)(dw)^{(-r+1)/2}
=f(dw)^{(r-1)/2}\,dw,
\]
\[
\theta(q_0)(dw)^{(-r+1)/2}
=e^{\sqrt{-1}w}(dw)^{(r-1)/2}dw.
\]
We obtain
\[
|\theta(q)|_{h_0}
=O(|\theta(q_0)|_{h_0}).
\]
Hence, there exists $C^{(3)}_2>0$ such that
the following holds on
$\{\Image(w)>R\}$:
\begin{equation}
 \label{eq;20.8.6.200}
 |\theta(q)|_{h_0}
 \leq
  C^{(3)}_2\Image(w)^{-1}.
\end{equation}
By the Hitchin equation for
$(\hyperk_{\cnum,r},\theta(q_0),h_0)$,
there exists $C^{(3)}_3>0$
such that the following holds
on $\{\Image(w)>R\}$:
\[
 |R(h_0)|_{h_0}\leq C^{(3)}_3\Image(w)^{-2}.
\]
We set
$F(h_0,\theta(q)):=
R(h_0)+[\theta(q),\theta(q)^{\dagger}_{h_0}]$.
There exists
$C^{(3)}_4$ such that
the following holds on
$\{\Image(w)>R\}$:
\begin{equation}
\label{eq;20.8.6.250}
 |F(h_0,\theta(q))|_{h_0}
 \leq C_4^{(3)}\Image(w)^{-2}.
\end{equation}

\subsubsection{Comparison with auxiliary metric}

We explain how to obtain an estimate
for the difference between a harmonic metric
and an auxiliary metric.
(We shall study a general case
in \S\ref{subsection;20.8.6.50}.)

Let $Y$ be a relatively compact open subset of $X(R,\delta)$
such that $\del Y$ is smooth,
and that
\begin{multline}
\label{eq;20.8.6.301}
 Y\cap \{2R-1<\Image w<4R+1\}
 =\\
 \{2R-1<\Image w<4R+1,\,|\Re(w)|<a_Y\}
\end{multline}
for some $a_Y>2$.
According to Proposition \ref{prop;20.5.29.20},
there exists $h_Y\in\Harm(q_{|Y})$
such that $h_{Y|\del Y}=h_{0|\del Y}$.
Let $s_Y$ be the automorphism of
$\hyperk_{Y,r}$
determined by $h_Y=h_{0|Y}s_Y$.
In the following $C^{(4)}_i$
are positive constants which are independent of $Y$.

By Proposition \ref{prop;20.4.15.11}
and (\ref{eq;20.8.6.110}),
there exists $C^{(4)}_1$
such that
the following holds
on $\{2R<\Image w<4R,\,|\Re(w)|<a_Y-1\}$:
\begin{equation}
\label{eq;20.8.6.201}
 |\theta(q)|_{h_Y}\leq C_1^{(4)}.
\end{equation}
Note that
$|\theta(q)|_{h_0}
=|\theta(q)|_{h_Y}$
on $\del Y$.
By Proposition \ref{prop;20.6.13.32},
(\ref{eq;20.8.6.110}),
and (\ref{eq;20.8.6.200}),
there exists $C^{(4)}_2$,
such that the following holds
on $\{2R<\Image w<4R,\,a_Y-1<|\Re(w)|<a_Y\}$:
\begin{equation}
\label{eq;20.8.6.202}
 |\theta(q)|_{h_Y}\leq C_2^{(4)}.
\end{equation}
By (\ref{eq;20.8.6.111}),
(\ref{eq;20.8.6.210}),
(\ref{eq;20.8.6.201}),
(\ref{eq;20.8.6.202})
and Corollary \ref{cor;20.4.17.10},
there exists $C^{(4)}_3>0$
such that the following holds
on $\{\Image(w)=3R\}\cap Y$:
\begin{equation}
\label{eq;20.8.6.220}
 \log\Tr(s_Y)\leq
 C_3^{(4)}|w|^{1-\epsilon}.
\end{equation}
There exists a complex number $\gamma$
such that
$\Re(\gamma w^{1-\epsilon})>0$
on $\{\Image(w)>R\}$.
There exists $C^{(4)}_4>0$
such that the following holds
on $\{\Image(w)=3R\}\cap Y$:
\begin{equation}
 \log\Tr(s_Y)\leq
 C_4^{(4)}\Re(\gamma w^{1-\epsilon}).
\end{equation}
By (\ref{eq;20.8.16.6}) and
(\ref{eq;20.8.6.250}),
there exists $C^{(4)}_5>0$
such that the following holds on
$\{\Image(w)>R\}$:
\[
-\del_z\del_{\zbar} \log\Tr(s_Y)
\leq C^{(4)}_{5}(\Image w)^{-2}.
\]
Hence, there exists $C^{(4)}_6>0$
such that
\begin{equation}
\label{eq;20.8.6.251}
 -\del_z\del_{\zbar}\Bigl(
 \log\Tr(s_Y)
 -C^{(4)}_3\Re(\gamma w^{1-\epsilon})
 -C^{(4)}_6\log(\Image(w))
 -\log r
  \Bigr)
  \leq 0.
\end{equation}
On $Y\cap \{\Image(w)=3R\}$,
we obtain
\begin{multline}
 \log\Tr(s_Y)
 -C^{(4)}_3\Re(\gamma w^{1-\epsilon})
 -C^{(4)}_6\log(\Image(w))
 -\log r
 \\
 \leq
  \log\Tr(s_Y)
  -C^{(4)}_3\Re(\gamma w^{1-\epsilon})
  \leq 0.
\end{multline}
On $\del(Y)\cap \{\Image(w)\geq 3R\}$,
we obtain
\begin{multline}
 \log\Tr(s_Y)
 -C^{(4)}_3\Re(\gamma w^{1-\epsilon})
 -C^{(4)}_6\log(\Image(w))
 -\log r
 \leq \\
  \log\Tr(s_Y)-\log r
  \leq 0.
\end{multline}
Hence, we obtain
$\log\Tr(s_Y)
-C^{(4)}_3\Re(\gamma w^{1-\epsilon})
-C^{(4)}_6\log(\Image(w))
-\log r\leq 0$
on $Y\cap \{\Image(w)\geq 3R\}$.
There exists
$C^{(4)}_7$ such that
the following holds
on $Y\cap \{\Image(w)\geq 3R\}$:
\begin{equation}
\label{eq;20.8.6.300}
 \log\Tr(s_Y)\leq
 C^{(4)}_7\Re(\gamma w^{1-\epsilon})=:\phi.
\end{equation}

\subsubsection{Construction of harmonic metrics}

Let us explain an outline of the proof of
the third claim of Theorem \ref{thm;20.6.14.10}
under the simplified setting in \S\ref{subsection;24.1.5.20}.
(We shall study a general case in \S\ref{subsection;20.8.6.31}.)

Suppose that $R>10$.
Let $\veca\in\nbigp$.
We take $h_0\in\Harm(q_0)$ such that $\veca(h_0)=\veca$.
Let $Y_i$ be a smooth exhaustive family  of
$X(R,\delta)$
such that
$Y_i\cap \{2R-1<\Image w<4R+1\}
=\{2R-1<\Image w<4R+1,\,\,|\Re w|<a_{Y_i}\}$
for $a_{Y_i}>2$.
We obtain $h_{Y_i}\in\Harm(q_{|Y_i})$
such that
$h_{Y_i|\del Y_i}=h_{0|\del Y_i}$.
According to Proposition \ref{prop;20.6.15.30},
we may assume that
the sequence $\{h_{Y_i}\}$  is convergent to
$h_{\infty}\in\Harm(q)$.
By (\ref{eq;20.8.6.300}),
we obtain
$\log\Tr(h_{Y_i}h_0^{-1})\leq \phi$
on $Y_i\cap\{\Image(w)\geq 3R\}$,
independently of $i$.
Hence, we obtain
$\log\Tr(h_{Y_{\infty}}h_0^{-1})\leq \phi$,
which implies
$\veca(h_{\infty})=\veca(h_0)=\veca$.

\subsubsection{Comparison with harmonic metrics}

Let us explain an outline of the proof of
Proposition \ref{prop;20.6.14.11}
under the simplified setting in \S\ref{subsection;24.1.5.20}.
(We shall study a general case
in \S\ref{subsection;20.8.6.40}.)

Suppose that $R>10$.
Let $Y$ be a relatively compact open subset
of $X(R,\delta)$
such that $\del Y$ is smooth and
that it satisfies (\ref{eq;20.8.6.301}).
Let $h\in\Harm(q)$.
Let $h_Y\in\Harm(q_{|Y})$ such that
$h_{Y|\del Y\cap\{\Image(w)>R\}}=h_{|\del Y\cap\{\Image(w)>R\}}$.
We obtain the automorphism $s_Y$
determined by $h_{Y}=h_{|Y}\cdot s_Y$.
In the following,
constants $C^{(5)}_i$ are independent of
$h$ and $Y$.

By (\ref{eq;20.8.6.110}),
there exists $C^{(5)}_1>0$
such that
the following holds
on $\{\Image (w)>2R\}$:
\begin{equation}
\label{eq;20.8.6.400}
 |\theta(q)|_h\leq C^{(5)}_1\Image(w)^{-1}.
\end{equation}
Because
$h_{Y|\del Y\cap\{\Image(w)>R\}}
=h_{|\del Y\cap\{\Image(w)>R\}}$,
the following holds
on $\del Y\cap \{\Image(w)>2R\}$:
\begin{equation}
 \label{eq;20.8.6.401}
 |\theta(q)|_{h_Y}\leq C^{(5)}_1\Image(w)^{-1}.
\end{equation}

By Proposition \ref{prop;20.4.15.11}
and (\ref{eq;20.8.6.110}),
there exists $C^{(5)}_2$
such that
the following holds
on $\{2R<\Image w<4R,\,|\Re(w)|\leq a_Y-1\}$:
\begin{equation}
\label{eq;20.8.6.403}
 |\theta(q)|_{h_Y}\leq C_2^{(5)}.
\end{equation}
By Proposition \ref{prop;20.6.13.32},
(\ref{eq;20.8.6.110})
and (\ref{eq;20.8.6.401}),
there exists $C^{(5)}_3$
such that the following holds
on $\{2R<\Image w<4R,\,a_Y-1\leq |\Re(w)|\leq a_Y\}$:
\begin{equation}
\label{eq;20.8.6.402}
 |\theta(q)|_{h_Y}\leq C_3^{(5)}.
\end{equation}
By (\ref{eq;20.8.6.111}),
(\ref{eq;20.8.6.400}),
(\ref{eq;20.8.6.403}),
(\ref{eq;20.8.6.402}),
and Corollary \ref{cor;20.4.17.10},
there exists $C^{(5)}_4>0$
such that the following holds
on $\{\Image(w)=3R\}\cap Y$:
\begin{equation}
\label{eq;20.8.6.405}
 \log\Tr(s_Y)\leq
 C_4^{(5)}
 \Re(\gamma w^{1-\epsilon}).
\end{equation}
Because $\log\Tr(s_Y)$ is subharmonic
by (\ref{eq;20.8.16.6}),
we obtain the following on
$Y\cap\{\Image(w)\geq 3R\}$:
\[
 \log\Tr(s_Y)
 \leq C_5^{(5)}\Re(\gamma w^{1-\epsilon})
 +\log r.
\]

\subsection{Preliminary for the proof of
Theorem \ref{thm;20.6.14.10} and Proposition \ref{prop;20.6.14.11}}

We prepare some notation and lemmas
for the proof of Theorem \ref{thm;20.6.14.10}
and Proposition \ref{prop;20.6.14.11}
in \S\ref{subsection;24.1.6.10}.

\subsubsection{Setting}
\label{subsection;20.6.24.10}

Take $R_0>10$ and $0<\delta_0<\pi/2$.
Let $f$ be a holomorphic function on $X(R_0,\delta_0)$,
which induces a section of
$\gbigb_{\projtilde^1_{\infty}}$
on $\Xbar(R_0,\delta_0)$.
We also assume the following.
\begin{itemize}
 \item $\nbigz(f)\cap \bigl(
       \openopen{-\delta_0}{0}
       \cup
       \openopen{\pi}{\pi+\delta_0}
       \bigr)=\emptyset$.
 \item For $\theta\in \openopen{0}{\pi}\setminus\nbigz(f)$,
       $\gminia(f,\theta)$ is of the form:
       \[
       \gminia(f,\theta)=
       \sqrt{-1}w
       +\sum_{0<\kappa<1}
       \gminia(f,\theta)_{\kappa}w^{\kappa}.
       \]
       Note that $\gminia(f,\theta)_{\kappa}$ are constant
       on each connected component of
       $\openopen{0}{\pi}\setminus\nbigz(f)$.
\end{itemize}
For each $\theta\in \openopen{-\delta_0}{\pi+\delta_0}$,
there exist a finite subset
$\nbigi(f,\theta)\subset
 \bigoplus_{\kappa>0} \cnum w^{\kappa}$
and a sector
\[
 W(f,\theta)=
\bigl\{|w|>R(\theta),|\arg(w)-\theta|<L(\theta)
\bigr\}
\]
such that $f$ is expressed as
$f=\sum_{\gminia\in\nbigi(f,\theta)} f_{\theta,\gminia}$
on $W(f,\theta)$,
where $f_{\theta,\gminia}$ has a single growth order $\gminia$,
and $R(\theta)$ and $L(\theta)$ are positive numbers.
We assume that the expression
$f=\sum_{\gminia\in\nbigi(f,\theta)}f_{\theta,\gminia}$
is reduced,
i.e.,
for two distinct elements $\gminia,\gminib\in\nbigi(f,\theta)$,
neither $\gminia\prec_{\theta}\gminib$
nor $\gminib\prec_{\theta}\gminia$ holds.
(See Lemma \ref{lem;20.8.18.1}.)

\begin{lem}
\label{lem;20.9.9.1}
If $0<\theta<\pi$,
we obtain
$\deg(\gminia-\sqrt{-1}w)<1$
for any $\gminia\in\nbigi(f,\theta)$.
For each
$\gminib=\sum_{\kappa>0} \gminib_{\kappa}w^{\kappa}
\in \nbigi(f,0)$,
the following holds.
\begin{itemize}
 \item We have $\deg(\gminib)\geq 1$.
       If $\deg(\gminib)>1$,
       we obtain $-\sqrt{-1}\gminib_{\deg(\gminib)}>0$.
       If $\deg(\gminib)=1$,
       we obtain $-\sqrt{-1}\gminib_{\deg(\gminib)}\geq 1$.
\end{itemize}
For each
$\gminic=\sum_{\kappa>0} \gminic_{\kappa}w^{\kappa}
\in \nbigi(f,\pi)$,
the following holds.
\begin{itemize}
 \item We have $\deg(\gminic)\geq 1$.
       If $\deg(\gminic)>1$,
       we obtain
\[
       -\sqrt{-1}\gminic_{\deg(\gminic)}
       e^{\pi\sqrt{-1}\deg(\gminic)}<0.
\]
       If $\deg(\gminic)=1$,
       we obtain $\sqrt{-1}\gminic_{\deg(\gminic)}\leq -1$.
\end{itemize}
\end{lem}
\pf
Let $0<\theta<\pi$
and $\gminia\in \nbigi(f,\theta)$.
We have the expression:
\[
 \gminia=\sum_{0<\kappa\leq \deg(\gminia)}\gminia_{\kappa}w^{\kappa}.
\]
Take $\theta_-<\theta<\theta_+$
such that
$|\theta-\theta_{\pm}|$ is sufficiently small.
We have the expressions
\[
\gminia(f,\theta_{\pm})
=\sqrt{-1}w+\sum_{0<\kappa<1}\gminia(f,\theta_{\pm})_{\kappa}w^{\kappa}.
\]
If $\deg(\gminia)>1$,
we obtain
$\Re(\gminia_{\deg(\gminia)}e^{\sqrt{-1}\theta_{\pm}\deg(\gminia)})<0$
from $\gminia\prec_{\theta_{\pm}}\gminia(f,\theta_{\pm})$.
It implies that
\[
 \Re(\gminia_{\deg(\gminia)}e^{\sqrt{-1}\theta\deg(\gminia)})<0,
\]
and hence we obtain
$\gminia\prec_{\theta}\gminia(f,\theta_{\pm})$,
which contradicts that
$\gminia,\gminia(f,\theta_{\pm})\in \nbigi(f,\theta)$
and that the expression $f=\sum f_{\theta,\gminia}$ is reduced.
If $\deg(\gminia)<1$,
we obtain
$\gminia(f,\theta_{\pm})\prec_{\theta}\gminia$,
which contradicts that
$\gminia,\gminia(f,\theta_{\pm})\in \nbigi(f,\theta)$
and that the expression $f=\sum f_{\theta,\gminia}$ is reduced.
Hence, we obtain $\deg(\gminia)=1$.
Because
$\gminia,\gminia(f,\theta_{\pm})\in\nbigi(f,\theta)$,
we obtain
$\Re(\gminia_{1}e^{\sqrt{-1}\theta})
=\Re(\sqrt{-1}e^{\sqrt{-1}\theta})$.
If $\gminia_1\neq \sqrt{-1}$,
then
either one of
$\Re(\gminia_{1}e^{\sqrt{-1}\theta_+})
>\Re(\sqrt{-1}e^{\sqrt{-1}\theta+})$
or 
$\Re(\gminia_{1}e^{\sqrt{-1}\theta_-})
>\Re(\sqrt{-1}e^{\sqrt{-1}\theta-})$
holds.
It contradicts 
$\gminia\prec_{\theta_{\pm}}
\gminia(f,\theta_{\pm})$.
Therefore, we obtain
the first claim
$\deg(\gminia-\sqrt{-1}w)<1$.

Let us study the second claim.
We have the expansion of $\gminib\in \nbigi(f,0)$:
\[
 \gminib=\sum_{0<\kappa}\gminib_{\kappa}w^{\kappa}.
\]
Choose a sufficiently small $\theta_+>0$.
Note that
\[
\gminia(f,\theta_+)
=\sqrt{-1}w+\sum_{0<\kappa<1}
\gminia(f,\theta_+)_{\kappa}w^{\kappa}\in\nbigi(f,0). 
\]
If $\deg(\gminib)<1$,
we obtain
$\gminia(f,\theta_+)\prec_{\theta_+}\gminib$,
which contradicts the choice of $\gminia(f,\theta_+)$.
Hence, we obtain 
$\deg(\gminib)\geq 1$.
Because $\gminib,\gminia(f,\theta_+)\in\nbigi(f,0)$,
we obtain
$\Re(\gminib_{\deg(\gminib)})=0$.
Note that
$\gminib\prec_{\theta_+}\gminia(f,\theta_+)$.
Hence, if $\deg(\gminib)>1$,
we obtain
$\Re(\gminib_{\deg\gminib}e^{\sqrt{-1}\deg(\gminib)\theta_+})<0$.
If $\deg(\gminib)=1$,
we obtain
$\Re\bigl((\gminib_{\deg\gminib}-\sqrt{-1})
e^{\sqrt{-1}\deg(\gminib)\theta_+}\bigr)\leq 0$.
Therefore,
$\gminib_{\deg(\gminib)}$
satisfies the desired condition.
We can check the third claim
in a similar way.
\hfill\qed
\vspace{.1in}
 
For each $\gminia=\sum\gminia_kw^{\kappa}\in\nbigi(f,\theta)$,
we set $\supp(\gminia):=
\{\kappa\in\real_{>0}\,|\,\gminia_{\kappa}\neq 0\}$,
and we set
\[
 S(f)=\bigcup_{-\delta_0<\theta<\pi+\delta_0}
 \bigcup_{\gminia\in\nbigi(f,\theta)}
  \supp(\gminia)
  \subset\real.
\]
Note that $S(f)$ is a finite subset of $\real$.
We choose $0<\epsilon_0<1/10$
satisfying the following condition:
\begin{equation}
 \epsilon_0<
  \min\bigl\{
  |\kappa_1-\kappa_2|/10\,\big|\,
  \kappa_1,\kappa_2\in S,\,\,
  \kappa_1\neq\kappa_2
\bigr\}.
\end{equation}

\subsubsection{Another coordinate $g$ and some estimates}
\label{subsection;20.8.6.30}

Let $A$ be a complex number such that
$\Re(\sqrt{-1}A)>0$
and
$\Re(\sqrt{-1}A e^{\sqrt{-1}(1-\epsilon_0)\pi})>0$.
We put $g=w+Aw^{1-\epsilon_0}$
on $\{|w|>R_0,\,\,|\arg(w)-\pi/2|<\pi\}$.
The following lemma is easy to see.
\begin{lem}
 \label{lem;20.8.6.20}
For any $0<L_2<L_1<L_0<\pi$,
there exist $R_0<R_1<R_2$
and an open subset
\[
\{|w|>R_2,|\arg(w)-\pi/2|<L_2\}
 \subset
 \nbigu
 \subset
  \{|w|>R_0,|\arg(w)-\pi/2|<L_0\},
\]
such that
 $g$ induces an isomorphism
$\nbigu
 \simeq \bigl\{
 \zeta\in\cnum\,\big|\,
 |\zeta|>R_1,\,\,|\arg(\zeta)-\pi/2|<L_1
 \bigr\}$.
\end{lem}
\pf
We set
$\nbigu_0:=\{|w|>R_0,\,|\arg(w)-\pi/2|<L_0\}$.
Let $\nbigubar_0$ denote the closure of $\nbigu_0$
in $\projtilde^1_{\infty}$.
Note that
$g$ induces a continuous map
$\nbigubar_0\lrarr \projtilde^1_{\infty}$,
which is also denoted by $g$.
For any $Q\in\varpi_{\infty}^{-1}(\infty)\cap
\nbigubar_0$,
we obtain $g(Q)=Q$.
Moreover,
for the real coordinate systems
$(|w|^{-1},\arg(w))$,
the map $g$ is $C^1$ around $Q$,
and the derivative at $Q$ is the identity.
Then, the claim of the lemma follows from
the inverse function theorem.
\hfill\qed

\begin{lem}
\label{lem;24.1.6.1}
There exist $C_{i}>0$ $(i=1,2)$ such that
\[
  |f|\leq C_{1}e^{-\Image(g)}
\]
on $\{0\leq \Image(g),\,\,|g|>C_2\}$. 
\end{lem}
\pf
It follows from Lemma \ref{lem;20.6.6.1}
and Lemma \ref{lem;20.8.15.20} below.
\qed
 
\begin{lem}
\label{lem;20.8.15.3}
 There exists $B_0>0$ 
 such that
 $g$ induces a holomorphic isomorphism
 \[
 \bigl\{w\in X(R_0,\delta_0)\,\bigr|\,\Image(g(w))\geq B_0\bigr\}
\simeq
\bigl\{
\zeta\in\cnum\,\big|\,\Image(\zeta)\geq B_0
\bigr\}.
 \]
\end{lem}
\pf
Take $\delta_0<\delta_1<\delta_2<\pi/2$.
There exist $R_0<R_1<R_2$ and
an open subset
$X(R_2,\delta_0)
\subset
\nbigu\subset
 X(R_0,\delta_2)$
such that
$g$ induces an isomorphism
$\nbigu\simeq
X(R_1,\delta_1)$.
In particular,
$g$  induces an injection
on $X(R_2,\delta_0)$.
There exists $R_3>R_2$
such that
\[
g^{-1}\Bigl(
g\bigl(X(R_0,\delta_0)
\setminus X(R_2,\delta_0)
\bigr)
\Bigr)
\cap X(R_3,\delta_0)
=\emptyset.
\]
We obtain
$g^{-1}\bigl(
g(X(R_3,\delta_0))
\bigr)
\cap X(R_0,\delta_0)=
X(R_3,\delta_0)$.
By Lemma \ref{lem;20.8.6.20},
there exists $R_4>0$
such that
$X(R_4,0)\subset g(X(R_3,\delta_0))$.
Then, any $B_0>R_4$ satisfies the desired condition.
\hfill\qed

\subsubsection{Some sectors}
\label{subsection;24.1.6.20}

For any $R>0$, $B\in\real$ and $0<L<\pi/4$,
we set
\begin{multline}
S_0(g,R,B,L):=
\\
 \bigl\{
  w\in X(R_0,\delta_0)\,\big|\,|g(w)|>R,\,
  \Image(g(w))>B,\,
  \arg(g(w))<L
  \bigr\},
\end{multline}
\begin{multline}
S_{\pi}(g,R,B,L):=
\\
 \bigl\{
  w\in X(R_0,\delta_0)\,\big|\,|g(w)|>R,\,
  \Image(g(w))>B,\,
  \arg(g(w))>\pi-L
  \bigr\}.
\end{multline}
Let $B_0$ be a positive constant as in Lemma \ref{lem;20.8.15.3}.
By Lemma \ref{lem;24.1.6.1},
there exists $C>0$ such that
\begin{equation}
\label{eq;24.1.6.3}
 |f|\leq Ce^{-\Image g}
\end{equation}
on $\{w\in X(R_0,\delta_0)\,|\,\Image g(w)\geq B_0\}$.
Let $R^{(0)}>0$, $B^{(0)}>B_0$ and $0<L^{(0)}<\pi/4$
such that
$R^{(0)}\sin(L^{(0)})>10B^{(0)}$ and
\begin{equation}
\label{eq;20.9.11.4}
S_0(g,R^{(0)},B^{(0)},L^{(0)})\subset W(f,0),
\quad
S_{\pi}(g,R^{(0)},B^{(0)},L^{(0)})
\subset
W(f,\pi),
\end{equation}
\begin{equation}
\label{eq;20.9.11.5}
 \nbigz(f)\cap \{0<\theta<2L^{(0)}\}=\emptyset,
 \quad
\nbigz(f)\cap \{\pi-2L^{(0)}<\theta<\pi\}=\emptyset.
\end{equation}

\begin{lem}
\label{lem;20.6.24.2}
 For any $B^{(0)}_1>B^{(0)}$
 and $0<L^{(0)}_1<L^{(0)}$,
 there exist $R^{(0)}_1>R^{(0)}$
 and a subset
 $\nbigc_0
 \subset S_0(g,R^{(0)},B^{(0)},L^{(0)})$
 such that the following holds.
\begin{itemize}
 \item There exist $C^{(0)}_i>0$ $(i=1,2)$ such that
       \[
       |f|\geq C^{(0)}_1
       \exp\Bigl(
       \Re\bigl(\gminia(f,\theta_0)\bigr)
       -C^{(0)}_{2}|g|^{\epsilon_0/2}
       \Bigr)
       \]
       on $S_0(g,R^{(0)},B^{(0)},L^{(0)})\setminus\nbigc_0$,
       where $\theta_0>0$ is sufficiently small. 
 \item Let $\nbigd$ be any connected component of $\nbigc_0$
       such that
\[
       \nbigd\cap
        S_0(g,R^{(0)}_{1},B^{(0)}_1,L^{(0)}_1)\neq\emptyset. 
\]
       Then, $\nbigd$ is relatively compact
       in $S_0(g,R^{(0)},B^{(0)},L^{(0)})$.
 \item There exists an increasing sequence 
       of positive numbers $T_{0,i}$ such that
       (i) $\lim T_{0,i}=\infty$,
       (ii) $\nbigc_0\cap \{|g(w)|=T_{0,i}\}
       =\emptyset$.
\end{itemize}
\end{lem}
\pf
Note that
$\gminia(f,\theta_0)$ is contained in $\nbigi(f,0)$.
Consider
\[
 \gminib\in\nbigi(f,0)\setminus\{\gminia(f,\theta_0)\}.
\]
According to Lemma \ref{lem;20.9.9.1},
there exist the following three cases;
(a) $\deg(\gminib)>1$ and $-\sqrt{-1}\gminib_{\deg\gminib}>0$,
(b) $\deg(\gminib)=1$ and
    $-\sqrt{-1}\gminib_1>1$,
(c) $\deg(\gminib)=1$ and $-\sqrt{-1}\gminib_1=1$.
    
According to the first claim of
Lemma \ref{lem;20.8.15.20} below,
if either (a) or (b) is 
satisfied,
there exist
$0<L^{(0)}_2(\gminib)\leq L^{(0)}$,
$C_1(\gminib)>0$ and $C'_1(\gminib)>0$
such that
$\Re(\gminib-\gminia(f,\theta_0))\leq
 -C_1(\gminib)|g|^{\epsilon_0}+C'_1(\gminib)$ 
holds on
$S_0(g,R^{(0)},B^{(0)},L^{(0)}_2(\gminib))$.
Because
(a)
$\Re(\gminib_{\deg(\gminib)}
 e^{\sqrt{-1}\deg(\gminib)\theta})<0$
or
(b) $\Re((\gminib_{\deg(\gminib)}-\sqrt{-1})
 e^{\sqrt{-1}\deg(\gminib)\theta})<0$ 
holds for any $0<\theta\leq L^{(0)}$,
there exist
$0<C_2(\gminib)<C_1(\gminib)$
and
$C'_2(\gminib)>C'_1(\gminib)$
such that
$\Re(\gminib-\gminia(f,\theta_0))\leq
-C_2(\gminib)|g|^{\epsilon_0}+C'_2(\gminib)$
holds on
$S_0(g,R^{(0)},B^{(0)},L^{(0)})
\cap
\{\arg(g(w))>L^{(0)}_2(\gminib)/2\}$.
Therefore,
we obtain
$\Re(\gminib-\gminia(f,\theta_0))\leq
-C_2(\gminib)|g|^{\epsilon_0}+C_2'(\gminib)$
holds on
$S_0(g,R^{(0)},B^{(0)},L^{(0)})$.

In the case (c),
by the third claim of Lemma \ref{lem;20.8.15.20},
there exist
$0<L^{(0)}_2(\gminib)<L^{(0)}$,
$C_1(\gminib)>0$
and $C_1'(\gminib)>0$
such that
either
$\Re(\gminib-\gminia(f,\theta_0))\leq
 -C_1(\gminib)|g|^{\epsilon_0}+C_1'(\gminib)$
or 
 $\Re(\gminia(f,\theta_0)-\gminib)\leq
 -C_1(\gminib)|g|^{\epsilon_0}+C_1'(\gminib)$
holds on
$S_0(g,R^{(0)},B^{(0)},L^{(0)}_2(\gminib))$.
Note that
$\gminib\prec_{\theta}\gminia(f,\theta_0)$
for $0<\theta\leq L^{(0)}$
by Lemma \ref{lem;20.9.11.3}
and (\ref{eq;20.9.11.5}).
Hence,
$\Re(\gminib-\gminia(f,\theta_0))\leq
 -C_1(\gminib)|g|^{\epsilon_0}+C_1'(\gminib)$ holds,
indeed.
Moreover,
because 
$\gminib\prec_{\theta}\gminia(f,\theta_0)$
for $0<\theta\leq L^{(0)}$,
and because
$\deg(\gminia(f,\theta)-\gminib)>2\epsilon_0$,
there exist
$0<C_2(\gminib)<C_1(\gminib)$
and $C_2'(\gminib)>C_1'(\gminib)$
such that
$\Re(\gminib-\gminia(f,\theta_0))\leq
-C_2(\gminib)|g|^{\epsilon_0}+C_2'(\gminib)$
holds on
$S_0(g,R^{(0)},B^{(0)},L^{(0)})
\cap
 \{\arg(w)>L^{(0)}_2(\gminib)/2\}$.
 Hence, we obtain
$\Re(\gminib-\gminia(f,\theta_0))\leq
-C_2(\gminib)|g|^{\epsilon_0}+C_2'(\gminib)$
on 
$S_0(g,R^{(0)},B^{(0)},L^{(0)})$.

In all,
there exist $C_3>0$ and $C_3'>0$
such that
 \[
  \Re(\gminib-\gminia(f,\theta_0))
  \leq
  -C_3|g|^{\epsilon_0}+C_3'
 \]
 on
 $S_0(g,R^{(0)},B^{(0)},L^{(0)})$
 for any $\gminib\in\nbigi(f,0)\setminus\{\gminia(f,\theta_0)\}$.
 Hence, there exists $C_i>0$ $(i=4,5)$ such that
 \[
 \Bigl|
  \sum_{\gminib\neq\gminia(f,\theta_0)} f_{\gminib}
  \Bigr|
  \leq
  C_4\exp\Bigl(\Re\bigl(\gminia(f,\theta_0)\bigr)
  -C_5|g|^{\epsilon_0}\Bigr).
 \]
Then, we obtain the claim of the lemma
by applying Lemma \ref{lem;20.6.24.1}
with $\rho=\epsilon_0/2$. 
\hfill\qed

\vspace{.1in}

The following lemma is similar to
Lemma \ref{lem;20.6.24.2}.
\begin{lem}
\label{lem;20.8.4.30}
 For any $B^{(0)}_1>B^{(0)}$,
 and $0<L^{(0)}_1<L^{(0)}$,
 there exist $R^{(0)}_1>R^{(0)}$,
 and a subset
 $\nbigc_{\pi}
 \subset S_{\pi}(g,R^{(0)},B^{(0)},L^{(0)})$
 such that the following holds.
\begin{itemize}
 \item There exists $C^{(0)}_i>0$ $(i=1,2)$ such that
       \[
       |f|\geq C^{(0)}_1
       \exp\Bigl(
       \Re\bigl(\gminia(f,\theta_{\pi})\bigr)
       -C^{(0)}_{2}|g|^{\epsilon_0/2}
       \Bigr)
       \]
       on $S_{\pi}(g,R^{(0)},B^{(0)},L^{(0)})\setminus\nbigc_{\pi}$,
       where $\pi-\theta_{\pi}>0$ is sufficiently small.
 \item Let $\nbigd$ be any connected component of $\nbigc_{\pi}$
       such that
\[
       \nbigd\cap
       S_{\pi}(g,R^{(0)}_{1},B^{(0)}_1,L^{(0)}_1)\neq\emptyset. 
\]
       Then, $\nbigd$ is relatively compact
       in $S_{\pi}(g,R^{(0)},B^{(0)},L^{(0)})$.
 \item There exists an increasing sequence
       of positive numbers $T_{\pi,i}$ such that
       (i) $\lim T_{\pi,i}=\infty$,
       (ii) $\nbigc_{\pi}\cap
       \{|g(w)|=T_{\pi,i}\}
       =\emptyset$.
       \hfill\qed
\end{itemize}
\end{lem}

\subsubsection{Relatively compact open subsets}
\label{subsection;24.1.6.21}

We introduce a condition for
a relatively compact open subset $Y$ in $X(R_0,\delta_0)$,
which we shall use in
\S\ref{subsection;20.8.6.50}--\S\ref{subsection;20.8.5.20}.

Let $B_0$ be a constant as in Lemma \ref{lem;20.8.15.3}.
We take $R^{(0)}>0$,
$B^{(0)}>B_0$
and $0<L^{(0)}<\pi/4$
as in \S\ref{subsection;24.1.6.20}.
We set $B_1^{(0)}=2B^{(0)}$
and $L^{(0)}_1=\frac{1}{2}L^{(0)}$,
and let $R_1^{(0)}>R^{(0)}$
be as in Lemma \ref{lem;20.6.24.2} and Lemma \ref{lem;20.8.4.30}.
We take $B_0<B^{(0)}_{-1}<B^{(0)}$,
$0<R^{(0)}_{-1}<R^{(0)}$,
and $L^{(0)}<L^{(0)}_{-1}<\pi/4$
such that
\[
S_0(g,R^{(0)}_{-1},B^{(0)}_{-1},L^{(0)}_{-1})
\subset W(f,0),
\quad
S_{\pi}(g,R^{(0)}_{-1},B^{(0)}_{-1},L^{(0)}_{-1})
\subset W(f,\pi).
\]
We also assume
$R^{(0)}_{-1}\sin(L^{(0)})>10B^{(0)}$.

\begin{condition}
 \label{condition;20.6.15.10}
 \mbox{{}}
 \begin{itemize}
 \item The boundary $\del Y$ is smooth.      
 \item There exists $T_{0,i(0)}>2R^{(0)}_{-1}$ such that
\[
       Y\cap S_0(g,R_{-1}^{(0)},B_{-1}^{(0)},L_{-1}^{(0)})
       =S_0(g,R_{-1}^{(0)},B_{-1}^{(0)},L_{-1}^{(0)})\cap
       \{|g(w)|\leq T_{0,i(0)}\}.
\]
  \item There exists $T_{\pi,i(\pi)}>2R^{(0)}_{-1}$ such that
\[
Y\cap S_{\pi}(g,R_{-1}^{(0)},B_{-1}^{(0)},L_{-1}^{(0)})
       =S_{\pi}(g,R_{-1}^{(0)},B_{-1}^{(0)},L_{-1}^{(0)})\cap
       \{|g(w)|\leq T_{\pi,i(\pi)}\}. 
\]
  \item $\bigl\{w\in X(R_0,\delta_0)\,\big|\,
	 \Image g(w)\geq B_{-1}^{(0)},\,
	 |g(w)|\leq 2R_{-1}^{(0)}
	 \bigr\}\subset Y$.
 \end{itemize}
\end{condition}
We shall use an exhaustive family
satisfying Condition \ref{condition;20.6.15.10}
to construct harmonic metrics with the desired property
in Theorem \ref{thm;20.6.14.10}
and Theorem \ref{thm;20.6.9.30}.

\subsubsection{Appendix: A choice of complex coordinate}

Let $R>0$ and $0<\delta<\pi/2$.
Let $\gminia_i$ $(i=1,\ldots,n(1))$ be
mutually distinct holomorphic functions
on $X(R,\delta)$
expressed as a finite sum
\[
 \gminia_i=
 \sqrt{-1}w+\sum_{0<\kappa<1}\gminia_{i,\kappa}w^{\kappa}.
\]
Let $\gminib_i$ $(i=1,\ldots,n(2))$ be holomorphic functions
on $X(R,\delta)$
expressed as a finite sum
\[
 \gminib_i=\sum_{0<\kappa\leq \rho(\gminib_i)}
 \gminib_{i,\kappa} w^{\kappa},
\]
where
one of the following holds.
\begin{itemize}
 \item $\rho(\gminib_i)>1$ and
       $-\sqrt{-1}\gminib_{i,\rho(\gminib_i)}>0$.
 \item $\rho(\gminib_i)=1$ and
       $-\sqrt{-1}\gminib_{i,1}>1$.
\end{itemize}
Let $\gminic_i$ $(i=1,\ldots,n(3))$ be holomorphic functions
on $X(R,\delta)$ expressed as a finite sum
\[
 \gminic_i=\sum_{0<\kappa\leq \rho(\gminic_i)}
  \gminic_{i,\kappa}w^{\kappa},
\]
where one of the following holds.
\begin{itemize}
 \item $\rho(\gminic_i)>1$ and
       $-\sqrt{-1}\gminic_{i,\rho(\gminic_i)}
       e^{\pi\sqrt{-1}\rho(\gminic_i)}<0$.
 \item $\rho(\gminic_i)=1$ and
      $\sqrt{-1}\gminic_{i,1}<-1$.
\end{itemize}

We set
\[
 S=\{0,1\}\cup \bigcup_i\{\kappa\,|\,\gminia_{i,\kappa}\neq 0\}
 \cup\bigcup_i\{\kappa\,|\,\gminib_{i,\kappa}\neq 0\}
 \cup \bigcup_i\{\kappa\,|\,\gminic_{i,\kappa}\neq 0\}.
\]
Take $\epsilon>0$
such that
$2\epsilon<|\kappa_1-\kappa_2|$
for any two distinct elements
$\kappa_1,\kappa_2\in S$.

Let $A$ be a complex number such that
\[
\Re(\sqrt{-1}A)>0,
\quad
\Re(\sqrt{-1}Ae^{\sqrt{-1}(1-\epsilon)\pi})>0.
\]
We set $g(w)=w+Aw^{1-\epsilon}$
on $X(R,\delta)$.
In the proof of Lemma \ref{lem;20.8.15.3},
we observed that
there exists
$R_0>R$ such that
(i) $g^{-1}\bigl(g(X(R_0,\delta))\bigr)
 =X(R_0,\delta)$,
(ii)
$g$ induces an isomorphism
$X(R_0,\delta)\simeq g(X(R_0,\delta))$.
According to Lemma \ref{lem;20.8.6.20},
there exists $R_0'>R_0$
such that
$X(R_0',\delta/2)\subset
g(X(R_0,\delta))$.

\begin{lem}
 \label{lem;20.6.6.1}
There exist
positive constants $R_j>R'_0$ $(j=1,2)$
and $C_1$ such that the following
conditions are satisfied:
\begin{itemize}
 \item  $\Re(\gminia_p)-\Re(\sqrt{-1}g)<-C_1|w|^{1-\epsilon}$
	for any $p$ on
	$X(R_1,0)=\{|w|\geq R_1,\,\,\,\Image(w)\geq 0\}$.
 \item $\bigl\{w\in X(R,\delta)\,\big|\,
        |g(w)|>R_2,\,\,\Image g(w)\geq 0\bigr\}
       \subset X(R_1,0)$.
\end{itemize}
\end{lem}
\pf
Note that
$\sqrt{-1}g-\gminia_p
=\sqrt{-1}Aw^{1-\epsilon}
+\sum_{0<\kappa<1-\epsilon}
\alpha_{p,\kappa}w^{\kappa}$
for $\alpha_{p,\kappa}\in\cnum$.
Because $\{(p,\kappa)\,|\,\alpha_{p,\kappa}\neq 0\}$
is finite,
there exist $R_1$ and $C_1$
as in the first condition.

By our choice of $R_0<R_0'$,
the set
\begin{multline}
\nbigu_1:=\bigl\{w\in X(R,\delta)\,\big|\,
|g(w)|>R'_0,\,\,\Image g(w)\geq 0\bigr\}
\simeq \\
\bigl\{
\zeta\in\cnum\,\big|\,
|\zeta|>R'_0,\,\Image\zeta\geq 0
\bigr\}
\end{multline}
is connected.
Clearly, we obtain
$\nbigu_1\cap
 X(R,0)\neq\emptyset$.
Because
$\Image g(w)=\Image(w)-\Re(\sqrt{-1}Aw^{1-\epsilon})$,
the intersection
\[
\bigl\{w\in X(R,\delta)\,\big|\,
|g(w)|>R'_0,\Image g(w)\geq 0\bigr\}
\cap
\{w\in X(R,\delta)\,|\,\Image(w)=0\}
\]
is empty.
Hence, if $R_2>R_0'$ is sufficiently large,
we obtain
$\bigl\{w\in X(R,\delta)\,\big|\,
|g(w)|>R_2,\,\,\Image g(w)\geq 0\bigr\}
\subset
X(R_1,0)$.
\hfill\qed

\begin{lem}
 \label{lem;20.8.15.20}
There exist
positive constants
$L$, $R_3$ and $C_j$ $(j=2,3)$ such that the following holds:
\begin{itemize}
 \item  $\Re(\gminib_i)-\Re(\gminia_p)
	<-C_2|g|^{\rho(\gminib_i)-\epsilon}$
	for any $i$ and $p$
	on
	\[
	\nbigu_0:=\bigl\{w\in X(R,\delta)\,\big|\,
	\Re g(w)>0,\,\,|g(w)|\geq R_3,\,\,
        0\leq \arg g(w)\leq L\bigr\}.
	\]
 \item $\Re(\gminic_i)-\Re(\gminia_p)
       <-C_2|g|^{\rho(\gminic_i)-\epsilon}$
       for any $i$ and $p$
       on
       \[
       \nbigu_{\pi}:=\bigl\{w\in X(R,\delta)\,\big|\,
       \Re g(w)<0,\,\,|g(w)|\geq R_3,\,\,
       \pi-L\leq \arg g(w)\leq \pi\bigr\}.
       \]
 \item For $p\neq q$,
       we obtain
       $C_3|g|^{\epsilon}\leq
       \bigl|\Re(\gminia_p-\gminia_q)\bigr|$
       on $\nbigu_0\cup\nbigu_{\pi}$.
\end{itemize}
\end{lem}
\pf
We use the following estimate
on $X(R_0,\delta)$ for some $\epsilon'>0$:
\begin{equation}
\label{eq;20.8.15.10}
w-\bigl(g-Ag^{1-\epsilon}\bigr)
=O(|w|^{1-\epsilon-\epsilon'}).
\end{equation}

Let us study the first claim.
To simplify the notation,
we set $\rho(i):=\rho(\gminib_i)$.
If $\rho(i)>1$,
we obtain the following expression
\[
 \gminib_i-\gminia_p
=\gminib_{i,\rho(i)}w^{\rho(i)}
 +\sum_{0<\kappa<\rho(i)}
  \gminib_{(i,p),\kappa}w^{\kappa}.
\]
By (\ref{eq;20.8.15.10}),
there exists $\delta(i,p)>0$
such that
\[
\gminib_i-\gminia_p
-\Bigl(
\gminib_{i,\rho(i)}
g^{\rho(i)}
-\rho(i)\cdot
\gminib_{i,\rho(i)}
Ag^{\rho(i)-\epsilon}
\Bigr)
=O(|g|^{\rho(i)-\epsilon-\delta(i,p)}).
\]
By the assumption,
there exists $b>0$ such that
$\gminib_{i,\rho(i)}=\sqrt{-1}b$.
We obtain
\begin{multline}
\Re\Bigl(
 \gminib_{i,\rho(i)}
g^{\rho(i)}
-\rho(i)\cdot
\gminib_{i,\rho(i)}
Ag^{\rho(i)-\epsilon}
\Bigr) 
= \\
-b|g|^{\rho(i)}
\sin\bigl(\rho(i)\arg(g)\bigr)
-b\rho(i)|g|^{\rho(i)-\epsilon}
\Re\bigl(
\sqrt{-1}Ae^{\sqrt{-1}(\rho(i)-\epsilon)\arg(g)}
\bigr).
\end{multline}
There exists $L_1(i,p)>0$
such that
\[
\sin\bigl(\rho(i)\arg(g)\bigr)\geq 0,
\quad
\Re\bigl(
\sqrt{-1}Ae^{\sqrt{-1}(\rho(i)-\epsilon)\arg(g)}
\bigr)>0
\]
on
$\bigl\{
0\leq \arg(g)\leq L_1(i,p)
\bigr\}$.
Hence,
there exist $R_2'>0$ and $C_3(i,p)>0$
such that
\begin{equation}
\label{eq;20.8.15.11}
 \Re(\gminib_i-\gminia_p)
 \leq -C_3(i,p)|g|^{\rho(i)-\epsilon}
\end{equation}
on 
$\bigl\{
w\in X(R,\delta)\,\big|\,
|g(w)|>R_2',\,\,
0\leq \arg(g)\leq L_1(i,p)
\bigr\}$.
If $\rho(i)=1$,
there exists $b>1$ such that
$\gminib_{i,1}=\sqrt{-1}b$.
We obtain
\begin{multline}
 \gminib_i-\gminia_p
 =
 (\gminib_{i,1}-\sqrt{-1})w
 +\sum_{0<\kappa<1}
 \gminib_{(i,p),\kappa}w^{\kappa}
 \\
 =\sqrt{-1}(b-1)w
  +\sum_{0<\kappa<1}
 \gminib_{(i,p),\kappa}w^{\kappa}.
\end{multline}
By the same argument,
we can prove that there exist
$L_1(i,p)>0$, $R_2'>0$ and $C_3(i,p)>0$
such that 
(\ref{eq;20.8.15.11})
holds on
$\bigl\{
w\in X(R,\delta)\,\big|\,
|g(w)|>R_2',\,\,
0\leq \arg(g)\leq L_1(i,p)
\bigr\}$
with $\rho(i)=1$.
Thus, we obtain the first claim.
The second claim can be proved similarly.

Let us prove the third claim.
For $p\neq q$,
we set $\rho(p,q):=\deg(\gminia_p-\gminia_q)>0$.
We obtain the following expression:
\[
 \gminia_p-\gminia_q
 =
 \sum_{0<\kappa\leq \rho(p,q)}
 \gminia_{(p,q),\kappa}w^{\kappa}.
\]
Here, $\gminia_{(p,q),\rho(p,q)}\neq 0$.
There exists $\delta(p,q)>0$
such that the following holds:
\begin{multline}
 \Re(\gminia_p-\gminia_q)
 -\Re\Bigl(
 \gminia_{(p,q),\rho(p,q)}
 \bigl(g^{\rho(p,q)}-\rho(p,q)Ag^{\rho(p,q)-\epsilon}
 \bigr)
 \Bigr)
\\
 =O\bigl(
 |g|^{\rho(p,q)-\epsilon-\delta(p,q)}
\bigr).
\end{multline}
If $\gminia_{p,q,\rho(p,q)}\not\in\sqrt{-1}\real$,
there exist $L_3(p,q)>0$ and $C_4(p,q)>0$
such that
\begin{equation}
\label{eq;20.8.15.12}
\Bigl|
 \Re\bigl(
  \gminia_{(p,q),\rho(p,q)}
  g^{\rho(p,q)}
  \bigr)
  \Bigr|
  >C_4(p,q)|g|^{\rho(p,q)}
\end{equation}
on $\bigl\{
w\in X(R_0,\delta)\,\big|\,
 0\leq \arg(g(w))\leq L_3(p,q)
 \bigr\}$.
If $\gminia_{p,q,\rho(p,q)}\in\sqrt{-1}\real\setminus\{0\}$,
there exists $b\neq 0$
such that $\gminia_{p,q,\rho(p,q)}=\sqrt{-1}b$.
Note that
{\small
\begin{multline}
 \Bigl|
 \Re\Bigl(
 \gminia_{(p,q),\rho(p,q)}
 \bigl(g^{\rho(p,q)}-\rho(p,q)Ag^{\rho(p,q)-\epsilon}
 \bigr)
 \Bigr)
 \Bigr|
 = 
 |b|\times \\
 \Bigl|
|g|^{\rho(p,q)}\sin\bigl(\rho(p,q)\arg g\bigr)
+\rho(p,q)|g|^{\rho(p,q)-\epsilon}
 \Re\bigl(\sqrt{-1}Ae^{\sqrt{-1}(\rho(p,q)-\epsilon)\arg(g)}
 \bigr)
 \Bigr|.
\end{multline}
}
Hence,
there exist $L_3(p,q)>0$ and $C_4(p,q)>0$
\[
 \Bigl|
 \Re\bigl(
  \gminia_{(p,q),\rho(p,q)}
  g^{\rho(p,q)}
  \bigr)
  \Bigr|
  >C_4(p,q)|g|^{\rho(p,q)-\epsilon}.
\]
on $\bigl\{
w\in X(R_0,\delta)\,\big|\,
 0\leq \arg(g(w))\leq L_3(p,q)
 \bigr\}$.
Then, we can easily deduce the third claim.
\hfill\qed

\subsection{Proof of Theorem \ref{thm;20.6.14.10}
and Proposition \ref{prop;20.6.14.11}}
\label{subsection;24.1.6.10}

Take $R_0>10$ and $0<\delta_0<\pi/2$.
Let $f$ be a holomorphic function on $X(R_0,\delta_0)$
as in \S\ref{subsection;20.6.24.10}.
We set $q=f\,(dw)^r$.
We first study $G_r$-invariant harmonic metrics
of $(\hyperk_{X(R_0,\delta_0)},\theta(q))$
in \S\ref{subsection;20.8.6.500}--\S\ref{subsection;20.8.6.40}.
Then, we shall derive
Theorem \ref{thm;20.6.14.10} and Proposition \ref{prop;20.6.14.11}
in \S\ref{subsection;20.8.5.20}.

\subsubsection{The associated parabolic weights}
\label{subsection;20.8.6.500}

Let $\epsilon_0>0$ be as in \S\ref{subsection;20.6.24.10}.

\begin{prop}
\label{prop;20.6.9.22}
 Let $h\in\Harm(q)$.
 Then, there exists
 $\veca(h)\in\nbigp$
 satisfying the following estimate
 on any
 $\bigl\{w\in\cnum\,\big|\,|w|>R_1,
 |\arg(w)-\pi/2|<(1-\delta)\pi/2\bigr\}$
 for any $R_1>R_0$ and $\delta>0$:
\begin{equation}
\label{eq;20.6.14.100}
 \log
 \bigl|
  (dw)^{(r+1)/2-i}
 \bigr|_h
 =a_i(h)\cdot \Image(w)
 +O\Bigl(
 |w|^{1-\epsilon_0}
 \Bigr).
\end{equation}
\end{prop}
\pf
For $R>0$, $B>0$ and $0<L<\pi/4$,
let
$\Stilde(w,R,B,L)$
denote the following set:
\begin{multline}
 \bigl\{
  w\in\cnum\,\big|\,
   |w|>R,\,\Image(w)\geq B,\,\arg(w)<L
 \bigr\}
\cup \\
  \bigl\{
  w\in\cnum\,\big|\,
   |w|>R,\,\Image(w)\geq B,\,\arg(w)>\pi-L
 \bigr\}.
\end{multline}
By using the function $g=w+Aw^{1-\epsilon_0}$
as in \S\ref{subsection;20.8.6.30},
from now on we use $w$ to refer to $g$ for
simplicity.
By the results in \S\ref{subsection;20.8.6.30},
the following condition is satisfied.
\begin{condition}
 \label{condition;20.6.15.1}
 There exists $C^{(1)}_0>0$
 such that
 $|f|\leq C^{(1)}_0e^{-\Image(w)}$ on
 $\{|w|>R_0,\,\,\Image(w)\geq 0\}$.
 Moreover, 
 there exist
 $2R_0<R^{(1)}<R^{(1)}_1$,
 $2R_0< B^{(1)}<B^{(1)}_1$,
 $0<L^{(1)}_1<L^{(1)}$,
 and
 $\nbigc\subset
 \Stilde(w,R^{(1)},B^{(1)},L^{(1)})$
 such that the following conditions are satisfied.
\begin{itemize}
 \item There exists $C^{(1)}_i>0$ $(i=1,2)$ such that
       \[
      |f|\geq C^{(1)}_1\exp(-\Image(w)-C^{(1)}_2|w|^{1-\epsilon_0})
       \]
       on $\Stilde(w,R^{(1)},B^{(1)},L^{(1)})\setminus\nbigc$.
 \item Let $\nbigd$ be any connected component of
       $\nbigc$ such that
\[
      \nbigd\cap \Stilde(w,R^{(1)}_1,B^{(1)}_1,L^{(1)}_1)\neq\emptyset.
\]
       Then, $\nbigd$ is relatively compact in
       $\Stilde(w,R^{(1)},B^{(1)},L^{(1)})$.
\end{itemize}
\end{condition}

Let $\theta(q)$ denote the Higgs field of
$\hyperk_{X(R_0,\delta_0),r}$
associated with $q$.
For $w_1$ with $\Image(w_1)\geq 2R_0$,
$\{|w-w_1|<\Image(w_1)-R_0/2\}$ is contained in
$\{\Image(w)> R_0\}$.
There exists an isomorphism
\[
 \Phi_{w_1}:\{|\eta|<1\}
\simeq
\{|w-w_1|<\Image(w_1)-R_0/2\}
\]
given by
$\Phi_{w_1}(\eta)=w_1+(\Image(w_1)-R_0/2)\eta$.
Note that
there exists $C^{(1)}_3>0$
such that $|\Phi_{w_1}^{\ast}(q)|<C^{(1)}_3$
independently from $w_1$.
By applying Proposition \ref{prop;20.4.15.11}
to $\Phi_{w_1}^{\ast}(\hyperk_{X(R_0,\delta_0),r},\theta(q),h)$,
we obtain that
the sup norm of 
$\Phi_{w_1}^{\ast}(\theta(q))$
with respect to $\Phi_{w_1}^{\ast}(h)$
on $\{|\eta|<1/2\}$
is dominated by a constant independently from $w_1$.
Because
$\Phi_{w_1}^{\ast}(dw)=(\Image(w_1)-R_0/2)d\eta$,
we obtain $|\theta(q)|_h=O\bigl(\Image(w)^{-1}\bigr)$
on $\{\Image(w)\geq 2R_0\}$.
By the Hitchin equation,
we also obtain $|R(h)|_h=O\bigl(\Image(w)^{-2}\bigr)$
on $\{\Image(w)\geq 2R_0\}$.

\begin{lem}
\label{lem;20.6.24.101}
 The following estimates hold on
 $\Stilde(w,R^{(1)}_1,B^{(1)}_1,L^{(1)}_1)$
 for $i=1,\ldots,r$:
 \begin{equation}
\label{eq;20.6.24.20}
  \log
 \bigl|
  (dw)^{(r+1)/2-i}
 \bigr|_h
 =
 O\Bigl(
 \Image(w)+|w|^{1-\epsilon_0}
 \Bigr).
  \end{equation}
\end{lem}
\pf
We obtain the estimate
(\ref{eq;20.6.24.20})
on $\Stilde(w,R^{(1)},B^{(1)},L^{(1)})\setminus\nbigc$
by
the estimate $|\theta(q)|_h=O\bigl(\Image(w)^{-1}\bigr)$,
Condition \ref{condition;20.6.15.1}
and Corollary \ref{cor;20.4.21.2}.
Because $|R(h)|=O(\Image(w)^{-2})$ on $\{\Image(w)\geq 2R_0\}$,
there exist functions $\beta_i$ on $\{\Image(w)\geq 2R_0\}$
such that
$|\beta_i|=O\bigl(\log(\Image(w))\bigr)$
and
$\log\bigl|(dw)^{(r+1-2i)/2}\bigr|_h-\beta_i$
are harmonic functions on $\{\Image(w)\geq 2R_0\}$.
(See Lemma \ref{lem;20.4.21.3}.)
Then, by using Condition \ref{condition;20.6.15.1} again,
we obtain the estimate (\ref{eq;20.6.24.20})
on $\Stilde(w,R^{(1)}_1,B^{(1)}_1,L^{(1)}_1)$.
\hfill\qed

\vspace{.1in}

\begin{lem}
\label{lem;20.6.24.102}
 Let $K$ be any compact subset of
 $\openopen{0}{\pi}$,
 which we regard as a closed subset
 in $\varpi_{\infty}^{-1}(\infty)\cap\Xbar(R_0,\delta_0)$.
 Then, there exists a neighbourhood $\nbigu$
 of $K$ in $\Xbar(R_0,\delta_0)$
 such that the estimate
 {\rm(\ref{eq;20.6.24.20})} holds
 on $\nbigu\setminus\varpi_{\infty}^{-1}(\infty)$. 
\end{lem}
\pf
It is enough to study the case where $K$ consists of a point.
If $K\subset \openopen{0}{\pi}\setminus\nbigz(f)$,
we obtain the claim from Lemma \ref{lem;20.6.15.3}.
If $K\subset \nbigz(f)\cap\openopen{0}{\pi}$,
we obtain the claim from Lemma \ref{lem;20.6.14.2}
and Corollary \ref{cor;20.6.12.21}.
(See the proof of Theorem \ref{thm;20.6.9.21}
in \S\ref{subsection;20.6.24.100} for a more detailed argument.)
\hfill\qed

\vspace{.1in}

By Lemma \ref{lem;20.6.24.101}
and Lemma \ref{lem;20.6.24.102},
the estimate (\ref{eq;20.6.24.20})
holds on
$\bigl\{
 |w|>R^{(1)}_1,
\Image(w)\geq B^{(1)}_1
\bigr\}$.
We also have the following estimate
on $\{\Image(w)\geq 2R_0\}$:
\[
 \del\delbar\log
 \bigl|
  (dw)^{(r+1)/2-i}
  \bigr|_h
=O\bigl(\Image(w)^{-2}\bigr).
\]
By Proposition \ref{prop;20.4.22.11},
there exist $a_i(h)\in\real$ $(i=1,\ldots,r)$
such that
\[
 \log\bigl|(dw)^{(r+1-2i)/2}\bigr|_h
 -a_i(h)\Image(w)
 =O(|w|^{1-\epsilon_0})
\]
on $\{|w|>R^{(1)}_1,\,\Image(w)\geq B_1^{(1)}\}$.
Because the above $w$ is referring to $g$,
it means
\[
 \log\bigl|(dg(w))^{(r+1-2i)/2}\bigr|_h
 -a_i(h)\Image(g(w))
 =O(|g(w)|^{1-\epsilon_0})
\]
on $\{|g(w)|>R^{(1)}_1,\,\Image g(w)\geq B_1^{(1)}\}$.
By Lemma \ref{lem;20.8.6.20},
for any $\delta>0$,
there exists $R_1>0$ such that
$\{|w|>R_1,\,\,|\arg(w)-\pi/2|<(1-\delta)\pi/2\}
\subset
 \{|g(w)|>R^{(1)}_1,\Image g(w)\geq B_1^{(1)}\}$.
Using $dg/dw-1=O(|w|^{-\epsilon_0})$
and
$\Image g(w)-\Image w=O(|w|^{1-\epsilon_0})$,
we obtain
\[
 \log\bigl|(dw)^{(r+1-2i)/2}\bigr|_h
 -a_i(h)\Image(w)
 =O(|w|^{1-\epsilon_0})
\]
on $\{|w|>R_1,\,\,|\arg(w)-\pi/2|<(1-\delta)\pi/2\}$.
Because
\[
 \theta(q) (dw)^{(r+1-2i)/2}=(dw)^{(r+1-2(i+1))/2}\,dw
 \quad
 (i=1,\ldots,r-1)
\]
and $|\theta(q)|_h=O\bigl(\Image(w)^{-1}\bigr)$,
we obtain $a_i(h)\geq a_{i+1}(h)$.
We also obtain $a_{r}(h)\geq a_1(h)-1$
by using
\[
 \theta(q)(dw)^{(-r+1)/2}=f(dw)^{(r-1)/2}(dw)
\]
and the estimate of $f$ around any point of
$(\{\infty\}\times\openopen{0}{\pi})\setminus \nbigz(f)$
as in Condition \ref{condition;20.4.28.100}.
Thus, the proof of Proposition \ref{prop;20.6.9.22} is completed.
\hfill\qed

\subsubsection{Mutually boundedness}
\label{subsection;20.8.6.501}

\begin{prop}
\label{prop;20.6.9.23}
 Let $h_i\in\Harm(q)$ $(i=1,2)$
 such that $\veca(h_1)=\veca(h_2)$.
Then,
 for any $R_1>R_0$ and $0<\delta_1<\delta_0$,
 $h_1$ and $h_2$ are mutually bounded
 on $X(R_1,\delta_1)$.
\end{prop}
\pf
Take any $0<\delta_1<\delta_0$
and $R_1>R_0$.
It is enough to prove that
there exists $R'_1>R_1$
such that $h_1$ and $h_2$
are mutually bounded on $X(R'_1,\delta_1)$
because $X(R_1,\delta_1)\setminus X(R_1',\delta_1)$
is compact.

Let $s$ be determined by
$h_2=h_1\cdot s$.
Because $\det(s)=1$,
we obtain $\Tr(s)\geq r$.
By Proposition \ref{prop;20.6.12.20},
$\Tr(s)$ is bounded
on $X(R_1,\delta_1)\setminus X(R_1,\delta_3)$
for any $0<\delta_3<\delta_1$.

By the assumption $\veca(h_1)=\veca(h_2)$,
the following estimate holds on
$\bigl\{|w|>R_1,\,\,|\arg(w)-\pi/2|<(1-\delta)\pi/2
\bigr\}$
for any $\delta>0$:
\[
 \log\Tr(s)=O\bigl(|w|^{1-\epsilon_0}\bigr).
\]
By Lemma \ref{lem;20.6.14.2},
there exist $N>0$
and a subset $Z_1\subset \real_{>0}$ with
$\int_{Z_1}dt/t<\infty$
such that
\[
 \log\Tr(s)
 =O(|w|^{N})
\]
on
$\bigl\{|w|>R_1,\,\,
-\delta<\arg(w)<\delta,\,\,
|w|\not\in Z_1
\bigr\}$
for any $\delta>0$.
By Corollary \ref{cor;20.6.12.21},
we obtain that
\[
 \log\Tr(s)=O\bigl(|w|^{1-\epsilon_0}\bigr)
\]
on
$\bigl\{|w|>R_1,\,\,
-\delta_1<\arg(w)<\delta_1
\bigr\}$.
Similarly,
we obtain 
\[
 \log\Tr(s)=O\bigl(|w|^{1-\epsilon_0}\bigr)
\]
on $\{|w|>R_1,\,\,\pi-\delta_1<\arg(w)<\pi+\delta_1\}$.
If $0<\delta_4\leq\delta_1$
is sufficiently small,
we obtain that
$\log\Tr(s)$ is bounded
on $X(R_1,\delta_4)$
by Corollary \ref{cor;20.5.2.1}.
Because $\log\Tr(s)$ is bounded
on $X(R_1,\delta_1)\setminus X(R_1,\delta_4/2)$,
we obtain that
$\log\Tr(s)$ is bounded
on $X(R_1,\delta_1)$.
\hfill\qed

\subsubsection{Auxiliary metrics}
\label{subsection;20.8.6.502}

To introduce an auxiliary metric $h_0$,
we set $q_0:=e^{\sqrt{-1}g(w)}(dg(w))^r$
on $X(R_0,\delta_0)$.

\begin{lem}
 \label{lem;20.6.15.100}
For any $\veca\in\nbigp$,
there exists $h_0\in\Harm(q_0)$
 such that the following holds on
 $\nbigu_1:=
 \bigl\{w\in X(R_0,\delta_0)\,\big|\,\Image(g(w))\geq B_0\bigr\}$:
\[
\log
 \bigl|
  (dw)^{(r+1)/2-i}
 \bigr|_{h_0}
 =a_i\cdot\Image(w)
 +O\Bigl(
 |w|^{1-\epsilon_0}
 \Bigr).
\]
We also obtain
\[
 |\theta(q)|_{h_0}=
 O\Bigl(\bigl(\Image g(w)\bigr)^{-1}\Bigr),
 \quad
 |F(h_0,\theta(q))|_{h_0}=
 O\Bigl(\bigl(\Image g(w)\bigr)^{-2}\Bigr)
\]
on $\nbigu_1$.
 Here,
 $F(h_0,\theta(q))=R(h_0)+
 [\theta(q),\theta(q)^{\dagger}_{h_0}]$.
\end{lem}
\pf
We consider the map
$\Phi:\cnum\lrarr \cnum^{\ast}$
by
$\Phi(w)=e^{\sqrt{-1}g(w)/r}$.
We set
$q_1:=(-\sqrt{-1}rdz)^r$
on $\cnum^{\ast}$.
We obtain
$q_0=\Phi^{\ast}(q_1)$.
Note that
$dg/dw-1=O(|w|^{-\epsilon_0})$
and
$\Image(g(w))-y=O(|w|^{1-\epsilon_0})$.
Then, the first claim is reduced to
\cite{Toda-lattice}.
(See Proposition \ref{prop;20.6.24.200}.)
We obtain
$|\theta(q_0)|_{h_0}=
O\Bigl(\bigl(\Image g(w)\bigr)^{-1}\Bigr)$.
Because $|f|=O\bigl(|e^{\sqrt{-1}g(w)}|\bigr)$,
we obtain
$|\theta(q)|_{h_0}=O\bigl(|\theta(q_0)|_{h_0}\bigr)$.
By the Hitchin equation
$R(h_0)+[\theta(q_0),\theta(q_0)^{\dagger}_{h_0}]=0$,
we obtain
$R(h_0)=O\Bigl(\bigl(\Image g(w)\bigr)^{-2}\Bigr)$.
Hence, we obtain
$|F(h_0,\theta(q))|_{h_0}
=O\Bigl(\bigl(\Image g(w)\bigr)^{-2}\Bigr)$.
Thus, we obtain the second claim.
\hfill\qed

\subsubsection{Comparison with auxiliary metrics}
\label{subsection;20.8.6.50}

Let $h_0$ be as in Lemma \ref{lem;20.6.15.100}.
For any relatively compact open subset $Y\subset X(R_0,\delta_0)$
satisfying Condition \ref{condition;20.6.15.10},
there uniquely exists $h_Y\in\Harm(q_{|Y})$ 
such that $h_{Y|\del Y}=h_{0|\del Y}$,
according to Proposition \ref{prop;20.5.29.20}.

\begin{lem}
\label{lem;20.8.15.2}
There exists $C^{(2)}_0>0$
which is independent of $Y$,
such that the following holds
on the region
$S_0(g,R^{(0)},B^{(0)},L^{(0)})\cap Y$:
\begin{equation}
|\theta(q)|_{h_Y}
\leq C^{(2)}_0\Image(g(w))^{-1}.
\end{equation}
 \end{lem}
\pf
Because $h_{Y|\del Y}=h_{0|\del Y}$,
there exists $C^{(2)}_1>0$,
which is independent of $Y$,
such that
the following holds
on $\del Y\cap S_0(g,R_{-1}^{(0)},B_{-1}^{(0)},L_{-1}^{(0)})$:
\[
|\theta(q)|_{h_Y}
\leq 
C^{(2)}_1(\Image(g(w))^{-1}).
\]

We set
{\small
\begin{multline}
 D:=\\
 \frac{1}{10}
 \min\bigl\{
  B^{(0)}-B^{(0)}_{-1},
  R^{(0)}-R^{(0)}_{-1},
  R^{(0)}_{-1}\sin(L^{(0)}_{-1}-L^{(0)}),
  R^{(0)}_{-1}\sin(L^{(0)}_{-1})-B^{(0)}
 \bigr\}.
\end{multline}
}

Take any $w_1\in S_0(g,R^{(0)},B^{(0)},L^{(0)})\cap Y$
such that
\[
 B^{(0)}\leq \Image(g(w_1))\leq B^{(0)}+D.
\]
We set
\begin{multline}
 \nbigb(w_1):=\bigl\{
 w\in X(R_0,\delta_0)\,\big|\,
 |g(w)-g(w_1)|<D
 \bigr\}
 \simeq
\\
 \bigl\{\zeta\in\cnum\,\big|\,
 |\zeta-g(w_1)|<D
 \bigr\}.
\end{multline}
Note that
$\nbigb(w_1)\subset S_0(g,R^{(0)}_{-1},B^{(0)}_{-1},L^{(0)}_{-1})$.
By applying Proposition \ref{prop;20.6.13.32}
to
$h_{Y|Y\cap \nbigb(w_1)}$,
we obtain that there exists
$C_2^{(2)}>0$,
which is independent of $Y$ and $w_1$ as above,
such that
the following holds on
$Y\cap
\bigl\{w\in X(R_0,\delta_0)\,\big|\,
|g(w)-g(w_1)|<D/2
\bigr\}$:
\[
 \bigl|\theta(q)\bigr|_{h_Y}
 \leq C_2^{(2)}.
\]

Take any $w_1\in S_0(g,R^{(0)},B^{(0)},L^{(0)})$
such that
$\Image g(w_1)>D+B^{(0)}$
and that
$\arg g(w_1)<L^{(0)}_{-1}/2$.
We set
$\rho(w_1):=(\Image g(w_1)-B^{(0)})/2$.
Then,
$\nbigb(w_1):=
\bigl\{
w\in X(R_0,\delta_0)\,\big|\,
|g(w)-g(w_1)|<\rho(w_1)
\bigr\}$
is contained in the following set:
\begin{multline}
\label{eq;20.8.15.1}
S_0(g,R^{(0)}_{-1},B^{(0)}_{-1},L^{(0)}_{-1})
\cup \\
\bigl\{
w\in X(R_0,\delta_0)\,\big|\,
\Image g(w)\geq B^{(0)}_{-1},
|g(w)|\leq 2R^{(0)}_{-1}
\bigr\}.
\end{multline}
Let $\Psi_{w_1}$ be the isomorphism
$\bigl\{\zeta\in\cnum\,\big|\,
|\zeta-g(w_1)|<\rho(w_1)
 \bigr\}
\simeq
\{|\eta|<1\}$
determined by
$\Psi_{w_1}(\zeta)
 =
\rho(w_1)^{-1}
 \bigl(\zeta-g(w_1)\bigr)$.
Note that
there exists $C^{(2)}_3\geq 1$,
which is independent of $Y$ and $w_1$ as above,
such that
\[
(C^{(2)}_3)^{-1}\leq
\frac{
\bigl|
(\Psi_{w_1}\circ g)^{\ast}(d\eta)/dw
\bigr|}
{\Image g(w_1)-B^{(0)}}
\leq
C^{(2)}_3.
\]
There exists $q_{w_1}$ on
$\bigl\{
|\eta|<1
\bigr\}$
such that
$(\Psi_{w_1}\circ g)^{\ast}(q_{w_1})
=q$
on $\nbigb(w_1)$.
By (\ref{eq;24.1.6.3}),
there exists $C^{(2)}_4>0$,
which is independent of $Y$ and $w_1$ as above,
such that
$|q_{w_1}|\leq C^{(2)}_4$.
We also obtain
\[
h_{w_1}\in\Harm\bigl(
q_{w_1|\{|\eta|<1\}\cap (\Psi_{w_1}\circ g)(Y)}
\bigr)
\]
such that
$(\Psi_{w_1}\circ g)^{\ast}(h_{w_1})
=h_Y$ on
$Y\cap \nbigb(w_1)$.
By applying Proposition \ref{prop;20.6.13.32}
to $h_{w_1|\{|\eta|<1\}\cap (\Psi_{w_1}\circ g)(Y)}$,
we obtain that there exists
$C^{(2)}_5>0$,
which is independent of $Y$ and $w_1$ as above,
such that the following holds
on 
$\{|\eta|<1/2\}\cap (\Psi_{w_1}\circ g)(Y)$:
\[
\bigl|
 \theta(q_{w_1})
 \bigr|_{h_{w_1}}
 \leq
 C^{(2)}_5.
\]
Hence, there exists $C^{(2)}_{6}>0$,
which is independent of $Y$ and $w_1$ as above,
such that the following holds
on $\bigl\{
 |g(w)-g(w_1)|<\rho(w_1)/2
\bigr\}\cap Y$:
\[
 \bigl|
 \theta(q)
 \bigr|_{h_Y}
 \leq C^{(2)}_6
\bigl(\Image g(w)\bigr)^{-1}.
\]

Take any $w_1\in S_0(g,R^{(0)},B^{(0)},L^{(0)})$
such that
$\arg g(w_1)>L^{(0)}/2$.
We set
\[
\rho(w_1):=
 |g(w_1)|
 \sin\bigl(
 (L^{(0)}_{-1}-\arg w_1)/2
 \bigr).
\]
Then,
$\nbigb(w_1)=
\Bigl\{
w\in X(R_0,\delta_0)\,\big|\,
 |g(w)-g(w_1)|<
\rho(w_1)
\Bigr\}$
is contained in
(\ref{eq;20.8.15.1}).
Let $\Psi_{w_1}$ be the isomorphism
$\bigl\{\zeta\in\cnum\,\big|\,
|\zeta-g(w_1)|<
\rho(w_1)
 \bigr\}
\simeq
\{|\eta|<1\}$
determined by
$\Psi_{w_1}(\zeta)
=\rho(w_1)^{-1}(\zeta-g(w_1))$.
Note that
there exists $C^{(2)}_7\geq 1$,
which is independent of $Y$ and $w_1$
such that
\[
(C^{(2)}_7)^{-1}\leq
\frac{
\bigl|
(\Psi_{w_1}\circ g)^{\ast}(d\eta)/dw
\bigr|}
{|g(w_1)|}
\leq
C^{(2)}_7.
\]
There exists $q_{w_1}$ on
$\bigl\{
|\eta|<1
\bigr\}$
such that
$(\Psi_{w_1}\circ g)^{\ast}(q_{w_1})
=q$
on $\nbigb(w_1)$.
By (\ref{eq;24.1.6.3}),
there exists $C^{(2)}_8>0$,
which is independent of $Y$ and $w_1$,
such that
$|q_{w_1}|\leq C^{(2)}_8$.
We also obtain
$h_{w_1}\in\Harm\bigl(
q_{w_1|\{|\eta|<1\}\cap (\Psi_{w_1}\circ g)(Y)}
\bigr)$
such that
$(\Psi_{w_1}\circ g)^{\ast}(h_{w_1})
=h_Y$ on
$Y\cap \nbigb(w_1)$.
By applying Proposition \ref{prop;20.6.13.32}
to $h_{w_1|\{|\eta|<1\}\cap (\Psi_{w_1}\circ g)(Y)}$,
we obtain that there exists
$C^{(2)}_9>0$,
which is independent of $Y$ and $w_1$,
such that the following holds
on 
$\{|\eta|<1/2\}\cap (\Psi_{w_1}\circ g)(Y)$:
\[
\bigl|
 \theta(q_{w_1})
 \bigr|_{h_{w_1}}
 \leq
 C^{(2)}_9.
\]
Hence, there exists $C^{(2)}_{10}>0$
which is independent of $Y$ and $w_1$
such that the following holds
on $\bigl\{
 |g(w)-g(w_1)|<
 \rho(w_1)/2
\bigr\}\cap Y$:
\[
 \bigl|
 \theta(q)
 \bigr|_{h_Y}
 \leq C^{(2)}_{10}|g(w)|^{-1}.
\]
Note that there exists
$C^{(2)}_{11}\geq 1$,
such that
$(C^{(2)}_{11})^{-1}|g(w)|
\leq\Image g(w)
\leq
C^{(2)}_{11}|g(w)|$
on
$S_0(g,R^{(0)},B^{(0)},L^{(0)})
\cap \{\arg(g(w))>L^{(0)}/2\}$.
Thus, we obtain the claim of Lemma \ref{lem;20.8.15.2}.
\hfill\qed

\begin{lem}
 \label{lem;20.6.9.10}
 There exists an $\real_{\geq 0}$-valued
 harmonic function $\phi$ on
 $\{\Image(g(w))\geq B_0\}$
 such that the following holds.
\begin{itemize}
 \item
      $\phi=O\bigl(|g(w)|^{1-\epsilon_0}\bigr)$.
 \item For any $Y$ satisfying Condition
  {\rm\ref{condition;20.6.15.10}},
      let $s_Y$ be the automorphism of
       $\hyperk_{Y,r}$ determined by
       $h_Y=h_{0|Y}\cdot s_Y$.
       Then, we obtain
       $\log \Tr(s_Y)\leq \phi$ on
       $Y\cap \{\Image(g(w))\geq B_1^{(0)}\}$.
\end{itemize}
\end{lem}
\pf
We recall that
$\bigl|F(h_0,\theta(q))\bigr|_{h_0}=
O\bigl((\Image g(w))^{-2}\bigr)$
by Lemma \ref{lem;20.6.15.100}.
Hence,
there exists an $\real_{\geq 0}$-valued function $\phi_1$
on $\{\Image g(w)\geq B_0\}$
such that
(i) $\phi_1=O\bigl(\log(1+\Image g(w))\bigr)$,
(ii) $\sqrt{-1}\Lambda\delbar\del \phi_1=
 \bigl|F(h_0,\theta(q))\bigr|$.

By Corollary \ref{cor;20.4.17.10},
Lemma \ref{lem;20.6.24.2}
and Lemma \ref{lem;20.8.15.2},
there exist $C^{(2)}_i$ $(i=20,21)$,
which are independent of $Y$,
such that the following holds
on $Y\cap S_0(g,R^{(0)},B^{(0)},L^{(0)})\setminus\nbigc_{0}$:
\begin{equation}
 \log\Tr(s_{Y})
 \leq
 -2\Re\gminia(f,\theta_0)
  +C^{(2)}_{20}|g(w)|^{\epsilon_0}
  +C^{(2)}_{21}.
\end{equation}
Here, $\theta_0>0$ is sufficiently small.
There exists $\beta\in\cnum^{\ast}$
such that
$\Re(\beta g(w)^{\epsilon_0})>0$
on $\{\Image g(w)\geq B_0\}$.
Then, there exist
$C^{(2)}_i$ $(i=22,23)$
which are independent of $Y$,
such that the following inequality holds
on $Y\cap S_0(g,R^{(0)},B^{(0)},L^{(0)})\setminus\nbigc_{0}$:
\begin{equation}
 \label{eq;20.8.4.100}
 \log\Tr(s_{Y})-\phi_1
 \leq
 -2\Re\gminia(f,\theta_0)
  +C^{(2)}_{22}\Re(\beta g(w)^{\epsilon_0})
  +C^{(2)}_{23}.
\end{equation}
Let $\nbigd$ be a connected component of
$\nbigc_{0}$ such that
\[
 \nbigd\cap S_0(g,R^{(0)}_{1},B^{(0)}_1,L^{(0)}_1)\neq\emptyset.
\]
Then, $\nbigd$ is relatively compact
in $Y\cap S_0(g,R^{(0)},B^{(0)},L^{(0)})$.
Because $\log\Tr(s_Y)-\phi_1$ is subharmonic
on $Y\cap S_0(g,R^{(0)},B^{(0)},L^{(0)})$,
the estimate (\ref{eq;20.8.4.100}) holds
on $\nbigd$.
Hence, (\ref{eq;20.8.4.100}) holds
on the closure of $Y\cap S_0(g,R^{(0)}_{1},B^{(0)}_1,L_1^{(0)})$.
Similarly,
there exist $C^{(2)}_{i}>0$ $(i=24,25)$
which are independent of $Y$
such that the following holds
on the closure of
$Y\cap S_{\pi}(g,R^{(0)}_{1},B_1^{(0)},L^{(0)}_1)$:
\begin{equation}
 \label{eq;20.8.4.101}
 \log\Tr(s_{Y})-\phi_1
 \leq
 -2\Re\gminia(f,\theta_{\pi})
  +C^{(2)}_{24}\Re(\beta g(w)^{\epsilon_0})
  +C^{(2)}_{25}.
\end{equation}
Here $\pi-\theta_{\pi}>0$ is sufficiently small.

As in Proposition \ref{prop;20.6.15.30},
there exists $C^{(2)}_{26}>0$
which is independent of $Y$ such that
the following holds
on $\bigl\{w\in X(R_0,\delta_0)\,\big|\,
\Image g(w)=B^{(0)}_1,\,\,|g(w)|\leq 2R^{(0)}_{-1}\bigr\}
\cap Y$:
\begin{equation}
\label{eq;20.8.5.2}
 \log\Tr(s_Y)-\phi_1
 \leq C^{(2)}_{26}.
\end{equation}
There exists $\gamma\in\cnum^{\ast}$
such that
$\Re(\gamma g(w)^{1-\epsilon_0})>0$
on
\[
 \{w\in X(R_0,\delta_0)\,|\,\Image g(w)\geq B_0\}.
\]
There exists $C^{(2)}_{27}>0$
which is independent of $Y$
such that the following holds
on $Y\cap\{\Image g(w)=B^{(0)}_1\}$:
\begin{equation}
 \label{eq;20.8.5.3}
  \bigl|
  \Re\gminia(f,\theta_0)
  \bigr|
  +\bigl|
  \Re\gminia(f,\theta_{\pi})
  \bigr|
 +\bigl|\Re(\beta g(w)^{\epsilon_0})\bigr|
 \leq C^{(2)}_{27}\Re(\gamma g(w)^{1-\epsilon_0}).
\end{equation}
By (\ref{eq;20.8.4.100}),
 (\ref{eq;20.8.4.101}),
 (\ref{eq;20.8.5.2})
 and  (\ref{eq;20.8.5.3}),
there exist
$C^{(2)}_{i}>0$ $(i=28,29)$,
which are independent of $Y$,
such that
the following holds
on $Y\cap \{\Image g(w)=B^{(0)}_1\}$:
\[
  \log\Tr(s_{Y})-\phi_1
 \leq
  C^{(2)}_{28}\Re(\gamma g(w)^{1-\epsilon_0})
 +C^{(2)}_{29}=:\phi_2.
\]
 
We obtain the following on
$Y\cap \{\Image g(w)\geq B_0\}$:
\[
-\del_{w}\del_{\wbar}
\bigl(
\log\Tr(s_Y)-\phi_1-\phi_2
\bigr)\leq 0.
\]
Note that $\log\Tr(s_Y)=\log r$ on $\del Y$,
and hence
\[
 \log\Tr(s_Y)-\phi_1-\phi_2
 \leq
 \log\Tr(s_Y)
 \leq
 \log r
\]
on $\del Y\cap \{\Image g(w)\geq B^{(0)}_1\}$.
We also obtain
\[
 \log\Tr(s_Y)-\phi_1-\phi_2
 \leq
 0\leq \log r
\]
on $Y\cap \{\Image g(w)=B^{(0)}_1\}$.
Hence, we obtain
\[
 \log\Tr(s_Y)-\phi_1-\phi_2
 \leq \log r
\]
on $\del\bigl(Y\cap \{\Image g(w)\geq B^{(0)}_1\}\bigr)$.
Hence, we obtain
$\log\Tr(s_Y)\leq
 \log r+\phi_1+\phi_2$
 on $Y\cap\{\Image g(w)\geq B^{(0)}_1\}$.
 Because $\phi_1=O\bigl(\log(1+|g(w)|)\bigr)$,
 there exists $C^{(2)}_{30}>0$
 such that
 $\phi_1+\phi_2+\log r
 \leq C^{(2)}_{30}\Re(\gamma g(w)^{1-\epsilon_0})=:\phi$
 on $\{\Image g(w)\geq B_0\}$.
Thus, we obtain Lemma  \ref{lem;20.6.9.10}.
 \hfill\qed

\subsubsection{Construction of harmonic metrics}
\label{subsection;20.8.6.31}

\begin{prop}
\label{prop;20.6.9.24}
 For any $\veca\in\nbigp$,
 there exists $h\in\Harm(q)$
 such that $\veca(h)=\veca$.
\end{prop}
\pf
For $\veca\in\nbigp$,
we take an auxiliary metric $h_0$
as in Lemma \ref{lem;20.6.15.100}.
Let $Y_1\subset Y_2\subset\cdots$
be a smooth exhaustive family of $X(R_0,\delta_0)$
satisfying Condition \ref{condition;20.6.15.10}.
According to Proposition \ref{prop;20.5.29.20},
we obtain a sequence
$h_{Y_i}\in\Harm(q_{|Y_i})$
such that $h_{Y_i|\del Y_i}=h_{0|\del Y_i}$.
According to Proposition \ref{prop;20.5.29.1},
by taking a subsequence,
we may assume that
the sequence $h_{Y_i}$ is convergent to
$h_{\infty}\in\Harm(q)$.
Lemma  \ref{lem;20.6.15.100}
and Lemma \ref{lem;20.6.9.10}
imply $\veca(h_{\infty})=\veca$.
Thus, the proof of Proposition \ref{prop;20.6.9.24}
is completed.
\hfill\qed

\subsubsection{Comparison with harmonic metrics}
\label{subsection;20.8.6.40}

Let $B_0<B^{(0)}_{-1}<B^{(0)}_0<B^{(0)}_1$
be as in \S\ref{subsection;24.1.6.21}.
We can prove the following proposition
by the argument of Lemma \ref{lem;20.6.9.10}.

\begin{prop}
\label{prop;20.6.15.31}
 There exist
 a harmonic function $\phi$ on
 $\{w\in X(R_0,\delta_0)\,|\,\Image g(w)\geq  B_0\}$
 such that the following holds:
 \begin{itemize}
  \item $\phi=O\bigl(|g(w)|^{1-\epsilon_0}\bigr)$.
  \item  Let $h\in\Harm(q_{|\{\Image g(w)> B_0\}})$.
Let $Y$ be a relatively compact subset of
$X(R_0,\delta_0)$
satisfying Condition {\rm\ref{condition;20.6.15.10}}.
 Suppose that $h_Y\in\Harm(q_{|Y})$
 satisfies
\[
 h_{Y|\del Y\cap \{\Image(g(w))> B_0\}}
 =h_{|\del Y\cap\{\Image g(w)> B_0\}}.
\]
 Let $s_Y$ be the automorphism of
 $\hyperk_{Y\cap \{\Image g(w)>B_0\},r}$
 determined by
	 $h_{Y| Y\cap\{\Image g(w)>B_0\}}
	 =h_{|Y\cap\{\Image g(w)>B_0\}}\cdot s_Y$.
	 Then, we obtain $\log\Tr(s_Y)\leq \phi$
	 on
	 $Y\cap \{\Image (g(w))>B^{(0)}_1\}$.
 \end{itemize}
\end{prop}
\pf
In this proof,
constants are independent of $h$ and $Y$.
Let $h$, $h_Y\in\Harm(q_{|Y})$
and $s_Y$ be as in the statement.
By Proposition \ref{prop;20.4.15.11},
there exists a constant $C_1>0$
such that the following holds
on $\{w\in X(R_0,\delta_0)\,|\,\Image g(w)\geq B^{(0)}_{-1}\}$:
\[
 |\theta(q)|_{h}
 \leq C_1\bigl(
  \Image g(w)
  \bigr)^{-1}.
\]
By the same argument as
the proof of Lemma \ref{lem;20.8.15.2},
we can prove that
there exists a constant $C_2>0$
such that
the following inequality holds on
$Y\cap S_0(g,R^{(0)},B^{(0)},L^{(0)})$:
\[
 |\theta(q)|_{h_Y}
 \leq
 C_2\bigl(\Image g(w)\bigr)^{-1}.
\]
By Corollary \ref{cor;20.4.21.2}
and Lemma \ref{lem;20.6.24.2},
there exist constants $C_i>0$ $(i=3,4)$
such that the following inequality holds
on $Y\cap S_0(g,R^{(0)},B^{(0)},L^{(0)})\setminus\nbigc_{0}$:
\begin{equation}
\label{eq;20.8.6.40}
 \log\Tr(s_{Y})
 \leq
 -2\Re\gminia(f,\theta_0)
  +C_3\Re(\beta g(w)^{\epsilon_0})
  +C_4.
\end{equation}
Here, $\theta_0>0$ is sufficiently small,
and $\beta$ is a complex number
such that
$\Re(\beta g(w)^{\epsilon_0})>0$
on $\{\Image g(w)\geq B_0\}$.
Note that $\log\Tr(s_Y)$ is subharmonic
by (\ref{eq;20.8.16.6}).
Hence,
the inequality (\ref{eq;20.8.6.40}) holds
on the closure of
$Y\cap S_0(g,R^{(0)}_1,B^{(0)}_1,L^{(0)}_1)$.
Similarly,
there exist constants $C_i>0$  $(i=5,6)$
such that the following holds
on the closure of
$Y\cap S_{\pi}(g,R^{(0)}_{1},B^{(0)}_1,L^{(0)}_1)$:
\begin{equation}
 \log\Tr(s_{Y})
 \leq
 -2\Re\gminia(f,\theta_{\pi})
  +C_5\Re(\beta g(w)^{\epsilon_0})
  +C_6.
\end{equation}
Here, $\pi-\theta_{\pi}>0$ is sufficiently small.
As in Proposition \ref{prop;20.6.15.30},
there exists a constant $C_7>0$
such that
$\log\Tr(s_Y)\leq C_7$
on $\{\Image g(w)=B^{(0)}_1,\,\,|g(w)|\leq 2R^{(0)}_{-1}\}
\cap Y$.
Therefore,
there exist positive constants $C_i$ $(i=8,9)$
such that the following holds
on $Y\cap \{\Image g(w)=B^{(0)}_1\}$:
\[
  \log\Tr(s_{Y})
 \leq
  C_8\Re(\gamma g(w)^{1-\epsilon_0})
 +C_9=:\phi_0.
\]
Here $\gamma$ is a complex number
such that
$\Re(\gamma g(w)^{1-\epsilon_0})>0$
on $\{\Image g(w)\geq B_0\}$.

We set $\phi=\phi_0+\log (r)$.
Note that $\log\Tr(s_Y)-\phi$ is a subharmonic function.
By the construction,
we obtain
$\log\Tr(s_Y)-\phi\leq 0$
on $\del(Y\cap\{\Image g(w)\geq B_1^{(0)}\})$,
and hence 
on $Y\cap\{\Image g(w)\geq B_1^{(0)}\}$.
Thus, we obtain Proposition \ref{prop;20.6.15.31}.
\hfill\qed

\subsubsection{Proof of Theorem \ref{thm;20.6.14.10}
and Proposition \ref{prop;20.6.14.11}}
\label{subsection;20.8.5.20}

We consider the map $\Phi:X(R,\delta)\lrarr \cnum$
defined by $w=-\sqrt{-1}\alpha z^{-\rho}$.
If $R$ is sufficiently large
and $\delta$ is sufficiently small,
we obtain
$\Phi:X(R,\delta)\lrarr V$,
and it is a holomorphic embedding.
By applying Proposition \ref{prop;20.6.9.22}
and Proposition \ref{prop;20.6.9.23}
to $\Phi^{\ast}(q)=f\cdot (dw)^r$,
we obtain the first two claims of
Theorem \ref{thm;20.6.14.10}.

Let $\veca\in\nbigp$.
By Proposition \ref{prop;20.6.9.24},
there exists $h_{0,\veca}\in\Harm(\Phi^{\ast}(q))$
such that $\veca(h_{0,\veca})=\veca$.
It induces
$h_{0,\veca}\in
\Harm(q_{|\Phi(X(R,\delta))})$.
We extend it to a $G_r$-invariant Hermitian metric
of $\hyperk_{V\setminus\varpi^{-1}(0),r}$ such that
$\det(h_{1,\veca})=1$.
Let $\{X_i\}$ be a smooth exhaustive family of $V$
such that
each
$\Phi^{-1}(X_i)$ satisfies
Condition \ref{condition;20.6.15.10}.
Let $h_{i}\in\Harm(q)$
such that $h_{i|\del X_i}=h_{1,\veca|\del X_i}$.
According to Proposition \ref{prop;20.6.15.30},
we may assume that the sequence $h_i$  is convergent,
and we obtain $h_{\infty}\in\Harm(q)$
as the limit of a subsequence of $h_i$.
By Proposition \ref{prop;20.6.15.31},
$h_{\infty}$ satisfies $\veca_I(h_{\infty})=\veca$.
Thus, we obtain the first claim of Theorem \ref{thm;20.6.14.10}.

For the proof of Proposition \ref{prop;20.6.14.11},
we use the notation in \S\ref{subsection;24.1.6.21}.
We set $\nbigu_{I,1}=\Phi(\{\Image(g(w))>B_0\})$
and $\nbigu_{I,2}=\Phi(\{\Image(g(w))>B_1^{(0)}\})$.
Let $\{K_i\}$ be a smooth exhaustive family of
$V\setminus\varpi^{-1}(0)$
such that each $\Phi^{-1}(K_i)$  satisfies
Condition \ref{condition;20.6.15.10}.
Let $\phi_I$ denote the function on
$\nbigu_{I,1}$ induced by $\phi$
in Proposition \ref{prop;20.6.15.31}.
Then, we obtain 
Proposition \ref{prop;20.6.14.11}
from Proposition \ref{prop;20.6.15.31}.
\hfill\qed

\subsection{Refined estimate in an easy case}

Let $f$ be a section of $\gbigb$
on $V\subset\cnumtilde$.
Let $I$ be an interval of $V\cap\varpi^{-1}(0)$
such that $\Ibar\subset V$
and that $I$ is special with respect to $f$.
We assume the following.
\begin{condition}
 There exist a neighbourhood $\nbigu$ of $\Ibar$ in $V$,
a finite sum
$\gminia_I=
 \sum_{0<\kappa\leq\rho}\gminia_{I,\kappa}z^{\kappa}$
 $(\gminia_{I,\rho}\neq 0)$,
 $\alpha\in\cnum^{\ast}$,
 $\ell_I\in\real$ and
 $C>1$
 such that
$\bigl|
 e^{-\gminia_I}z^{-\ell_I}
 (f-\alpha e^{\gminia_I}z^{\ell_I})
 \bigr|
 \leq C|z|^{\epsilon}$
 for $\epsilon>0$
 on $\nbigu\setminus\varpi^{-1}(0)$.
\end{condition}

Recall $\deg(\gminia_I)=\rho$.
We set $q=f\,(dz/z)^r$
and
 \[
 \gminid_I(q):=\frac{\ell_I}{\deg(\gminia_I)}+r.
 \]
 
Let $h\in\Harm(q)$.
We obtain $\veca_I(h)\in\nbigp$
as in Theorem \ref{thm;20.6.14.10}.
In this easier case,
we can obtain the estimate on the behaviour of $h$
up to boundedness.
Let $\veck_I(h)$ be the tuple of integers
obtained from $\veca_I(h)$ as in \S\ref{subsection;20.6.24.200}.
We shall prove the following proposition
in \S\ref{subsection;20.6.26.1}--\ref{subsection;20.6.26.2}.

\begin{prop}
\label{prop;20.6.26.3}
The following estimates holds
as $-\Re(\gminia_I(z))+\gminid_I(q)\log|\gminia_I(z)|\to\infty$:
\begin{multline}
 \log\Bigl|
 (dz)^{(r+1-2i)/2}
 \Bigr|_h
-a_{I,i}(h)\Bigl(-\Re(\gminia_I(z))+\gminid_I(q)\log|\gminia_I(z)|
 \Bigr)
\\
 +(\deg(\gminia_I)+1)\left(
 \frac{r+1}{2}-i
 \right)\log|z|
\\
 -\frac{k_{I,i}}{2}
 \log\Bigl(
 -\Re(\gminia_I(z))+\gminid_I(f)\log|\gminia_I(z)|
 \Bigr)
 =O(1).
\end{multline}
\end{prop}

\subsubsection{Normalization}
\label{subsection;20.6.26.1}

Let $\nbigu$ denote an open neighbourhood of $\Ibar$.
We define the map
$F:\nbigu\setminus\varpi^{-1}(0)\lrarr \cnum$
by
\[
 F(z)=-\sqrt{-1}
\Bigl(
 \gminia_I(z)
 -\gminid_I(q)\log\bigl(-\sqrt{-1}\gminia_I(z)\bigr)
\Bigr).
\]
Let $\varpi_{\infty}:\projtilde^1_{\infty}\lrarr\proj^1$
denote the oriented real blowing up at $\infty$.
We regard $\cnum$ as an open subset of 
$\projtilde^{1}_{\infty}$.
Let $w$ denote the standard coordinate system on $\cnum$.
By using the polar decomposition of $w$,
we identify
$\projtilde^1_{\infty}\setminus\{0\}$
with
$(\real_{>0}\cup\{\infty\})\times S^1$.

\begin{lem}
There exist a neighbourhood $\nbigu'$ of $\Ibar$ in $\nbigu$
and a neighbourhood $\nbigv$
of $\{\infty\}\times\closedclosed{0}{\pi}$
in $\projtilde^1_{\infty}$
such that $F$ induces a homeomorphism
$\nbigu'\simeq\nbigv$.
\end{lem}
\pf
It is easy to see that
$F$ extends to a continuous map
$\nbigu\lrarr\projtilde^1_{\infty}$,
which is also denoted by $F$.
Moreover, we can easily check that
$F$ induces a homeomorphism of
a neighbourhood of $\Ibar$ in $\varpi^{-1}(0)$
and a neighbourhood of $\{\infty\}\times\closedclosed{0}{\pi}$
in $\varpi_{\infty}^{-1}(\infty)$.

We set $\rho:=\deg(\gminia_I)$.
On a neighbourhood of $\Ibar$,
we use the real coordinate system
given by
$(|z|^{\rho},\arg(z))$.
On a neighbourhood of $\{\infty\}\times\closedclosed{0}{\pi}$,
we use the real coordinate system
given by
$(|w|^{-1},\arg(w))$.
Then, it is easy to check that
$F$ is of $C^1$-class,
and the tangent map at each point of
$\Ibar$ is an isomorphism.
Then, the claim follows from
the inverse function theorem.
\hfill\qed

\vspace{.1in}
Then, $(F^{-1})^{\ast}(q)$ is expressed as follows:
\[
 (F^{-1})^{\ast}(q)
 =\beta e^{\sqrt{-1}w}(1+\upsilon(w))\,(dw)^r.
\]
Here,
$\upsilon$ is a holomorphic function
satisfying
$\upsilon(w)=O(|w|^{-\epsilon})$
for some $\epsilon>0$,
and 
$\beta$ is a non-zero complex number.

\subsubsection{Comparison with the model metric}
\label{subsection;20.6.26.2}

Let $\delta>0$ and $R>0$.
Let $f_0$ be a nowhere vanishing holomorphic function
on $X(R,\delta)$
such that
$\log|f_0|=-\Image w+O(1)$.
  We set
  $q_0:=f_0\,(dw)^r$.
According to Theorem \ref{thm;20.6.14.10},
for any $h_0\in\Harm(q_0)$,
there exist $\veca(h_0)\in\nbigp$
and $0<\rho<1$
such that the following estimates hold
on $\bigl\{|w|>R,\,\,\big|\arg(w)-\pi/2|<(1-\delta_1)\pi/2\bigr\}$
for any $0<\delta_1<1$:
\[
 \log\bigl|(dw)^{(r+1-2i)/2}\bigr|_{h_0}
 -a_i(h_0)\Image(w)
 =O\bigl(|w|^{\rho}\bigr).
\]
Note that
$\log|f_0|=O(1)$
on $\bigl\{w\in X(R+1,\delta)\,\big|\,|\Image(w)|<A\bigr\}$
for any $A>0$
in this case.
By Corollary \ref{cor;20.6.11.1},
we obtain that
$\log\bigl|(dw)^{(r+1-2i)/2}\bigr|_{h_0}$
is bounded
on $\bigl\{w\in X(R+1,\delta)\,\big|\,|\Image(w)|<A\bigr\}$
for any $A>0$.
By Proposition \ref{prop;20.4.22.11},
we obtain the following estimate
on $\{w\in X(R+1,\delta)\,|\,\Image(w)\geq 1\}$:
\[
 \log\bigl|(dw)^{(r+1-2i)/2}\bigr|_{h_0}
 -a_i(h_0)\Image(w)
 =O\bigl(\log(1+\Image(w))\bigr).
\]

 Let $f_1$ be a holomorphic function
 on $X(R,\delta)$
 such that there exist $C_{1}>0$ and $\kappa>0$
 such that
 \[
  |f_0-f_1|\leq C_{1}|w|^{-\kappa}|f_0|.
 \]
 We obtain $\log|f_1|=-\Image(w)+O(1)$.
We set $q_1:=f_1\,(dw)^r$.
Similarly,
for $h_1\in\Harm(q_1)$,
there exists
$\veca(h_1)\in\nbigp$
such that the following estimate holds
on $\bigl\{w\in X(R+1,\delta)\,\big|\,\Image(w)\geq 1\bigr\}$:
\[
 \log\bigl|(dw)^{(r+1-2i)/2}\bigr|_{h_1}
 -a_i(h_1)\Image(w)
 =O\bigl(\log(1+\Image(w))\bigr).
\]
  
 \begin{prop}
\label{prop;20.4.28.41}
  Suppose that there exist $h_0\in\Harm(q_0)$
  and $h_1\in\Harm(q_1)$
  such that
  $\veca(h_0)=\veca(h_1)$.
  Then, $h_0$ and $h_1$ are mutually bounded
 on $\bigl\{w\in X(R+1,\delta)\,\big|\,\Image(w)\geq 1\bigr\}$.
 \end{prop}
 \pf
 We set $e_i:=(dw)^{(r+1-2i)/2}$ for $i=1,\ldots,r$.
 Note that
 $\theta(q_1) e_i=\theta(q_0)e_i$ $(i=1,\ldots,r)$,
 $\theta(q_1)e_{r}=\theta(q_0)e_{r}+(f_1-f_0)e_1\,dw$
 and $\theta(q_0)(e_{r})=f_0\,e_1\,dw$.
 Because $|f_1-f_0|\leq C_{1}|w|^{-\kappa}|f_0|$,
 there exists $C_{2}>0$ such that
 $|\theta(q_1)-\theta(q_0)|_{h_0}
 \leq C_2|w|^{-\kappa}|\theta(q_0)|_{h_0}$.
 Hence,
 on $\bigl\{w\in X(R+1,\delta)\,\big|\,\Image(w)\geq 1\bigr\}$,
 we obtain
 $|\theta(q_1)-\theta(q_0)|_{h_0}
 =O\bigl((1+\Image(w))^{-1-\kappa}\bigr)$
 and
 $|\theta(q_1)|_h=O\bigl((1+\Image(w))^{-1}\bigr)$.
 Because
 $R(h_0)+\bigl[\theta(q_0),\theta(q_0)_{h_0}^{\dagger}\bigr]=0$,
 we obtain
 \[
 \bigl|
  R(h_0)+\bigl[\theta(q_1),\theta(q_1)_{h_0}^{\dagger}\bigr]
  \bigr|_{h_0}
  =O\bigl((1+\Image(w))^{-2-\kappa}\bigr).
 \]
 Let $s$ be the automorphism of
 $\hyperk_{X(R,\delta),r}$
 determined by $h_1=h_0s$.
 By (\ref{eq;20.8.16.6}),
 we obtain the following on
 $\bigl\{w\in X(R+1,\delta)\,\big|\,\Image(w)\geq 0\bigr\}$:
\[
 -2\del_z\del_{\zbar}\log\Tr(s)\leq
 \bigl|
 R(h_0)+\bigl[\theta(q_1),\theta(q_1)_{h_0}^{\dagger}\bigr]
 \bigr|_{h_0}
 =O\bigl((1+\Image(w))^{-2-\kappa}\bigr).
\]
Note that $\Tr(s)$ is bounded on
$\bigl\{w\in X(R+1,\delta)\,\big|\,
 |\Image(w)|<A
 \bigr\}$.
 By Proposition \ref{prop;20.4.22.11},
 there exists $a\in\real$
 such that
 $|\log\Tr(s)-a\Image(w)|$ is bounded
 on
 $\bigl\{
  w\in X(R+1,\delta)\,\big|\,\Image(w)\geq 0
 \bigr\}$.
Because $\veca(h_0)=\veca(h_1)$,
we obtain $a=0$,
which implies
the boundedness of
$\log\Tr(s)$.
\hfill\qed

\vspace{.1in}

Let $\upsilon(w)$ be a holomorphic function on $X(R,\delta)$
such that $\upsilon(w)=O(|w|^{-\kappa})$.
Let $\alpha$ be a non-zero complex number.
We set $q=\alpha(1+\upsilon) e^{\sqrt{-1}w}(dw)^r$.
Let $h\in\Harm(q)$.
We obtain $\veca(h)\in\nbigp$
and a tuple of integers $\veck(h)$
as in \S\ref{subsection;20.6.24.200}.

\begin{cor}
\label{cor;20.6.26.4}  
 We obtain the following estimates
 for $i=1,\ldots,r$
 on $\bigl\{w\in X(R,\delta)\,\big|\,\Image(w)\geq 0\bigr\}$:
\[
\log \bigl|
  (dw)^{(r+1-2i)/2}
\bigr|_h
 -a_i(h)\Image(w)
 -\frac{k_i(h)}{2}\log\bigl(1+\Image(w)\bigr)
 =O(1).
\]
 \end{cor}
\pf
In the case $\upsilon=0$,
the claim follows from Proposition \ref{prop;20.6.24.200}
and Theorem \ref{thm;20.6.14.10}.
The general case follows from Proposition \ref{prop;20.4.28.41}.
\hfill\qed

\vspace{.1in}
We obtain Proposition \ref{prop;20.6.26.3}
by applying the normalization in \S\ref{subsection;20.6.26.1}
and Corollary \ref{cor;20.6.26.4}.
\hfill\qed

\section{Global results}

\subsection{Holomorphic $r$-differentials with multiple growth orders
on a punctured disc}

\subsubsection{General case}
\label{subsection;20.7.13.30}

Let $U$ be a neighbourhood of $0$ in $\cnum$.
Let $\varpi:\Utilde\lrarr U$
denote the oriented real blowing up.
Let $f$ be a section of $\gbigb_{\Utilde}$ on $\Utilde$.
Let $\nbigs(f)$ denote the set of the intervals
which are special with respect to $f$.
For any $I\in \nbigs(f)$,
we fix a branch of $\log z$ around $I$.
Then,
$\rho(I)\in\real_{>0}$
and $\alpha_I\in\cnum^{\ast}$ are determined as
\[
\deg(\gminia(f,\theta)-\alpha_Iz^{-\rho(I)})<\rho(I)
\]
for $\theta\in I\setminus\nbigz(f)$.
We set $q=f\,(dz)^r$.
We set $\nbigs(q):=\nbigs(f)$.
Let $U_0\subset U$ be a relatively compact neighbourhood of $0$
in $U$.
We set
$U^{\ast}=U\setminus\{0\}$ and
$U_0^{\ast}=U_0\setminus\{0\}$.

\begin{thm}
\mbox{{}}\label{thm;20.6.9.30}
\begin{itemize}
 \item 
 For any $h\in\Harm(q)$ and for any $I\in \nbigs(q)$,
 there exist $\veca_I(h)\in\nbigp$ and $\epsilon>0$
 such that the following estimates hold
 as $|z|\to 0$
 on $\bigl\{|\arg(\alpha _Iz^{-\rho(I)})-\pi|<(1-\delta)\pi/2
 \bigr\}$
for any $\delta>0$:
\[
 \log\bigl|
  (dz)^{(r+1)/2-i}
 \bigr|_h
+a_i(h)\Re\bigl(\alpha_I z^{-\rho(I)}\bigr)
=O\bigl(|z|^{-\rho(I)+\epsilon}\bigr).
 \]
\item
 If $h_i\in\Harm(q)$  $(i=1,2)$
 satisfy $\veca_I(h_1)=\veca_I(h_2)$
 for any $I\in \nbigs(q)$,
     then $h_1$ and $h_2$ are mutually bounded
     on $U_0^{\ast}$.
\item
     For any $(\veca_I)_{I\in \nbigs(q)}\in \prod_{I\in \nbigs(q)}\nbigp$,
     there exists $h\in\Harm(q)$
     such that $\veca_I(h)=\veca_I$.
\end{itemize}
\end{thm}
\pf
We obtain the first claims from
Theorem \ref{thm;20.6.14.10}.
We obtain the second claim from
Proposition \ref{prop;20.6.12.20},
Theorem \ref{thm;20.6.9.21}
and Theorem \ref{thm;20.6.14.10}.

Let us prove the third claim.
Let $(\veca_I)_{I\in \nbigs(q)}\in \prod_{I\in \nbigs(q)}\nbigp$.
For each $I\in \nbigs(q)$,
let $V_I$ be a relatively compact neighbourhood of
$\Ibar$ in $\Utilde$.
There exist
relatively compact neighbourhoods $\nbigu_{I,i}$ $(i=1,2)$,
a function $\phi_I$,
and a smooth exhaustive family $\{K_{I,i}\}$,
as in Proposition \ref{prop;20.6.14.11}.
By Theorem \ref{thm;20.6.14.10},
there exists a $G_r$-invariant Hermitian metric $h_0$ of
$\hyperk_{U\setminus\{0\},r}$
such that $\det(h_0)=1$
and that
$h_{0|\nbigu_{I,1}}\in\Harm(q_{|\nbigu_{I,1}})$
with $\veca_I(h_0)=\veca_I$
for $I\in \nbigs(q)$.
Let $\{X_i\}$ be a smooth exhaustive family for $U\setminus\{0\}$
such that
(i) $\coprod_{I}K_{I,i}\subset X_i$,
(ii) $K_{I,i}\cap \nbigu_{I,1}=X_i\cap\nbigu_{I,1}$
for any $I$.
We obtain $h_i\in\Harm(q_{|X_i})$
such that $h_{i|\del X_i}=h_{0|\del X_i}$.
We may assume that
the sequence $\{h_i\}$ is convergent,
and we obtain $h_{\infty}\in\Harm(q)$ as the limit.
By Proposition \ref{prop;20.6.14.11},
we obtain $\veca_I(h_{\infty})=\veca_I$
for $I\in \nbigs(q)$.
\hfill\qed

\subsubsection{Refined estimate in the nowhere vanishing case}

Let us state a refined result
in the case where $q$ is nowhere vanishing.

\begin{lem}
 Let $q$ be a holomorphic $r$-differential with multiple growth orders
 on $(U,0)$.
 If $q$ is nowhere vanishing,
 then there exist an integer $\ell(q)$
 and a meromorphic function $\gminia$ on $(U,0)$
 such that $q=z^{\ell(q)}e^{\gminia(z)}(dz/z)^r$.
\end{lem}
\pf
We describe $q$ as
$q=\beta (dz/z)^r$ for a holomorphic function $\beta$.
There exists an integer $\ell$
such that
$\gminia=\log(z^{-\ell}\beta)$
is well defined as a single-valued holomorphic function
on $U^{\ast}$.
Because $\beta$ induces a section of $\gbigb$,
we obtain $\Re(\gminia)=O\bigl(|z|^{-\rho}\bigr)$
for some $\rho>0$.
Because $\Re(\gminia)$ is harmonic,
there exists $\rho_1>0$ such that
$\del_x\Re(\gminia)=O(|z|^{-\rho_1})$
and
$\del_y\Re(\gminia)=O(|z|^{-\rho_1})$,
where $(x,y)$ is the real coordinate system
induced by $z=x+\sqrt{-1}y$.
By the Cauchy-Riemann equation,
we obtain that
$\del_x\Image(\gminia)=O(|z|^{-\rho_1})$
and
$\del_y\Image(\gminia)=O(|z|^{-\rho_1})$.
Hence, there exists $\rho_2>0$
such that 
$\Image(\gminia)=O(|z|^{-\rho_2})$.
Thus, we obtain that $\gminia$ is meromorphic at $z=0$.
\hfill\qed

\vspace{.1in}
For the description $q=z^{\ell(q)}e^{\gminia(z)}(dz/z)^r$,
we obtain the integer $m(q)$
such that
$z^{m(q)}\gminia(z)$ is holomorphic at $z=0$,
and $\alpha:=(z^{m(q)}\gminia(z))_{z=0}\neq 0$.
We set $\gminid(q):=\frac{\ell(q)}{m(q)}+r$.
We also set
\[
T(q):=\{0\leq \theta<2\pi\,|\,-m(q)\theta+\arg(\alpha)\equiv \pi,
\mod 2\pi
 \}.
\]
 Note that special intervals with respect to $q$
 are 
 $\bigl\{
  |\theta-\theta_0|<\pi/2m(q)
  \bigr\}$
  for $\theta_0\in T(q)$.

\begin{prop}
\label{prop;20.4.25.30}
  Let $h\in\Harm(q)$.
For any $\theta_0\in T(q)$,
 we obtain the tuple of real numbers
\[
 \veca(h,\theta_0)
 =(a_0(h,\theta_0),\ldots,a_{r-1}(h,\theta_0))
  \in\nbigp
\]
 determined by the following estimates
for $i=1,\ldots,r$
 as
 $-\Re\bigl(
 \gminia(z)
 \bigr)
 +\gminid(q)\log|\gminia(z)|\to\infty$
 on $\bigl\{|\arg(z)-\theta_0|<\pi/2m(q)\bigr\}$:
  \begin{multline}
 \log\Bigl|
 (dz)^{(r+1)/2-i}
 \Bigr|_h
-a_i(h,\theta_0)
 \Bigl(
-\Re\bigl(
 \gminia(z)
 \bigr)
\\
+\gminid(q)\log|\gminia(z)|
 \Bigr)
 +(m(q)+1)\left(\frac{r+1}{2}-i\right)\log|z|
 \\
 =O\left(
 \log\Bigl(
 -\Re(\gminia(z))
+\gminid(q)\log|\gminia(z)|
 \Bigr)
 \right).
\end{multline}
>From $\veca(h,\theta_0)\in\nbigp$,
we determine the integers $k_i(h,\theta_0)$ $(i=0,\ldots,r-1)$
as in {\rm\S\ref{subsection;20.6.24.200}}.
Then,
we obtain the following estimates
as
$-\Re\bigl(
 \gminia(z)
 \bigr)
 +\gminid(q)\log|\gminia(z)|\to\infty$
 on $\bigl\{|\arg(z)-\theta_0|<\pi/2m(q)\bigr\}$:
\begin{multline}
 \log\Bigl|
 (dz)^{(r+1)/2-i}
 \Bigr|_h
 -a_i(h,\theta_0)
 \Bigl(
 -\Re\bigl(
 \gminia(z)
 \bigr)
+\gminid(q)\log|\gminia(z)|
 \Bigr)
\\
 +(m(q)+1)\left(\frac{r+1}{2}-i\right)\log|z|
 \\
 -\frac{k_i(h,\theta_0)}{2}
 \log\Bigl(
 -\Re(\gminia(z))
 +\gminid(q)\log|\gminia(z)|
 \Bigr)
 =O\bigl(1\bigr).
\end{multline}
\end{prop}
\pf
It follows from Proposition \ref{prop;20.6.26.3}.
\hfill\qed
 
\subsection{The case of punctured Riemann surfaces}

\subsubsection{Statement}

Let $X$ be a Riemann surface.
Let $D$ be a finite subset of $X$.
For each $P\in D$,
let $(X_P,z_P)$ be a holomorphic coordinate neighbourhood
such that $z_P(P)=0$.
Set $X_P^{\ast}:=X_P\setminus\{P\}$.
Let $\varpi_P:\Xtilde_P\to X_P$ denote the oriented blow up
along $P$.

Let $q$ be an $r$-differential on $X$.
We obtain the expression
$q_{|X_P^{\ast}}=f_P\,(dz_P/z_P)^r$.
Let $D_{\mero}$ be the points $P$ of $D$
such that $f_P$ are meromorphic at $P$.
Let $D_{>0}$ denote the set of $P$
such that $f_P$ is holomorphic at $P$
and that $f_P(P)=0$.
We set $D_{\leq 0}:=D_{\mero}\setminus D_{>0}$.

For each $P\in D_{>0}$,
we have the description
$q_P=w_P^{m_P}f_P(dw_P/w_P)^r$,
where $f_P$ is holomorphic at $P$
such that $f_P(P)\neq 0$.
Let 
$\nbigp(q,P)$
denote the set of 
$\vecb=(b_1,b_2,\ldots,b_{r})\in\real^r$
satisfying
\[
b_1\geq b_2\geq \cdots \geq b_{r}\geq b_1-m_P,
\quad
\sum b_i+r(r+1)/2=0.
\]

We set $D_{\ess}:=D\setminus D_{\mero}$.
We assume that
for any $P\in D_{\ess}$,
$f_P$ is a section of $\gbigb_{\Xtilde_P}$.
At each $P\in D_{\ess}$,
we obtain the set of the intervals $\nbigs(q,P)$
in $\varpi_P^{-1}(P)$ which are special
with respect to $q_{|X_P^{\ast}}$.

Let $h\in\Harm(q)$.
For any $P\in D_{\ess}$ and $I\in \nbigs(q,P)$,
we obtain
$\veca_{P,I}(h)\in\nbigp$
as in Theorem \ref{thm;20.6.9.30}.
For any $P\in D_{>0}$,
we obtain
$\vecb(h)\in \nbigp(q,P)$
as in Proposition \ref{prop;20.7.1.20}.
Thus, we obtain the map
\begin{equation}
\label{eq;20.6.6.10}
 \Harm(q)\lrarr
 \prod_{P\in D_{\ess}}
 \prod_{I\in \nbigs(q,P)}
 \nbigp
 \times
 \prod_{P\in D_{>0}}
 \nbigp(q,P).
\end{equation}

We introduce a boundary condition at infinity of $X$
when $X$ is non-compact.
Because
$(K_{X\setminus D}^{(r+1-2i)/2})^{-1}\otimes
 K_{X\setminus D}^{(r+1-2(i+1))/2}\simeq K_{X\setminus D}^{-1}$,
 the restrictions of
 $h_{|K_{X\setminus D}^{(r+1-2i/2)}}$
 and $h_{|K_{X\setminus D}^{(r+1-2(i+1))/2}}$
 induce a K\"ahler metric $g(h)_i$ of $X\setminus D$
$(i=1,\ldots,r-1)$.
 \begin{df}
We say that $h\in\Harm(q)$ is complete at infinity of $X$
if there exists a relatively compact open neighbourhood $N$ of $D$
such that
$g(h)_{i|X\setminus N}$ $(i=1,\ldots,r-1)$
are complete.
Let $\Harm(q;D,\rc)$ denote the set of
$h\in\Harm(q)$ which is complete at infinity of $X$.
\end{df}

\begin{lem}
\label{lem;20.10.6.20}
 Let $N$ be any relatively compact open neighbourhood
of $D$ in $X$.
\begin{itemize}
 \item For any $h\in\Harm(q;D,\rc)$,
       $g(h)_{i|X\setminus N}$ $(i=1,\ldots,r-1)$
       are complete.
       Moreover,
       $g(h)_{i|X\setminus N}$ $(i=1,\ldots,r-1)$
       are mutually bounded,
       and $|q|_{g(h)_i|X\setminus N}$ are bounded.
 \item Any $h_j\in\Harm(q;D,\rc)$ $(j=1,2)$ are mutually
       bounded on $X\setminus N$.
\end{itemize}
\end{lem}
\pf
The first claim follows from
\cite[Proposition 3.27]{Note0}.
The second follows from
\cite[Proposition 3.29]{Note0}.
\hfill\qed

\vspace{.1in}

We obtain the map
\begin{equation}
\label{eq;20.7.2.30}
 \Harm(q;D,\rc)\lrarr
 \prod_{P\in D_{\ess}}
 \prod_{I\in \nbigs(q,P)}
 \nbigp
 \times
 \prod_{P\in D_{>0}}
 \nbigp(q,P).
\end{equation}

We set
$\Harm^{\real}(q;D,\rc):=
\Harm^{\real}(q)\cap
\Harm(q;D,\rc)$.
Let $\nbigp^{\real}(q,P)$ denote the set of
$\vecb\in\nbigp(q,P)$
such that
$b_i+b_{r+1-i}=-r-1$.
Let $\nbigp^{\real}$ denote the set of
$\veca\in\nbigp$
such that $a_i+a_{r+1-i}=0$.
As the restriction of 
(\ref{eq;20.6.6.10})
and 
(\ref{eq;20.7.2.30}),
we obtain the following maps:
\begin{equation}
\label{eq;20.7.6.10}
 \Harm^{\real}(q)\lrarr
 \prod_{P\in D_{\ess}}
 \prod_{I\in \nbigs(q,P)}
 \nbigp^{\real}
 \times
 \prod_{P\in D_{>0}}
 \nbigp^{\real}(q,P),
\end{equation}
\begin{equation}
\label{eq;20.7.6.11}
 \Harm^{\real}(q;D,\rc)\lrarr
 \prod_{P\in D_{\ess}}
 \prod_{I\in \nbigs(q,P)}
 \nbigp^{\real}
 \times
 \prod_{P\in D_{>0}}
 \nbigp^{\real}(q,P).
\end{equation}

\begin{thm}
 \label{thm;20.6.26.10}
 The maps {\rm(\ref{eq;20.7.2.30})}
 and {\rm(\ref{eq;20.7.6.11})}
 are bijective.
In particular,
if $X$ is compact,
 the maps {\rm(\ref{eq;20.6.6.10})} and
 {\rm(\ref{eq;20.7.6.10})}
 are bijective.
\end{thm}

The theorem follows from
Lemma \ref{lem;20.7.5.31},
Lemma \ref{lem;20.7.5.32}
and Lemma \ref{lem;20.6.29.2}
below.

\begin{rem}
If the zero set of $q$ is finite,
we obtain a refined estimate
around each $P\in D_{\ess}$
as in Proposition {\rm\ref{prop;20.4.25.30}}.
\end{rem}

\subsubsection{Uniqueness}

We set $\nbigs(q):=\coprod_{P\in D_{\ess}}\nbigs(q,P)$.

\begin{lem}
\label{lem;20.7.5.31}
 Suppose that $h_j\in\Harm(q;D,\rc)$ $(j=1,2)$
 satisfy $\vecb_P(h_1)=\vecb_P(h_2)$
 for any $P\in D_{>0}$
 and $\veca_I(h_1)=\veca_I(h_2)$
 for $I\in \nbigs(q)$.
 Then, we obtain $h_1=h_2$.
\end{lem}
\pf
As in Theorem \ref{thm;20.9.23.100},
there uniquely exists
$h^{\rc}\in\Harm(q)$
such that $g(h^{\rc})_i$ are complete on $X\setminus D$.
Recall that $g(h^{\rc})_i$ are mutually bounded,
and that the Gaussian curvature of $g(h^{\rc})_i$
are bounded from below \cite[Lemma 3.15]{Note0}.

Let $g$ be a K\"ahler metric of $X$
such that
(i) the Gaussian curvature is bounded from below,
(ii) the following condition is satisfied for 
a relatively compact neighbourhood $N$ of $D$.
\begin{condition}
\label{condition;20.7.5.40}
 $g_{|X\setminus N}$
 is mutually bounded with the metrics
 $g(h^{\rc})_i$ $(i=1,\ldots,r-1)$.
\end{condition}

Let $s$ be the automorphism of
$\hyperk_{X\setminus D,r}$
determined by
$h_2=h_1\cdot s$.
Let $\Lambda$ denote the adjoint
of the multiplication of the K\"ahler form
associated with $g$.
By (\ref{eq;20.8.16.2}),
we obtain the following on $X\setminus D$:
\begin{multline}
 \label{eq;20.7.5.22}
  \sqrt{-1}\Lambda\delbar\del\Tr(s)
  =-\bigl|
  [\theta(q),s]s^{-1/2}
  \bigr|^2_{h_1,g}
  -\bigl|
   \delbar_E(s)s^{-1/2}
   \bigr|^2_{h_1,g}
\\
 \leq
   -\bigl|[\theta(q),s]s^{-1/2}
   \bigr|^2_{h_1,g}.
\end{multline}
By
Proposition \ref{prop;20.7.1.25},
Proposition \ref{prop;20.7.1.20},
Theorem \ref{thm;20.6.9.30},
and Lemma \ref{lem;20.10.6.20},
we obtain the boundedness of $\Tr(s)$.
Hence, the inequality
(\ref{eq;20.7.5.22}) holds across $D$.
(See \cite[Lemma 2.2]{s2}.)
It particularly implies that
$\Tr(s)$ is subharmonic on $X$.

Let $N$ be a relatively compact neighbourhood of $D$.
If
\[
 \sup_{Q\in X\setminus D}\Tr(s)(Q)
=\sup_{Q\in N\setminus D}\Tr(s)(Q),
\]
the maximum principle for subharmonic functions
implies that
$\Tr(s)$ is constant on $X\setminus D$,
and hence
$[\theta(q),s]=0$ on $X\setminus D$.
Together with $\det(s)=1$,
we obtain $s=\id$.

Let us study the case
$\sup_{Q\in X\setminus D}\Tr(s)(Q)
>\sup_{Q\in N\setminus D}\Tr(s)(Q)$.
By Omori-Yau maximum principle \cite[Lemma 3.2]{Note0},
there exists
a sequence $Q_{\ell}\in X\setminus N$ such that
\[
-2\sqrt{-1}\Lambda\delbar\del\Tr(s)(Q_{\ell})\leq \ell^{-1},
\quad
\Tr(s)(Q_{\ell})>\sup_{Q\in X\setminus D}\Tr(s)(Q_{\ell})-\ell^{-1}.
\]
By (\ref{eq;20.7.5.22}),
there exists $C_1>0$ such that
the following holds for any $\ell$:
\[
 \bigl|
 [s,\theta(q)]s^{-1/2}
 \bigr|^2_{h_1,g}(Q_{\ell})
 \leq C_{1}\ell^{-1}.
\]
Let $z_{\ell}$ be a holomorphic coordinate
around $Q_{\ell}$
such that $|dz_{\ell}|^2_{g}(Q_{\ell})=2$.
By Condition \ref{condition;20.7.5.40},
there exists $B>0$
such that
$\bigl|(dz_{\ell})^{(r+1-2i)/2}\bigr|_{h_1}
\cdot
\bigl|(dz_{\ell})^{(r+1-2(i+1))/2}\bigr|_{h_1}^{-1}
\leq B$
for any $\ell$.
Hence, by \cite[Lemma 3.9]{Note0}
there exists $C_{2}>0$ and $\ell_2$
such that
$\Tr(s)(Q_{\ell})\leq r(1+C_2\ell^{-1/2})$
for any $\ell>\ell_2$.
Hence, we obtain
$\sup_{Q\in X\setminus D}\Tr(s)(Q)\leq r$.
Because $\det(s)=1$,
we obtain $\Tr(s)\geq r$ for any $Q$.
Therefore, $\Tr(s)$ is constantly $r$,
and we obtain $s=\id$.
\hfill\qed

\begin{lem}
 \label{lem;20.7.5.32}
 Suppose that $h\in\Harm(q;D,\rc)$
 satisfies
 $\vecb_P(h)\in\nbigp^{\real}(q,P)$ for any $P\in D_{>0}$
 and $\veca_I(h)\in\nbigp^{\real}$ for any $I\in \nbigs(q)$.
Then, we obtain $h\in\Harm^{\real}(q;D,\rc)$. 
\end{lem}
\pf
We obtain $h^{\lor}\in\Harm(q)$
as the dual of $h$
by using the natural identification of
$\hyperk_{X\setminus D,r}$
with its dual.
Clearly, $h^{\lor}$ is also complete at infinity of $X$.
Because $\vecb_P(h)\in\nbigp^{\real}(q,P)$ $(P\in D_{>0})$
and $\veca_I(h)\in\nbigp^{\real}$ $(I\in \nbigs(q))$,
we obtain
$\vecb_P(h)=\vecb_P(h^{\lor})$
and
$\veca_I(h)=\veca_I(h^{\lor})$.
By Lemma \ref{lem;20.7.5.31},
we obtain
$h=h^{\lor}$.
\hfill\qed

\subsubsection{Existence}
   
\begin{lem}
 \label{lem;20.6.29.2}
 For any $\vecb_P\in\nbigp(q,P)$ $(P\in D_{>0})$
 and $\veca_I\in\nbigp$ $(I\in \nbigs(q))$,
 there exists $h\in\Harm(q;D,\rc)$
 such that
 $\vecb_P(h)=\vecb_P$ for any $P\in D_{>0}$
 and that
 $\veca_I(h)=\veca_I$ for any $I\in \nbigs(q)$.
\end{lem}
\pf
Let $h^{\rc}$ be as in the proof of Lemma \ref{lem;20.7.5.31}.
For any $(\veca_I)_{I\in \nbigs(q)}$
and $(\vecb_P)_{P\in D_{>0}}$,
by Proposition \ref{prop;20.7.6.20}
and Theorem \ref{thm;20.6.9.30},
there exists a $G_r$-invariant Hermitian metric $h_0$
of $\hyperk_{X\setminus D,r}$
such that the following holds.
\begin{itemize}
 \item There exists a relatively compact neighbourhood $N_1$
       of $D$ such that
       $h_{0|X\setminus N_1}=h^{\rc}_{|X\setminus N_1}$.
 \item There exists a relatively compact neighbourhood
       $N_2\subset N_1$ of $D$ such that
       $h_{0|N_2\setminus D}\in\Harm(q_{|N_2\setminus D})$
       with $\vecb_P(h_{0})=\vecb_P$ for $P\in D_{>0}$,
        and $\veca_I(h_0)=\veca_I$ for $I\in\nbigs(q)$.
 \item $\det(h_0)=1$.
\end{itemize}
Then, according to Proposition \ref{prop;20.6.29.11},
there exists $h\in \Harm(q;D,h^{\rc})$
which is mutually bounded with $h_0$.
Thus, we obtain Lemma \ref{lem;20.6.29.2},
and the proof of Theorem \ref{thm;20.6.26.10}
is completed.
\hfill\qed

\subsection{Examples on $\cnum$}
\label{subsection;20.7.1.10}

Let $\varpi_{\infty}:\projtilde^1_{\infty}\lrarr\proj^1$
be the oriented real blow up at $\infty$.
We identify
$\varpi_{\infty}^{-1}(\infty)$ with $S^1$
by the polar decomposition $z=|z|e^{\sqrt{-1}\theta}$.

\subsubsection{Classification 
  in the case of  $\gamma(z) e^{\gminia(z)}(dz)^r$}
\label{subsection;20.4.27.40}

Let $\rho$ be a positive integer.
Let $\alpha$ be a non-zero complex number.
Let $\Lambda(\alpha,\rho)$ denote
the set of the connected components of
\[
 \bigl\{
 e^{\sqrt{-1}\theta}\in S^1\,\big|\,
 \Re(\alpha e^{\sqrt{-1}\rho\theta})<0
 \bigr\}.
\]
Note that $|\Lambda(\alpha, \rho)|=\rho$. 
Let $\gminia$ be a polynomial of the form
\[
 \gminia=\alpha z^{\rho}+\sum_{j=1}^{\rho-1}\gminia_jz^j.
\]
Let $\gamma(z)$ be a non-zero polynomial.
We set $q=\gamma(z)e^{\gminia}(dz)^r$.
Clearly,
$\Lambda(\alpha,\rho)$ is the set of the intervals
special with respect to $q$.
For any $h\in\Harm(q)$,
we obtain
$\veca_{I}(h)\in\nbigp$ $(I\in\Lambda(\alpha,\rho))$.
The following proposition is a special case of
Theorem \ref{thm;20.6.26.10}.

\begin{prop}
\label{prop;20.7.6.2}
 The map 
 $\Harm(q)\lrarr \nbigp^{\Lambda(\alpha,\rho)}$
 is bijective.
 \hfill\qed
\end{prop}

\begin{rem}
 We can obtain a refined estimate around $\infty$
 as in Proposition {\rm\ref{prop;20.4.25.30}}. 
\end{rem}

\subsubsection{The case of
   $\sum_{i=1}^n \gamma_{i}(z)e^{\gminia_i(z)}\,(dz)^r$}

Let $n\geq 2$.
Let $\gminia_i(z)$ $(i=1,2,\ldots,n)$
be mutually distinct polynomials such that
$\gminia_i(0)=0$.
If $\gminia_i\neq 0$,
we set $\rho_i:=\deg(\gminia_i)$,
and let $\alpha_i\neq 0$ denote the coefficients of
the top term of $\gminia_i$,
i.e.,
\[
 \gminia_i(z)
 =\alpha_iz^{\rho_i}+\sum_{j=1}^{\rho_i-1}\gminia_{i,j}z^j.
\]
If $\gminia_i=0$,
we set $\rho_i=0$ and $\alpha_i=0$.
Let $\gamma_i(z)$ $(i=1,\ldots,n)$ be non-zero polynomials.
We set
$q=\sum_{i=1}^n\gamma_i(z)e^{\gminia_i(z)}(dz)^r$.
The following proposition is a special case of
Theorem \ref{thm;20.6.26.10}.
\begin{prop}
  \mbox{{}}\label{prop;20.7.6.3}
\begin{itemize}
 \item Suppose that
there exist $\rho>0$ and $\alpha\in\cnum^{\ast}$
	    such that $\rho_i=\rho$
	    and $\alpha_i/\alpha\in\real_{>0}$ for any $i$.
Then, $\Lambda(\alpha,\rho)$ is the set of the intervals
which are special with respect to $q$.
Hence, for any $h\in\Harm(q)$,
      we obtain $\veca_I(h)\in\nbigp$ $(I\in\Lambda(\alpha,\rho))$,
      which induces a bijection
      $\Harm(q)\lrarr \nbigp^{\Lambda(\alpha,\rho)}$.
 \item
      Otherwise, we obtain $\Harm(q)=\{h^{\rc}(q)\}$.
\end{itemize} 
\end{prop}
\pf
Let us prove the first claim.
We set
$\kappa_0:=\max\{\alpha_i/\alpha\}$
and $\kappa_1:=\min\{\alpha_i/\alpha\}$.
Let $f$ be the entire function determined by
$q=f(dz)^r$.
We set
$\nbigz(f)=
\bigl\{
 \theta\in S^1=\real/2\pi\seisuu\,\big|\,
 \Re(\alpha e^{\sqrt{-1}\theta})=0 
\bigr\}$.
Let $\theta\in S^1\setminus\nbigz(q)$.
There exists $\gminia_i$
such that $\gminia_i=\gminia(f,\theta)$
(See Definition \ref{df;20.9.11.10} for $\gminia(f,\theta)$.)
Note that $\Re(\alpha e^{\sqrt{-1}\rho\theta})\neq 0$.
If $\Re(\alpha e^{\sqrt{-1}\rho\theta})>0$,
we obtain $\alpha_i/\alpha=\kappa_0$.
If $\Re(\alpha e^{\sqrt{-1}\rho\theta})<0$,
we obtain $\alpha_i/\alpha=\kappa_1$.
Then, the first claim is clear.

Let us prove the second claim.
Note that $q$ is not meromorphic
at $\infty$ under the assumption that $n\geq 2$.
Suppose that there exists a special interval
$I\subset S^1$ with respect to $f$.
There exists $\rho>0$
such that the length of $I$ is $\pi/\rho$.
Because there exists $\gminia_i$
such that $\rho=\rho_i$,
$\rho$ is a positive integer.
Suppose that there exists
$\gminia_j$ such that $\rho_j>\rho$.
Then, there exists
$\theta\in I\setminus \nbigz(f)$ such that
$\Re(\alpha_je^{\sqrt{-1}\rho_j\theta})>0$.
It contradicts that $I$ is negative.
Suppose that there exists $\rho_j<\rho$.
For any $\theta\in I\setminus\nbigz(f)$,
there exists $\gminia_i$ such  that $\rho_i=\rho$
and $\gminia_i=\gminia(f,\theta)$.
But, we obtain
$\Re(\alpha_js^{\rho_j}e^{\sqrt{-1}\rho_j\theta})
>\Re(\alpha_is^{\rho_i}e^{\sqrt{-1}\rho_i\theta})$
for any sufficiently large $s$,
which contradicts $\gminia_i=\gminia(f,\theta)$.
Therefore, we obtain $\rho_i=\rho$ for any $i$.
There exists $\alpha\in\cnum^{\ast}$ such that
for any $\theta\in I\setminus\nbigz(f)$
we obtain $\deg(\gminia(f,\theta)-\alpha z^{\rho})<\rho$.
Suppose that there exists $\gminia_j$
such that $\alpha_j/\alpha\not\in\real_{>0}$.
Then, there exists $\theta\in I$
such that
$\Re(\alpha_j e^{\sqrt{-1}\rho\theta})>0$.
It contradicts that $I$ is negative.
Hence, we obtain
$\alpha_j/\alpha\in\real_{>0}$ for any $j$.
\hfill\qed

\begin{cor}
Let $P$ be any non-zero polynomial.
and let $Q$ be any non-constant polynomial.
We set
$q_{\cos}=P\cos(Q)(dz)^r$,
$q_{\sin}=P\sin(Q)(dz)^r$ ,
$q_{\cosh}=P\cosh(Q)(dz)^r$
and
$q_{\sinh}=P\sinh(Q)(dz)^r$.
Because
$P\cos(Q)=(Pe^{\sqrt{-1}Q}+Pe^{-\sqrt{-1}Q})/2$,
we obtain
$\Harm(q_{\cos})=\{h^{\rc}(q_{\cos})\}$.
Similarly,
we obtain
$\Harm(q_{\sin})=\{h^{\rc}(q_{\sin})\}$,
$\Harm(q_{\cosh})=\{h^{\rc}(q_{\cosh})\}$
and 
 $\Harm(q_{\sinh})=\{h^{\rc}(q_{\sinh})\}$.
\hfill\qed 
 \end{cor}

\subsubsection{Airy function}

Let $\Gamma$ be a continuous map
$\real\lrarr \cnum$
such that
$\Gamma(t)=t(1+\sqrt{-3})/2$ for $t<-1$
and
$\Gamma(t)=t(-1+\sqrt{-3})/2$ for $t>1$.
Recall that an Airy function is defined as
\begin{equation}
\label{eq;20.7.2.40}
 \Ai(z):=\int_{\Gamma}
 e^{zt-t^3/3}dt.
\end{equation}
(For example, see \cite[\S22]{Wasow}.)
It is an entire function.
Recall that it is a solution of
the linear differential equation
$\del_z^2u-zu=0$,
and hence it induces a section of $\gbigb_{\projtilde^1}$.
On $|\arg(z)|<\pi$,
it satisfies the following estimate
(for example, see \cite[\S23]{Wasow}):
\begin{equation}
\label{eq;20.6.26.10}
e^{\frac{2}{3}z^{3/2}}z^{1/4}
\Bigl(
 \Ai(z)
 -\frac{z^{-1/4}}{2\sqrt{\pi}}e^{-\frac{2}{3}z^{3/2}}
\Bigr)
=O(|z|^{-3/2}).
\end{equation}
Let $\gamma(z)$ be any non-zero polynomial.
We set $q=\gamma(z)\Ai(z)\,(dz)^r$.
Because of (\ref{eq;20.6.26.10}),
the special interval with respect to $q$
is $I=\{-\pi/3<\theta<\pi/3\}$.
We obtain the following proposition
from Theorem \ref{thm;20.6.26.10}.
 \begin{prop}
\label{prop;20.7.6.4}
  For any $h\in\Harm(q)$,
  we obtain $\veca_I(h)\in\nbigp$,
  which induces a bijection
 $\Harm(q)\lrarr \nbigp$.
 \hfill\qed
 \end{prop}

 \begin{rem}
More generally,
let $f$ be a solution of a differential equation
 $(\del_z^n+\sum_{j=0}^{n-1}a_j(z)\del_z^j)f=0$,
 where $a_j(z)$ are polynomials.
  Then, $f$ induces a section of $\gbigb_{\projtilde^1_{\infty}}$
  according to the classical asymptotic analysis.
 (For example, see {\rm\cite[\S II.1]{Majima}}.)
  Hence, Theorem {\rm\ref{thm;20.6.26.10}} allows us
  to classify $\Harm(f\,(dz)^r)$
  in terms of the asymptotic expansion of $f$.
 \end{rem}

\section{A generalization}

\subsection{Preliminary}

\subsubsection{Meromorphic extensions}

Let $X$ be a Riemann surface
with a discrete subset $S$.
Let $\nbigo_X(\ast S)$ denote the sheaf of
meromorphic functions on $X$ which may have poles
along $S$.
Let $\iota:X\setminus S\lrarr X$
denote the inclusion.
For any locally free $\nbigo_{X\setminus S}$-module
$\nbigm$,
we obtain the $\iota_{\ast}\nbigo_{X\setminus S}$-module
$\iota_{\ast}\nbigm$.
\begin{df}
A meromorphic extension of $\nbigm$ on $(X,S)$
is a locally free $\nbigo_X(\ast S)$-submodule
$\nbigmtilde\subset \iota_{\ast}\nbigm$
such that 
$\iota^{\ast}\nbigmtilde=\nbigm$.
\end{df}

\subsubsection{Meromorphic extensions of
   cyclic Higgs bundles on a punctured disc}

Let $U$ be a neighbourhood of $0$ in $\cnum$.
We set $U^{\circ}:=U\setminus\{0\}$.
Let $\gbigl_i$ $(i=1,\ldots r)$
be locally free $\nbigo_{U^{\circ}}$-modules
of rank one.
Let $\psi_i:\gbigl_i\lrarr\gbigl_{i+1}\otimes K_{U^{\circ}}$
$(i=1,\ldots,r-1)$
and $\psi_r:\gbigl_r\lrarr\gbigl_1\otimes K_{U^{\circ}}$
be $\nbigo_{U^{\circ}}$-morphisms.
We obtain
a holomorphic section
$q_{\leq r-1}:=\psi_{r-1}\circ\cdots\circ\psi_1$
of
$\Hom(\gbigl_1,\gbigl_r)\otimes
K_{U^{\circ}}^{\otimes (r-1)}$.
We assume the following.
\begin{itemize}
\item The zero set of $q_{\leq r-1}$ is finite.
\end{itemize}

We set $\gbigv:=\bigoplus_{i=1}^r\gbigl_i$.
Let $\theta$ be the Higgs field of $\gbigv$
induced by $\psi_i$ $(i=1\ldots,r)$.

\begin{prop}
\label{prop;19.10.23.100}
For any meromorphic extension 
$\widetilde{\det(\gbigv)}$ of $\det(\gbigv)$ on $(U,0)$,
there uniquely exist meromorphic extensions 
$\gbigltilde_i$
of $\gbigl_i$ $(i=1,\ldots,r)$ on $(U,0)$
such that 
(i)
 $\theta(\gbigltilde_i)\subset\gbigltilde_{i+1}
 \otimes K_U$ $(i=1,\ldots,r-1)$,
 (ii) $\det(\gbigvtilde)=\widetilde{\det(\gbigv)}$,
 where $\gbigvtilde=\bigoplus\gbigltilde_i$.
\end{prop}
\pf
By shrinking $U$,
we may assume that $q_{\leq r-1}$ is nowhere vanishing
on $U^{\circ}$.
We obtain the morphisms
$\beta_i:\gbigl_i\lrarr\gbigl_{i+1}$ $(i=1,\ldots,r-1)$
by $\psi_i=\beta_i\,(dz/z)$.

Because $U^{\circ}$ is Stein,
there exists a global frame $u^{\star}_1$ of $\gbigl_1$.
We set
$u^{\star}_{i}:=\beta_{i-1}\circ\cdots\beta_1(u^{\star}_1)$
on $U^{\circ}$ for $i=2,\ldots,r$.
We obtain the meromorphic extensions
$\gbigltilde^{\star}_i=\nbigo_U\cdot u^{\star}_i$
of $\gbigl_i$.

Let $v_0$ be a frame of $\widetilde{\det(\gbigv)}$.
Let $f$ be the holomorphic function on $U^{\circ}$
determined by
$u^{\star}_1\wedge\cdots\wedge
 u^{\star}_r=f\cdot v_0$.
 There exist an integer $\ell$
 and  a holomorphic function
$g$ on $U^{\circ}$
such that $f=z^{\ell}e^g$.
We set $u_i:=e^{-g/r}u_i^{\star}$
and  we obtain meromorphic extensions
$\gbigltilde_i:=\nbigo_{U}(\ast 0)u_i$.
We set $\gbigvtilde=\bigoplus\gbigltilde_i$.
Then, by the construction,
we obtain
$\theta(\gbigltilde_i)
\subset
\gbigltilde_{i+1}\otimes K_{U}$
for $i=1,\ldots,r-1$,
and 
$\det(\gbigvtilde)=\widetilde{\det(\gbigv)}$.

Let $\gbigltilde^{\sharp}_{i}$ be
meromorphic extensions of $\gbigl_i$
such that
(i) $\theta\gbigltilde^{\sharp}_{i}
\subset\gbigltilde^{\sharp}_{i+1}\otimes K_U$ $(i=1,\ldots,r-1)$,
(ii) $\det(\gbigvtilde^{\sharp})=
 \widetilde{\det(\gbigv)}=\det(\gbigvtilde)$,
 where $\gbigvtilde^{\sharp}=\bigoplus\gbigltilde^{\sharp}_i$.
There exist frames
$u^{\sharp}_{i}$ of $\gbigltilde^{\sharp}_{i}$
such that 
$\theta(u^{\sharp}_{i})=u^{\sharp}_{i+1} dz/z$ for $i=1,\ldots,r-1$.

Let $\gamma$ be the holomorphic function on $U^{\circ}$
determined by
$u^{\sharp}_{1}=\gamma u_{1}$.
Then, we obtain
$u^{\sharp}_{i}=\gamma u_i$ for $i=1,\ldots,r$.
Because both
$u^{\sharp}_{1}\wedge\cdots\wedge
 u^{\sharp}_r$
and 
$u_{1}\wedge\cdots\wedge
 u_r$
are sections of
$\widetilde{\det(\gbigv)}$,
$\gamma^r$ is meromorphic at $0$.
Hence, we obtain that $\gamma$ is meromorphic,
i.e.,
$\gbigltilde_i=\gbigltilde_i^{\sharp}$.
\hfill\qed

\subsection{Cyclic Higgs bundles}

Let $X$ be a Riemann surface
with a finite subset $D\subset X$.
Assume that $X\setminus D$ is an open Riemann surface,
i.e.,
$X$ is open,
or $X$ is compact and $D\neq\emptyset$.
(See Remark \ref{rem;20.8.6.10}.)
For each $P\in D$,
let $(X_P,z_P)$ be a holomorphic coordinate neighbourhood
such that $z_P(P)=0$.
Set $X_P^{\ast}:=X_P\setminus\{P\}$.

Let $r\geq 2$.
Let $L_i$ $(i=1,\ldots,r)$ be 
holomorphic line bundles on $X\setminus D$.
Let
$\psi_i:L_i\lrarr
L_{i+1}\otimes K_{X\setminus D}$
and 
$\psi_r:L_r\lrarr
 L_1\otimes K_{X\setminus D}$
be non-zero morphisms.
We set $E:=\bigoplus L_i$.
Let $\theta$ be the cyclic Higgs field 
of $E$ induced by $\psi_i$
$(i=1,\ldots,r)$.

We obtain the holomorphic section  
$q:=\psi_{r}\circ\cdots\circ\psi_1=(-1)^{r-1}\det(\theta)$
of $K_{X\setminus D}^{\otimes r}$,
and the holomorphic section
$q_{\leq r-1}:=
 \psi_{r-1}\circ\cdots\circ\psi_1$
of $\Hom(L_1,L_r)\otimes K_{X\setminus D}^{\otimes (r-1)}$.
We assume the following.
\begin{itemize}
\item
     The zero set of $q_{\leq r-1}$ is finite.
\item $q$ is not constantly $0$.
\item Let $f_P$ $(P\in D)$ be holomorphic functions
       on $X_P^{\ast}$ obtained as 
       $q_{|X_P^{\ast}}=f_P\,(dz_P/z_P)^r$.
      Then, $f_P$ have multiple growth orders at $P$.
\end{itemize}
Let $D_{\mero}$ denote the set of $P\in D$
such that $f_P$ is meromorphic at $P$.
We put $D_{\ess}:=D\setminus D_{\mero}$.
For $P\in D_{\mero}$,
we describe
$f_P=z_P^{m_P}\beta_P$,
where $\beta_P$ is a nowhere vanishing holomorphic function on $X_P$.
Let $D_{>0}$ denote the set of $P\in D_{\mero}$
such that $m_P>0$.
We set $D_{\leq 0}:=D_{\mero}\setminus D_{>0}$.

\subsubsection{Flat metrics on the determinant and local frames}

Let $h_{\det(E)}$ be a Hermitian metric of
$\det(E)$
such that the Chern connection of
$\bigl(\det(E),h_{\det(E)}\bigr)$
is flat.
Note that
such a metric exists
because $\det(E)\simeq\nbigo_{X\setminus D}$
under the assumption that $X\setminus D$ is non-compact.
(For example, see \cite[Theorem 30.3]{Forster-book}.)

 \begin{lem}
For each $P\in D$,
there exist frames
 $\vecv_P=(v_{P,i}\,|\,i=1,\ldots,r)$
 of $L_{i|X_P^{\ast}}$
 and a real number $c(\vecv_P)$
 such that the following condition is satisfied.
\begin{itemize}
 \item $\theta(v_{P,i})=v_{P,i+1}(dz_P/z_P)$ for $i=1,\ldots,r-1$.
 \item We set $\omega_P=v_{P,1}\wedge\cdots\wedge v_{P,r}$.
       Then, we obtain
       $|z_P|^{c(\vecv_P)}\cdot |\omega_P|_{h_{\det(E)}}=1$.
\end{itemize}
 \end{lem}
\pf
There exists a frame $w_P$ of
$\det(E)_{|X_P^{\ast}}$
such that
$|w_P|_{h_{\det(E)}}=|z_P|^{d_P}$ for a real number $d_P$.
By Proposition \ref{prop;19.10.23.100},
there exists a frame
$v'_{P,i}$ of $L_{i|X_P^{\ast}}$
$(i=1,\ldots,r)$
such that
 $\theta(v'_{P,i})=v'_{P,i+1}(dz_P/z_P)$ for $i=1,\ldots,r-1$,
 and
 $v'_{P,1}\wedge\cdots\wedge v'_{P,r}
 =z_P^{n_P}\beta_P\cdot w_P$
 for an integer $n_P$
 and a nowhere vanishing holomorphic function $\beta_P$
 on $X_P$.
 By fixing an $r$-th root $\beta_P^{1/r}$ of $\beta_P$,
 and by setting $v_{P,i}:=\beta_P^{-1/r}v_{P,i}'$,
 we obtain the claim of the lemma.
\hfill\qed

\vspace{.1in}

For $P\in D_{>0}$, 
let $\nbigp(q,P,\vecv_P)$ denote the set of
$\vecb=(b_i)\in\real^r$
satisfying the following conditions:
\[
 \sum_{i=1}^r b_i=c(\vecv_P),\quad
 b_i\geq b_{i+1}\,\,(i=1,\ldots,r-1),
\quad
 b_r\geq b_{1}-m_{P}.
\]

\begin{rem}
If
 $(E,\theta)
 =(\hyperk_{X\setminus D,r},\theta(q))$
 for an $r$-differential $q$, 
 we usually choose
 $h_{\det(E)}=1$
 and
 \[
 v_{P,i}=z_P^{i}(dz_P)^{(r+1-2i)/2}
 \quad
 (i=1,\ldots,r),
 \]
 for which we obtain
 $c(\vecv_P)=-r(r+1)/2$.
\end{rem}

\subsubsection{Harmonic metrics and the associated parabolic weights}

We consider the $G_r$-action on $L_i$
by $a\bullet u_i=a^iu_i$,
which induces a $G_r$-action on $E$.
For any open subset $Y\subset X\setminus D$,
let $\Harm^{G_r}((E,\theta)_{|Y},h_{\det(E)})$
denote the set of $G_r$-invariant harmonic metrics $h$
of $(E,\theta)_{|Y}$
such that $\det(h)=h_{\det(E)|Y}$.

\begin{prop}
\label{prop;20.7.13.41}
 Let $h\in \Harm^{G_r}(E,\theta,h_{\det(E)})$.
\begin{itemize}
 \item For any $P\in D_{>0}$,
       there exists $\vecb_P(h)\in\nbigp(q,P,\vecv_P)$
       determined by the following condition.
       \[
       b_{P,i}(h)=
       \inf\Bigl\{
       b\in\real\,\Big|\,
      \mbox{\rm $|z_P|^b|v_{P,i}|_h$ is bounded}
       \Bigr\}.
       \]
 \item For any $P\in D_{\ess}$ and any $I\in \nbigs(q,P)$,
       there exist $\veca_I(h)\in\nbigp$ and $\epsilon>0$
       such that the following estimates hold
       as $|z_P|\to 0$ on
       $\bigl\{
        |\arg(\alpha_Iz_P^{-\rho(I)})-\pi|<(1-\delta)\pi/2
       \bigr\}$ for any $\delta>0$:
       \[
       \log|v_i|_h+a_{I,i}(h)\Re(\alpha_Iz_P^{-\rho(I)})
       =O\bigl(|z_P|^{-\rho(I)+\epsilon}\bigr)
       \]
       Here, $\alpha_I$ and $\rho(I)$ are determined for $I$
       and $q=(-1)^{r-1}\det(\theta)$
       as in {\rm\S\ref{subsection;20.7.13.30}}.
\end{itemize}
\end{prop}
\pf
It is enough to study the case
$D=\{P\}$ and $X=X_P$.
There exists a holomorphic line bundle
$\det(E)^{1/r}$ on $X_P^{\ast}$
with an isomorphism
\begin{equation}
\label{eq;20.7.13.40}
 (\det(E)^{1/r})^{\otimes r}\simeq \det(E).
\end{equation}
There exists a flat metric
$h_{\det(E)^{1/r}}$ of $\det(E)^{1/r}$
such that
$h_{\det(E)^{1/r}}^{\otimes\,r}
=h_{\det(E)}$
under the isomorphism (\ref{eq;20.7.13.40}).
By using the frames $\vecv_{P}$
and $z_P^i(dz_P)^{(r+1-2i)/2}$,
we obtain an isomorphism
$(E,\theta)\otimes \det(E)^{-1/r}
\simeq
 (\hyperk_{X_P^{\ast},r},\theta(q))$.
Then, the claims of the proposition are reduced to
Proposition \ref{prop;20.7.1.20}
and Theorem \ref{thm;20.6.9.30}.
\hfill\qed

\begin{prop}
\label{prop;20.7.13.110}
 Suppose that $h_1,h_2\in\Harm^{G_r}(E,\theta,h_{\det(E)})$
 satisfy $\vecb_P(h_1)=\vecb_P(h_2)$ for any $P\in D_{>0}$
 and $\veca_I(h_1)=\veca_I(h_2)$ for any $I\in\nbigs(q)$.
 Then, for any relatively compact neighbourhood $N$ of $D$ in $X$,
 $h_1$ and $h_2$ are mutually bounded
 on $N\setminus D$.
\end{prop}
\pf
By using the local isomorphisms
in the proof of Proposition \ref{prop;20.7.13.41},
the claim is reduced to
Proposition \ref{prop;20.7.1.25},
Proposition \ref{prop;20.7.1.20}
and Theorem \ref{thm;20.6.9.30}.
\hfill\qed

\subsubsection{Completeness at infinity}
\label{subsection;20.10.6.10}

Let $Z(q_{\leq r-1})$ denote the zero set of $q_{\leq r-1}$,
which is assumed to be finite.
We set $\Dtilde=D\cup Z(q_{\leq r-1})$.
Note that $\psi_i$ $(i=1,\ldots,r-1)$ induce
isomorphisms
$(L_{i+1}/L_i)_{|X\setminus\Dtilde}
\simeq K_{X\setminus\Dtilde}^{-1}$.
Hence,
for any $h\in \Harm^{G_r}(E,\theta,h_{\det(E)})$,
the Hermitian metric
$h_{|L_{i+1}}\otimes h_{|L_i}^{-1}$
$(i=1,\ldots,r-1)$ on $L_{i+1}/L_i$
induce K\"ahler metrics $g(h)_i$
on $X\setminus \Dtilde$.

\begin{df} 
$h\in\Harm^{G_r}(E,\theta,h_{\det(E)})$
is called complete at infinity of $X$
if there exists a relatively compact open subset
$N$ of $\Dtilde$
such that 
$g(h)_{i|X\setminus N}$ are complete.
Let $\Harm^{G_r}(E,\theta,h_{\det(E)};D,\rc)$
denote the set of
$h\in\Harm^{G_r}(E,\theta,h_{\det(E)})$
which are complete at infinity of $X$.
\end{df}

Because $X\setminus\Dtilde$ is an open Riemann surface,
there exists an isomorphism
$\det(E)_{|X\setminus\Dtilde}\simeq\nbigo_{X\setminus\Dtilde}$.
Hence,
there exists a holomorphic line bundle
$\det(E)^{1/r}$
on $X\setminus\Dtilde$
with an isomorphism
$\rho:(\det(E)^{1/r})^{\otimes r}\simeq \det(E)_{|X\setminus\Dtilde}$.
We say that
such $(\det(E)^{1/r}_1,\rho_1)$
and $(\det(E)^{1/r}_2,\rho_2)$ are isomorphic
if
there exists an isomorphism
$\kappa:\det(E)^{1/r}_1\simeq \det(E)^{1/r}_2$
such that
$\rho_2\circ\kappa^{\otimes\,r}=\rho_1$.

\begin{lem}
\label{lem;20.7.13.120}
If we choose $(\det(E)^{1/r},\rho)$ appropriately,
there exist isomorphisms
 $\varphi_i:
 L_{i|X\setminus\Dtilde}\simeq
 K_{X\setminus\Dtilde}^{(r+1-2i)/2}
 \otimes\det(E)^{1/r}$
 $(i=1,\ldots,r)$
such that the following conditions are satisfied:
\begin{itemize}
  \item The following diagram is commutative.
	\[
	 \begin{CD}
	  E_{|X\setminus\Dtilde}
	  @>{\bigoplus\varphi_i}>>
	  \hyperk_{X\setminus \Dtilde,r}
	   \otimes\det(E)^{1/r}\\
	  @V{\theta}VV @V{\theta(q)}VV \\
  	  E_{|X\setminus\Dtilde}
	   \otimes K_{X\setminus\Dtilde}
	  @>{\bigoplus\varphi_i}>>
	  \hyperk_{X\setminus\Dtilde,r}
  	   \otimes\det(E)^{1/r}	  
	  \otimes K_{X\setminus\Dtilde}.
	 \end{CD}
	\]
  \item The induced morphism
\[
	\det\Bigl(\bigoplus\varphi_i\Bigr):
	\det(E_{|X\setminus\Dtilde})
	\simeq
	\det\bigl(\hyperk_{X\setminus\Dtilde,r}\otimes
	\det(E)^{1/r}
	\bigr)
	=(\det(E)^{1/r})^{\otimes\,r}
\]
	is equal to $\rho^{-1}$.
 \end{itemize}
Such $(\det(E)^{1/r},\rho)$ is unique up to isomorphisms.
Such an isomorphism $\bigoplus\varphi_i$
is unique up to the multiplication of
an $r$-th root of $1$.
\end{lem}
\pf
We take
a holomorphic line bundle
$\det(E)_0^{1/r}$
with an isomorphism
$\rho_0:(\det(E)_0^{1/r})^{\otimes r}\simeq\det(E)$.
There exists an isomorphism of holomorphic line bundles
$\nu_1:L_{1|X\setminus \Dtilde}
\simeq
 K_{X\setminus\Dtilde}^{\otimes(r-1)/2}
 \otimes\det(E)^{1/r}_{0|X\setminus\Dtilde}$.
 We obtain the holomorphic isomorphisms
$\nu_i:L_{i|X\setminus\Dtilde}
\simeq
K_{X\setminus\Dtilde}^{\otimes(r+1-2i)/2}
\otimes\det(E)^{1/r}_{0|X\setminus\Dtilde}$ $(i=2,\ldots,r)$
such that
the following diagram is commutative:
\[
 \begin{CD}
  \bigoplus_{i=1}^rL_{i|X\setminus\Dtilde}
  @>{\bigoplus\nu_i}>>
  \hyperk_{X\setminus\Dtilde,r}
  \otimes\det(E)^{1/r}_{0}\\
  @V{\theta}VV @V{\theta(q)}VV \\
  \bigoplus_{i=1}^rL_{i|X\setminus\Dtilde}
  \otimes K_{X\setminus\Dtilde}
  @>{\bigoplus\nu_i}>>
  \hyperk_{X\setminus\Dtilde,r}
  \otimes\det(E)_0^{1/r}
   \otimes K_{X\setminus\Dtilde}.
 \end{CD}
\]
The induced isomorphism
$\det\bigl(\bigoplus\nu_i\bigr):
\det(E)_{|X\setminus\Dtilde}\simeq
(\det(E)_{0}^{1/r})^{\otimes r}$
is equal to
$\beta\rho_0^{-1}$,
where $\beta$
is a nowhere vanishing holomorphic function.

Let $X\setminus\Dtilde=\bigcup_{\lambda\in\Lambda} U_{\lambda}$
be a locally finite covering
by simply connected open subsets $U_{\lambda}$.
There exists an $r$-th root
$\beta^{1/r}_{\lambda}$ of $\beta_{|U_{\lambda}}$.
If $U_{\lambda}\cap U_{\mu}\neq\emptyset$,
then
$\gamma_{\lambda,\mu}:=
\beta^{1/r}_{\lambda}\cdot (\beta^{1/r}_{\mu})^{-1}$
induces a locally constant function
$U_{\lambda}\cap U_{\mu}\lrarr G_r\subset\cnum^{\ast}$.
We obtain
$\gamma_{\lambda,\mu}\cdot\gamma_{\mu,\lambda}=1$.
If $U_{\lambda}\cap U_{\mu}\cap U_{\nu}\neq\emptyset$,
then 
$\gamma_{\lambda,\mu}\cdot\gamma_{\mu,\nu}\cdot\gamma_{\nu,\lambda}=1$.
We obtain a holomorphic line bundle $\nbigl$ on $X\setminus\Dtilde$
by gluing holomorphic line bundles
$\nbigo_{U_{\lambda}}\cdot e_{\lambda}$ $(\lambda\in \Lambda)$
on $U_{\lambda}$
via the relation
$e_{\lambda}=\gamma_{\lambda,\mu}e_{\mu}$
on $U_{\lambda}\cap U_{\mu}$.
Because $\gamma_{\lambda,\mu}^r=1$,
we obtain an isomorphism
$\rho_{\nbigl}:
\nbigl^{\otimes\,r}\simeq \nbigo_{X\setminus\Dtilde}$
by $e_{\lambda}^{\otimes r}\longmapsto 1$.

For each $\lambda$,
we obtain an isomorphism
$\varphi_{i,\lambda}:=
\beta_{\lambda}^{-1/r}\nu_{i|U_{\lambda}}\otimes e_{\lambda}:
L_{i|U_{\lambda}}
\simeq\bigl(
 K_{X\setminus\Dtilde}^{(r+1-2i)/2}
 \otimes \det(E)^{1/r}_0
 \bigr)_{|U_{\lambda}}
 \otimes \nbigl_{|U_{\lambda}}$.
By the construction, the morphisms
$\varphi_{i,\lambda}$ $(\lambda\in \Lambda)$
determine
an isomorphism
$\varphi_{i}:
L_{i|X\setminus\Dtilde}
\simeq
K_{X\setminus\Dtilde}^{(r+1-2i)/2}
\otimes\det(E)_0^{1/r}\otimes\nbigl$.
By setting
$(\det(E)^{1/r},\rho)=
(\det(E)_0^{1/r}\otimes\nbigl,
\rho_0\otimes\rho_{\nbigl})$,
we obtain the claim of the lemma.
\hfill\qed

\vspace{.1in}
In the following,
we choose
$(\det(E)^{1/r},\rho)$
as in Lemma \ref{lem;20.7.13.120}.
Let $h_{\det(E)^{1/r}}$
be the flat metric of $\det(E)^{1/r}$
determined by
$h_{\det(E)^{1/r}}^{\otimes r}=h_{\det(E)|X\setminus\Dtilde}$.
We obtain the bijection
\[
 \Upsilon:\Harm(q_{|X\setminus\Dtilde})
\simeq
\Harm^{G_r}\bigl(
(E,\theta)_{|X\setminus\Dtilde},h_{\det(E)|X\setminus\Dtilde}
\bigr)
\]
determined by
$\Upsilon(h)=h\otimes h_{\det(E)^{1/r}}$
under the isomorphism
in Lemma {\rm\ref{lem;20.7.13.120}}.

\begin{prop}
\label{prop;20.10.6.1}
 Let $N$ be any relatively compact open neighbourhood of
 $\Dtilde$ in $X$.
 \begin{itemize}
  \item For any
	$h_1,h_2\in \Harm^{G_r}(E,\theta,h_{\det(E)};D,\rc)$,
       the restrictions
       $h_{1|X\setminus N}$ and $h_{2|X\setminus N}$
	are mutually bounded.
  \item For any $h\in\Harm^{G_r}(E,\theta,h_{\det(E)};D,\rc)$,
	the metrics $g(h)_{i|X\setminus N}$
	$(i=1,\ldots,r-1)$ are mutually bounded,
	and $|q_{|X\setminus N}|_{g(h)_i}$ are bounded.
 \end{itemize}
\end{prop}
\pf
We may assume that the boundary of $N$ in $X$
is smooth and compact.
We obtain the first claim
by applying \cite[Proposition 3.29]{Note0}
to $\Upsilon^{-1}(h_i)_{|X\setminus N}$.
For any $h\in\Harm^{G_r}(E,\theta,h_{\det(E)})$,
we obtain
$g(h)_{i|X\setminus \Dtilde}
=g(\Upsilon^{-1}(h_{|X\setminus \Dtilde}))_i$
$(i=1,\ldots,r-1)$.
Hence, we obtain the second claim
from  \cite[Proposition 3.27]{Note0}.
\hfill\qed

\vspace{.1in}
Note that
we may naturally regard $|\psi_i|_h^2$
as real sections of
$K_X\otimes \overline{K_X}$.
Note that $g(h)_i=|\psi_i|_h^2$ $(i=1,\ldots,r-1)$.
The second claim of Proposition \ref{prop;20.10.6.1}
is reworded as follows.
\begin{cor}
\label{cor;20.10.6.2}
 Let $h\in\Harm^{G_r}(E,\theta,h_{\det(E)};D,\rc)$.
For any relatively compact open neighbourhood $N$ of
 $\Dtilde$ in $X$,
the metrics
$(|\psi_i|_h^2)_{|X\setminus N}$ $(i=1,\ldots,r-1)$
are mutually bounded,
and
 $\bigl(
 |\psi_r|_h^2/|\psi_i|_h^2\bigr)_{|X\setminus N}$
 $(i=1,\ldots,r-1)$
are bounded.
\hfill\qed
\end{cor}

\subsubsection{Existence and uniqueness}

By Proposition \ref{prop;20.7.13.41},
we obtain the following map:
\begin{equation}
 \label{eq;20.7.13.101}
 \Harm^{G_r}(E,\theta,h_{\det(E)})
 \lrarr
 \prod_{P\in D_{>0}}\nbigp(q,P,\vecv_P)
 \times
 \prod_{I\in\nbigs(q)}\nbigp.
\end{equation}
We obtain the following map
as the restriction of (\ref{eq;20.7.13.101}):
\begin{equation}
\label{eq;20.7.13.100}
 \Harm^{G_r}(E,\theta,h_{\det(E)};D,\rc)
 \lrarr
 \prod_{P\in D_{>0}}\nbigp(q,P,\vecv_P)
 \times
 \prod_{I\in\nbigs(q)}\nbigp.
\end{equation}

\begin{thm}
 \label{thm;20.8.5.121}
The map {\rm(\ref{eq;20.7.13.100})} is a bijection.
In particular, if $X$ is compact,
 the map {\rm(\ref{eq;20.7.13.101})} is a bijection.
\end{thm}
\pf
Suppose that
$h_1,h_2\in\Harm^{G_r}(E,\theta,h_{\det(E)};D,\rc)$
satisfy
$\vecb_P(h_1)=\vecb_P(h_2)$ for any $P\in D_{>0}$
and
$\veca_I(h_1)=\veca_I(h_2)$ for any $I\in\nbigs(q)$.
By Proposition \ref{prop;20.7.13.110}
and Proposition \ref{prop;20.10.6.1},
we obtain that $h_1$ and $h_2$ are mutually bounded.
By applying the argument in the proof of Theorem \ref{thm;20.6.26.10},
we obtain $h_1=h_2$.

Let $\vecb_P$ $(P\in D_{>0})$
and $\veca_I$ $(I\in\nbigs(q))$.
Let $N$ be a relatively compact open neighbourhood of $D$ in
$X\setminus Z(q_{\leq r-1})$.
By using the argument in the proof of Proposition \ref{prop;20.7.13.41},
we can construct
\[
 h_N\in \Harm((E,\theta,h_{\det(E)})_{|X\setminus Z(q_{\leq r-1})})
\]
such that
$\vecb_P(h_N)=\vecb_P$ for $P\in D_{>0}$
and 
$\veca_I(h_N)=\veca_I$ for $I\in\nbigs(q)$.
Let $N_1$ be a relatively compact neighbourhood of
$D$ in $N$.
Let $N_2$ be a relatively compact neighbourhood of
$Z(q_{\leq r-1})$.
We set
$h^{\rc}_{E,\theta}:=\Upsilon(h^{\rc})$,
where $h^{\rc}$ denotes the unique complete metric
in $\Harm(q_{|X\setminus \Dtilde})$.
(See \S\ref{subsection;20.10.6.10} for $\Upsilon$.)
Then, there exists a $G_r$-invariant
Hermitian metric $h_0$ of $E$
such that
(i)
$h_{0|N_1\setminus D}=
 h_{N|\setminus N_1\setminus D}$,
(ii)
$h_{0|X\setminus (N\cup N_2)}=
 h^{\rc}_{E,\theta|X\setminus(N\cup N_2)}$,
(iii) $\det(h_0)=h_{\det(E)}$.
By Proposition \ref{prop;20.6.29.11},
we obtain
$h\in \Harm^{G_r}(E,\theta,h_{\det(E)})$
which is mutually bounded with $h_0$,
which satisfies
$\vecb_P(h)=\vecb_P$ $(P\in D_{>0})$
and
$\veca_I(h)=\veca_I$ $(I\in \nbigs(q))$.
\hfill\qed

\begin{rem}
\label{rem;20.8.6.10}
This type of classification is well known
in the case where $X$ is compact and $D$ is empty.
Indeed, under the assumption that $q\neq 0$,
$(E,\theta)$ is stable with respect to the $G_r$-action.
Hence, for a prescribed Hermitian metric $h_{\det(E)}$ of $\det(E)$,
there uniquely exists a $G_r$-invariant Hermitian metric $h$ of $E$
such that
 (i) $h$ is Hermitian-Einstein in the sense that
 the trace-free part of $F(h,\theta)^{\bot}$ is $0$,
(ii) $\det(h)=h_{\det(E)}$.
If moreover $c_1(E)=0$,
we may choose a flat metric $h_{\det(E)}$ of $\det(E)$,
then such $h$ is a harmonic metric.
\end{rem}

We obtain the following theorem
as a consequence of Theorem \ref{thm;20.8.5.121}
and Corollary \ref{cor;20.10.6.2}.

\begin{thm}
\label{thm;20.10.8.3}
Suppose $D$ is empty.
Let $N\subset X$ be a relatively compact open set
containing all zeros of $q_{\leq r-1}$.  
There uniquely exists a metric
$h\in\Harm^{G_r}(E,\theta,h_{\det(E)})$
such that the metrics $|\psi_i|_{h}^2$ $(i=1,\cdots, r-1)$
are complete on $X\setminus N$. 
Moreover, on $X\setminus N$, $|\psi_i|^2_{h}$ $(i=1,\ldots,r-1)$
are mutually bounded, 
and $\frac{|\psi_r|^2_{h}}{|\psi_i|^2_{h}}$
$(i=1,\cdots, r-1)$ are bounded.
\hfill\qed
\end{thm}

\subsubsection{Boundedness at infinity}

Let $g$ be a complete K\"ahler metric of $X$.
Let $N$ be a relatively compact open neighbourhood of $\Dtilde$.
We say that $h\in\Harm^{G_r}(E,\theta,h_{\det(E)})$
is bounded at infinity of $X$ with respect to $g$
if there exists $0<\delta=\delta(N)<1$
such that on $X\setminus N$
\begin{equation}
\label{eq;20.10.8.1}
  \delta\leq |\psi_i|_{h,g}^2\leq \delta^{-1} \quad(i=1,\ldots,r-1).
\end{equation}
In other words,
$g(h)_i$ $(i=1,\ldots,r-1)$ are mutually bounded
with $g$ on $X\setminus N$.

\begin{lem}
 \label{lem;20.10.8.2}
Suppose that there exists
$h\in\Harm^{G_r}(E,\theta,h_{\det(E)})$
which is bounded at infinity of $X$.
Then, 
$h$ is complete at infinity of $X$.
Moreover,
$|q|_g$ is bounded on $X\setminus N$,
and hence $|\psi_r|_{h,g}$ is bounded on $X\setminus N$.
\end{lem}
\pf
The first claim is clear by definition.
The boundedness of $|q|_g$ on $X\setminus N$
follows from Proposition \ref{prop;20.10.6.1}.
\hfill\qed

\vspace{.1in}
We obtain the following uniqueness
up to boundedness.
\begin{cor}
Suppose that $h_i\in\Harm^{G_r}(E,\theta,h_{\det(E)})$ $(i=1,2)$
are bounded at infinity of $X$.
Then,
$h_1$ and $h_2$ are mutually bounded on $X\setminus N$.
If moreover $D$ is empty,
we obtain $h_1=h_2$ on $X$.
\end{cor}
\pf
The first claim follows from
Proposition \ref{prop;20.10.6.1}
and Lemma \ref{lem;20.10.8.2}.
The second claim follows
from Theorem \ref{thm;20.10.8.3}
and Lemma \ref{lem;20.10.8.2}.
\hfill\qed

\vspace{.1in}
A bounded solution does not necessarily exist
for a prescribed complete K\"ahler metric.
\begin{cor}
Let $g'$ be another complete K\"ahler metric of $X$.
If there exist
$h,h'\in\Harm^{G_r}(E,\theta,h_{\det(E)})$
such that $h$ (resp. $h'$)
is bounded at infinity of $X$ with respect to $g$ (resp. $g'$),
then $g$ and $g'$ are mutually bounded on $X$.
\end{cor}
\pf
By Proposition \ref{prop;20.10.6.1}
and Lemma \ref{lem;20.10.8.2},
$h$ and $h'$ are mutually bounded on $X\setminus N$.
It implies that
$|g(h)_i|^2$ $(i=1,\ldots,r-1)$
are mutually bounded with $g$ and $g'$ on $X\setminus N$.
Hence,
we obtain that $g$ and $g'$ are mutually bounded on $X$.
\hfill\qed

\vspace{.1in}

\begin{prop}
\label{prop;20.10.8.4}
 Let $g$ be a complete hyperbolic metric of $X$.
Assume that $|q|_g$  is bounded  on $X\setminus N$.
Then, $h\in\Harm^{G_r}(E,\theta,h_{\det(E)})$
is bounded at infinity of $X$ with respect to $g$
if and only if
$h$ is complete at infinity.
\end{prop}
\pf
The ``only if'' part of the claim is clear by definition.
Let us prove that the ``if'' part  of the claim.
Let $K$ be a compact neighbourhood of
$\Dtilde$ in $N$ with smooth boundary.
Note that $X\setminus K$ is hyperbolic.

\begin{lem}
Let $g_{X\setminus K}$ be a complete hyperbolic metric
of $X\setminus K$.
\begin{itemize}
 \item $g\leq g_{X\setminus K}$ on $X\setminus K$.
 \item $g$ and $g_{X\setminus K}$
       are mutually bounded on $X\setminus N$.
\end{itemize}
\end{lem}
\pf
We may assume that the Gaussian curvature
of $g_{X\setminus K}$ and $g$ are constantly $-2$.
Recall that they correspond
to complete solutions of the Toda equation (\ref{eq;20.7.2.1})
with $r=2$ and $q=0$
on $X\setminus K$ and $X$, respectively.
(See also \S\ref{subsection;20.10.10.2}.)
The first claim follows from
\cite[Theorem 1.7]{Note0}.
The second claim of the lemma follows from
\cite[Proposition 3.29]{Note0}.
\hfill\qed

\vspace{.1in}

Because $|q|_g$ is bounded on $X$,
$|q_{|X\setminus K}|_{g_{X\setminus K}}$
is bounded.
Let $h^{\rc}_{X\setminus K}\in\Harm(q_{|X\setminus K})$
be the complete solution on $X\setminus K$.
According to \cite[Theorem 1.8]{Note0},
$g(h^{\rc}_{X\setminus K})_i$ $(i=1,\ldots,r-1)$
are mutually bounded with $g_{X\setminus K}$.
Hence,
$g(h^{\rc}_{X\setminus K})_{i|X\setminus N}$
$(i=1,\ldots,r-1)$
are mutually bounded with $g_{|X\setminus N}$.

Suppose that $h\in\Harm^{G_r}(E,\theta,h_{\det(E)})$
is complete at infinity.
By \cite[Proposition 3.29]{Note0},
$h^{\rc}_{X\setminus K|X\setminus N}$
and $\Upsilon^{-1}(h_{|X\setminus N})$
are mutually bounded.
(See \S\ref{subsection;20.10.6.10} for $\Upsilon$.)
Hence,
$g(h)_{i|X\setminus N}$ $(i=1,\ldots,r-1)$
are mutually bounded with
$g_{|X\setminus N}$,
i.e.,
$h$ is bounded at infinity
with respect to $g$.
\hfill\qed

\vspace{.1in}
We obtain the following theorem
as a consequence of
Theorem \ref{thm;20.10.8.3},
Lemma \ref{lem;20.10.8.2},
and Proposition \ref{prop;20.10.8.4}.

\begin{thm}
Suppose that $D$ is empty.
Let $g$ be a complete hyperbolic metric of $X$.
If $q$ is bounded with respect to $g$,
there uniquely exists a bounded solution
$h\in \Harm^{G_r}(E,\theta,h_{\det(E)})$.
Conversely, if there exists such a bounded solution $h$,
$q$ is bounded with respect to $g$.
\hfill\qed
\end{thm}

\begin{rem}
According to Proposition {\rm\ref{prop;20.10.8.4}},
if $g$ is complete hyperbolic
and if $|q|_g$ is bounded on $X\setminus N$,
we may replace the condition
``completeness at infinity'' in
Theorem {\rm\ref{thm;20.8.5.121}}
with
the condition
``boundedness at infinity with respect to $g$''.
 \end{rem}

\section{Appendix}

\subsection{The case of trivial $r$-differentials}

\subsubsection{Rank $2$ case}

Let $X$ be a Riemann surface with a finite subset $D$.
Suppose that
$X\setminus D$ has a K\"ahler metric $g_{X\setminus D}$
with infinite volume
such that 
the Gaussian curvature of
the Riemannian metric $\Re(g_{X\setminus D})$
is constantly $-k$ $(k>0)$.

\begin{prop}
\label{prop;20.10.7.20}
 For any $\veca=(a_P)_{P\in D}\in\real^D_{\geq 0}$,
 there exists a K\"ahler metric $g_{\veca}$ of $X\setminus D$
 whose Gaussian curvature is constantly $-k$,
 such that the following holds.
 \begin{itemize}
  \item For any neighbourhood $N$ of $D$,
	$g_{\veca|X\setminus N}$ and $g_{X\setminus D|X\setminus N}$
	are mutually bounded.
  \item Let 
        $(X_P,z_P)$ be
	a coordinate neighbourhood
	around $P$ such that $z_P(P)=0$.
	If $a_P>0$, then
	$g_{\veca|X_P\setminus\{P\}}|z_P|^{-2(a_P-1)}$
	is mutually bounded with $dz_P\,d\zbar_P$
	around $P$.
	If $a_P=0$, then
	$g_{\veca|X_P\setminus\{P\}}$
	is mutually bounded with
	$|z_P|^{-2}(\log|z_P|)^{-2}dz_P\,d\zbar_P$
	around $P$.
 \end{itemize}
\end{prop}
\pf
It is enough to study the case where $k=2$.
Recall that in this paper,
for a K\"ahler metric $g$ on $X\setminus D$,
we use the Hermitian metric on $K_X^{\ell/2}$
denoted as $g^{-\ell/2}$,
determined as follows.
\begin{itemize}
 \item
If $g=g_0\,dz\otimes d\zbar$ for a local holomorphic coordinate $z$,
then
$|(dz)^{\ell/2}|^2_{g^{-\ell/2}}=(\frac{g_0}{2})^{-\ell/2}$.
\end{itemize}
Let $\omega_g$ denote the associated K\"ahler form,
which is locally $\frac{\sqrt{-1}}{2}g_0dz\wedge d\zbar$.

We set
$\hyperk_{X\setminus D,2}:=
 K_{X\setminus D}^{1/2}
 \oplus
 K_{X\setminus D}^{-1/2}$.
Let $0_{X\setminus D,2}$ denote the quadratic differential
on $X\setminus D$
which is constantly $0$.
It induces the Higgs field $\theta(0_{X\setminus D,2})$
of $\hyperk_{X\setminus D,2}$.
We consider the action of $G_2=\{\pm 1\}$
on $K_{X\setminus D}^{\pm 1/2}$
given by $a\bullet x=a^{\pm 1}x$,
which induces a $G_2$-action on $\hyperk_{X\setminus D,2}$.
A K\"ahler metric $g$ of $X\setminus D$
induces a Hermitian metric $h^{(1)}(g)=g^{-1/2}\oplus g^{1/2}$.

\begin{lem}
$h^{(1)}(g)$ is a harmonic metric
if and only if
the Gaussian curvature of $\Re(g)$ is
constantly $-2$.
\end{lem}
\pf
Let $R(g)$ denote the curvature of
the Chern connection of $(K^{-1}_{X\setminus D},g)$.
Then, by a direct calculation,
we can check that
$h^{(1)}$ is harmonic if and only if
$\sqrt{-1}R(g)+2\omega_g=0$,
which implies that the Gaussian curvature of $g$ is $-2$.
\hfill\qed

\vspace{.1in}

Let $\Harm(0_{X\setminus D,2})$ denote the set of
$G_2$-invariant harmonic metrics $h$ of
$(\hyperk_{X\setminus D,2},\theta(0_{X\setminus D,2}))$
such that $\det(h)=1$.
Note that $h$ is decomposed as
$h=h_{|K_{X\setminus D}^{1/2}}\oplus h_{|K_{X\setminus D}^{-1/2}}$,
and
it is equal to $h^{(1)}(g)$
for a K\"ahler  metric $g$ of $X\setminus D$
whose Gaussian curvature is $-2$.

Let us give some preliminaries on
hyperbolic metrics 
on $X_P^{\ast}:=X_P\setminus\{P\}$ $(P\in D)$.
We may assume that
$X_P=\{|z_P|<1/2\}$.
It is well known and easy to check that
$g_{P,1}:=2(1-|z_P|^2)^{-2}dz_P\,d\zbar_P$
is a hyperbolic metric of $X_P$ normalized as
$\sqrt{-1}R(g_{P,1})+2\omega_{g_{P,1}}=0$.
For $a>0$,
we set as follows on $X_P^{\ast}$:
\[
 g_{P,a}=\frac{2dz_P^a\,d\zbar_P^a}{(1-|z_P|^{2a})^2}
 =\frac{2a^2|z_P|^{2(a-1)}\,dz_P\,d\zbar_P}{(1-|z_P|^{2a})^2}.
\]
Because $g_{P,a}$ is locally obtained as the pull back of $g_{P,1}$
by the map $z_P\longmapsto z_P^a$,
it is a hyperbolic metric on $X_P^{\ast}$ normalized as
$\sqrt{-1}R(g_{P,a})+2\omega_{g_{P,a}}=0$.
We also set
\[
 g_{P,0}:=\frac{2dz_P\,d\zbar_P}{|z_P|^2(\log|z_P|^2)^2}.
\]
It is easy to check
$\sqrt{-1}R(g_{P,0})+2\omega_{g_{P,0}}=0$.

\begin{lem}
\label{lem;20.8.20.1}
Let $g_P$ be any hyperbolic metric of $X_P^{\ast}$
normalized as
$\sqrt{-1}R(g_{P})+2\omega_{g_{P}}=0$.
\begin{itemize}
 \item There exists $a\geq 0$
       such that $g_P$ and $g_{P,a}$ are mutually bounded
       around $P$.
 \item Let $X_{P,1}$ be a relatively compact neighbourhood
       of $P$ in $X_P$.
       Then, the volume of
       $X_{P,1}\setminus\{P\}$
       is finite with respect to $g_{P}$.      
\end{itemize}
\end{lem}
\pf
Let $h\in\Harm(0_{X_P^{\ast},2})$.
Because $\theta(0_{X_P^{\ast},2})$ is nilpotent,
the harmonic bundle
$(\hyperk_{X_P^{\ast},2},\theta(0_{X_P^{\ast},2}),h)$
is tame.
As in the proof of Proposition \ref{prop;20.7.1.20},
we obtain the locally free
$\nbigo_{X_P}$-modules
\[
 \nbigp^{h}_a\hyperk_{X_P^{\ast},2}
=\nbigp^{h}_aK^{1/2}_{X_P^{\ast}}
\oplus
 \nbigp^{h}_aK^{-1/2}_{X_P^{\ast}}\quad(a\in\real),
\]
and we can prove that
there exists $b\in\real$
such that 
$z_P(dz_P)^{1/2}$ and $z_P^2(dz_P)^{-1/2}$
are sections of
$\nbigp^{h}_b\hyperk_{X_P^{\ast},2}$.
We obtain the numbers
\[
 b_i(h):=
  \inf\bigl\{b\in\real\,\big|\,
  \mbox{$|z_P|^{b}\cdot |z_P^i(dz_P)^{(3-2i)/2}|_h$ is bounded}
  \bigr\}
  \quad(i=1,2).
\]
Because $\theta(0_{X_P^{\ast},2})$ is nilpotent,
we obtain
\[
 \theta(0_{X_P^{\ast},2})
=O\bigl(\bigl|\log|z_P|\bigr|^{-1}\bigr)dz_P/z_P
\]
with respect to $h$,
according to \cite{s2}.
Because
$\theta(0_{X_P^{\ast},2})(z_P(dz_P)^{1/2})
=z_P^2(dz_P)^{-1/2}\cdot (dz_P/z_P)$,
we obtain $b_1(h)\geq b_2(h)$.
Because $\det(h)=1$,
we obtain $b_1(h)+b_2(h)+3=0$.
Moreover,
according to the norm estimate of Simpson \cite{s2},
if $h_j\in\Harm(0_{X_P^{\ast},2})$ $(j=1,2)$
satisfy
$b_i(h_1)=b_i(h_2)$ $(i=1,2)$
then $h_1$ and $h_2$
are mutually bounded.

For a hyperbolic metric $g_{P}$
of $X_P^{\ast}$
normalized as
$\sqrt{-1}R(g_P)+2\omega_{g_P}=0$,
we obtain
$h^{(1)}(g_P)\in\Harm(0_{X_P^{\ast},2})$
induced by $g_P$.
For any $(b_1,b_2)\in\real^2$ such that
$b_1\geq b_2$ and $b_1+b_2+3=0$,
there exists $a\in\real_{\geq 0}$
such that
$b_1=\frac{a}{2}-\frac{3}{2}$
and $b_2=-\frac{a}{2}-\frac{3}{2}$.
Then, we can directly check that
$b_i\bigl(h^{(1)}(g_{P,a})\bigr)=b_i$.
Then, we obtain the first claim.
The second claim can be checked
in the case of $g_{P,a}$ $(a\in\real_{\geq 0})$.
\hfill\qed

\vspace{.1in}

Take any $\veca\in\real^D_{\geq 0}$.
We obtain
$h_{P,a_P}=h^{(1)}(g_{P,a_P})\in\Harm(0_{X_P^{\ast},2})$
corresponding to $g_{P,a_P}$.
Let $h_{X\setminus D}\in\Harm(0_{X\setminus D,2})$
on $X\setminus D$
corresponding to
the hyperbolic metric $g_{X\setminus D}$.
Let $X_{1,P}$ $(P\in D)$
be relatively compact neighbourhoods of $P$
in $X_P$.
There exists a $G_2$-invariant Hermitian metric $h_0$
such that
(i) $h_0=h_{P,a_P}$ on $X_{1,P}$,
(ii) $h_0=h_{X\setminus D}$
on $X\setminus\coprod_{P\in D} X_{P}$,
(iii) $\det(h_0)=1$.
Note that the support of
$F(h_0)=
R(h_0)
+[\theta(0_{X\setminus D,2}),
 \theta(0_{X\setminus D,2})^{\dagger}_{h_0}]$
is compact.
We also note that
$K^{-1/2}_{X\setminus D}$ is the unique proper Higgs subbundle
of $\bigl(
\hyperk_{X\setminus D,2},\theta(0_{X\setminus D,2})
\bigr)$.

Let $h_1$ denote the Hermitian metric of
$K^{-1/2}_{X\setminus D}$
induced by $h_0$.
Because $\rank K^{-1/2}_{X\setminus D}=1$,
we obtain $F(h_1)=R(h_1)$.
Let us prove that
\begin{equation}
 \label{eq;20.8.16.1}
  \int_{X\setminus D}
  \sqrt{-1}R(h_1)=-\infty.
\end{equation}
Let $h_2$ denote the Hermitian metric of
$K^{-1/2}_{X\setminus D}$
induced by the hyperbolic metric $g_{X\setminus D}$.
By the assumption, we have
\[
\int_{X\setminus D}
  \sqrt{-1}R(h_2)=
-2\int_{X\setminus D}\omega_{g_{X\setminus D}}
=-\infty.
\]
By the construction,
$h_2=h_1$ holds on $X\setminus \coprod X_P$.
By Lemma \ref{lem;20.8.20.1},
the volume of
$\coprod_{P\in D}(X_P\setminus D)$
is finite with respect to $g_{X\setminus D}$.
We obtain
\[
 \int_{X\setminus \coprod X_P}
 \sqrt{-1}R(h_1)=
 \int_{X\setminus \coprod X_P}
 \sqrt{-1}R(h_2)=-\infty.
\]
Because $h_{1|X_{1,P}}$ $(P\in D)$ are induced by
$g_{P,a_P}$,
we obtain
\[
 \int_{X_{1,P}\setminus D}
 \sqrt{-1}R(h_1)
 =\int_{X_{1,P}\setminus D}
 -2\omega_{g_{P,\veca}}<0.
\]
Because $\coprod_{P\in D}(X_{P}\setminus X_{1,P})$ is compact,
we obtain
\[
 \sum_{P\in D}
  \int_{X_P\setminus X_{1,P}}
  \sqrt{-1}R(h_1)<\infty.
\]
In all,
we obtain (\ref{eq;20.8.16.1}).
As a result,
$(\hyperk_{X\setminus D,2},\theta(0_{X\setminus D,2}),h_0)$
is analytically stable.
By Proposition \ref{prop;20.6.11.10},
there exists $h_{\veca}\in\Harm(0_{X\setminus D,2})$
such that
$h_{\veca}$ and $h_0$ are mutually bounded on $X\setminus D$.
The K\"ahler metric $g_{\veca}$
on $K_{X\setminus D}^{-1}$ induced by $h_{\veca}$
satisfies the desired conditions.
\hfill\qed

\begin{rem}
If there exists a complete K\"ahler metric
$g_{X\setminus D}$ 
of $X\setminus D$ with finite volume,
then there exist a compact Riemann surface $\Xbar$
and an open embedding $X\setminus D\lrarr \Xbar$
whose complement is a finite subset.
(For example, see {\rm\cite{Eberlein, Siu-Yau}},
where much more general results are proved.)
Therefore, we may apply
the theory of tame harmonic bundles  
{\rm \cite{s2}} due to Simpson
to study a problem as in Proposition {\rm\ref{prop;20.10.7.20}}.
The stability condition provides us with a constraint
on $\veca$.
\end{rem}

\subsubsection{Rank $r$ case}

Let $r$ be a positive integer larger than $2$.
Let $0_{X\setminus D,r}$ denote the $r$-differential on
$X\setminus D$ which is constantly $0$.
Let $h\in\Harm(0_{X\setminus D,r})$.
Because
$\theta(0_{X\setminus D,r})$ is nilpotent,
the harmonic bundle
$(\hyperk_{X\setminus D,r},\theta(0_{X\setminus D,r}),h)$
is tame.
As in the proof of Proposition \ref{prop;20.7.1.20},
from $(\hyperk_{X_P^{\ast},r},\theta(0_{X_P^{\ast},r}),h_{|X_P^{\ast}})$,
we obtain the locally free
$\nbigo_{X_P}$-modules
$\nbigp^h_a\hyperk_{X_P^{\ast},r}
=\bigoplus_{i=1}^r
\nbigp_a^hK^{(r+1-2i)/2}_{X_P^{\ast}}$ $(a\in\real)$,
and we can prove that
there exists $b\in\real$
such that 
$z_P^i(dz_P)^{(r+1-2i)/2}$ $(i=1,\ldots,r)$
are sections of
$\nbigp^{h}_b\hyperk_{X_P^{\ast},r}$.
For $i=1,\ldots,r$,
we obtain the numbers
\[
 b_{P,i}(h):=
  \inf\bigl\{b\in\real\,\big|\,
  \mbox{$|z_P|^{b}\cdot |z_P^i(dz)_P^{(r+1-i)/2}|_h$ is bounded}
  \bigr\}.
\]
Let $\vecb_P(h)$ denote the tuple
$(b_{P,1}(h),\ldots,b_{P,r}(h))$.
Because $\theta(0_{X\setminus D,r})$ is nilpotent,
we obtain
$|\theta(0_{X\setminus D,r})|_h=O\bigl((-\log|z_P|)^{-1}\bigr)dz_P/z_P$
with respect to $h$.
Because
\[
\theta(0_{X_P^{\ast},r})(z_P^i(dz_P)^{(r+1-2i)/2})
=z_P^{i+1}(dz_P)^{(r+1-2(i+1))/2}\cdot (dz_P/z_P)
\]
for $i=1,\ldots,r-1$,
we obtain $b_{P,i}(h)\geq b_{P,i+1}(h)$ $(i=1,\ldots,r-1)$.
Because $\det(h)=1$,
we obtain $\sum b_{P,i}(h)=-r(r+1)/2$.

For simplicity,
we assume that the Gaussian curvature of $g_{X\setminus D}$ is $-2$.
Note that
$(\hyperk_{X\setminus D,r},\theta(0_{X\setminus D,r}))$
is isomorphic to
the $(r-1)$-th symmetric product of
$(\hyperk_{X\setminus D,2},\theta(0_{X\setminus D,2}))$.
Let $h_r^{(1)}(g_{X\setminus D})
\in\Harm(0_{X\setminus D,r})$
denote the harmonic metric
induced by $h^{(1)}(g_{X\setminus D})$.

\begin{prop}
Suppose that $\vecb_P=(b_{P,1},\ldots,b_{P,r})\in\real^r$ $(P\in D)$
satisfy
\[
 b_{P,1}\geq\cdots\geq b_{P,r},
 \quad
 \sum b_{P,j}=-\frac{r(r+1)}{2}.
\]
 Then, there exists $h\in\Harm(0_{X\setminus D,r})$
 such that
 (i)  $\vecb_P(h)=\vecb_P$ $(P\in D)$,
 (ii) For any neighbourhood $N$ of $D$,
 $h$ and $h_r^{(1)}(g_{X\setminus D})$
 are mutually bounded on $X\setminus N$.
\end{prop}
\pf
Let us explain an outline of the proof.
\begin{lem}
\label{lem;20.8.20.2}
 There exists
 $h_{P,\vecb_P}\in\Harm(0_{X_P^{\ast},r})$
 such that
 $\vecb(h_{P,\vecb})=\vecb_P$.
\end{lem}
\pf
Let $Y$ be a compact Riemann surface
whose genus $g(Y)$ is larger than
$10\sum_i|b_{P,i}|+10r^2$.
Let $Q$ be a point of $Y$.
Let $(Y_Q,z_Q)$ be a holomorphic coordinate neighbourhood
of $Q$ in $Y$ such  that $z_Q(Q)=0$.
For any $a\in\real$,
$\nbigp_{a}\bigl(\hyperk_{Y_Q,r}(\ast Q)\bigr)$
denote the locally free
$\nbigo_{Y_Q}$-module
obtained as follows:
\[
 \nbigp_{a}\bigl(\hyperk_{Y_Q,r}(\ast Q)\bigr)
 =
 \bigoplus_{i=1}^r
  \nbigo_{Y_Q}([a-b_{P,i}])
  z_Q^i(dz_Q)^{(r+1-2i)/2}.
\]
Here, for $c\in\real$,
$[c]$ denotes $\max\{n\in\seisuu\,|\,n\leq c\}$.
From
$\nbigp_{a}\bigl(\hyperk_{Y_Q,r}(\ast Q)\bigr)$
and $\hyperk_{Y\setminus Q,r}$,
we obtain locally free
$\nbigo_Y$-modules
$\nbigp_{a}\bigl(\hyperk_{Y,r}(\ast Q)\bigr)$.
Thus, we obtain a filtered bundle
$\nbigp_{\ast}(\hyperk_{Y,r}(\ast Q))$
on $(Y,Q)$.
We can easily check
$\theta(0_{Y\setminus Q,r})
 \nbigp_a(\hyperk_{Y,r}(\ast Q))
 \subset
  \nbigp_a(\hyperk_{Y,r}(\ast Q))
  \otimes K_Y(Q)$.
We obtain
\[
  \deg\Bigl(
  \nbigp_{\ast}\bigl(
   K^{(r+1-2i)/2}_{Y}(\ast Q)
   \bigr)
  \Bigr)
  =\frac{r+1-2i}{2}(2g(Y)-2)
  -(b_{P,i}+i).
\]
For any $1\leq j\leq r$,
we obtain the following:
\begin{multline}
 \sum_{i=j}^r
 \deg\Bigl(\nbigp_{\ast}\bigl(
  K^{(r+1-2i)/2}_Y(\ast Q)
  \bigr)
  \Bigr)
\\
 =-\frac{(j-1)(r+1-j)}{2}
  (2g(Y)-2)
 -\sum_{i=j}^r(b_{P,i}+i).
\end{multline}
Hence, for any $1<j\leq r$,
we obtain
\[
 \sum_{i=j}^r
 \deg\Bigl(\nbigp_{\ast}\bigl(
  K^{(r+1-2i)/2}_Y(\ast Q)\bigr)
  \Bigr)
  < 0=\deg\Bigl(
  \nbigp_{\ast}\bigl(
  \hyperk_{Y,r}(\ast Q)\bigr)
   \Bigr).
\]
If a non-zero subbundle
$E\subset \hyperk_{Y\setminus Q,r}$
satisfies
$\theta(0_{Y\setminus Q,r})E
\subset E\otimes K_{Y\setminus Q}$,
then there exists $1\leq j\leq r$
such that
$E=\bigoplus_{i=j}^rK_{Y\setminus Q}^{(r+1-2i)/2}$.
Therefore,
the regular filtered Higgs bundle
$\bigl(\nbigp_{\ast}\hyperk_{Y,r}(\ast Q),
\theta(0_{Y\setminus Q},r)\bigr)$
is stable.
According to \cite{s2},
there uniquely exists a harmonic metric $h_{Y,\vecb_P}$
of $(\hyperk_{Y\setminus Q,r},\theta(0_{Y\setminus Q,r}))$
such that $\det(h_{Y,\vecb_P})=1$
and that
$\hyperk_{Y\setminus Q,r}$
with $h_{Y,\vecb_P}$
induces $\nbigp_{\ast}\bigl(
\hyperk_{Y,r}(\ast Q)\bigr)$.
By the uniqueness,
we obtain that $h_{Y,\vecb_P}$ is $G_r$-invariant.

We embed $X_P$ into $Y$
by using $z_P$ and $(Y_Q,z_Q)$.
Then, we can construct $h_{P,\vecb_P}$
as the pull back of $h_{Y,\vecb_P}$.
\hfill\qed

\vspace{.1in}

Let $h_0$ be a $G_r$-invariant Hermitian metric
of $\hyperk_{X\setminus D,r}$
satisfying the following conditions.
\begin{itemize}
 \item $h_0=h_r^{(1)}(g_{X\setminus D})$
      on $X\setminus\coprod_{P\in D}X_P$.
 \item There exist relatively compact neighbourhoods
       $X_{1,P}$ $(P\in D)$ in $X_P$
       such that
       $h_{0}=h_{P,\vecb_P}$ on $X_{1,P}\setminus\{P\}$.
 \item $\det(h_0)=1$.
\end{itemize}
Note that the support of $F(h_0)$ is compact.

\begin{lem}
 $(\hyperk_{X\setminus D,r},\theta(0_{X\setminus D,r}),h_0)$
 is analytically stable
 with respect to the $G_r$-action.
\end{lem}
\pf
For a $G_r$-invariant Hermitian metric
$h$ of $\hyperk_{X\setminus D,r}$,
let $h_{\geq j}$ $(1\leq j\leq r)$
denote the induced Hermitian metric of 
$\bigoplus_{i=j}^rK^{(r+1-2i)/2}_{X\setminus D}$.
Because the volume of $X\setminus D$
with respect to $g_{X\setminus D}$ is infinite,
we obtain
the following for $1<j\leq r$:
\begin{multline}
 \int_{X\setminus D}
 \sqrt{-1}\Lambda
 \Tr F\bigl(h_r^{(1)}(g_{X\setminus D})_{\geq j}\bigr)
 = \\
 \int_{X\setminus D}
 \sqrt{-1}\Lambda
 \Tr R(h_r^{(1)}\bigl(g_{X\setminus D})_{\geq j}\bigr)
=-\infty.
\end{multline}
By Lemma \ref{lem;20.8.20.1},
and by
$h_{0|X\setminus \coprod X_P}
=h_r^{(1)}(g_{X\setminus D})_{|X\setminus \coprod X_P}$,
we obtain
\begin{equation}
\label{eq;20.8.20.10}
 \int_{X\setminus \coprod X_P}
 \sqrt{-1}\Lambda
 \Tr F\bigl((h_0)_{\geq j}\bigr)
=\int_{X\setminus \coprod X_P}
 \sqrt{-1}\Lambda
 \Tr R\bigl((h_0)_{\geq j}\bigr)
 =-\infty.
\end{equation}
On $X_{1,P}$,
$h_{0}=h_{P,\vecb_P}$ holds by the construction.
Because $h_{P,\vecb_P}$ is harmonic,
we obtain $F(h_{P,\vecb_P})=0$.
Let $\pi_j$ denote the projection of
$\hyperk_{X_{1,P},r}$
onto
$\bigoplus_{i=j}^rK^{(r+1-2i)/2}_{X_{1,P}^{\ast}}$.
It is the orthogonal projection
with respect to $h_0$.
By the Chern-Weil formula \cite[Lemma 3.2]{s1},
we obtain
\begin{equation}
\label{eq;20.8.20.11}
\int_{X_{1,P}}
 \sqrt{-1}\Lambda \Tr F\bigl((h_0)_{\geq j}\bigr)
 =-\int_{X_{1,P}}\bigl|
  (\delbar_E+\theta)(\pi_j)
 \bigr|^2_{h_0}\leq 0.
\end{equation}
Because $X_{P}\setminus X_{1,P}$ is relatively compact,
we obtain
\begin{equation}
\label{eq;20.8.20.12}
 \int_{X_P\setminus X_{1,P}}
  \sqrt{-1}\Lambda\Tr F\bigl((h_0)_{\geq j}\bigr)
  <\infty.
\end{equation}
By (\ref{eq;20.8.20.10}),
(\ref{eq;20.8.20.11})
and
(\ref{eq;20.8.20.12}),
we obtain
\[
 \int_{X\setminus D}
 \sqrt{-1}\Lambda
 \Tr F\bigl((h_0)_{\geq j}\bigr)=-\infty<0.
\]
If a non-zero subbundle $E$
of $\hyperk_{X\setminus D,r}$
satisfies
$\theta(0_{X\setminus D,r})E\subset
E\otimes K_{X\setminus D}$,
there exists $1\leq j\leq r$
such that
$E=\bigoplus_{i=j}^rK^{(r+1-2i)/2}_{X\setminus D}$.
Therefore,
$(\hyperk_{X\setminus D,r},\theta(0_{X\setminus D,r}),h_0)$
is analytically stable
with respect to the $G_r$-action.

\hfill\qed

\vspace{.1in}

By Proposition \ref{prop;20.6.11.10},
there exists $h\in\Harm(0_{X\setminus D,r})$
such that
$h$ and $h_0$ are mutually bounded.
Then, $h$ satisfies the desired conditions.
\hfill\qed

\subsection{Existence of harmonic metrics
  in the potential theoretically hyperbolic case}

\subsubsection{Solvability of the Poisson equation and
   the existence of harmonic metrics}

Let $G$ be a compact Lie group
with a character $\kappa:G\lrarr S^1$.
Let $X$ be an open Riemann surface
equipped with a $G$-action.
Let $g_X$ be any $G$-invariant K\"ahler metric of $X$.
Let $\Lambda$ denote the adjoint of the multiplication
of the K\"ahler form associated with $g_X$.
Let $(E,\delbar_E,\theta)$ be a $(G,\kappa)$-homogeneous
Higgs bundle on $X$.
Let $h_0$ be a $G$-invariant Hermitian metric of $E$.

A function on $X$ is called locally bounded
if its restriction to any compact set is bounded.

\begin{condition}
 Assume that there exists a $G$-invariant
 $\real_{\geq 0}$-valued
 locally bounded function $\alpha$
 on $X$ such that
\[
 \sqrt{-1}\Lambda\delbar\del\alpha\geq
  \bigl|\Lambda F(h_0)\bigr|_{h_0}
\]
 in the sense of distributions.
 Note that this condition is independent of
 the choice of $g_X$.
\end{condition}

\begin{prop}
\label{prop;20.7.12.1}
There exists a $G$-invariant harmonic metric $h$ of
the Higgs bundle $(E,\delbar_E,\theta)$
 such that
\begin{equation}
\label{eq;20.5.31.1}
 \log\Tr(h\cdot h_0^{-1})\leq \alpha +\log\rank(E),
 \quad
 \log\Tr(h_0\cdot h^{-1})\leq \alpha +\log\rank(E).
\end{equation}
\end{prop}
 \pf
 Let $X_1\subset X_2\subset\cdots$
 be a smooth exhaustive sequence of $X$.
 Let $h_i$ be a harmonic metric of
 $(E,\delbar_E,\theta)_{|X_i}$
 such that $h_{i|\del X_i}=h_{0|\del X_i}$.
 Let $s_i$ be the automorphism of $E_{|X_i}$
determined by
$h_{i}=h_{0|X_i}s_i$.
By (\ref{eq;20.8.16.6}),
we obtain
\[
 \sqrt{-1}\Lambda\delbar\del\log\Tr(s_i)
 \leq
 \Bigl(
 |F(h_0)|_{h_0}
 \Bigr)_{|X_i}.
\]
By the condition, we obtain the following on $X_i$:
\[
  \sqrt{-1}\Lambda
  \delbar\del
  \Bigl(
  \log\Tr(s_i)
  -\alpha
 \Bigr)
 \leq 0.
\]
Note that
$s_{i|\del X_i}=\id_{E_i|\del X_i}$.
Hence, we obtain
\[
 \log\Tr(s_i)
 \leq
 \alpha_{|X_i}
 +\log (\rank E)
 -\min_{\del X_i}(\alpha)
 \leq
  \alpha_{|X_i}
 +\log (\rank E).
\]
Similarly,
we obtain
\[
 \log\Tr(s_i^{-1})
 \leq
  \alpha_{|X_i}
 +\log (\rank E).
\]
Then, we obtain the claim of the proposition.
\hfill\qed

 \begin{cor}
If $\det(h_0)$ is flat,
there exists a $G$-invariant harmonic metric $h$
 of $(E,\delbar_E,\theta)$
  satisfying {\rm(\ref{eq;20.5.31.1})}
  and $\det(h)=\det(h_0)$.
  Moreover, if $\alpha$ is bounded,
  and if  $\int_X |F(h_0)|<\infty$,
  then
 $(\delbar_E+\theta)(h\cdot h_0^{-1})$
 is $L^2$.
 \hfill\qed
 \end{cor}

\begin{rem}
In {\rm\cite{Ni1}},
the solvability of the Poisson equation
is applied to the solvability of
the Hermitian-Einstein equation.
Proposition {\rm\ref{prop;20.7.12.1}}
is an analogue in the context of
harmonic metrics on Riemann surfaces.
\end{rem}
 
\subsubsection{Potential theoretically hyperbolic
Riemann surfaces}

  Let $X$ be a potential theoretically hyperbolic Riemann surface,
  i.e., there exists a non-constant non-positive subharmonic function
  on $X$.
  It is well known that there exists a Green function
  $\ttG(p,q):X\lrarr \real_{<0}$
  characterized as follows
  (see \cite[\S7]{Varolin}):
\begin{itemize}
 \item $\ttG(p,\bullet)$ is harmonic on $X\setminus\{p\}$.
 \item Let $(X_p,z_p)$ be a holomorphic coordinate neighbourhood
       such that $z_p(p)=0$. Then,
       $\ttG(p,\bullet)-\log|z_p|^2$ is harmonic on $X_p$.
 \item If $H$ is any other subharmonic function satisfying
       the above two conditions, then
       $\ttG(p,\bullet)\geq H$.
\end{itemize}
 Recall that for any compact support $2$-form $\psi$ on $X$
 we obtain
\[
 \sqrt{-1}\del\delbar\Bigl(
  \frac{1}{2\pi}\int \ttG(p,q)\psi(q)
  \Bigr)=\psi.
\]
 \begin{lem}
  Let $g$ be any K\"ahler  metric of $X$.
  The volume form is denoted by $\vol$.
 For a function $v:X\lrarr \real_{\geq 0}$,
suppose that
$|\ttG(p,q) (v\cdot\vol)(q)|$ is integrable on $X$ for any $p$.
Then, by setting
 $\alpha(v)(p):=-\frac{1}{2\pi}\int \ttG(p,q)(v\vol)(q)$,
we obtain $\alpha(v)\geq 0$
and 
 $\sqrt{-1}\Lambda\delbar\del\alpha(v) =v$.
 \hfill\qed
 \end{lem}

The following is clear by the construction of
the Green function (see \cite[\S7]{Varolin}).
\begin{lem}
 Let $K$ be any compact subset of $X$.
 Let $N$ be any relatively compact neighbourhood of $K$
 in $X$.
 Then,
 $\ttG(p,q)$ is bounded on
 $K\times (X\setminus N)$.
 \hfill\qed
\end{lem}

\begin{cor}
 Let $\psi$ be a compact support $2$-form on $X$.
 Then, there exists a bounded function $\alpha$ on $X$
 such that
 $\del\delbar\alpha=\psi$.
 \hfill\qed
\end{cor}

Let $(E,\delbar_E,\theta)$ be a
$(G,\kappa)$-homogeneous Higgs bundle on $X$.
Let $h_0$ be a $G$-invariant Hermitian metric of $E$.
 Suppose that
\[
  \int \ttG(p,q)\bigl(|F(h_0)|_{h_0}\bigr)\vol(q)<\infty
\]
 for any $p$.

 \begin{prop}
  There exists a $G$-invariant harmonic metric $h$ of
  $(E,\delbar_E,\theta)$ such that
\[
\max\Bigl\{
 \log\bigl|h\cdot h_0^{-1}\bigr|_{h_0},\,\,
 \log\bigl|h_0\cdot h^{-1}
  \bigr|_{h_0}
  \Bigr\}
  \leq
  \alpha(|F(h_0)|_{h_0})+\log\rank(E).
\]
  If $\det(h_0)$ is flat,
  $h$ also satisfies
 $\det(h)=\det(h_0)$.
\hfill\qed
 \end{prop}

\begin{cor}         
 Suppose that the support of $F(h_0)$ is compact.
  Then, there exists a $G$-invariant harmonic metric $h$ of
  $(E,\delbar_E,\theta)$
  such that $h$ and $h_0$ are mutually bounded.
  If $\det(h_0)$ is flat,
  it also satisfies
  (i) $\det(h)=\det(h_0)$,
  (ii) $(\delbar_E+\theta)(h\cdot h_0^{-1})$ is $L^2$.
 \hfill\qed
\end{cor}

\subsubsection{Upper half plane}

Let us state a consequence
in the case $X=\hyperh=\{z\in\cnum\,|\,\Image(z)>0\}$.
Let $(x,y)$ be the real coordinate system
determined by $z=x+\sqrt{-1}y$.

\begin{prop}
 Let $h_0$ be a $G$-invariant Hermitian metric of $E$
 such that
 $|F(h_0)|=O((1+y)^{-2})$.
 Then, there exists a $G$-invariant harmonic metric $h$
 of $(E,\delbar_E,\theta)$
 such that
 \[
 \max\Bigl\{
 \log\bigl|
  h\cdot h_0^{-1}
 \bigr|_{h_0},\,
\log\bigl|
  h_0\cdot h^{-1}
 \bigr|_{h_0}
 \Bigr\}
 =O\Bigl(\log(2+y)
  \Bigr).
 \]
 If $\det(h_0)$ is flat,
 we obtain
 $\det(h)=\det(h_0)$.
 If 
 $|F(h_0)|=O((1+y)^{-2-\epsilon})$
 for a positive constant $\epsilon$,
 then $h$ and $h_0$ are mutually bounded.
\end{prop}
\pf
For a function $f$ on $\hyperh$
such that $f=O\bigl((1+y)^{-2}\bigr)$,
we set
\[
 u(x,y)=-\frac{1}{4\pi}
  \int_{\hyperh}\log\Bigl(
  \frac{(x-\xi)^2+(y-\eta)^2}{(x-\xi)^2+(y+\eta)^2}
    \Bigr)
  f(\xi,\eta)\,d\xi\,d\eta.
\]
Then, $-(\del_x^2+\del_y^2)u=f$.
If $f$ is positive,
then $u$ is also positive on $\hyperh$.
The rest follows from Lemma \ref{lem;20.4.21.3}
and Proposition \ref{prop;20.7.12.1}.
\hfill\qed

\subsection{Uniqueness of harmonic metrics
  in the potential theoretically parabolic case}

  Let $X$ be a potential theoretically parabolic Riemann surface,
  i.e., $X$ is non-compact,
and any non-positive subharmonic function on $X$ is constant.
 Let $(E,\delbar_E,\theta)$   be a Higgs bundle on $X$.
\begin{prop}
\label{prop;20.7.5.20}  
 Let $h_i$ $(i=1,2)$ be harmonic metrics of
 the Higgs bundle $(E,\delbar_E,\theta)$
 which are mutually bounded.
 Then, there exists a decomposition
\[
 (E,\delbar_E,\theta)=
 \bigoplus_{j=1}^m(E_j,\delbar_{E_j},\theta_j)
\]
 such that
 (i) the decomposition is orthogonal with respect to
 both $h_1$ and $h_2$,
 (ii) there exist positive numbers $a_j$
 such that
 $h_{1|E_j}=a_j h_{2|E_j}$. 
 \end{prop}
\pf
Let $s$ be the automorphism of $E$
determined by
$h_2=h_1s$.
Then,
$\log\Tr(s)$ 
is a bounded subharmonic function on $X$.
Because $X$ is assumed to be
potential theoretically parabolic,
$\log\Tr(s)$ is constant.
According to (\ref{eq;20.8.16.3}),
it implies that $(\delbar_E+\theta)s=0$,
and hence we obtain the desired decomposition.
\hfill\qed

\vspace{.1in}
The following corollary of the proposition
has been well known after \cite{s1,s2}.

\begin{cor}
\label{cor;20.7.6.1}
 Let $X$ be the complement of a finite subset $D$
of a compact Riemann surface $\Xbar$.
Then, the claim of Proposition {\rm\ref{prop;20.7.5.20}} holds.
\end{cor}
\pf
Let $f$ be an $\real_{\leq 0}$-valued subharmonic function
on $X$.
According to \cite[Theorem 3.3.25]{Noguchi-Ochiai},
it uniquely extends to a subharmonic function on $\Xbar$.
By the maximum principle,
we obtain that $f$ is constant,
i.e., $X$ is potential theoretically parabolic.
\hfill\qed

\begin{rem}
\label{rem;20.7.10.1}
 It is also easy and standard to prove the claim of
 the corollary directly.
 Indeed, in the proof of Proposition {\rm\ref{prop;20.7.5.20}},
 without any assumption on $X$,
 we obtain the boundedness of
 $\log \Tr(s)$
 and $\sqrt{-1}\Lambda\delbar\del\log\Tr(s)\leq 0$ on $X$.
 If $X=\Xbar\setminus D$,
 $\sqrt{-1}\Lambda\delbar\del\log\Tr(s)\leq 0$
 holds on $\Xbar$
 according to {\rm\cite[Lemma 2.2]{s2}}.
 Hence, $\log\Tr(s)$ is constant,
 and we obtain $(\delbar_E+\theta)s=0$.
\end{rem}

\end{document}